\documentclass{elsart}
\usepackage{amssymb}
\usepackage{color}
\usepackage{graphicx}
\usepackage{subfigure}
\usepackage{amsmath}
\usepackage{amsfonts}
\usepackage{epstopdf}

\usepackage{morefloats}

\renewcommand{\qed}{\hfill \nobreak \ifvmode \relax \else
      \ifdim\lastskip<1.5em \hskip-\lastskip
      \hskip1.5em plus0em minus0.5em \fi \nobreak
      \vrule height0.75em width0.5em depth0.25em\fi}

\def\pd#1#2{\frac{\partial #1}{\partial #2}}

\renewcommand{\vec}[1]{\mbox{\boldmath \small $#1$}}
\def\dd#1#2{\frac{\partial #1}{\partial #2}}

 \usepackage[ 
              unicode=true,
             pdfstartview=FitH,
             CJKbookmarks=true,
             bookmarksnumbered=true,
             bookmarksopen=true,
             colorlinks=true, 
             citecolor=magenta,
             linkcolor=blue,
             linktocpage=true,
             ]{hyperref}       

\usepackage{multirow}
\allowdisplaybreaks
    \numberwithin{equation}{section}

    \newtheorem{Example}{Example} [section]
        \def\DG{DG methods}
        \def\CDG{central DG methods}

  \newtheorem{Remark}{Remark}[section]

\begin{document}

\begin{frontmatter}

\title{
Runge-Kutta discontinuous Galerkin methods
for the special relativistic magnetohydrodynamics}
  \author{Jian Zhao}
\ead{everease@163.com}
\address{HEDPS, CAPT \& LMAM, School of Mathematical Sciences, Peking University,
Beijing 100871, P.R. China}
\author[label2]{Huazhong Tang}
\thanks[label2]{Corresponding author. Tel:~+86-10-62757018;
Fax:~+86-10-62751801.}
\ead{hztang@math.pku.edu.cn}
\address{HEDPS, CAPT \& LMAM, School of Mathematical Sciences, Peking University,
Beijing 100871, P.R. China; School of Mathematics and Computational Science,
 Xiangtan University, Hunan Province, Xiangtan 411105, P.R. China}
 \date{\today{}}

\maketitle

\begin{abstract}
 This paper develops   $P^K$-based non-central and central Runge-Kutta discontinuous Galerkin (DG) methods with WENO limiter for the one- and two-dimensional special relativistic magnetohydrodynamical (RMHD) equations, $K=1,2,3$.
The non-central DG methods are locally divergence-free, while the central DG are
``exactly'' divergence-free  but have to find
two approximate solutions defined on mutually dual meshes.
The adaptive WENO limiter   first identifies
 the ``troubled'' cells by using a modified TVB minmod function,
and then uses the  WENO technique to locally reconstruct a new
polynomial of degree $(2K+1)$ inside the ``troubled'' cells
replacing the DG solution by  based on the cell average values of the DG solutions in the neighboring cells as well as the original cell averages of the ``troubled'' cells.
The WENO limiting procedure does not destroy  the locally or ``exactly'' divergence-free property of magnetic field
and is only employed for finite ``troubled'' cells so that the computational cost can be as little as possible.
Several test problems in one and two dimensions are solved by using
our non-central and central Runge-Kutta DG methods with WENO limiter. The numerical results demonstrate
that our methods are stable, accurate, and robust
in resolving complex wave structures.

\end{abstract}

\end{frontmatter}

\section{Introduction}

Relativistic hydrodynamics (RHD) or relativistic magnetohydrodynamics (RMHD)
play major roles in astrophysics, nuclear physics, plasma physics and other fields.
They are necessary in situations where the local velocity of the flow is
close to the light speed in vacuum or where the local internal energy
density is comparable (or larger) than the local rest mass density
of the fluid.
For example, in the formation of neutron stars and black holes,
and the high-speed jet of physical phenomena, the relativistic effect can not be neglected
so that  the  RHDs or  RMHDs are needed.
The dynamics of  RMHD system requires solving highly
nonlinear equations so that the analytic treatment of practical
RMHD problems is extremely difficult.
Numerical simulation has become an important way in studying RHDs and RMHDs.
 In the past few decades, significant progress
is made and several numerical methods have been developed to investigate
the RMHD equations.

 The pioneering numerical work may date back to
the finite difference code via artificial viscosity for the
spherically symmetric general RHD equations in the Lagrangian coordinate
\cite{May-White1966,May-White1967}.
Wilson first
attempted to solve multi-dimensional RHD equations in the
Eulerian coordinate by using the finite difference method
with the artificial viscosity technique  \cite{Wilson:1972}.
After that, a lot of modern shock-capturing methods in non relativistic hydrodynamics were extended
to the RHD and RMHD equations.
For example, some latest works are the adaptive moving mesh methods
 \cite{HeAdaptiveRHD,HeAdaptiveRMHD},   second-order
 generalized problem schemes \cite{YangGRPRHD1D,YangGRPRHD2D,Wu-Tand-SISC2016} and third-order GRP scheme
 \cite{WuGRPRHD3}, the locally evolution Galerkin method \cite{WuEGRHD},
 approximate Riemann solvers based on the local linearization \cite{GodunovRMHD,Koldoba:2002},
 TVD scheme \cite{tvdRMHD},
 HLL scheme \cite{Zanna:2003}, HLLC schemes \cite{MignoneHLLCRMHD,Honkkila:2007},
 the kinetic scheme \cite{QamarKinetic2004}, and adaptive mesh refinement methods \cite{Anderson:2006,Host:2008} etc.
Recently three physical-constraints-preserving (PCP)
schemes were developed for the special RHD equations.
They are the high-order accurate PCP finite difference weighted essentially non-oscillatory
(WENO) schemes and discontinuous Galerkin (DG) methods proposed in \cite{Wu-Tand-JCP2015,Wu-Tang-RHD2016,QinShuYang-JCP2016}.
Moreover,  the set of admissible states  and the PCP schemes of
the ideal RMHDs was also studied
for the first time in \cite{Wu-Tang-RMHD2016}, where  the importance of divergence-free fields
was revealed in achieving PCP methods especially.

The magnetohydrodynamics mainly studies the interaction between magnetic field and conducting fluid.
Compared with the RHD equations,
the RMHD equations not only contain
more number of the equations with   complicated  forms,
but also an additional divergence-free  constraint of magnetic field.
Violating such constraint leads
to nonphysical plasma transport orthogonal to the magnetic field.
Some special techniques must be employed to preserve that constraint.
Up to now, three popular approaches have been suggested
in the context of  MHD  equations.
They are the eight-wave formulation of the MHD equations
\cite{Powell:1994}, the constrained transport  method
 \cite{Evans:1988}, and the the projection method \cite{Brackbill:1980}.
 In the DG framework, locally divergence-free  DG  methods \cite{DGMHDLi}
and ``exactly'' divergence-free central DG methods \cite{CentralMHD} have been developed.
All of the above works are on the non-relativistic MHD equations and
some of them have been extended to the relativistic case.

The DG  methods have been rapidly developed  in recent decades and become a kind of important methods in computational fluid dynamics. It is easy to achieve high order accuracy, be suitable for parallel computing, and adapt to complex boundary.
The DG method was first developed by Reed and Hill \cite{ReedFirstDG} to solve scalar, steady-state linear hyperbolic equation, but it had not been widely used.
A major development of the DG methods was carried out by Cockburn and coworkers in a series of papers \cite{CockburnDG1,CockburnDG2,CockburnDG3,CockburnDG4,CockburnDG5}, where
the DG spatial approximation was combined with explicit Runge-Kutta time discretization to develop
the Runge-Kutta DG methods, and a general framework of DG methods was established  for the nonlinear equation or system.
The Runge-Kutta DG  methods have gotten a wide range of research and application,
 such as the Euler equations \cite{TangRKDG,momentlimiter,Remade2003},
Maxwell equations \cite{DGMaxwell}, nonlinear Dirac equations \cite{ShaoTangDiracDG} etc.
Moreover, the DG methods have also been used to solve other partial differential
equations, such as convection-diffusion type equation or system \cite{BassiDGNS,CockLDG} and
Hamilton-Jacobi equation  \cite{hu1999DGHamilton,AnaDGHJ,CDGHamiltonLi} etc.
The readers are referred to the review article \cite{CockDGReview}.

The \CDG{} \cite{LiuCentralDG}
were developed by combing the \DG{} and central scheme \cite{LiuCentralSch} and  found two approximate solutions
defined  on mutually dual meshes. Although
 two approximate solutions are redundant, the numerical flux may be avoided
due to  the use of the solution on the dual grid to calculate the
flux at the cell interface. It is one of the advantages of the central scheme.
Because the \CDG{} can be regarded as a variant of traditional/no-central \DG{}, they keep  many advantages of no-central \DG{}, such as
compact stencil and  parallel implementation etc.
Moreover, the \CDG{} allow a larger CFL number than no-central \DG{}
and may reduce numerical oscillations for some problems.
Up to now, the \CDG{} have also been used to solve
the Euler equations \cite{LiuCentralDG},
 the ideal magneto-hydrodynamical equations
\cite{CentralMHD,arbiCDGMHD},   the special
relativistic hydrodynamical equations
\cite{ZhaoTang2016A,Wu-Tang-RHD2016}  and so on.

When the strong discontinuity appears in the solution,
the numerical oscillations in the DG solutions should be suppressed
after each Runge-Kutta inner stage or after some complete Runge-Kutta
steps by using the nonlinear limiter, which  is a commonly used technique
in the modern shock-capturing methods for hyperbolic conservation laws.
The commonly used limiter  is  the minmod limiter, which limits
the slope of solution such that the values of limited solution in the cell
falls in the certain interval determined by the cell average values of neighboring cells.
It has good robustness but is  only first-order accurate near extreme points.
 Cockburn et al. gave modified TVB minmod limiter and applied it to the non-central DG methods.
 The modified TVB minmod limiter does not limit the solution near extreme points by
choosing a parameter $M$, thus it does not destory the accuracy  of DG methods near the extreme point.
In general, for nonlinear equation, the parameter $M$ is dependent on the problem. Moreover, for $K\ge3$,
the accuracy of $P^K$-based DG methods may still be destroyed because
 more than three of the higher order moments will be set to zero by the  modified TVB minmod limiter.
Besides that commonly used limiter,
some other limiters are porposed, such as
the moment based limiters \cite{momentlimiter}  and its improvement
 \cite{Krivodonova2007} etc. Those limiters may suppress numerical oscillations
  near the discontinuity, however, the accuracy of \DG{}
  can be reduced in   some region.

In the modern shock-capturing methods, the ENO (non-oscillatory essentially)  and WENO   methods have been widely used \cite{Shu-SIEV2009}
and are more robust than the slope limiter methodology,
especially for high order schemes. An attempt has been made to use an ENO or WENO methodology
as limiter for the DG methods \cite{QiuShuWENOlimiter,QiuShuWENOlimiterUn,QiuShuWENOlimiter3D}.
The WENO limiter  first identifies  the ``troubled'' cells  by using a modified TVB minmod function,
and then a new  polynomial inside the ``troubled'' cells is locally reconstructed to
replace the DG solution by using the  WENO technique and the cell average values of
the DG solutions in the neighboring cells as well as the original cell averages of the ``troubled'' cells.
Because the WENO limiter is only employed for finite ``troubled'' cells,
the computational cost can be as little as possible.

The aim of this paper develops  $P^K$-based non-central and central Runge-Kutta discontinuous
Galerkin (DG) methods with WENO limiter for the one- and two-dimensional special RMHD equations,
$K=1,2,3$. The former is locally divergence-free, while the latter is ``exactly'' divergence-free.
It
is organized as follows. Section \ref{chap:rmhddg} introduces the special RMHD equations
and  calculation of eigenvalues.
 Sections \ref{sec:rmhdlocaldg} and \ref{sec:rmhdcdg} give  $P^K$-based locally and ``exactly'' divergence-free DG methods  with WENO limiter for the special RMHD equations, respectively, $K=1,2,3$.
 Section \ref{sec:rmhdnum} conducts several numerical experiments
to demonstrate the accuracy and efficiency of the proposed DG methods.
 Conclusions are given in Section \ref{Section-conclusion}.

\section{Relativistic magnetohydrodynamical equations}
\label{chap:rmhddg}
This section introduces the relativistic magnetohydrodynamical (RMHD) equations
and calculation of the eigenvalues  for the Jacobian matrix.

The RMHDs is investigating the interaction between magnetic field and conducting fluid.
In the covariant form,  the four-dimensional space-time RMHD equations may
be written as follows \cite{AnileRFMF,Anile:Pennisi:1987}
\begin{equation}
  \label{eq:RMHDcov-form}
    \begin{cases}
      \partial_\alpha (\rho u^\alpha) = 0,   \\
      \partial_\alpha \Big((\rho h+|\vec b|^2) u^\alpha u^\beta
      -b^\alpha b^\beta+ p_{tot} g^{\alpha \beta} \Big) = 0,\\
       \partial_\alpha (u^\alpha b^\beta-u^\beta b^\alpha)=0,
    \end{cases}
\end{equation}
which stand for the laws of local baryon number conservation and
 energy-momentum conservation, and the induction equation for the magnetic field.
In \eqref{eq:RMHDcov-form},
the Greek indices $\alpha$ and $\beta$ run from 0 to 3,
$\partial_\alpha =\partial_{x^\alpha}$ denotes the covariant derivative with
$x^\alpha=(ct,x_1,x_2,x_3)^T$,
 $u^\alpha=\gamma(c,v_1,v_2,v_3)^T$ stands for
the four-velocity vector,
 $\gamma=1/\sqrt{1-|\vec
  v|^2/c^2}$ is the Lorentz factor with the fluid velocity vector $\vec v:=(v_1,v_2,v_3)^T$,
  $g^{\alpha\beta}$ denotes the metric tensor, which is restricted to
  to the Minkowski tensor in this paper, i.e.
  $\big( g^{\alpha \beta}\big)_{4\times4} =
\mbox{diag}\{-1,1,1,1\}$,
the relativistic enthalpy $h$ is defined by
$$h=1+\frac{e}{c^2}+\frac{p}{\rho c^2},$$
where $e$ is the specific internal energy and related to other thermodynamic
or state variables such as the temperature, the pressure,
the volume, or the internal energy etc. by the equation of state
(EOS). The simplest EOS is the ideal gas law given by
 \begin{equation}\label{RMHDgamma-law}
 p  = {(\Gamma-1)\rho e},
 \end{equation}
where $\Gamma$ is the adiabatic index (also ratio of specific heats).
In comparison to the RHD equations, 
the RMHD equations \eqref{eq:RMHDcov-form}  involve  the four-magnetic field vector
$$ b^\alpha=\gamma\left(\frac{\vec v\cdot \vec B}{c},
B_1/\gamma^2+v_1\frac{\vec v\cdot \vec
  B}{c^2},B_2/\gamma^2+v_2\frac{\vec v\cdot \vec
  B}{c^2},B_3/\gamma^2+v_3\frac{\vec v\cdot \vec B}{c^2}\right)^T,$$
where $\vec B=(B_1,B_2,B_3)$ is the magnetic field, the total pressure $p_{tot}$
consists of the gas pressure $p$
and magnetic pressure $p_m$, that is, $p_{tot}=p+p_m$, here the magnetic pressure
$p_m=\frac{1}{2}b^\alpha b_\alpha$,  $b_\alpha:=g_{\alpha\beta}b^\beta$,
and $g_{\alpha\beta}$ is the inverse of metric matrix $g^{\alpha\beta}$.

Throughout this paper, units in which the speed of light is equal
to one will be used so that
\begin{align*}
&x^\alpha=(t,x_1,x_2,x_3)^T,\quad u^\alpha=\gamma(1,v_1,v_2,v_3)^T,\\
 &b^\alpha=\gamma\big(\vec v\cdot \vec B, B_1/\gamma^2+v_1(\vec v\cdot
\vec B),B_2/\gamma^2+v_2(\vec v\cdot \vec B),B_3/\gamma^2+v_3(\vec
v\cdot \vec B)\big)^T.
\end{align*}

It is obvious from \eqref{eq:RMHDcov-form} that in modeling flow at speeds where
relativistic effects become important, space and time become
intrinsically coupled and the governing equations of the ideal RMHDs
become more complicated. Nonetheless,
it is still possible to write \eqref{eq:RMHDcov-form}
into a first-order system of time evolution equations
 in some fixed reference frame (i.e. the so-called lab frame)
as follows
 \begin{equation}
   \nabla \cdot \textbf B=0,
   \label{eq:divfree}
   \end{equation}
\begin{equation}\label{eqn:RMHDconeqn}
    \displaystyle\frac{\partial \vec{U}}{\partial t} +
    \sum^3_{i=1}\frac{\partial \vec{F}_i(\vec{U})}{\partial x_i}=0,
 \end{equation}
where $\vec{U}$ is the conservative vector, $\vec{F}_i$ is the flux in the $x_i$ direction, $i=1,~2,~3$,
and their expressions  are given by
 \begin{align*}
  \vec{U} =& \Big(D, m_1, m_2, m_3,B_1,B_2,B_3, E\Big)^T,\\
  \vec{F}_1 =& \Big(Dv_1, m_1 v_1 -B_1 b_1/\gamma+ p_{tot},  m_2
  v_1-B_1 b_2/\gamma,  \\
  & m_3 v_1-B_1 b_3/\gamma, 0, B_2 v_1-B_1 v_2, B_3
  v_1-B_1v_3,m_1 \Big)^T,\\
  \vec{F}_2 =& \Big(Dv_2, m_1 v_2 -B_2 b_1/\gamma,  m_2
  v_2-B_2 b_2/\gamma+p_{tot},  \\
  & m_3 v_2-B_2 b_3/\gamma, B_1 v_2-B_2 v_1, 0, B_3
  v_2-B_2v_3,m_2 \Big)^T,\\
  \vec{F}_3 =& \Big(Dv_3, m_1 v_3 -B_3 b_1/\gamma,  m_2
  v_3-B_3 b_2/\gamma,  \\
  & m_3 v_3-B_3 b_3/\gamma+p_{tot}, B_1 v_3-B_3 v_1, B_2 v_3-B_3 v_2, 0,m_3 \Big)^T
   ,\label{eqn:rmhdcondef}
 \end{align*}
here $D=\rho \gamma$, $m_i=(\rho h \gamma^2+|\vec B|^2) v_i-(\vec{v}\cdot
\vec{B}) B_i$, and $E=Dh\gamma-p_{tot}+|\vec B|^2$
denote the mass, $x_i$-momentum, and energy densities in the lab frame, respectively.
Eq. \eqref{eq:divfree} is a divergence-free constraint on the magnetic field,
and the solutions of \eqref{eqn:RMHDconeqn} should satisfy such constraint at any time.
Numerically preserving such condition
is very non-trivial  but important for the robustness of numerical scheme, and has to be
respected. In physics, numerically incorrect magnetic field topologies may lead to nonphysical plasma
transport orthogonal to the magnetic field. The condition \eqref{eq:divfree} is also very crucial for
the stability of induction equation. The existing numerical experiments in the non-relativistic
MHD case have also indicated that violating the divergence free condition of magnetic field may lead
to numerical instability and nonphysical or inadmissible solutions.
The flux $\vec F_i$ in \eqref{eqn:RMHDconeqn} cannot be explicitly expressed as
a function of $\vec U$
but can be cast into an explicit function of primitive variable vector $\vec
V:=(\rho,v_1,v_2,v_3,B_1,B_2,B_3,p)^T$.
Thus if giving the value of  $\vec U$, then one has to get the value of $\vec V$ in order to
 calculate $\vec F_i$.
Up to now, several approaches have been suggested to recover
the primitive variables from the conservative vector in the literature,
e.g. six numerical approaches discussed in \cite{nobleRMHDPrimitive}.
The approach in \cite{Koldoba:2002} is used here.
If introducing an auxiliary variable $\theta:=\rho h\gamma^2>0$ and denoting $\vec m:=(m_1,m_2,m_3)^T$,
then it is easy to prove the following identities
$$\vec v\cdot \vec B=\theta^{-1} (\vec m\cdot \vec B),$$
and
$$ \vec v=\frac{\vec m+ \theta^{-1} (\vec m\cdot \vec B) \vec
  B}{\theta+|\vec B|^2},$$
$$
p_m=\frac{b^\alpha b_\alpha}{2}=\frac{1}{2}\big(\frac{|\vec B|^2}{\gamma^2}+\frac{(\vec
  m\cdot \vec B)^2}{\theta^2}\big),
$$
where the Lorentz factor $\gamma$ is expressed as follows
$$\gamma=(1-|\vec v|^2)^{-\frac{1}{2}}=\big( 1-\frac{\theta^2 |\vec
  m|^2+2\theta(\vec m\cdot \vec B)^2+(\vec m\cdot \vec B)^2 |\vec
  B|^2}{\theta^2 (\theta+|\vec B|^2)^2}\big)^{-\frac{1}{2}},$$
and the density $\rho$ and pressure $p$ are written as
$$\rho=\frac{D}{\gamma},\quad p=\frac{\Gamma-1}{\Gamma}
(\frac{\theta}{\gamma^2}-\frac{D}{\gamma}).$$
Substituting the above equations into the energy density expression
$E=\rho h \gamma^2-p-p_m+|\vec B|^2$ gives a nonlinear equation with respect to
$\theta$ as follows
\begin{equation}\label{eq:RMHDCon2Pri}
\theta-\frac{\Gamma-1}{\Gamma}(\frac{\theta}{\gamma^2}-\frac{D}{\gamma})-\frac{1}{2}\big
(\frac{|\vec B|^2}{\gamma^2}+\frac{(\vec m\cdot \vec
B)^2}{\theta^2}\big )+|\vec B|^2-E=0,
\end{equation}
which may be solved by any standard root finding algorithm
such as the Newton¨CRaphson method.
Once $\theta$ is found to some
accuracy, the Lorentz factor $\gamma$,
 velocity $\vec v$,  density $\rho$,
and  gas pressure $p$ can be orderly calculated.

\begin{Remark}
It is not difficult to know that
  $\theta_{min}:=\theta^\star\le\theta<\Gamma E=:\theta_{max}$, where $\theta^\star$
  satisfies
  $\gamma(\theta^\star)=1+\epsilon$, and  $\epsilon$ is a small positive number.
 Thus in practical computations, the initial guess for the Newton-Raphson method of \eqref{eq:RMHDCon2Pri}
 may be chosen as
  $(\theta_{min}+\theta_{max})/2$, and the iteration number is generally less than 8.
 \end{Remark}

The characteristic structure of the RMHD equations was first studied in \cite{Anile:Pennisi:1987}.
The eigenvalues and eigenvectors of the Jacobian matrix  is needed in
our numerical methods for \eqref{eqn:RMHDconeqn}.
Here gives the calculation of eigenvalues.
Without loss of generality, consider $x_1$-split system.
Because the component of $\vec F_1$ corresponding to $B_1$ is zero,
$B_1=$const  and $x_1$-split system consists of the following seven equations
\begin{equation}
  \label{eq:RMHDxsplit}
    \displaystyle\frac{\partial \vec{U}}{\partial t} +\frac{\partial \vec{F}_1(\vec{U})}{\partial x_1}=0,
 \end{equation}
where
 \begin{align*}
  \vec{U} =& \Big(D, m_1, m_2, m_3,B_2,B_3, E\Big)^T,\\
  \vec{F}_1 =& \Big(Dv_1, m_1 v_1 -B_1 b_1/\gamma+ p_{tot},  m_2
  v_1-B_1 b_2/\gamma,  \\
  & m_3 v_1-B_1 b_3/\gamma, B_2 v_1-B_1 v_2, B_3
  v_1-B_1v_3,m_1 \Big)^T.
 \end{align*}
 The eigensystem of the Jacobian matrix $\partial \vec F_1/\partial \vec U$
 can be found by slightly modifying the eigen-system
for the one-dimensional  RMHD equations, which have seven
waves: two Alfv\'{e}n, two fast and two slow magnetosonic waves
(also magnetoacoustic wave), and an entropy wave, whose
the speeds are denoted by
$\lambda_1^{a,\pm}$,  $\lambda_1^{f,\pm}$, $\lambda_1^{s,\pm}$, and $\lambda_1^{e}$, separately.
They satisfy \cite{RMHDEigenAnton}
$$\lambda_1^{f,-}\le \lambda_1^{a,-}\le \lambda_1^{s,-} \le \lambda_1^e \le\lambda_1^{s,+}
\le \lambda_1^{a,+} \le \lambda_1^{f,+},$$
where the speed of entropy
wave  $\lambda_1^e$ is equal to $v_1$, speeds of two Alfv\'{e}n wave satisfy
the following quadratic equation
$$(\rho h+|\vec b|^2)\gamma^2(v_1-\lambda)^2-(b^1-b^0\lambda)^2=0,$$
which can be directly solved by the root formula of quadratic equation with one unknown,
 while speeds of four  magnetoacoustic wave satisfy the following quartic equation \cite{MignoneHLLCRMHD}
\begin{equation}\label{eq:RMHDeigvalue}
  \rho h(1-c_s^2) a^4=(1-\lambda^2)\big[( |\vec b|^2+\rho h c_s^2)
  a^2-c_s^2 (b^1-\lambda b^0)^2\big],
\end{equation}
where $a=\gamma(\lambda-v_1)$, and $c_s$ is the sound speed and becomes
$c_s=\sqrt{\Gamma p/(\rho h)}$  for the perfect gas.
In some cases, Eq.~\eqref{eq:RMHDeigvalue} may be simplified and solved.
In the following it is discussed in three cases.
\begin{itemize}
  \item If the fluid velocity $\vec v=\vec 0$, then Eq.~\eqref{eq:RMHDeigvalue}
  reduces to
\begin{equation*}
  \big(\rho h+|\vec b|^2\big) \lambda^4-\big( |\vec b|^2+\rho h c_s^2+B_1^2 c_s^2\big)
  \lambda^2+c_s^2  B_1^2=0,
\end{equation*}
which can be solved by the root formula of quadratic equation with one unknown
to get the value of $\lambda^2$, and so will the value of $\lambda$.
\item If the normal component of magnetic field is zero, that is, $B_1=0$,
then Eq.~\eqref{eq:RMHDeigvalue} degenerate into
\begin{equation}\label{eq:RMHDeigvalueStwo}
  a_2\lambda^2+a_1\lambda+a_0=0,
\end{equation}
where the coefficients are given by  $$a_2=\rho
h\big[c_s^2+\gamma^2(1-c_s^2)\big]+|\vec b|^2-c_s^2(\vec v\cdot\vec B)^2,$$
$$ a_1=-2\rho h \gamma^2 v_1 (1-c_s^2),$$
$$ a_0=\rho h\big[ -c_s^2+\gamma^2 v_1^2
(1-c_s^2)\big]-|\vec b|^2+c_s^2(\vec v\cdot \vec B)^2.$$
Eq.~\eqref{eq:RMHDeigvalueStwo} can also be
solved by the root formula of quadratic equation with one unknown.
\item In addition to the above two cases, one has to solve the original  quartic equation
 \eqref{eq:RMHDeigvalue}, which is solved by using the analytic formula in our work,
 certainly can also be solved by  some iterative method. 
  \end{itemize}

It has been seen that the expressions of eigenvalues of Jacobian matrix
of \eqref{eqn:RMHDconeqn} is more complicated than the non-relativistic
case,   the calculation of corresponding eigenvectors  is very complicated and may be
obtained by the matrix transformation \cite{RMHDEigenAnton}.

\section{Locally divergence-free DG methods}
\label{sec:rmhdlocaldg}
Because the divergence-free constraint \eqref{eq:divfree} of the magnetic field
is natural for one-dimensional RMHD equations,
no special numerical treatment is required for such constraint
and thus the one-dimensional RMHD equations may be directly solved
by using the non-central or central DG  methods described in
\cite{ZhaoTang2013,ZhaoThesis2014}, they are not repeated here.

This section   gives locally divergence-free DG methods
for the two-dimensional RMHD equations
 \begin{align}
   &\nabla \cdot \textbf B=0,\ \ (x,y)\in \Omega,
   \label{eq:2Ddivfree}
   \end{align}
\begin{equation}\label{eqn:RMHDconeqn2D}
    \displaystyle\frac{\partial \vec{U}}{\partial t} +
    \frac{\partial \vec{F}_1(\vec{U})}{\partial
      x}+\frac{\partial \vec{F}_2(\vec{U})}{\partial
      y}=0,\ \ (x,y)\in \Omega,
 \end{equation}
 on the rectangular mesh $\{ e_{j,k},\forall j,k\in \mathbb Z\}$, where
 \begin{align*}
  \vec{U} =& \Big(D, m_x, m_y, m_z,B_x,B_y,B_z, E\Big)^T,\\
  \vec{F}_1 =& \Big(Dv_x, m_x v_x -B_x b_x/\gamma+ p_{tot},  m_y
  v_x-B_x b_y/\gamma,  \\
  & m_z v_x-B_x b_z/\gamma, 0, B_y v_x-B_x v_y, B_z
  v_x-B_xv_z,m_x \Big)^T,\\
  \vec{F}_2 =& \Big(Dv_y, m_x v_y -B_y b_x/\gamma,  m_y
  v_y-B_y b_y/\gamma+p_{tot},  \\
  & m_z v_y-B_y b_z/\gamma, B_x v_y-B_y v_x, 0, B_z
  v_y-B_yv_z,m_y \Big)^T,\\
 \end{align*}
and $e_{j,k}=(x_{j-\frac{1}{2}},x_{j+\frac{1}{2}})\times (y_{k-\frac{1}{2}},y_{k+\frac{1}{2}})
$. 

The conservative  vector $\vec U$ is divided into two parts as follows
$$\vec R=(D,m_x,m_y,m_z,B_z,E)^T,\quad \vec Q=(B_x,B_y)^T,$$
and similarly the flux vector $\vec F_i$ is also written into $\vec {F}_i^R(\vec U)$ and
$\vec {F}^Q_i(\vec U)$. After that, the RMHD equations \eqref{eqn:RMHDconeqn2D} may be recast into
 \begin{align}\label{eq:RMHD2DR}
&    \displaystyle\frac{\partial \vec{R}}{\partial t} +
    \frac{\partial \vec{F}^R_1(\vec{U})}{\partial x}+\frac{\partial \vec{F}^R_2(\vec{U})}{\partial y}=0,
\\ \label{eq:RMHD2DQ}
 &   \displaystyle\frac{\partial \vec{Q}}{\partial t} +
    \frac{\partial \vec{F}^Q_1(\vec{U})}{\partial x}+\frac{\partial \vec{F}^Q_2(\vec{U})}{\partial y}=0.
 \end{align}
The aim of locally divergence-free DG methods is to find an approximation
$\vec U_h$ of $\vec U$,
where the approximations of $\vec R$ and $\vec Q$ are denoted by
$\vec R_h$ and $\vec Q_h$, respectively, but the spatial discretizations for $\vec R$ and $\vec Q$
are different, see below.

\subsection{Spatial approximation}

The divergence-free constraint of magnetic field is not  needed
in solving the system \eqref{eq:RMHD2DR} for $\vec R$ by using the DG methods,
thus
the DG spproximation of the system \eqref{eq:RMHD2DR} is similar to
that for the RHD equations, see \cite{ZhaoTang2013,ZhaoThesis2014}, and will be omitted here.

Let us discuss the spatial discretization of the governing equations for the dependent variable
$\vec Q$.
The  DG methods  find the approximate solution $\vec Q_h$ such that
for any $t$, $\vec Q_h$ belongs to the following finite dimensional
function space
$$
\vec M_h:=\left\{\vec v(\vec x)\in \vec W^K(e_{j,k}),~\text{}
  ~\vec x=(x,y) \in e_{j,k}, \forall j,k \right\},
$$
where
$$
\vec W^K(e_{j,k}):=\left\{\vec v=\big(v_1(\vec x),~v_2(\vec x)\big)^T \Big|
  ~v_1(\vec x),v_2(\vec x)\in \mathbb P^{K}(e_{j,k}), \frac{\partial
      v_1}{\partial x}+\frac{\partial
      v_2}{\partial y}=0\right\}.
$$
The basis functions of $\vec W^K(e_{j,k})$ may be obtained by
calculating the curl of any vector function whose components are basis functions
of $\mathbb P^{K+1}(e_{j,k})$, thus the dimension of  $\vec W^K(e_{j,k})$
is equal to $D_W=(K+1)(K+4)/2$.
The following lists a set of basis function of $\vec W^K(e_{j,k})$ for the case of $K=3$
\begin{align*}
&\vec \varphi^{(0)}_{j,k}(x,y)=\begin{pmatrix}0\\1\end{pmatrix},\quad  \vec \varphi^{(1)}_{j,k}(x,y)=\begin{pmatrix}1\\0\end{pmatrix},
\quad \vec \varphi^{(2)}_{j,k}(x,y)=\begin{pmatrix} 0\\
  \xi\end{pmatrix},\\
& \vec \varphi^{(3)}_{j,k}(x,y)=\begin{pmatrix}
  \eta\\ 0\end{pmatrix}, \quad \vec \varphi^{(4)}_{j,k}(x,y)=\begin{pmatrix}
  h^x_j\xi\\ -h^y_k \eta\end{pmatrix}, \quad \vec
\varphi^{(5)}_{j,k}=\begin{pmatrix}\eta^2-1/3\\0\end{pmatrix},\\
&\vec \varphi^{(6)}_{j,k}=\begin{pmatrix}0\\ \xi^2-1/3\end{pmatrix},
\quad
\vec \varphi^{(7)}_{j,k}=\begin{pmatrix} h^x_j (\xi^2-1/3)\\ -2 h^y_k \xi\eta\end{pmatrix},\quad \vec \varphi^{(8)}_{j,k}=\begin{pmatrix}
  -2  h^x_j \xi\eta \\ h^y_k(\eta^2-1/3)\end{pmatrix},\\
& \vec \varphi^{(9)}_{j,k}=\begin{pmatrix}
  \eta^3-\frac{3}{5}\eta\\0\end{pmatrix},\quad \vec \varphi^{(10)}_{j,k}=\begin{pmatrix}0\\ \xi^3-\frac{3}{5}\xi\end{pmatrix},\quad \vec \varphi^{(11)}_{j,k}=\begin{pmatrix} h^x_j (\xi^2-1/3)
  \eta\\ -  h^y_k \xi(\eta^2-1/3)\end{pmatrix},\\
& \vec \varphi^{(12)}_{j,k}=\begin{pmatrix}
  h^x_j \xi(\eta^2-1/3) \\ -h^y_k(\eta^3-\eta)/3\end{pmatrix},\quad \vec \varphi^{(13)}_{j,k}=\begin{pmatrix}
  h^x_j(\xi^3-\xi)/3 \\ - h^y_k(\xi^2-1/3)\eta\end{pmatrix},
\end{align*}
where $\xi=2(x-x_j)/h^x_j$, $\eta=2(y-y_k)/h^y_k$, $h^x_j=x_{j+\frac{1}{2}}-x_{j-\frac{1}{2}}$, and $h^y_k=y_{k+\frac{1}{2}}-y_{k-\frac{1}{2}}$.

Multiplying  Eq.~\eqref{eq:RMHD2DQ}  by the test function
$\vec w(\vec x) \in \vec W^K(e_{j,k})$, integrating it over the cell $e_{j,k}$, and using
the divergence theorem give
\begin{align}\nonumber
\frac{d}{dt}\int_{e_{j,k}} \vec Q(\vec x,t)\cdot \vec w(\vec x)~d\vec x
&+\int_{\partial e_{j,k}} ({\vec F}^Q_1 \cdot \vec w(\vec x)n_1+{\vec
  F}^Q_2 \cdot \vec w(\vec x) n_2~)~ds
\\ =&\int_{e_{j,k}} \vec F^Q  \cdot \nabla \vec w(\vec x) ~d\vec x,
\label{eq:RMHDEQ1}
\end{align}
where $\vec F^Q=(\vec F_1^Q,\vec F_2^Q)^T$, and
$(n_1,n_2)$ denote the outer normal vector
of cell boundary $\partial e_{j,k}$.

The approximate solution $\vec Q_h$ may be expressed as follows
\begin{equation}
\vec Q_h(\vec x,t)=\sum_{\ell=0}^{D_W-1} Q^{(\ell)}_{j,k}(t)\vec\varphi^{(\ell)}_{j,k}(\vec x)
=:\vec Q_{j,k}(\vec x,t),~\vec x\in e_{j,k}.
\end{equation}
If the solution $\vec Q$ in \eqref{eq:RMHDEQ1} and flux $\vec F^Q$ on the cell boundary
are replaced with the approximate solution $\vec Q_h$ and numerical flux $\hat{\vec F}^Q$
(e.g. the Lax-Friedrichs flux), respectively,
the test function is taken as the basis, and the integrals
are evaluated by using the Gaussian quadrature,
then the semi-discrete DG   methods may be given as follows
\begin{align}
\nonumber
\sum_{m=0}^{D_W-1} &
\Big(\int_{e_{j,k}} \vec \varphi^{(m)}_{j,k}(\vec x) \cdot \vec \varphi^{(\imath)}_{j,k}(\vec x) ~d\vec x\Big)
\frac{d  Q^{(m)}_{j,k}(t)} {dt}
\\
\nonumber
=& -|\partial e_{j,k}|\sum_{l=1}^{\tilde{q}}\omega_{l}
\Big(\hat{\vec F}_1^Q\big(\vec U(\tilde{\vec x}^G_l, t)\big)n_1+\hat{\vec F}_2^Q\big(\vec U(\tilde{\vec x}^G_l, t)\big)n_2\Big) \cdot
\vec \varphi_{j,k}^{(\imath)}(\tilde{\vec
  x}^G_{l})
 \\
&+|e_{j,k}|\sum_{\ell=1}^q\omega_\ell \vec F^Q \big(\vec
U_{j,k}(\vec x^G_\ell,t)\big) \cdot \nabla \vec \varphi^{(\imath)}_{j,k}(\vec x^G_\ell) ,~\imath=0,1,\hdots,D_W-1,
\label{eq:RMHDLoDGSemi}
\end{align}
where $\vec U_{j,k}=\big(\vec R_{j,k}^T, \vec Q_{j,k}^T\big)^T$.
The system \eqref{eq:RMHDLoDGSemi} may be considered as
a nonlinear system of ordinary differential equations with respect to
the degree of freedom or moments $Q^{(m)}_{j,k}$, and approximated
by using the explicit Runge-Kutta time discretization to give
a fully-discrete locally divergence-free DG methods.

\begin{Remark}
The approximate magnetic field $\vec Q_h$ obtained above
is divergence free in each cell, but it is not continuous in general across the cell boundary.
In view of this fact, such DG methods are called as locally divergence-free \cite{DGMHDLi}.

  \end{Remark}

\subsection{Adaptive WENO limiter}\label{subsection3.2}

The above DG methods may be directly used to solve the  problem whose solution is smooth or
only contains weak discontinuity,
but  the limiting procedure  is needed to
suppress the numerical oscillations in solving the problem containing the strong discontinuity.
In solving the system \eqref{eqn:RMHDconeqn2D}, the limiting procedure used in
the Runge-Kutta DG methods should ensure that the new approximate magnetic field
 does still satisfy the divergence-free constraint.
 This paper uses the WENO limiting procedure, whose implementation consists of two parts:
  identify the ``troubled'' cells and reconstruct the new approximate solution in the ``troubled'' cells.

\begin{description}
  \item[Step I: identify ``troubled'' cells]
It is the same as the way used for the RHD equations in \cite{ZhaoTang2013}
by identifying the cell $e_{j,k}$ in the $x$ and $y$ directions, respectively.
In each direction of the identification process,
the modified TVB minmod function 
    is applied to the characteristic variables in that direction,
    where the characteristic variables are calculated by using the eigenvectors
    corresponding to the system \eqref{eqn:RMHDconeqn2D},
for example,  the characteristic variables in   the $x$ direction is calculated by $\vec L_1 \vec U$,
here $\vec L_1$ denotes the left eigenvector matrix of the Jacobian matrix $ \partial{\vec F_1}/\partial{\vec U}$.
    If the cell $e_{j,k}$ is identified as a ``troubled'' cell in $x$ or $y$ direction
then the cell  $e_{j,k}$  is marked as $e_{j,k}^{tc}$, and go to \textbf{Step II};
otherwise the next cell is continuously checked.

\item[Srep II: reconstruct new approximate solution in ``troubled'' cell]
(i) Calculate the cell averages  of characteristic variables, denoted by $\{\vec W^{(0)}_{\imath,\ell}\}$, use the WENO technique for the characteristic variables $\vec W$
to reconstruct the new approximate values $\vec{W}^G_{m,p}$ of $\vec W$
at the Gaussian points $(x^G_m,y^G_p)$  in the ``troubled'' cell $e_{j,k}^{tc}$.
After that, calculate the  point value $\vec{U}^G_{m,p}$ of the conservative vector,
and divide it into two parts   $\vec R_{m,p}^G$ and  $\vec Q_{m,p}^G$.

(ii) 
Derive the new approximate solutions $\vec R_{j,k}^{WENO}(x,y,t_n)$ and
$\vec Q_{j,k}^{WENO}(x,y,t_n)$ by using the numerical integration and the approximate
values of solution at the Gaussian points as follows
$$
\vec R^{WENO}_{j,k}(x,y,t_n):=\vec R^{(0)}_{j,k} \phi_{j,k}^{(0)}(x,y)+\sum\limits_{\ell=1}^{K(K+3)/2} \vec R^{(\ell),WENO}_{j,k} \phi_{j,k}^{(\ell)}(x,y),\
\quad (x,y)\in e^{tc}_{j,k},
$$
where $\{\phi_{j,k}^{(\ell)}(x,y)\}$
are a set of basis functions of $\mathbb P^K(e_{j,k}^{tc})$,
satisfying the orthogonal properties, and the higher order moments are given by the following formula
\begin{align*}
\vec R_{j,k}^{(\ell),WENO} &:=\frac{1}{a_\ell}
\int_{e^{tc}_{j,k}} \vec R^{WENO}_{j,k}(x,y,t_n) \phi_{j,k}^{(\ell)}(x,y)~dxdy
\\
&\approx \frac{ h^x_j h^y_k}{a_\ell} \sum\limits_{m=1}^q \sum\limits_{p=1}^q
\omega_{m,p} \vec R^G_{m,p}
\phi_{j,k}^{(\ell)}(x^G_m,y^G_p),\quad 1\leq \ell\leq {K(K+3)/2},
\end{align*}
where $a_\ell=\int_{e^{tc}_{j,k}} \big(
\phi_{j,k}^{(\ell)}(x,y)\big)^2~dxdy$, the new approximate magnetic field is
given by
\begin{align*}
\vec
Q_{j,k}^{WENO}(x,y,t_n)=&Q^{(0)}_{j,k}\vec\varphi_{j,k}^{(0)}(x,y)+Q^{(1)}_{j,k}\vec\varphi_{j,k}^{(1)}(x,y)
\\
+&\sum\limits_{\ell=2}^{D_W-1} {  Q^{(\ell),WENO}_{j,k}}
\vec \varphi_{j,k}^{(\ell)}(x,y),\quad (x,y)\in e^{tc}_{j,k},
\end{align*}
where higher-order moments $Q^{(\ell),WENO}_{j,k},~l=2,\cdots,D_W-1$, satisfy the following
linear system
\begin{align*}
&\sum \limits_{\ell=2}^{D_W-1} \big(\int_{e^{tc}_{j,k}}\vec \varphi_{j,k}^{(\ell)}(x,y)
\cdot \vec \varphi_{j,k}^{(l)}(x,y) dxdy \big) Q_{j,k}^{(\ell),WENO}\\
&=  h^x_j h^y_k \sum\limits_{m=1}^q \sum\limits_{p=1}^q
\omega_{m,p} \vec Q^G_{m,p}
\cdot \vec \varphi_{j,k}^{(l)}(x^G_m,y^G_p),\quad 2\leq l\leq {D_W-1}.
\end{align*}
After getting $\vec R_{j,k}^{WENO}(x,y,t_n)$ and $\vec Q_{j,k}^{WENO}(x,y,t_n)$,
they are used to replace the original DG solutions in $e_{j,k}^{tc}$.
Up to now the solution in the
``troubled'' cell $e_{j,k}^{tc}$ is modified. Go to \textbf{Step I}
and check the next cell.
\end{description}

\section{``Exactly'' divergence-free  DG methods}
\label{sec:rmhdcdg}

This section introduces the ``exactly'' divergence-free
central DG methods for two-dimensional RMHD equations \eqref{eqn:RMHDconeqn2D}.
Similar to the locally divergence-free DG methods proposed in Section \ref{sec:rmhdlocaldg},
the conservative vector $\vec U$ in the central DG methods is still
written into two parts
 $$\vec
  R=(D,m_x,m_y,m_z,B_z,E),\quad \vec Q=(B_x,B_y),$$
and apply the different spatial approximations to  the governing equations for $\vec R$ and $\vec Q$.

\subsection{Spatial approximation}
For the sake of convenience, we first give the central DG methods with the explicit
Euler time discretization.
Divide the computational domain $\Omega$ into
two mutually dual meshes, denoted by $\{C_{j,k}\}$ and $\{D_{j+\frac{1}{2},k+\frac{1}{2}}\}$,
respectively, where $C_{j,k}=(x_{j-\frac{1}{2}},x_{j+\frac{1}{2}})\times
(y_{k-\frac{1}{2}},y_{k+\frac{1}{2}})$, $D_{j+\frac{1}{2},k+\frac{1}{2}}=(x_j,x_{j+1})\times
(y_k,y_{k+1})$.
The aim of central DG methods is to find two approximate approximate solutions $\vec U^C_h$
and $\vec U^D_h$, where
corresponding DG approximations of $\vec R$ is denoted by $\vec R^C_h$
and $\vec R^D_h$, while $\vec Q^C_h$ and $\vec Q^D_h$ are used to denote
 corresponding DG approximations of $\vec Q$.

Because   the divergence  free constraint \eqref{eq:2Ddivfree} is not considered in the
equations  \eqref{eq:RMHD2DR} for $\vec R$, Eq. \eqref{eq:RMHD2DR} may be directly
discretized in the central DG framework, see \cite{ZhaoThesis2014},
that is, we find two approximate solutions  $\vec R^C_h$ and $\vec R^D_h$
 such that their each component belongs to the following finite dimensional
 spaces
\begin{align}
  \nonumber
&\mathcal{V}^C:=\left\{v(\vec x)\in L^1(\Omega)\big|~  v(\vec x)\in
  \mathbb  P^{K}(C_{j,k}),\mbox{}~\vec x\in C_{j,k}\subset\Omega,
  \forall j,k \right\},\\
\nonumber
 &\mathcal{V}^D:=\left\{w(\vec x)\in L^1(\Omega)\big|~  w(\vec x)\in
  \mathbb  P^{K}(D_{j+1/2,k+1/2}),\mbox{}~\vec x\in D_{j+\frac{1}{2},k+\frac{1}{2}}\subset\Omega,
  \forall j,k \right\}.\end{align}
  If giving  the DG solutions $\vec U^{C,n}_h$ and $\vec U^{D,n}_h$
  at $t=t_n$, then the DG solution $\vec R^{C,n+1}_h$ at $t_{n+1}$ satisfies
\begin{align*}
&\int_{C_{j,k}} \vec R^{C,n+1}_h  v(\vec x) d\vec x
=\int_{C_{j,k}} \big (\theta \vec R^{D,n}_h+(1-\theta) \vec
R^{C,n}_h \big) v(\vec x)~d\vec x\\
&+ \Delta t_n\big(\int_{C_{j,k}} \vec F^R(\vec U^{D,n}_h) \cdot \nabla v(\vec
  x)~d\vec
x-\int_{\partial C_{j,k}}  \vec F^R(\vec U^{D,n}_h)\cdot \vec n v(\vec x)
 ds\big),~\forall v(\vec x)\in  \mathbb P^{K}(C_{j,k}),
  \end{align*}
  where $\Delta t_n=t_{n+1}-t_n$ denotes the practical time stepsize, $\theta=\Delta
  t_n/\tau_n$, and $\tau_n$ is the maximum timestep allowed by the CFL condition.
 Similarly, the approximate solution $\vec R^{D,n+1}_h$ on the mesh
  $\{D_{j+\frac{1}{2},k+\frac{1}{2}}\}$
  satisfies the following equation
\begin{align*}
&\int_{D_{j+\frac{1}{2},k+\frac{1}{2}}} \vec R^{D,n+1}_h  w(\vec x) d\vec x=
\int_{D_{j+\frac{1}{2},k+\frac{1}{2}}} \big (\theta \vec R^{C,n}_h+(1-\theta) \vec
R^{D,n}_h \big) w(\vec x)~d\vec x\\
&+ \Delta t_n\big(\int_{D_{j+\frac{1}{2},k+\frac{1}{2}}} \vec F^R(\vec U^{C,n}_h) \cdot \nabla w(\vec
  x)~d\vec
x-\int_{\partial D_{j+\frac{1}{2},k+\frac{1}{2}}}  \vec F^R(\vec
U^{C,n}_h)\cdot \vec n  w(\vec x) ds\big),~
  \end{align*}
  for all $w(\vec x)\in  \mathbb P^{K}(D_{j+\frac{1}{2},k+\frac{1}{2}})$.

Due to the need to consider the  divergence free constraint \eqref{eq:2Ddivfree},
the finite dimensional space for the approximation of $\vec Q$
is different from that for $\vec R$, and is denoted as follows
\begin{align*}
 \mathcal{M}^C:&=\{ \vec v(\vec x)~\Big| ~\nabla \cdot \vec v=0, \vec
   v(\vec x)
     \in \mathcal{W}^K(C_{j,k}), \text{}~\vec x\in C_{j,k},\forall
     j,k\},
     \\
  \mathcal{M}^D:&=\{ \vec w(\vec x)~\Big| ~\nabla \cdot \vec w=0, \vec
   w(\vec x)
     \in \mathcal{W}^K(D_{j+\frac{1}{2},k+\frac{1}{2}}), \text{}~\vec x\in D_{j+\frac{1}{2},k+\frac{1}{2}},\forall
     j,k\},
     \end{align*}
where
       $$ \mathcal{W}^K(\mathcal{T}):=\big[\mathbb{P}^K(\mathcal{T})\big]^2 \oplus span\big\{ \nabla
       \times (x^{K+1} y),\nabla \times ( xy^{K+1})\big\},$$
here $\mathcal{T}$ denotes
$C_{j,k}$ or $D_{j+\frac{1}{2},k+\frac{1}{2}}$,
the symbol ``$\oplus$'' denotes the direct sum.
The central DG methods for the $\vec Q$ equations are to find the approximate solution $\vec Q^{C}_h\in \mathcal{M}^C$  and
$\vec Q^{D}_h\in \mathcal{M}^D$ of $\vec Q$.

In the following, we provide how to get
the approximate magnetic field $\vec
Q^{C,n+1}_h=(B^{C,n+1}_{x,h},B^{C,n+1}_{y,h})^T$ on the mesh $\{C_{j,k}\}$
at time $t_{n+1}$.
It will consist of two steps:
 the approximate normal magnetic field is first gotten by solving the
magnetic equation on the cell boundary, and then used to reconstruct
the magnetic field in the cell.
Take the cell $C_{j,k}$ as an example in the following.

\noindent
{\tt Step (i)} Find the   approximate normal magnetic fields $b^x_{j\pm\frac{1}{2},k}(y)$
on the left and right boundaries $x=x_{j\pm \frac12}$
of cell $C_{j,k}$, and $b^y_{j,k\pm\frac{1}{2}}(x)$
on the top and bottom boundaries $y=y_{k\pm \frac12}$.
Because $B_x$ and $B_y$ satisfy
\begin{equation}\label{eq:magcdg1DBX}
\pd{B_x}{t}+\pd{G}{y}=0,
  \end{equation}
  \begin{equation}\label{eq:magcdg1DBY}
\pd{B_y}{t}-\pd{G}{x}=0,
  \end{equation}
  where $G=G(\vec U)=B_xv_y-B_yv_x$, they may be solved by using the one-dimensional central DG methods,
  that is, find $b^x_{j-\frac{1}{2},k}(y)\in
 \mathbb{P}^K\big((y_{k-\frac{1}{2}},y_{k+\frac{1}{2}})\big)$ satisfying
\begin{small}
\begin{align}\nonumber
&\int_{y_{k-\frac{1}{2}}}^{y_{k+\frac{1}{2}}} b^x_{j-\frac{1}{2},k}(y) \mu(y) dy=\int_{y_{k-\frac{1}{2}}}^{y_{k+\frac{1}{2}}}
\Big(\theta B_{x,h}^{D,n}(x_{j-\frac{1}{2}},y)+(1-\theta) B_{x,h}^{C,n}(x_{j-\frac{1}{2}},y)\Big) \mu(y)
dy\\
\nonumber
&+\Delta t_n \Big( \int_{y_{k-\frac{1}{2}}}^{y_{k+\frac{1}{2}}} G\big(\vec
U^{D,n}_h(x_{j-\frac{1}{2}},y)\big)\dd{\mu(y)}{y} dy-G^{D,n}_{j-\frac{1}{2},k+\frac{1}{2}}
\mu(y_{k+\frac{1}{2}})+G^{D,n}_{j-\frac{1}{2},k-\frac{1}{2}}\mu(y_{k-\frac{1}{2}})\Big),
\end{align}
\end{small}
 for any  $\mu(y)\in \mathbb{P}^K\big((y_{k-\frac{1}{2}},y_{k+\frac{1}{2}})\big)$,
 and find $b^y_{j,k-\frac{1}{2}}(x)\in
\mathbb P^K\big((x_{j-\frac{1}{2}},x_{j+\frac{1}{2}})\big)$ satisfying
\begin{small}
\begin{align}\nonumber
&\int_{x_{j-\frac{1}{2}}}^{x_{j+\frac{1}{2}}} b^y_{j,k-\frac{1}{2}}(x) \sigma(x) dx=\int_{x_{j-\frac{1}{2}}}^{x_{j+\frac{1}{2}}}
\Big(\theta B_{y,h}^{D,n}(x,y_{k-\frac{1}{2}})+(1-\theta) B_{y,h}^{C,n}(x,y_{k-\frac{1}{2}})\Big) \sigma(x)
dx\\
\nonumber
  &+\Delta t_n \Big(  \int_{x_{j-\frac{1}{2}}}^{x_{j+\frac{1}{2}}} -G\big(\vec
U^{D,n}_h(x,y_{k-\frac{1}{2}})\big)\dd{\sigma(x)}{x} dx+G^{D,n}_{j+\frac{1}{2},k-\frac{1}{2}}
\sigma(x_{j+\frac{1}{2}})-G^{D,n}_{j-\frac{1}{2},k-\frac{1}{2}}\sigma(x_{j-\frac{1}{2}})\Big),
\end{align}
\end{small}
for any
$\sigma(x)\in \mathbb P^K\big( (x_{j-\frac{1}{2}},x_{j+\frac{1}{2}})\big)$,
where $G^{D,n}_{j-\frac{1}{2},k+\frac{1}{2}}:=G\big(\vec
U^{D,n}_h(x_{j-\frac{1}{2}},y_{k+\frac{1}{2}})\big)$.
The normal magnetic fields $b^x_{j+\frac{1}{2},k}(y)$ and $b^y_{j,k+\frac{1}{2}}(x)$
may be solved by similar approach.

\noindent
{\tt Step (ii)}
Reconstruct the magnetic field $\vec Q^{C,n+1}_{h}|_{C_{j,k}}=\big(
B^{C,n+1}_{x,j,k},B^{C,n+1}_{y,j,k}\big)^T$ in the cell $C_{j,k}$, where $B^{C,n+1}_{x,j,k},B^{C,n+1}_{y,j,k}\in
\mathcal{W}^K(C_{j,k})$ and satisfy
\begin{align}\label{eq:magrecconsist}
\begin{aligned}
  B^{C,n+1}_{x,j,k}(x_{j+ \frac{1}{2}},y)=&b^x_{j+
    {\frac{1}{2},k}}(y), \quad y\in
  (y_{k-\frac{1}{2}},~y_{k+\frac{1}{2}}),\\
   B^{C,n+1}_{x,j,k}(x_{j- \frac{1}{2}},y)=&b^x_{j-
    {\frac{1}{2},k}}(y), \quad y\in
  (y_{k-\frac{1}{2}},~y_{k+\frac{1}{2}}),
  \\
  B^{C,n+1}_{y,j,k}(x,y_{k+\frac{1}{2}})=&b^y_{j,k+
    {\frac{1}{2}}}(x), \quad x\in
  (x_{j-\frac{1}{2}},~x_{j+\frac{1}{2}}),\\
   B^{C,n+1}_{y,j,k}(x,y_{k-\frac{1}{2}})=&b^y_{j,k-
    {\frac{1}{2}}}(x), \quad x\in
  (x_{j-\frac{1}{2}},~x_{j+\frac{1}{2}}),
\end{aligned}\end{align}
\begin{equation}
  \label{eq:magrecdiv}
\pd{B^{C,n+1}_{x,j,k}}{x}+\pd{B^{C,n+1}_{y,j,k}}{y}=0,\quad (x,y)\in C_{j,k}.
\end{equation}

Before giving the detailed expressions of
$B^{C,n+1}_{x,j,k}$ and $B^{C,n+1}_{y,j,k}$,
a necessary condition is first discussed for the solution of the system of equations
\eqref{eq:magrecconsist} and \eqref{eq:magrecdiv}.
If such the system exists the solution, then integrating \eqref{eq:magrecdiv} over the cell $C_{j,k}$ and
using the divergence theorem and Eq.~\eqref{eq:magrecconsist} gives
   \begin{equation}\label{eq:magrecnece}
   \int_{x_{j-\frac{1}{2}}}^{x_{j+\frac{1}{2}}}
   \big(b^y_{j,k-\frac{1}{2}}(x)-
   b^y_{j,k+\frac{1}{2}}(x)\big)
   dx+\int_{y_{k-\frac{1}{2}}}^{y_{k+\frac{1}{2}}}
   \big( b^x_{j-\frac{1}{2},k}(y)-
   b^x_{j+\frac{1}{2},k}(y)\big)  dy=0,
      \end{equation}
      which is a necessary condition  for the solution of the system of equations
      \eqref{eq:magrecconsist} and \eqref{eq:magrecdiv}.
It is proved in \cite{CentralMHD} that the normal magnetic field obtained by using the central DG methods
satisfies \eqref{eq:magrecnece}.

The following provides the expressions of $B^{C,n+1}_{x,j,k}$ and $B^{C,n+1}_{y,j,k}$
for $K=1,2,3$.

When $K=1$, the normal magnetic fields on the cell boundary are 
  $$ b^x_{j\pm \frac{1}{2},k}(y)=b^{x,(0)}_{j\pm \frac{1}{2},k}+b^{x,(1)}_{j\pm\frac{1}{2},k} \eta,$$
 $$b^y_{j,k\pm\frac{1}{2}}(x)=b^{y,(0)}_{j,k\pm \frac{1}{2}}+b^{y,(1)}_{j,k\pm\frac{1}{2}} \xi,$$
where $\xi=2(x-x_j)/h^x_j,~\eta=2(y-y_k)/h^y_k$,
and the reconstructed  magnetic fields in the cell are
\begin{align*}
B^{C,n+1}_{x,j,k}(x,y)&=a_0+a_1 \xi+a_2 \eta + a_3 (\xi^2-1/3)+a_4 \xi \eta,\\
B^{C,n+1}_{y,j,k}(x,y)&=b_0+b_1 \xi+b_2 \eta +b_3 \xi \eta+ b_4 (\eta^2-1/3),
\end{align*}
where the coefficients are given by
\begin{align}\label{eq:cdgp1mag}
\begin{aligned}
a_0=&\frac{1}{2}(b^{x,(0)}_{j+\frac{1}{2},k}+b^{x,(0)}_{j-\frac{1}{2},k})+\frac{h^x_j}{6h^y_k}(b^{y,(1)}_{j,k+\frac{1}{2}}  -b^{y,(1)}_{j,k-\frac{1}{2}}),
\\
b_0=&\frac{1}{2}(b^{y,(0)}_{j,k+\frac{1}{2}}+b^{y,(0)}_{j,k-\frac{1}{2}})+\frac{h^y_k}{6
  h^x_j}(b^{x,(1)}_{j+\frac{1}{2},k} -b^{x,(1)}_{j-\frac{1}{2},k}
),\\
 a_1=&\frac{1}{2}( b^{x,(0)}_{j+\frac{1}{2},k}-b^{x,(0)}_{j-\frac{1}{2},k}),\quad b_2=\frac{1}{2}(b^{y,(0)}_{j,k+\frac{1}{2}} - b^{y,(0)}_{j,k-\frac{1}{2}}),\\
 a_2=&\frac{1}{2}(b^{x,(1)}_{j+\frac{1}{2},k}+b^{x,(1)}_{j-\frac{1}{2},k}
 ),\quad b_1=\frac{1}{2}(b^{y,(1)}_{j,k+\frac{1}{2}}+b^{y,(1)}_{j,k-\frac{1}{2}} ),\\
 a_3=&-\frac{h^x_j}{4h^y_k}(b^{y,(1)}_{j,k+\frac{1}{2}}
 -b^{y,(1)}_{j,k-\frac{1}{2}}),\quad  b_4=-\frac{h^y_k}{4 h^x_j}(b^{x,(1)}_{j+\frac{1}{2},k} -b^{x,(1)}_{j-\frac{1}{2},k} ),\\
a_4=&\frac{1}{2}(b^{x,(1)}_{j+\frac{1}{2},k} -b^{x,(1)}_{j-\frac{1}{2},k}
), \quad b_3=\frac{1}{2}(b^{y,(1)}_{j,k+\frac{1}{2}}  -b^{y,(1)}_{j,k-\frac{1}{2}}).
\end{aligned}\end{align}
When $K=2$, the normal magnetic fields on the cell boundary are
  \begin{align*}
b^x_{j\pm\frac{1}{2},k}(y)&=b^{x,(0)}_{j\pm\frac{1}{2},k}+b^{x,(1)}_{j\pm\frac{1}{2},k} \eta+b^{x,(2)}_{j\pm\frac{1}{2},k} (\eta^2-\frac{1}{3}),\\
  b^y_{j,k\pm\frac{1}{2}}(x)&=b^{y,(0)}_{j,k\pm\frac{1}{2}}+b^{y,(1)}_{j,k\pm\frac{1}{2}}\xi+
  b^{y,(2)}_{j,k\pm\frac{1}{2}}(\xi^2-\frac{1}{3}),
  \end{align*}
and the reconstructed magnetic fields in the cell are given by
\begin{align*}
B^{C,n+1}_{x,j,k}(x,y) =& a_0+a_1 \xi+a_2 \eta + a_3 (\xi^2-1/3)+a_4 \xi
\eta\\
& +a_5
(\eta^2-1/3)
+a_6(\xi^3-\frac{3}{5}\xi)+a_7 \xi(\eta^2-\frac{1}{3}) ,\\
 B^{C,n+1}_{y,j,k}(x,y)=& b_0+b_1 \xi+b_2 \eta +b_3 (\xi^2-1/3)+
 b_4 \xi \eta\\
 &+ b_5 (\eta^2-1/3)+b_6
(\xi^2-\frac{1}{3})\eta+b_7 (\eta^3-\frac{3}{5}\eta),
\end{align*}
where the coefficients are
\begin{align*}
  a_0=&\frac{1}{2}(b^{x,(0)}_{j+\frac{1}{2},k}+b^{x,(0)}_{j-\frac{1}{2},k})+\frac{h^x_j}{6 h^y_k}(b^{y,(1)}_{j,k+\frac{1}{2}}-b^{y,(1)}_{j,k-\frac{1}{2}}),\\
b_0=&\frac{1}{2}(b^{y,(0)}_{j,k+\frac{1}{2}}+b^{y,(0)}_{j,k-\frac{1}{2}})+\frac{h^y_k}{6
  h^x_j}(b^{x,(1)}_{j+\frac{1}{2},k}-b^{x,(1)}_{j-\frac{1}{2},k}),\\
a_1=&\frac{1}{2}(b^{x,(0)}_{j+\frac{1}{2},k}-b^{x,(0)}_{j-\frac{1}{2},k})+\frac{h^x_j}{15h^y_k}(b^{y,(2)}_{j,k+\frac{1}{2}}-b^{y,(2)}_{j,k-\frac{1}{2}}),\\
b_2=&\frac{1}{2}(b^{y,(0)}_{j,k+\frac{1}{2}}-b^{y,(0)}_{j,k-\frac{1}{2}})+\frac{h^y_k}{15h^x_j}(b^{x,(2)}_{j+\frac{1}{2},k}-b^{x,(2)}_{j-\frac{1}{2},k}),\\
a_2=&\frac{1}{2}(b^{x,(1)}_{j+\frac{1}{2},k}+b^{x,(1)}_{j-\frac{1}{2},k}),
\quad
a_4=\frac{1}{2}(b^{x,(1)}_{j+\frac{1}{2},k}-b^{x,(1)}_{j-\frac{1}{2},k}),\\
b_1=&\frac{1}{2}(b^{y,(1)}_{j,k+\frac{1}{2}}+b^{y,(1)}_{j,k-\frac{1}{2}}),
\quad
b_4=\frac{1}{2}(b^{y,(1)}_{j,k+\frac{1}{2}}-b^{y,(1)}_{j,k-\frac{1}{2}}),
\\
a_5=&\frac{1}{2}(b^{x,(2)}_{j+\frac{1}{2},k}+b^{x,(2)}_{j-\frac{1}{2},k}),
\quad
a_7=\frac{1}{2}(b^{x,(2)}_{j+\frac{1}{2},k}-b^{x,(2)}_{j-\frac{1}{2},k}),\\
b_3=&\frac{1}{2}(b^{y,(2)}_{j,k+\frac{1}{2}}+b^{y,(2)}_{j,k-\frac{1}{2}}), \quad b_6=\frac{1}{2}(b^{y,(2)}_{j,k+\frac{1}{2}}-b^{y,(2)}_{j,k-\frac{1}{2}}).\\
a_3=&-\frac{h^x_j}{4h^y_k}(b^{y,(1)}_{j,k+\frac{1}{2}}-b^{y,(1)}_{j,k-\frac{1}{2}}),\quad
a_6=-\frac{h^x_j}{6 h^y_k}(b^{y,(2)}_{j,k+\frac{1}{2}}-b^{y,(2)}_{j,k-\frac{1}{2}}),\\
b_5=&-\frac{h^y_k}{4h^x_j}(b^{x,(1)}_{j+\frac{1}{2},k}-b^{x,(1)}_{j-\frac{1}{2},k}),\quad
b_7=-\frac{h^y_k}{6 h^x_j}(b^{x,(2)}_{j+\frac{1}{2},k}-b^{x,(2)}_{j-\frac{1}{2},k}).
\end{align*}
 When $K=3$, because the normal magnetic fields
 on the cell boundary
 are given by %
 \begin{align*}
   b^x_{j\pm\frac{1}{2},k}(y)&=b^{x,(0)}_{j\pm\frac{1}{2},k}+b^{x,(1)}_{j\pm\frac{1}{2},k} \eta+b^{x,(2)}_{j\pm\frac{1}{2},k} (\eta^2-\frac{1}{3})+b^{x,(3)}_{j\pm\frac{1}{2},k} (\eta^3-\frac{3}{5}\eta),\\
b^y_{j,k\pm\frac{1}{2}}(x)&=b^{y,(0)}_{j,k\pm\frac{1}{2}}+b^{y,(1)}_{j,k\pm\frac{1}{2}}
\xi+b^{y,(2)}_{j,k\pm\frac{1}{2}}(\xi^2-\frac{1}{3})+b^{y,(3)}_{j,k\pm\frac{1}{2}}(\xi^3-\frac{3}{5}\xi),
\end{align*}
the reconstructed magnetic fields in the cell  are 
\begin{align*}
B_{x,j,k}^{C,n+1}(x,y)=&a_0+a_1\xi+a_2\eta+a_3(\xi^2-\frac{1}{3})+a_4\xi\eta+a_5(\eta^2-\frac{1}{3})\\
&+a_6(\xi^3-\frac{3}{5}\xi)
+a_7(\xi^2-\frac{1}{3})\eta+a_8\xi(\eta^2-\frac{1}{3})\\
&+a_9(\eta^3-\frac{3}{5}\eta)
+a_{10}(\xi^4-\frac{6}{7}\xi^2+\frac{3}{35})+a_{11}\xi(\eta^3-\frac{3}{5}\eta),\\
B_{y,j,k}^{C,n+1}(x,y)=&b_0+b_1\xi+b_2\eta+b_3(\xi^2-\frac{1}{3})+b_4\xi\eta+b_5(\eta^2-\frac{1}{3})\\
&+b_6(\xi^3-\frac{3}{5}\xi)+b_7(\xi^2-\frac{1}{3})
\eta +b_8\xi(\eta^2-\frac{1}{3})\\
&+b_9(\eta^3-\frac{3}{5}\eta)+b_{10}\eta(\xi^3-\frac{3}{5}\xi)+b_{11}(\eta^4-\frac{6}{7}\eta^2+\frac{3}{35}).
\end{align*}
Some coefficients may be determined as follows
 \begin{align*}
    a_0 &=\frac{1}{2}(b^{x,(0)}_{j+\frac{1}{2},k}+b^{x,(0)}_{j-\frac{1}{2},k})+\frac{h^x_j}{6h^y_k}(b^{y,(1)}_{j,k+\frac{1}{2}}-b^{y,(1)}_{j,k-\frac{1}{2}}),\\
       b_0&
       =\frac{1}{2}(b^{y,(0)}_{j,k+\frac{1}{2}}+b^{y,(0)}_{j,k-\frac{1}{2}})+\frac{h^y_k}{6h^x_j}(b^{x,(1)}_{j+\frac{1}{2},k}-b^{x,(1)}_{j-\frac{1}{2},k}),\\
        a_1&=\frac{1}{2}(b^{x,(0)}_{j+\frac{1}{2},k}-b^{x,(0)}_{j-\frac{1}{2},k})+\frac{h^x_j}{15h^y_k}(b^{y,(2)}_{j,k+\frac{1}{2}}-b^{y,(2)}_{j,k-\frac{1}{2}}),\\
         b_2&=\frac{1}{2}(b^{y,(0)}_{j,k+\frac{1}{2}}-b^{y,(0)}_{j,k-\frac{1}{2}})+\frac{h^y_k}{15
         h^x_j}(b^{x,(2)}_{j+\frac{1}{2},k}-b^{x,(2)}_{j-\frac{1}{2},k}),\\
        a_3&=\frac{h^x_j}{h^y_k}\big(\frac{3}{70}(b^{y,(3)}_{j,k+\frac{1}{2}}-b^{y,(3)}_{j,k-\frac{1}{2}})-\frac{1}{4}(b^{y,(1)}_{j,k+\frac{1}{2}}-b^{y,(1)}_{j,k-\frac{1}{2}})\big),\\
       b_5&=\frac{h^y_k}{h^x_j}\big(\frac{3}{70}(b^{x,(3)}_{j+\frac{1}{2},k}-b^{x,(3)}_{j-\frac{1}{2},k})-\frac{1}{4}(b^{x,(1)}_{j+\frac{1}{2},k}-b^{x,(1)}_{j-\frac{1}{2},k})\big),\\
       a_4&=\frac{1}{2}(b^{x,(1)}_{j+\frac{1}{2},k}-b^{x,(1)}_{j-\frac{1}{2},k}),\quad
       b_4=\frac{1}{2}(b^{y,(1)}_{j,k+\frac{1}{2}}-b^{y,(1)}_{j,k-\frac{1}{2}}),\\
       a_5&=\frac{1}{2}(b^{x,(2)}_{j+\frac{1}{2},k}+b^{x,(2)}_{j-\frac{1}{2},k}),\quad
a_{8}=\frac{1}{2}(b^{x,(2)}_{j+\frac{1}{2},k}-b^{x,(2)}_{j-\frac{1}{2},k}),\\
b_3&=\frac{1}{2}(b^{y,(2)}_{j,k+\frac{1}{2}}+b^{y,(2)}_{j,k-\frac{1}{2}}),\quad b_{7}=\frac{1}{2}(b^{y,(2)}_{j,k+\frac{1}{2}}-b^{y,(2)}_{j,k-\frac{1}{2}}),\\
      a_6&=-\frac{h^x_j}{6h^y_k}(b^{y,(2)}_{j,k+\frac{1}{2}}-b^{y,(2)}_{j,k-\frac{1}{2}}),\quad a_{10}=-\frac{h^x_j}{ 8 h^y_k}(b^{y,(3)}_{j,k+\frac{1}{2}}-b^{y,(3)}_{j,k-\frac{1}{2}}),
      \\
      b_9&
      =-\frac{h^y_k}{6h^x_j}(b^{x,(2)}_{j+\frac{1}{2},k}-b^{x,(2)}_{j-\frac{1}{2},k}),\quad
      b_{11}=-\frac{h^y_k}{8h^x_j}(b^{x,(3)}_{j+\frac{1}{2},k}-b^{x,(3)}_{j-\frac{1}{2},k}),\\
       a_9&=\frac{1}{2}(b^{x,(3)}_{j+\frac{1}{2},k}+b^{x,(3)}_{j-\frac{1}{2},k}),\quad
a_{11}=\frac{1}{2}(b^{x,(3)}_{j+\frac{1}{2},k}-b^{x,(3)}_{j-\frac{1}{2},k}),\\
b_6&=\frac{1}{2}(b^{y,(3)}_{j,k+\frac{1}{2}}+b^{y,(3)}_{j,k-\frac{1}{2}}),\quad b_{10}=\frac{1}{2}(b^{y,(3)}_{j,k+\frac{1}{2}}-b^{y,(3)}_{j,k-\frac{1}{2}}),
      \end{align*}
while other four coefficients satisfy the following relations
\begin{align}
\begin{aligned}
   a_2+\frac{2}{3}a_7=\frac{1}{2}(b^{x,(1)}_{j+\frac{1}{2},k}+b^{x,(1)}_{j-\frac{1}{2},k}),\
  b_1+\frac{2}{3}b_8=\frac{1}{2}(b^{y,(1)}_{j,k+\frac{1}{2}}+b^{y,(1)}_{j,k-\frac{1}{2}}),\
  h_y^k a_7+h^x_j b_8=0.\end{aligned}
\label{eq:coeundefine}
\end{align}

The coefficients $a_2,~a_7,~b_1$, and $b_8$ cannot be uniquely determined by
\eqref{eq:magrecconsist} and \eqref{eq:magrecdiv}. In order to uniquely determine them,
an additional condition is needed.
In our computations, the coefficient $a_7$ at $t_{n+1}$  is obtained by solving
Eq. \eqref{eq:magcdg1DBX} with the central DG methods, i.e.
\begin{align*}
&\int_{C_{j,k}} B^{C,n+1}_{x,h}  v(\vec x) d\vec x
=\int_{C_{j,k}} \big (\theta  B^{D,n}_{x,h}+(1-\theta) B^{C,n}_{x,h} \big) v(\vec x)~d\vec x\\
&+ \Delta t_n\big(\int_{C_{j,k}}  G(\vec U^{D,n}_h)\dd{ v(\vec  x)}{y}~d\vec
x-\int_{\partial C_{j,k}}  G(\vec U^{D,n}_h) n_2 v(\vec x)
 ds\big),
  \end{align*}
 where $v(\vec x)$ is taken as  the basis function $(\xi^2-1/3)\eta$ corresponding to
$a_7$. After $a_7$ is gotten,
the coefficients $a_2,~b_8$, and $b_1$ may be calculated according to
\eqref{eq:coeundefine}.
Combing  $\vec R^{C,n+1}_h$ with
the approximate magnetic fields
 $\vec  Q^{C,n+1}_{h}$ on the mesh $\{C_{j,k}\}$ at $t_{n+1}$ derived by the above approach
 gives the fully approximate conservative vector $\vec
  U^{C,n+1}_h$. The approximate solution vector $\vec U^{D,n+1}_h$ may be similarly obtained,
but  is no longer repeated here.

 The time discretization in the above central DG methods is only first-order accurate,
 the higher-order accurate explicit  Runge-Kutta methods  may be used to replace
 the forward Euler   time discretization and improve the accuracy of central
  DG methods in time.

\begin{Remark}
 Because the normal magnetic fields obtained by the above central DG methods
 are continuous across the cell boundary,
the above central DG methods are  ``exactly'' divergence-free  \cite{CentralMHD}.
  \end{Remark}

  \begin{Remark}
The above methods need to calculate the flux at the corner of the rectangular cell
when the approximate normal magnetic field is solved on the cell boundary.
Naturally, the solutions on the dual mesh may be used.
    \end{Remark}

\subsection{Adaptive WENO limiter}

The limiting procedure is also needed for the ``exactly'' divergence-free central DG methods
when they are used to solve the RMHD problems with strong discontinuity.
The WENO limiting procedure is still independently applied to $\vec R_h$
and $\vec Q_h$. Especially, the limiting procedure for the approximate solutions
$\vec R^C_h$ and $\vec R^D_h$ is the same as that in the locally divergence-free DG
methods, see Section \ref{subsection3.2},
only but one needs to respectively identify the ``troubled'' cells
in two mutually dual meshes and reconstruct new WENO approximate solutions
replace $\vec R^C_h$ and $\vec R^D_h$ defined in  the ``troubled'' cells,
It will not be repeated here.

The above WENO limiting procedure cannot be directly applied to $\vec Q_h$,
otherwise the normal component of  limited magnetic field may be discontinuous across
the cell boundary.
In view of that $\vec Q_h$ is derived by the reconstruction based on $b^x(y)$ and $b^y(x)$,
a natural consideration is that
the WENO limiting procedure is first applied to $b^x(y)$ and $b^y(x)$,
then the limited magnetic field $b^x$ and $b^y$ are used to reconstruct
the new magnetic field within the cell.

Take $b^x_{j-\frac{1}{2},k}(y)$ as an example
to introduce the limiting procedure for the magnetic field.
In fact it is the same as that used for the one-dimensional
DG methods. For the sake of convenience,
use $\mathcal{I}_{j-\frac{1}{2},k}$ to
denote the boundary of cell $C_{j,k}$:
$x=x_{j-\frac{1}{2}},y_{k-\frac{1}{2}}\le y\le y_{k+\frac{1}{2}}$.

\begin{description}
   \item[Step I]:
Use the modified TVB minmod function to
check  whether the normal magnetic field $b^x_{j-\frac{1}{2},k}(y)$ on the
cell boundary $\mathcal{I}_{j-\frac{1}{2},k}$ is needed to be limited.
Calculate
\begin{align*}
\tilde{b}^x_k:=b^x_{j-\frac{1}{2},k}(y_{k+\frac{1}{2}})-b^{x,(0)}_{j-\frac{1}{2},k},\quad
\tilde{\tilde{b}}^x_k:=-b^x_{j-\frac{1}{2},k}(y_{k-\frac{1}{2}})+b^{x,(0)}_{j-\frac{1}{2},k},\end{align*}
then apply  the modified TVB minmod function to
$\tilde{b}^x_k$ and $\tilde{\tilde{b}}^x_k$ 
\begin{align*}
\tilde{b}^{x,mod}_k:=\tilde
m(\tilde{b}^x_k,\Delta_{+}{b}^{x,(0)}_k,\Delta_{-}{b}^{x,(0)}_k),\quad \tilde{\tilde{b}}^{x,mod}_k:=\tilde
m(\tilde{\tilde{b}}^x_k,\Delta_{+}{b}^{x,(0)}_k,\Delta_{-}{b}^{x,(0)}_k),
\end{align*}
where
\begin{align*}
\Delta_{+}{b}^{x,(0)}_k:=b^{x,(0)}_{j-\frac{1}{2},k+1}-b^{x,(0)}_{j-\frac{1}{2},k},\quad
\Delta_{-}{b}^{x,(0)}_k:=b^{x,(0)}_{j-\frac{1}{2},k}-b^{x,(0)}_{j-\frac{1}{2},k-1}.
\end{align*}
If $\tilde{b}^{x,mod}_k$ is different from $\tilde{b}^x_k$, or
$\tilde{\tilde{b}}^{x,mod}_k$ is different from $\tilde{\tilde{b}}^x_k$,
then mark $\mathcal{I}_{j-\frac{1}{2},k}$  as ``troubled'' cell boundary,
and go to \textbf{Step II}; otherwise check the approximate normal magnetic field
on the next cell boundary.

\item[Step II]:
Use the WENO technique to reconstruct the new normal
magnetic field on ``troubled'' cell boundary  $\mathcal{I}_{j-\frac{1}{2},k}$.
Using the cell averages $\{b^{x,(0)}\}$ of $b^x(y)$ on the neighboring cell boundary
of $\mathcal{I}_{j-\frac{1}{2},k}$ in the $y$ direction and
$2K+1$ order WENO reconstruction to get $b^{x,G}_m$,
the new approximation of $b^x$ at the Gaussian points $y^G_m$ within the interval
$(y_{k-\frac{1}{2}},y_{k+\frac{1}{2}})$,
  where $m=1,\cdots K+1$, then use numerical integration to give a new approximation of $b^x(y)$,
  for example,  the new approximate solution for  $K=2$ is given by
  $$b^{x,WENO}_{j-\frac{1}{2},k}(y)=b^{x,(0)}_{j-\frac{1}{2},k}+b^{x,WENO,(1)}_{j-\frac{1}{2},k} \phi_k^{(1)}(y)
  +b^{x,WENO,(2)}_{j-\frac{1}{2},k}\phi_k^{(2)}(y),$$
  where $\phi_k^{(1)}(y)=\eta,\phi_k^{(2)}(y)=\eta^2-1/3$, and
  $\eta=2(y-y_k)/{h^y_k}$,
  higher order moments may be determined by
$$
b^{x,WENO,(i)}_{j-\frac{1}{2},k}\int_{y_{k-\frac{1}{2}}}^{y_{k+\frac{1}{2}}}
\phi_k^{(i)}(y)\phi_k^{(i)}(y)dy   =h^y_k \sum\limits_{m=1}^{K+1} b^{x,G}_m
w_m^G \phi_k^{(i)}(y^G_m) ,\quad i=1,\cdots,K,
$$
where $w_m^G$ is weight corresponding to the point $y^G_m$.
Use the new normal magnetic field $b^{x,WENO}_{j-\frac{1}{2},k}(y)$
 to replace the old  normal magnetic field on $\mathcal{I}_{j-\frac{1}{2},k}$,
and   go to \textbf{Step I}.

\end{description}

\begin{Remark}
It is not difficult to know that the above WENO limiting procedure
does not change the cell average value of magnetic field over
the cell boundary so that the new  normal magnetic field
 satisfies the necessary condition \eqref{eq:magrecnece}.
  \end{Remark}

\section{Numerical Results}
\label{sec:rmhdnum}
This section uses our $P^K$-based locally and ``exactly'' divergence-free Runge-Kutta DG methods
with WENO limiter to solve several initial value problems or
initial-boundary-value problems of one- and two-dimensional RMHD equations in order
to demonstrate the accuracy and effectiveness of our Runge-Kutta DG methods.
Because two solutions of ``exactly'' divergence-free central Runge-Kutta DG methods on mutually dual meshes
are almost identical each other, only one of the central DG solutions
will be shown in the following.

\subsection{1D examples}

This section will solve a smooth problem
and three Riemann problems by the proposed DG methods.
Unless otherwise stated, the third-order accurate explicit TVD Runge-Kutta method
is used for the time discretization, see \cite{ZhaoThesis2014,ZhaoTang2013},
the CFL numbers of $P^1$-, $P^2$-, $P^3$-based non-central \DG{} are taken as
  $0.3,~0.2$, and $0.1$, respectively,  while
  the CFL numbers of $P^1$-, $P^2$-, $P^3$-based \CDG{} are chosen as
  $0.4,~0.3$, and $0.2$, respectively, and  $\theta=1$.
 The determination of time stepsize  may be found in \cite{ZhaoThesis2014},  
and the parameter $M$ in the modified TVB minmod function  is taken as $500$.

\begin{Example}[Smooth problem]\label{exRMHDSmooth1D}\rm
 This problem describes the periodic propagation of a sine wave within the domain
  $\Omega=[0,1]$ and
 is used to test the accuracy of non-central and central DG methods.
   The detailed initial data are
  \begin{align*}
& \rho(x,0)=1,\ \ v_x(x,0)=0,\ \ v_y(x,0)=0.1\sin(2\pi x),\\
& v_z(x,0) =0.1\cos(2\pi x),\ \ B_x(x,0)=1,\ \  B_y(x,0)=\kappa v_y(x,0),\\
& B_z(x,0) =\kappa v_z(x,0),\ \ p(x,0)=0.1,
\end{align*}
where $\kappa=\sqrt{1+\rho h\gamma^2}$, and corresponding exact solutions are
\begin{align*}
 &\rho(x,t)=1,\ \ v_x(x,t)=0,\ \
 v_y(x,t)=0.1\sin\big(2\pi(x+t/\kappa)\big),\\
& v_z(x,t)=0.1\cos\big(2\pi (x+t/\kappa)\big),\ \ B_x(x,t)=1,\ \ B_y(x,t)=\kappa v_y(x,t),\\
& B_z(x,t)=\kappa v_z(x,t),\ \ p(x,t)=0.1.
\end{align*}

In our computations, the adiabatic index $\Gamma=5/3$,
the computational domain $\Omega$ is divided into $N$ uniform cells and the periodic
conditions are specified, and the fourth-order Runge-Kutta method mentioned in \cite{ZhaoTang2013} is used
in order to ensure the accuracy in time. Table~\ref{tab:RMHDsmooth1D} shows
the $l^1$ errors in $B_y$ and orders at $t=1$ obtained by using the non-central
and central DG methods without or with WENO limiter in global.
It is seen that both non-central
and central DG methods may get the theoretical orders,
and the WENO limiter may not destroy the accuracy of DG methods.
  \end{Example}

    \begin{table}[!htbp]
  \centering
  \caption
  { $l^1$ errors in $B_y$ and orders at $t=1$ obtained by
the non-central and central DG methods.   The fourth-order accurate Runge-Kutta matheod and $N \times 2N$ uniform cells are used.}
\begin{tabular}{|c|c|c|c|c|c|c|c|c|c|}
  \hline
 \multirow{3}{20pt}{}
 &\multirow{2}{2pt}{}
 &\multicolumn{4}{|c|}{without limiter}&\multicolumn{4}{|c|}{with limiter in global}\\
 \cline{3-10}
 & & \multicolumn{2}{|c|}{non-central DG}&\multicolumn{2}{|c|}{central DG}&
 \multicolumn{2}{|c|}{non-central DG}&\multicolumn{2}{|c|}{central DG}\\
 \cline{2-10}
 & $N$ &  $l^1$ error& order &  $l^1$ error & order &  $l^1$ error& order &  $l^1$ error & order \\
\hline
\multirow{6}{20pt}{$P^{1}$}
&10& 1.15e-03& --& 9.21e-04& --& 9.94e-03& --& 2.42e-03& --\\
\cline{2-10}
&20&2.99e-04& 1.94&2.40e-04& 1.94&2.28e-03& 2.12&4.82e-04& 2.33\\
\cline{2-10}
&40&7.56e-05& 1.98&6.05e-05& 1.99&5.52e-04& 2.05&1.12e-04& 2.11\\
\cline{2-10}
&80&1.89e-05& 2.00&1.52e-05& 2.00&1.37e-04& 2.02&2.74e-05& 2.03\\
\cline{2-10}
&160&4.74e-06& 2.00&3.79e-06& 2.00&3.40e-05& 2.01&6.81e-06& 2.01\\
\cline{2-10}
&320&1.19e-06& 2.00&9.48e-07& 2.00&8.49e-06& 2.00&1.70e-06& 2.00\\
\hline
\multirow{6}{20pt}{$P^{2}$}
&10& 5.73e-05& --& 2.93e-05& --& 2.61e-03& --& 1.20e-03& --\\
\cline{2-10}
&20&7.21e-06& 2.99&3.83e-06& 2.94&4.41e-05& 5.88&4.68e-05& 4.68\\
\cline{2-10}
&40&9.09e-07& 2.99&4.83e-07& 2.98&5.47e-06& 3.01&1.99e-06& 4.55\\
\cline{2-10}
&80&1.14e-07& 2.99&6.06e-08& 3.00&1.01e-06& 2.43&1.31e-07& 3.92\\
\cline{2-10}
&160&1.43e-08& 2.99&7.59e-09& 3.00&1.37e-07& 2.89&1.29e-08& 3.35\\
\cline{2-10}
&320&1.80e-09& 3.00&9.49e-10& 3.00&1.75e-08& 2.97&1.50e-09& 3.10\\
\hline
\multirow{6}{20pt}{$P^{3}$}
&10& 2.40e-06& --& 1.28e-06& --& 5.91e-05& --& 3.43e-05& --\\
\cline{2-10}
&20&1.53e-07& 3.97&7.96e-08& 4.01&7.89e-07& 6.23&4.77e-07& 6.17\\
\cline{2-10}
&40&9.53e-09& 4.01&5.01e-09& 3.99&8.16e-08& 3.28&1.74e-08& 4.78\\
\cline{2-10}
&80&5.96e-10& 4.00&3.13e-10& 4.00&5.77e-09& 3.82&9.44e-10& 4.20\\
\cline{2-10}
&160&3.73e-11& 4.00&1.96e-11& 4.00&3.71e-10& 3.96&5.69e-11& 4.05\\
\cline{2-10}
&320&2.33e-12& 4.00&1.22e-12& 4.00&2.33e-11& 3.99&3.51e-12& 4.02\\
\hline
  \end{tabular}
  \label{tab:RMHDsmooth1D}
  \end{table}

What follows is that the non-central  and central DG methods are used to solve three
Riemann problems of 1D RMHD equations, whose exact solution are obtained by using
the approach provided in \cite{RMHDExcat}.

  \begin{Example}[Riemann problem 1]\label{exRMHDRMT1}\rm
The initial data of the first Riemann problem are
$$(\rho,v_x,v_y,v_z,B_x,B_y,B_z,p)=\begin{cases}
  (1,0,0,0,0.5,1,0,1),&x<0,\\
  (0.125,0,0,0,0.5,-1,0,0.1),&x>0,
  \end{cases}
  $$
  with the
adiabatic index $\Gamma=2$. It is an extension of Brio-Wu shock tube problem \cite{Brio:Wu:1988} in the
non-relativistic MHDs. Its solution consists of a left-moving fast rarefaction wave,
a slow compound wave, a contact discontinuity, a right-moving slow shock wave, and a right-moving fast rarefaction wave. There is an argument about the validity of the compound wave, which exists in the numerical results given by any shock capturing scheme, but does not appear in the exact solution obtained by the analytical calculation \cite{ZannaBucciantini:2002}.
    \end{Example}

\begin{figure}[!htbp]
    \centering{}
  \begin{tabular}{cc}

    \includegraphics[width=0.35\textwidth]{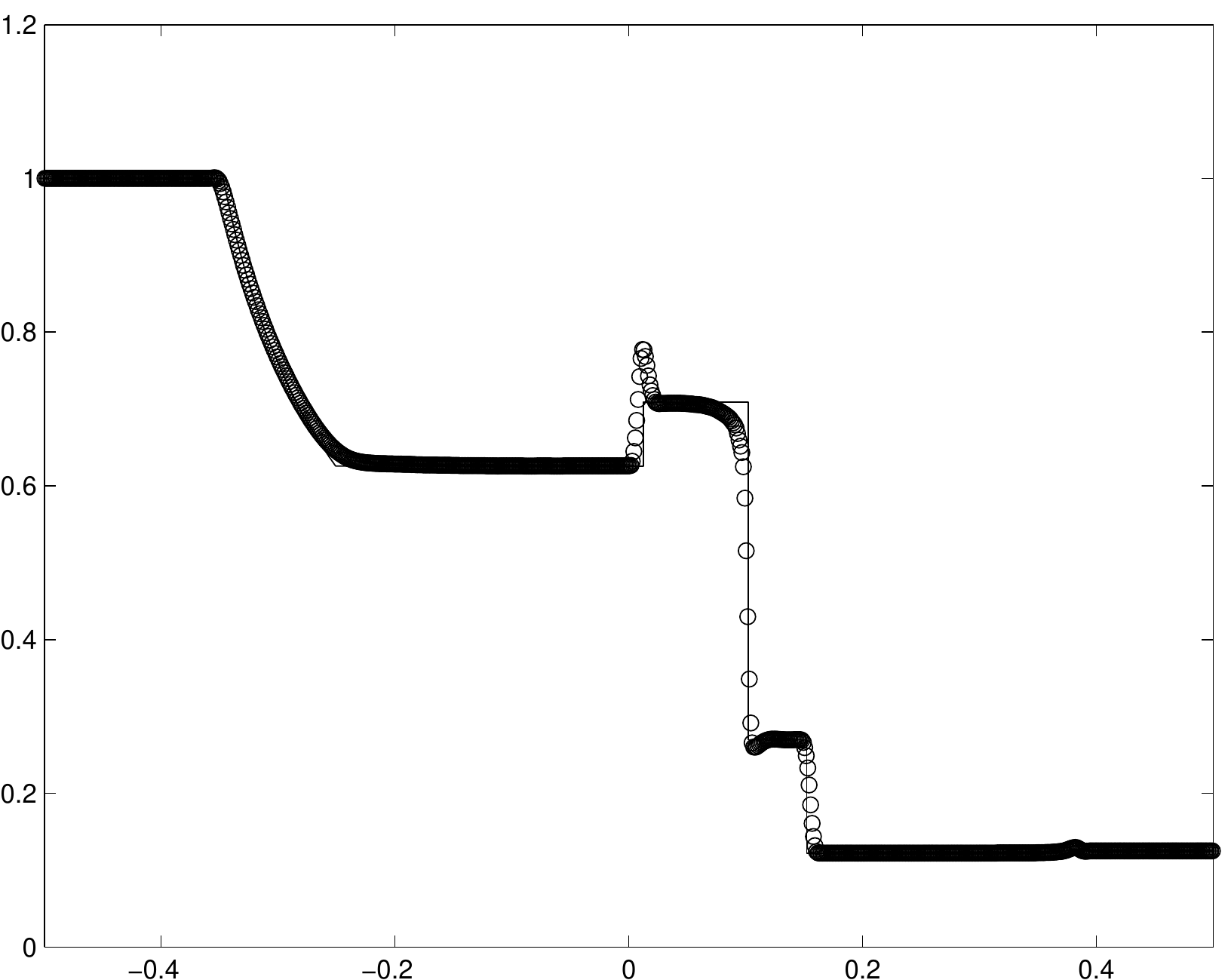}&
\includegraphics[width=0.35\textwidth]{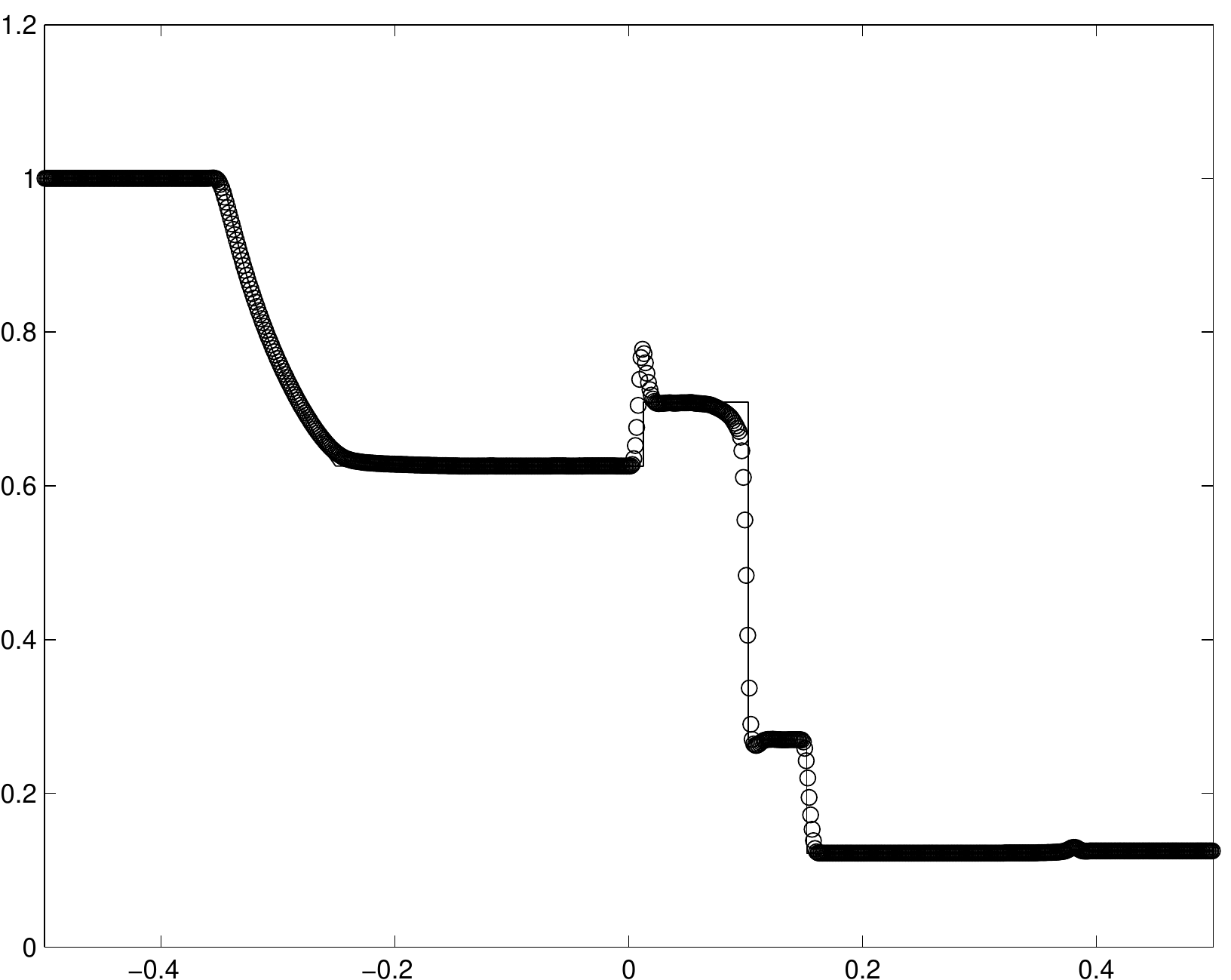}\\
\includegraphics[width=0.35\textwidth]{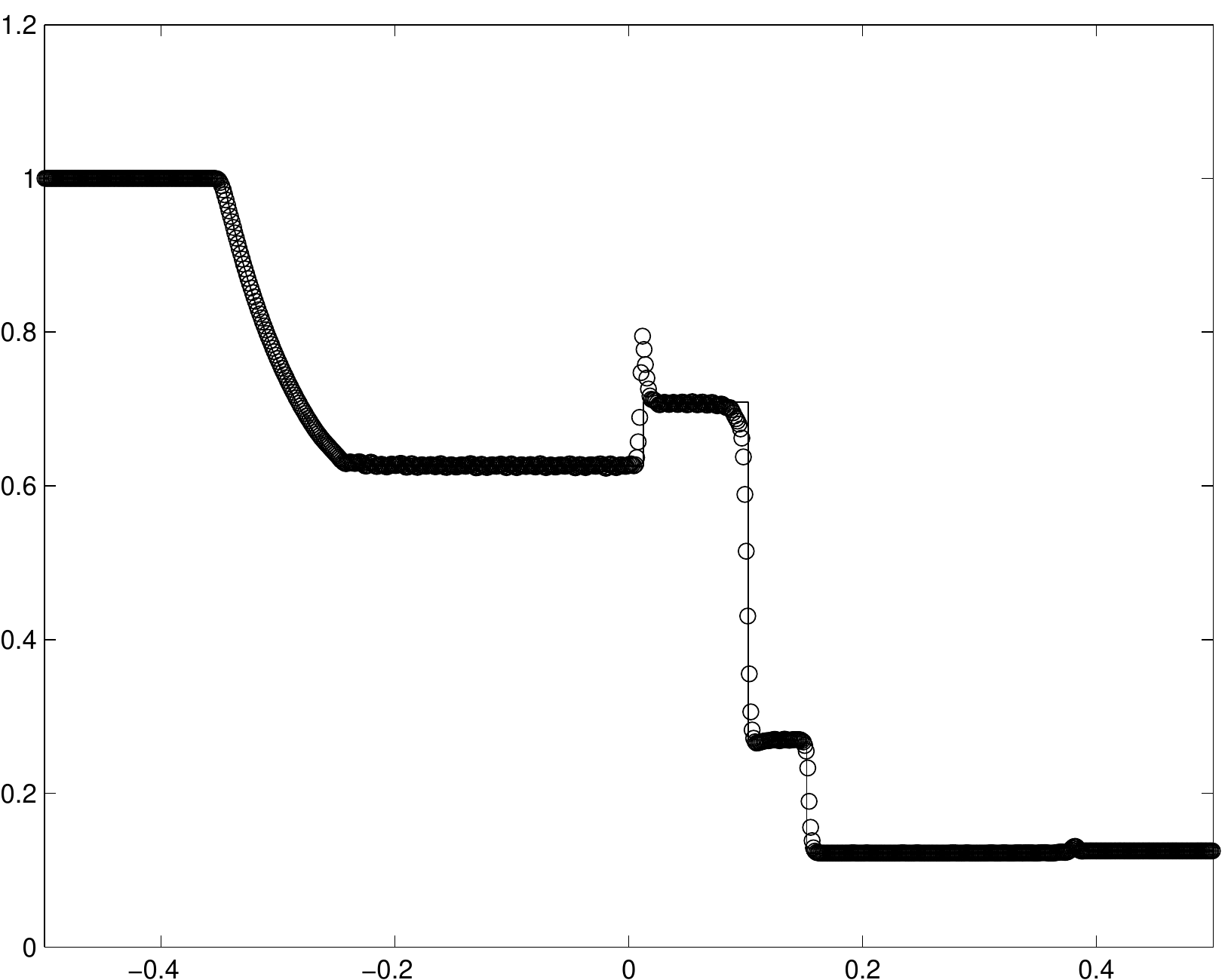}&
\includegraphics[width=0.35\textwidth]{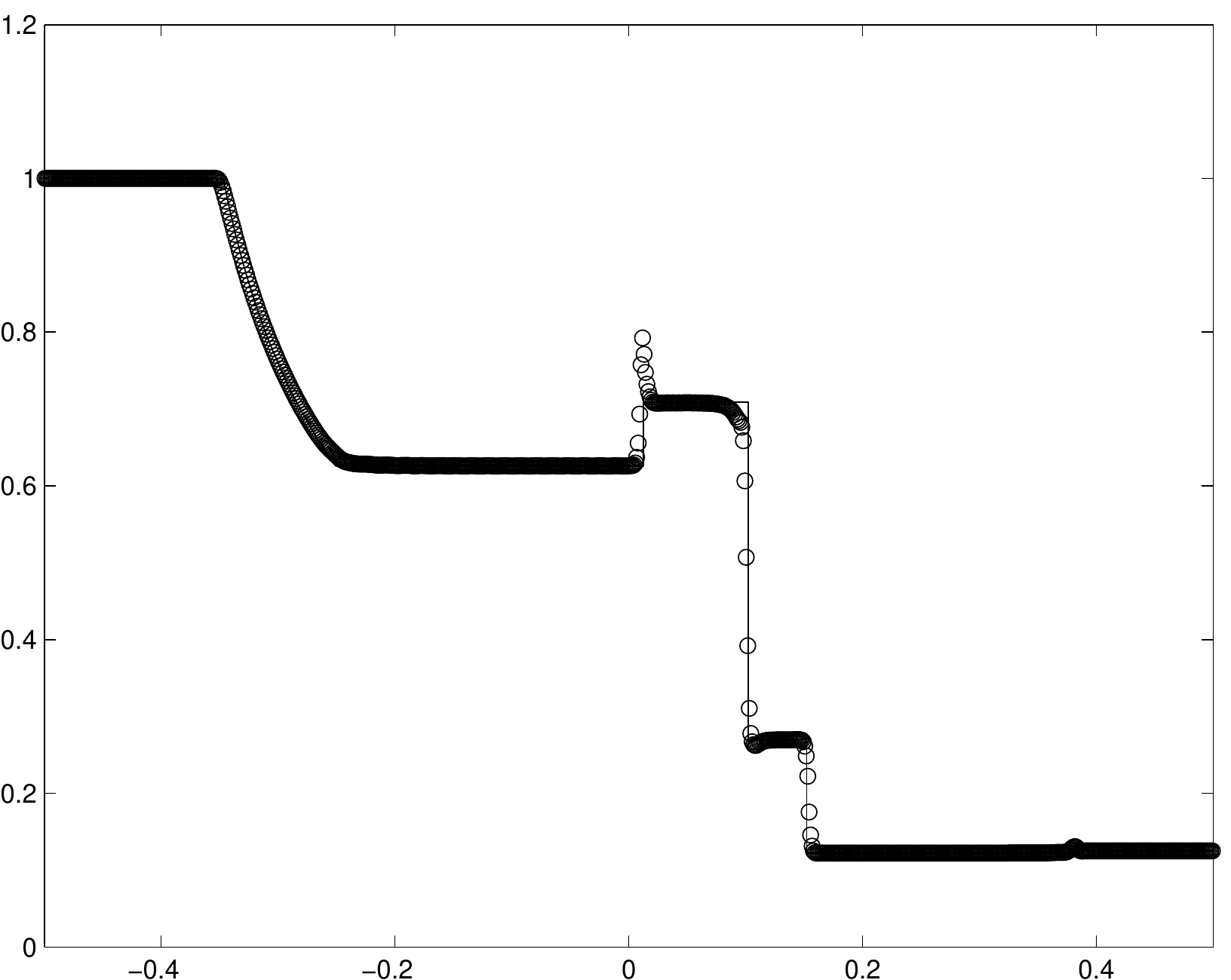}\\
\includegraphics[width=0.35\textwidth]{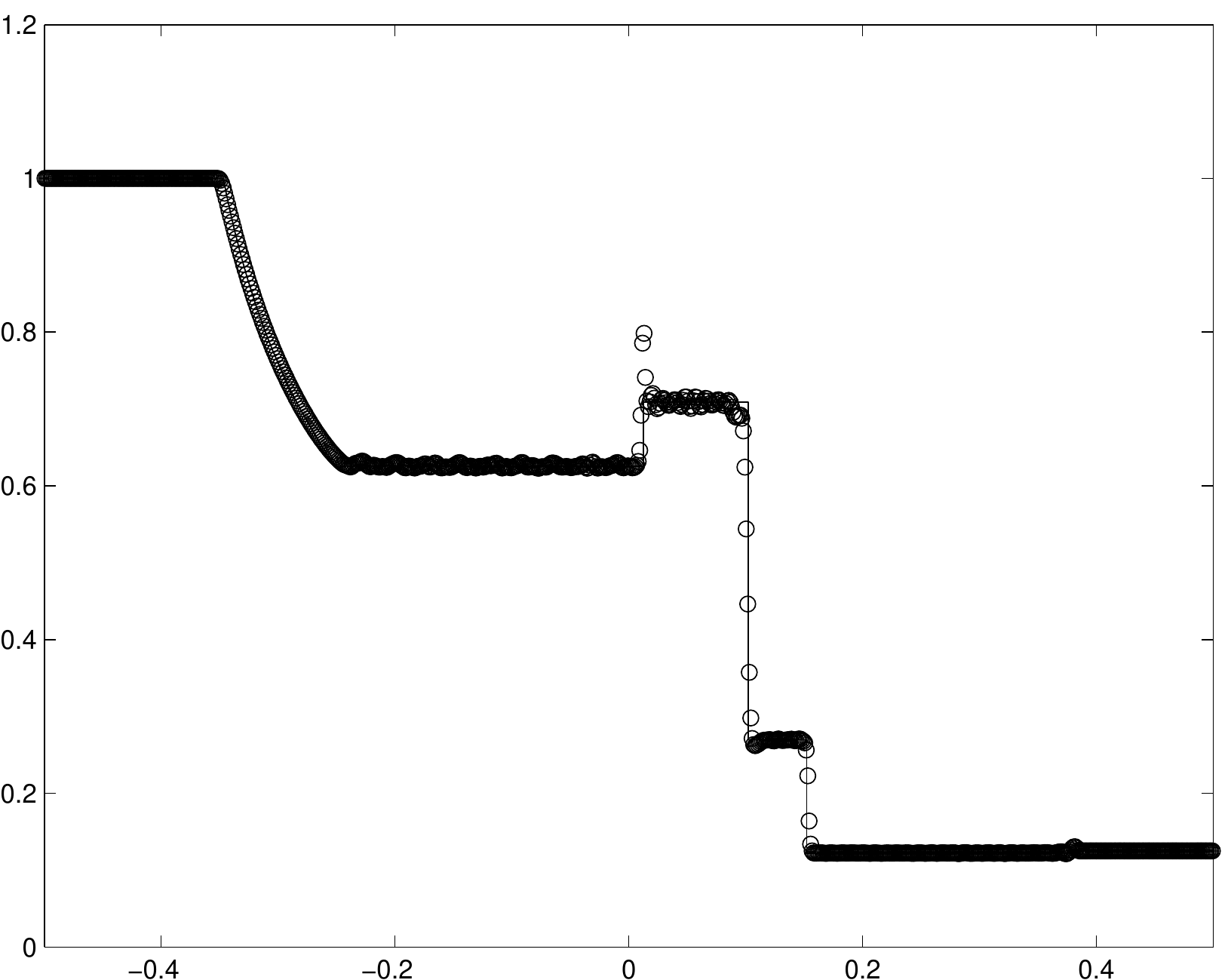}&
\includegraphics[width=0.35\textwidth]{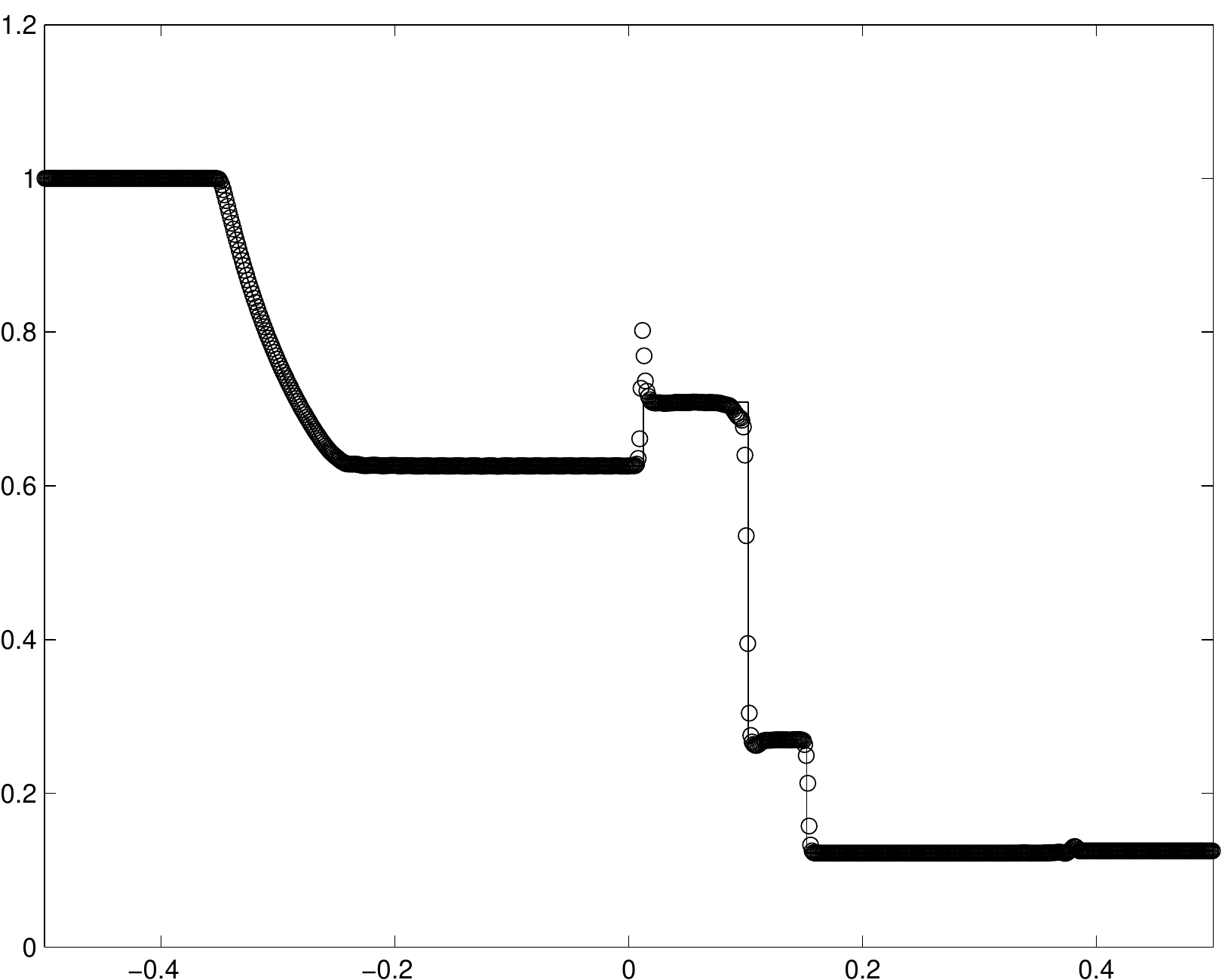}\\
    \end{tabular}
    \caption{Example ~\ref{exRMHDRMT1}£ºThe densities $\rho$ at $t=0.4$.
    The solid line denotes the exact solution, while the symbol ``$\circ$'' is numerical solution
   obtained with $800$ cells. Left: $P^K$-based non-central \DG{}; right: $P^K$-based \CDG{}.
      From top to bottom: $K=1,~2,~3$. }
    \label{fig:RMHDRMT1rho}
  \end{figure}

\begin{figure}[!htbp]
    \centering{}
  \begin{tabular}{cc}
    \includegraphics[width=0.35\textwidth]{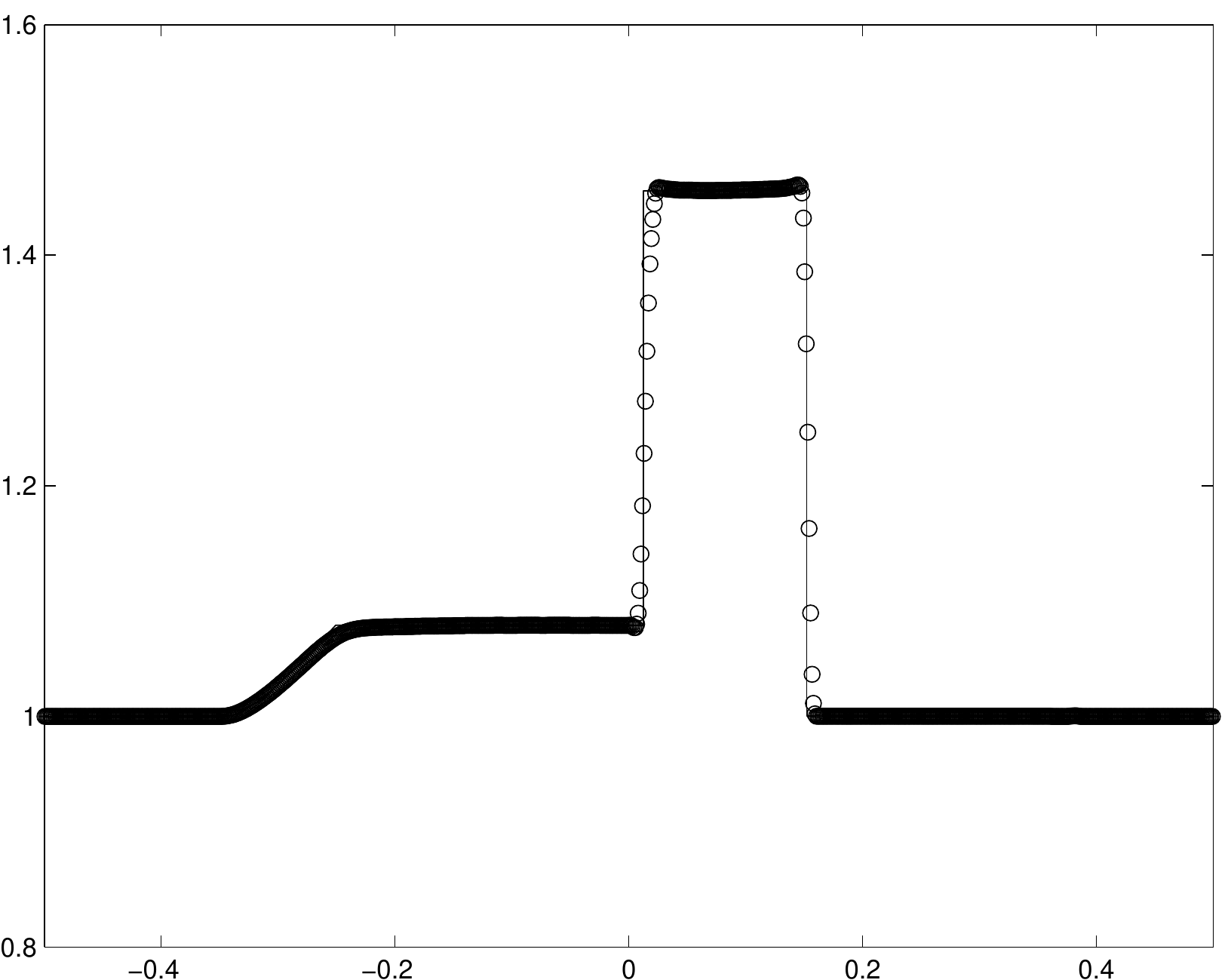}&
\includegraphics[width=0.35\textwidth]{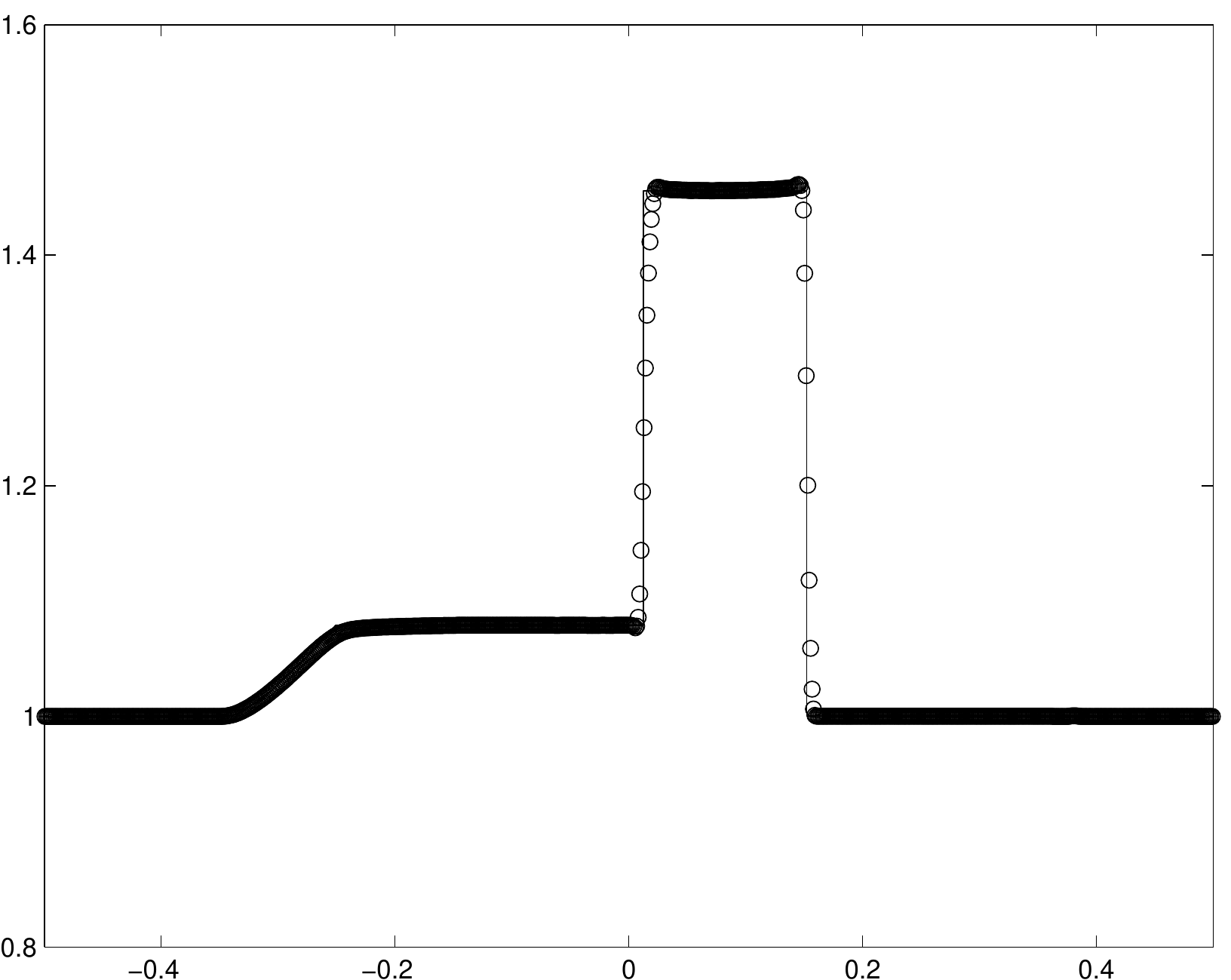}\\
\includegraphics[width=0.35\textwidth]{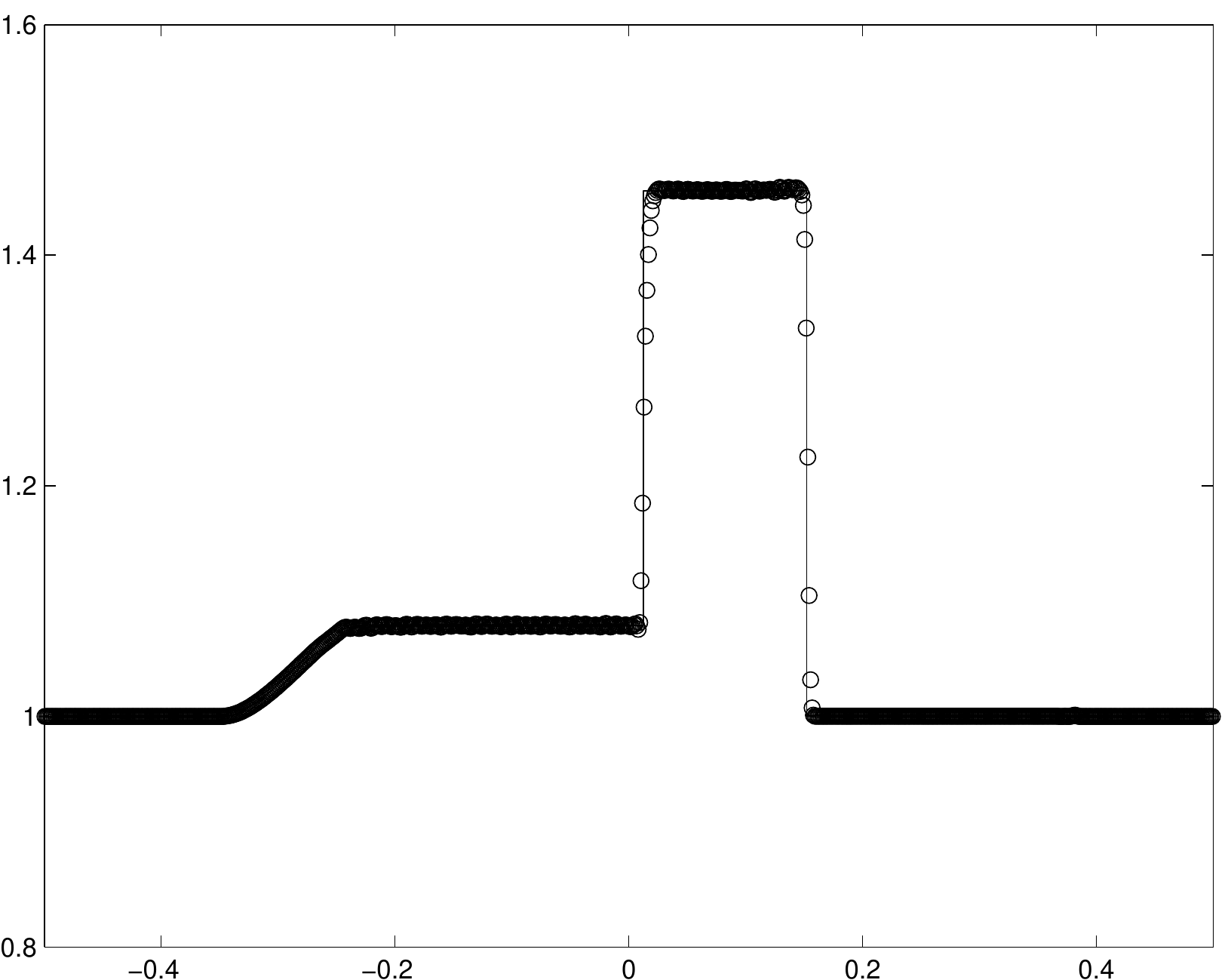}&
\includegraphics[width=0.35\textwidth]{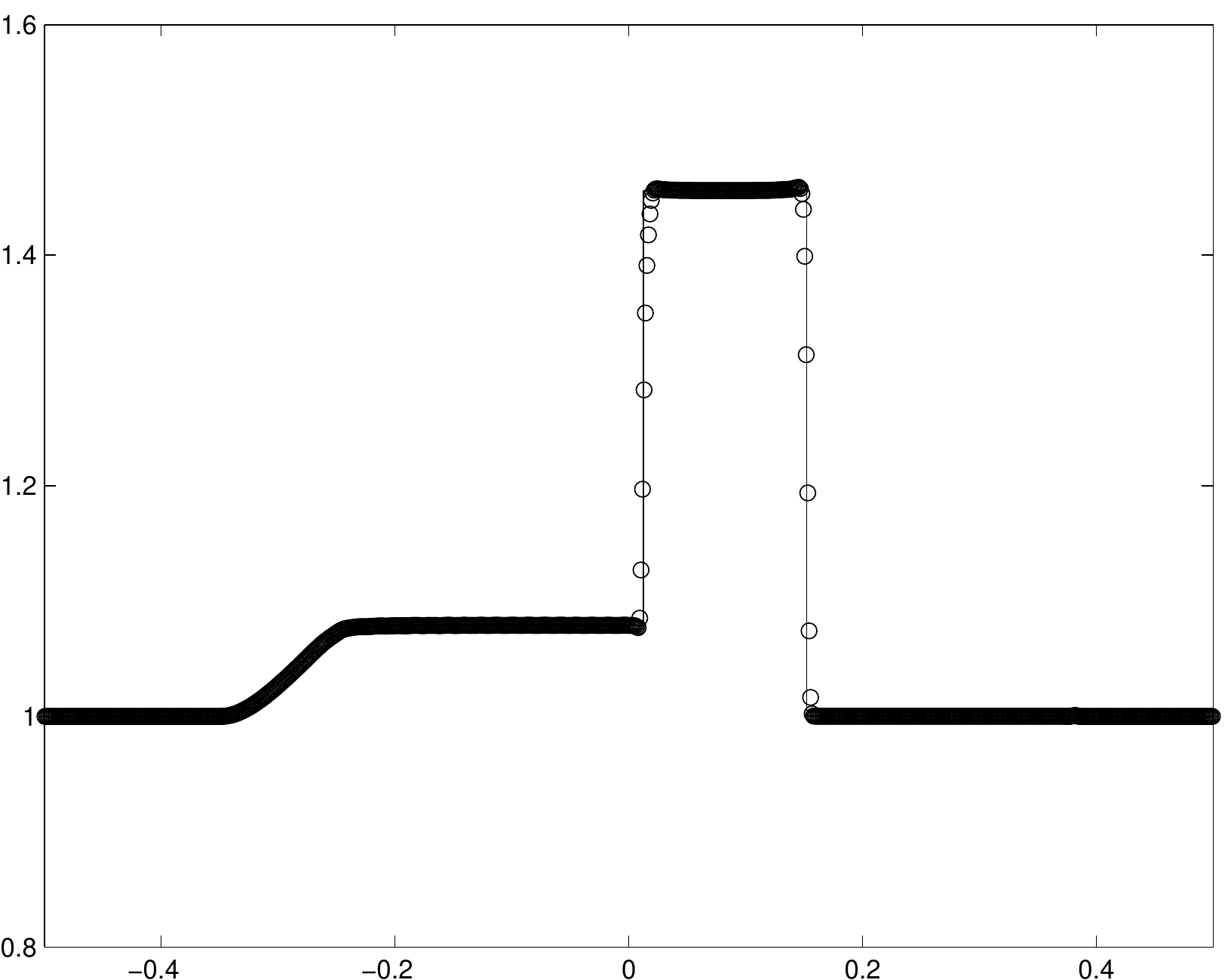}\\
\includegraphics[width=0.35\textwidth]{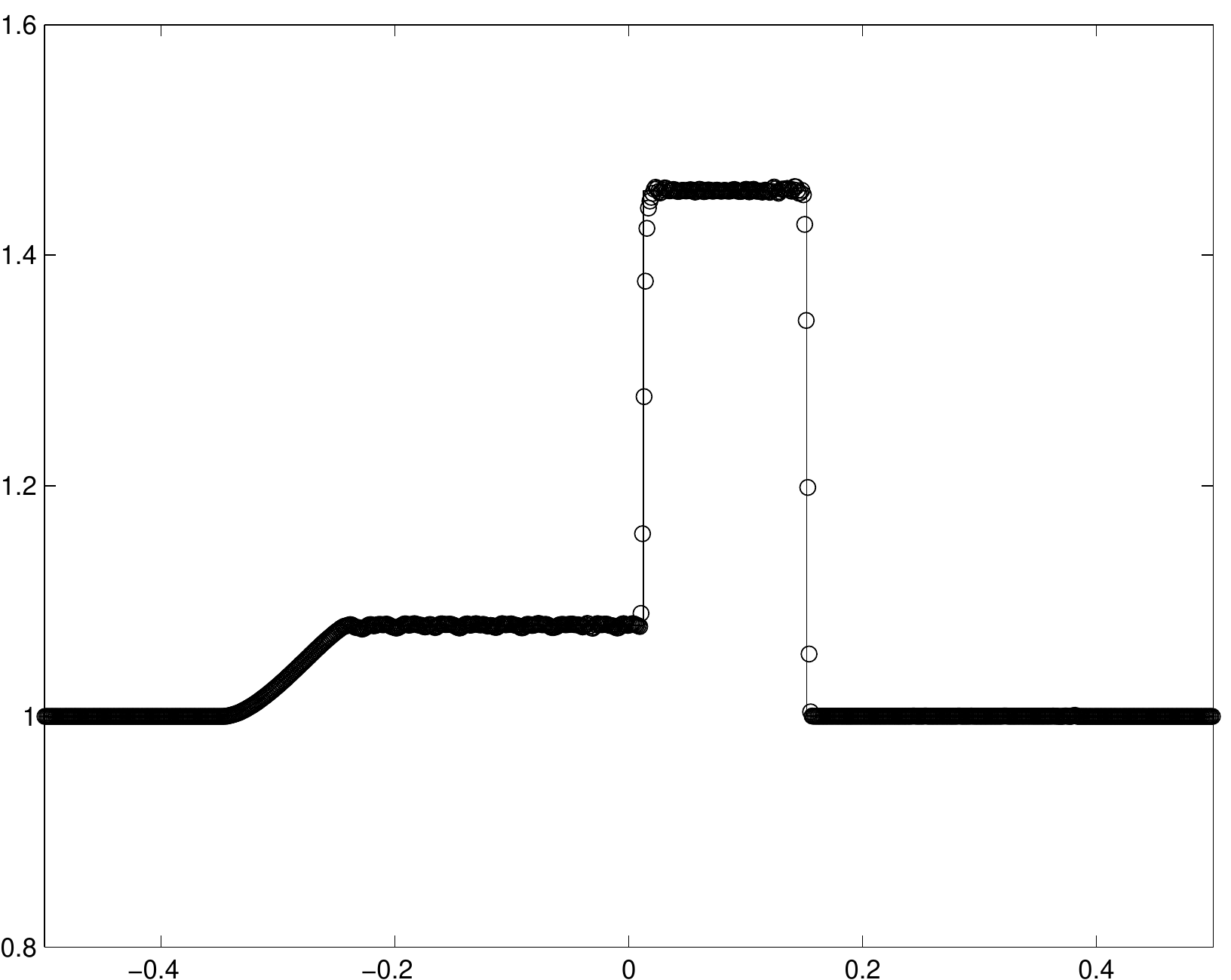}&
\includegraphics[width=0.35\textwidth]{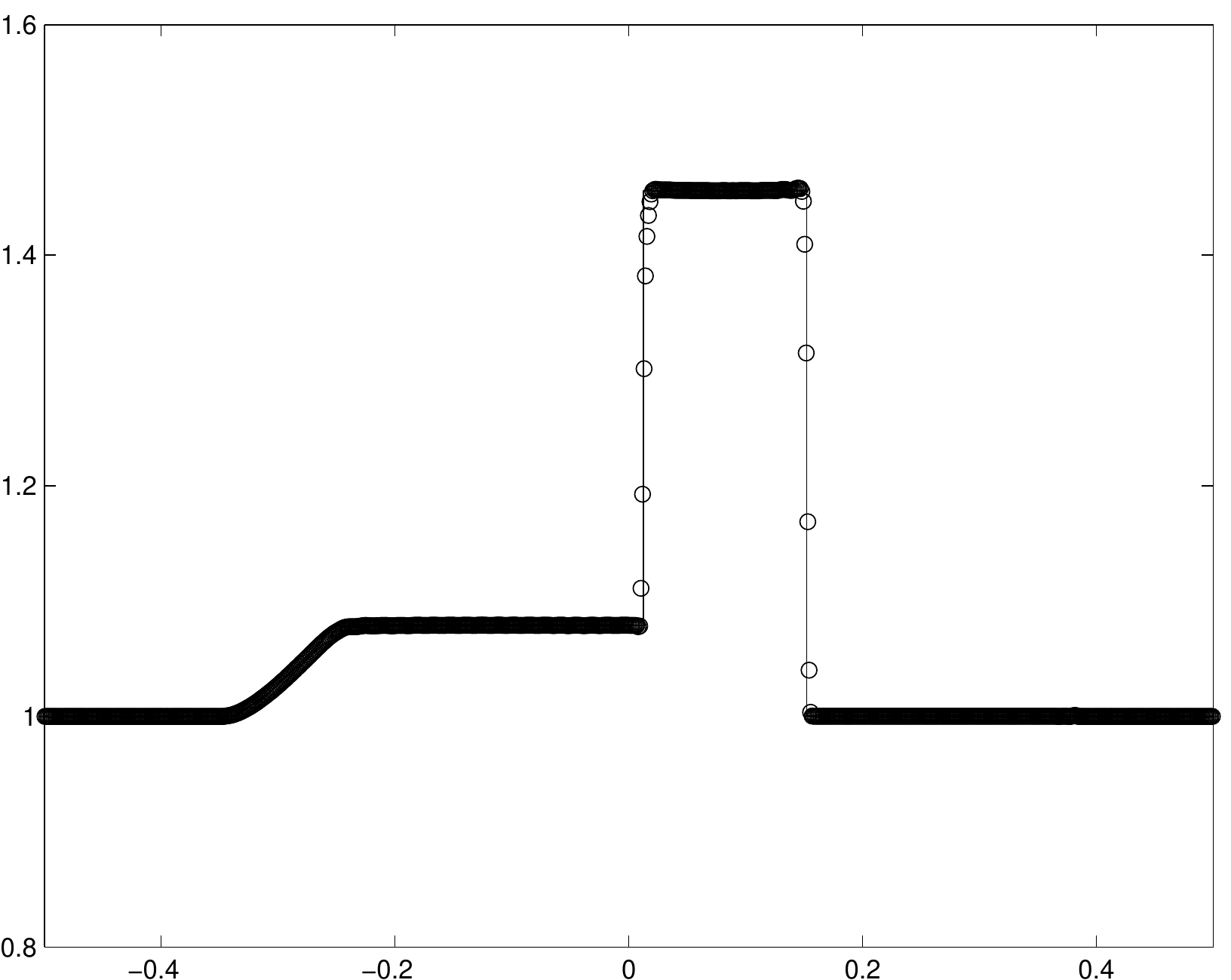}\\
    \end{tabular}
    \caption{Same as Fig.~\ref{fig:RMHDRMT1rho} except for the Lorentz factor $\gamma$. }
    \label{fig:RMHDRMT1gam}
  \end{figure}

\begin{figure}[!htbp]
    \centering{}
  \begin{tabular}{cc}
    \includegraphics[width=0.35\textwidth]{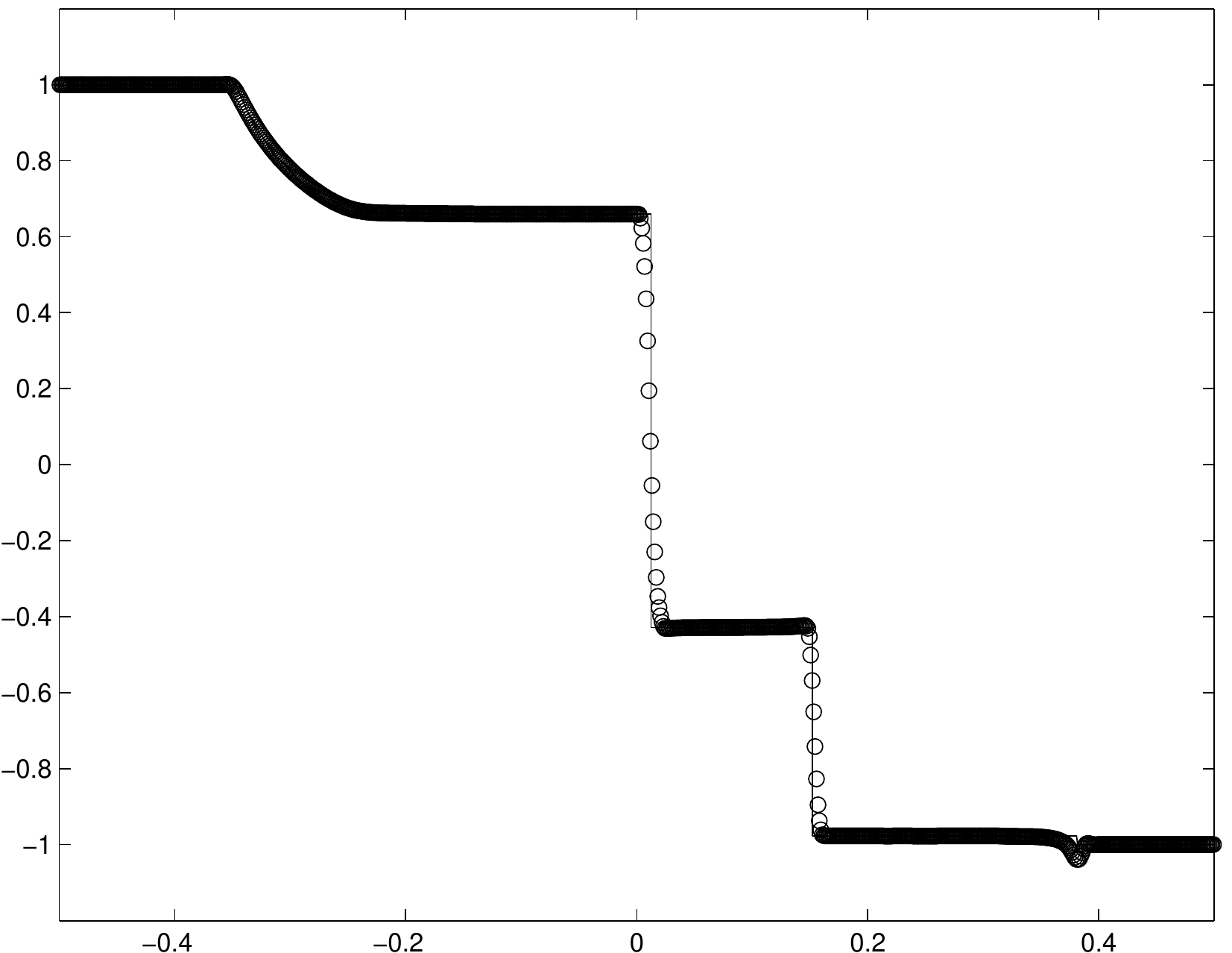}&
\includegraphics[width=0.35\textwidth]{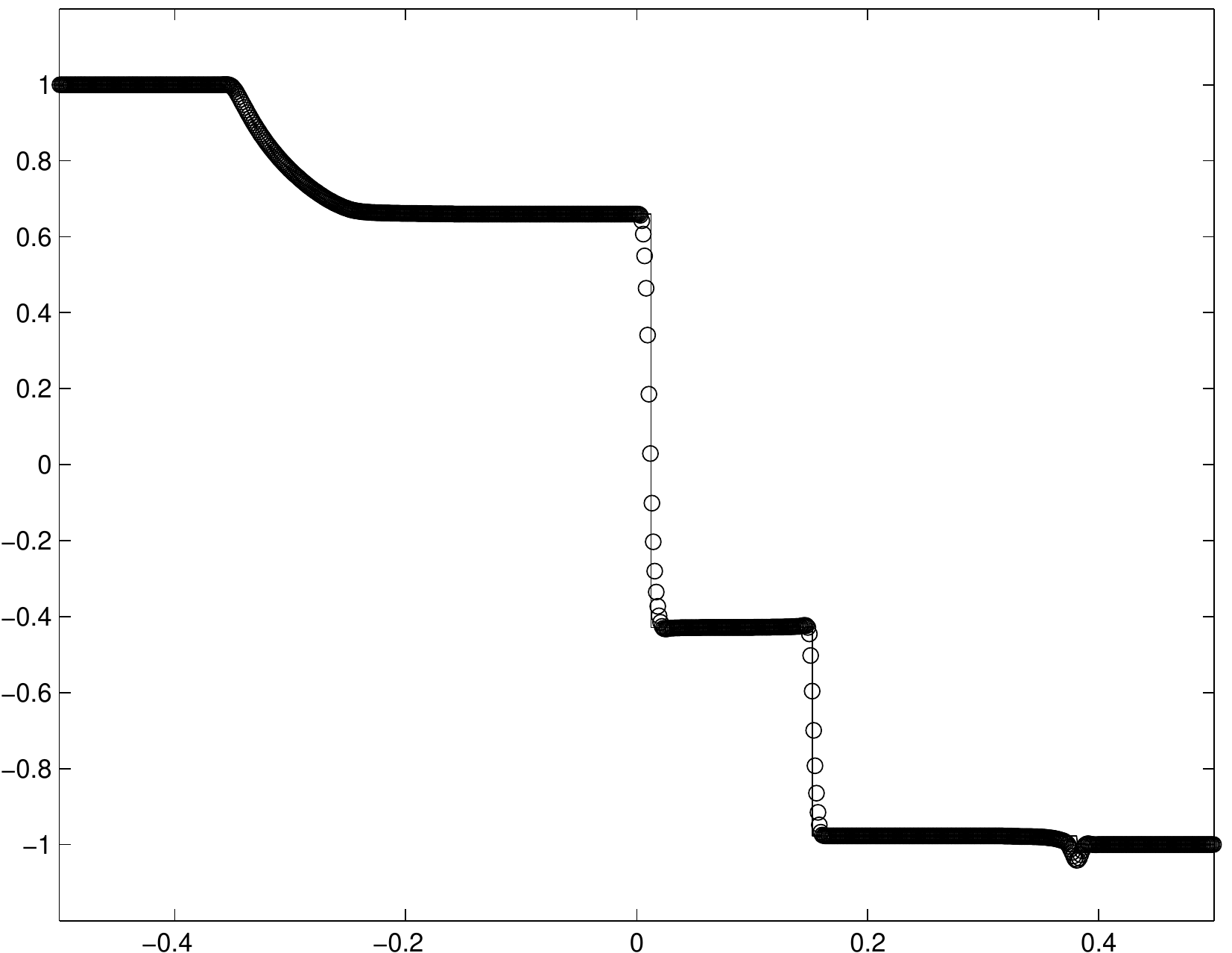}\\
\includegraphics[width=0.35\textwidth]{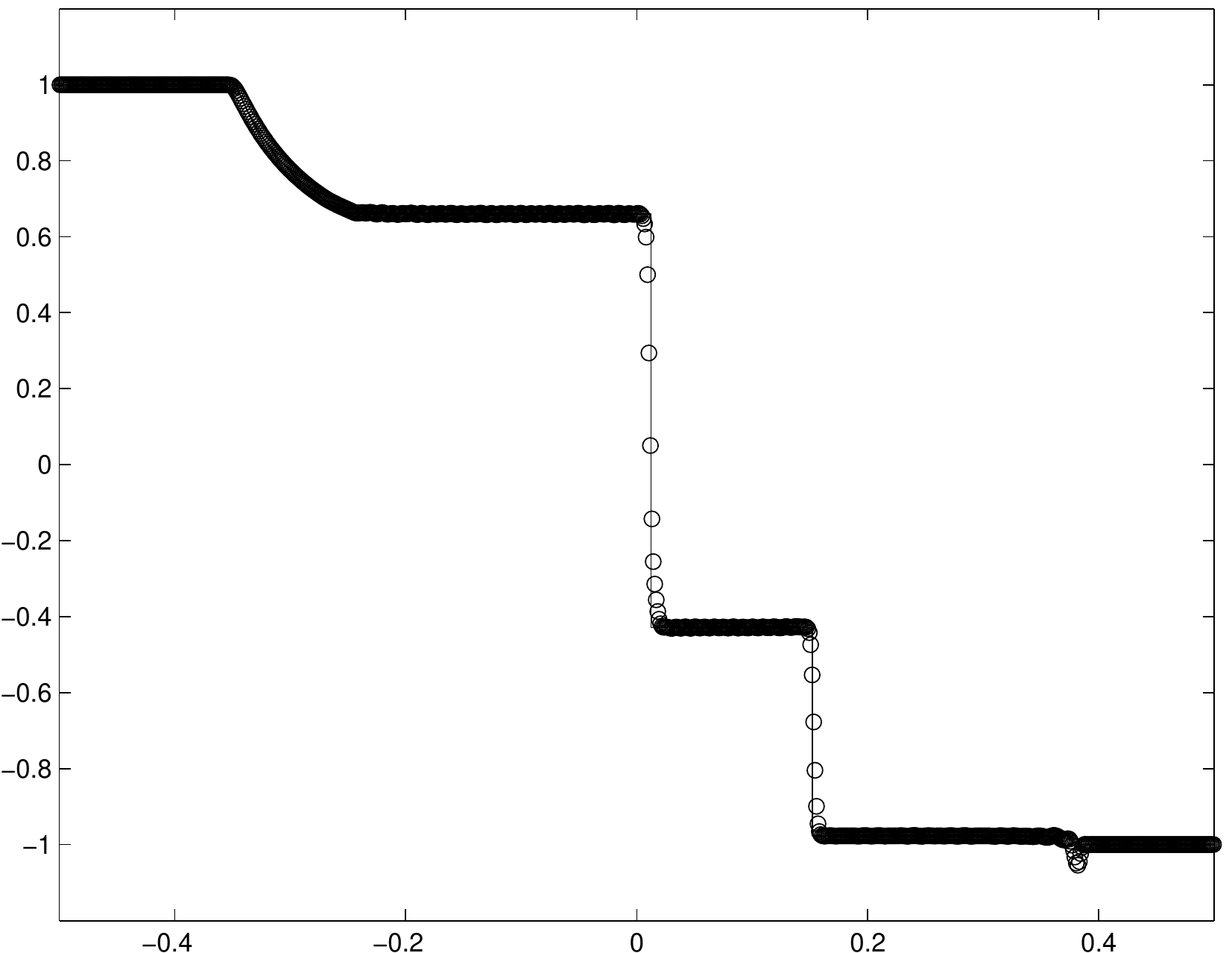}&
\includegraphics[width=0.35\textwidth]{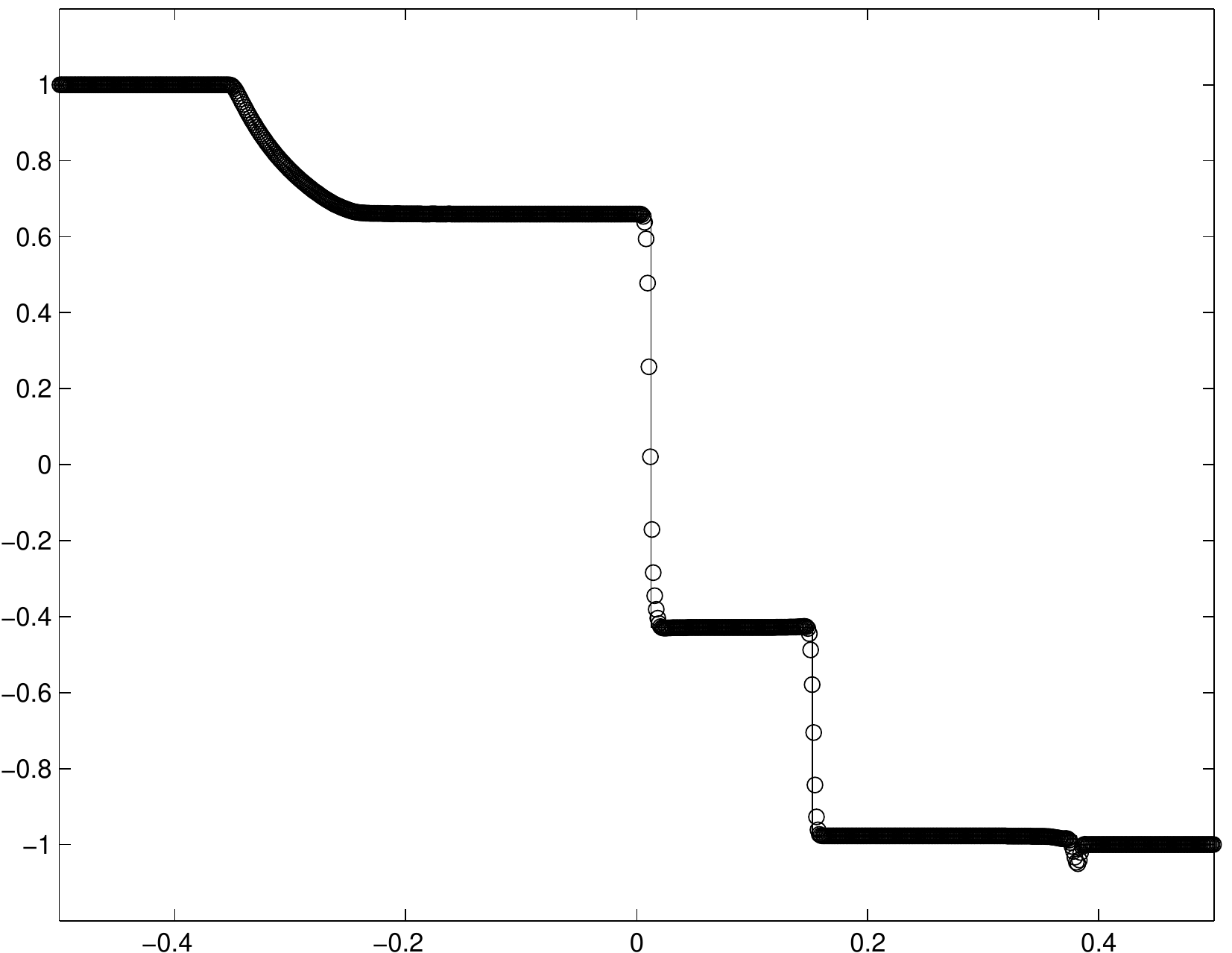}\\
\includegraphics[width=0.35\textwidth]{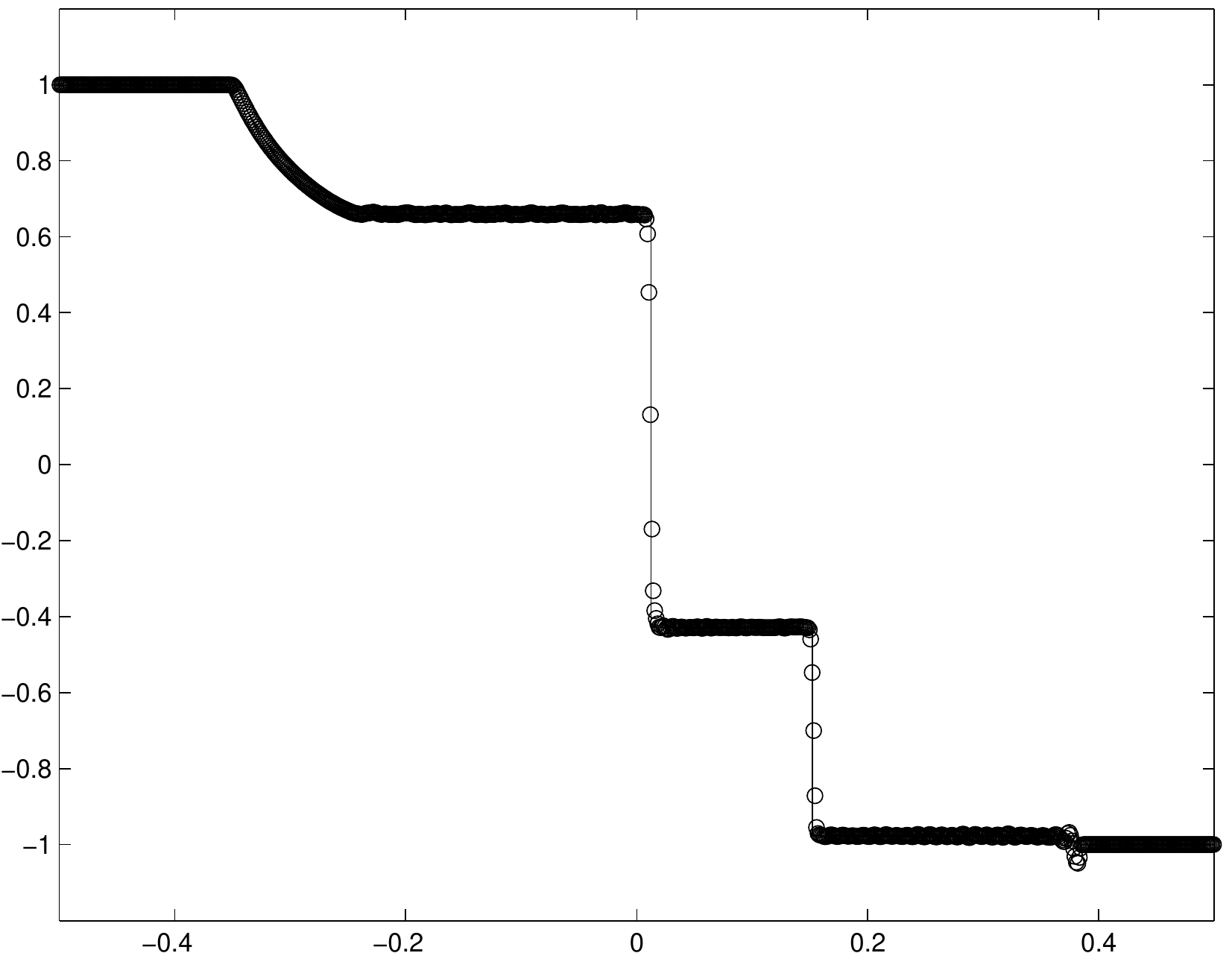}&
\includegraphics[width=0.35\textwidth]{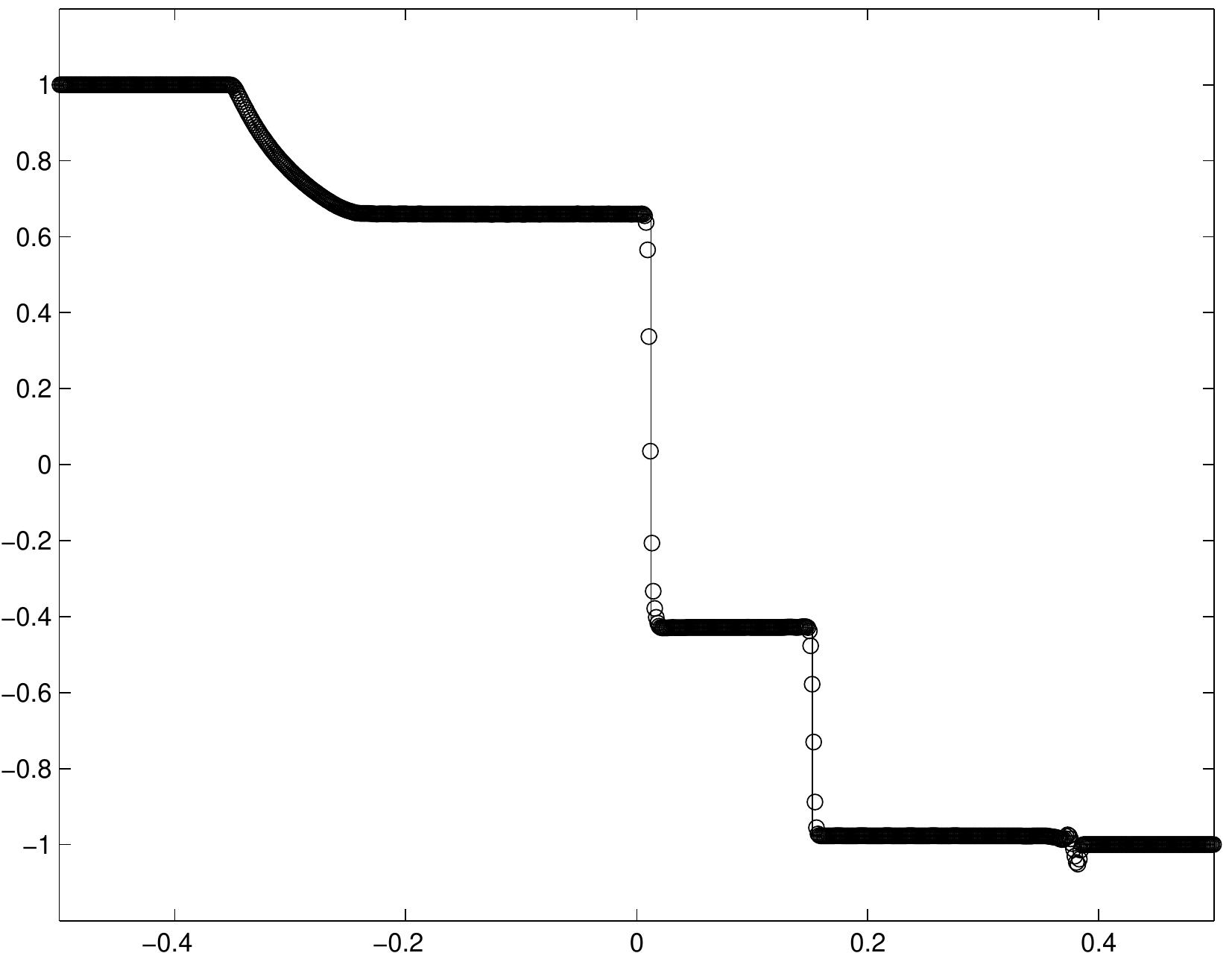}\\
    \end{tabular}
    \caption{Same as Fig.~\ref{fig:RMHDRMT1rho} except for $B_y$.}
    \label{fig:RMHDRMT1by}
  \end{figure}

\begin{figure}[!htbp]
    \centering{}
  \begin{tabular}{cc}
    \includegraphics[width=0.35\textwidth]{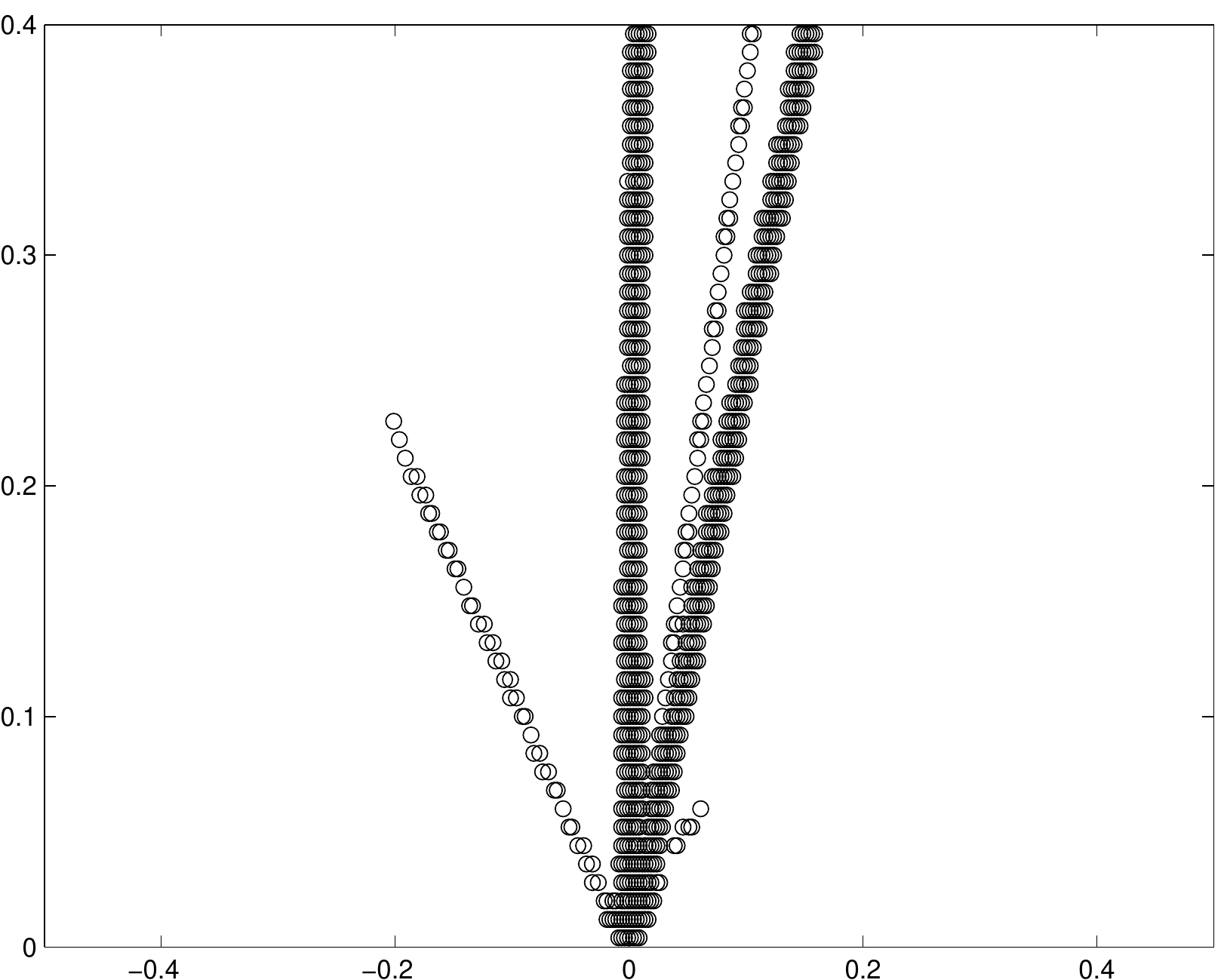}&
\includegraphics[width=0.35\textwidth]{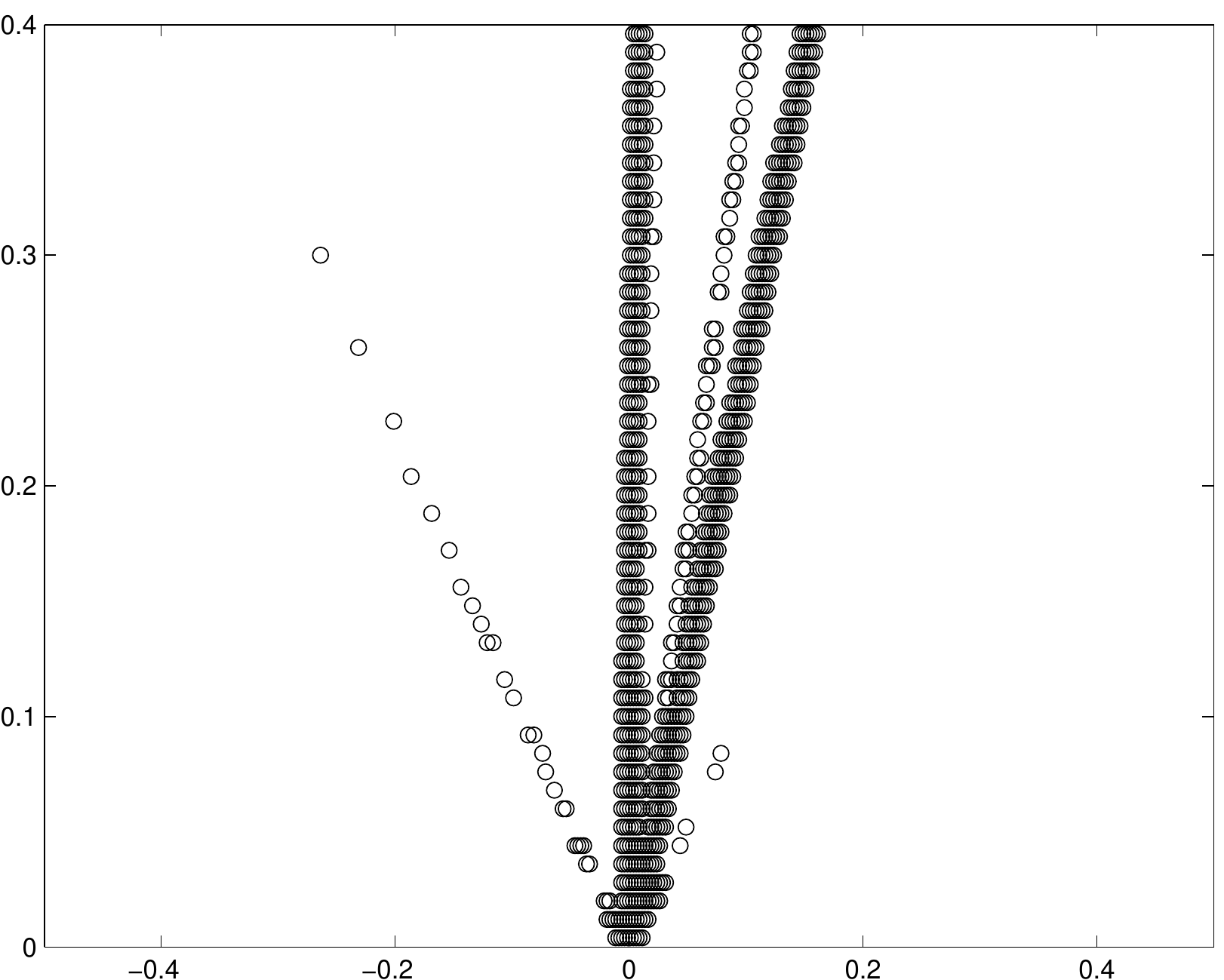}\\
\includegraphics[width=0.35\textwidth]{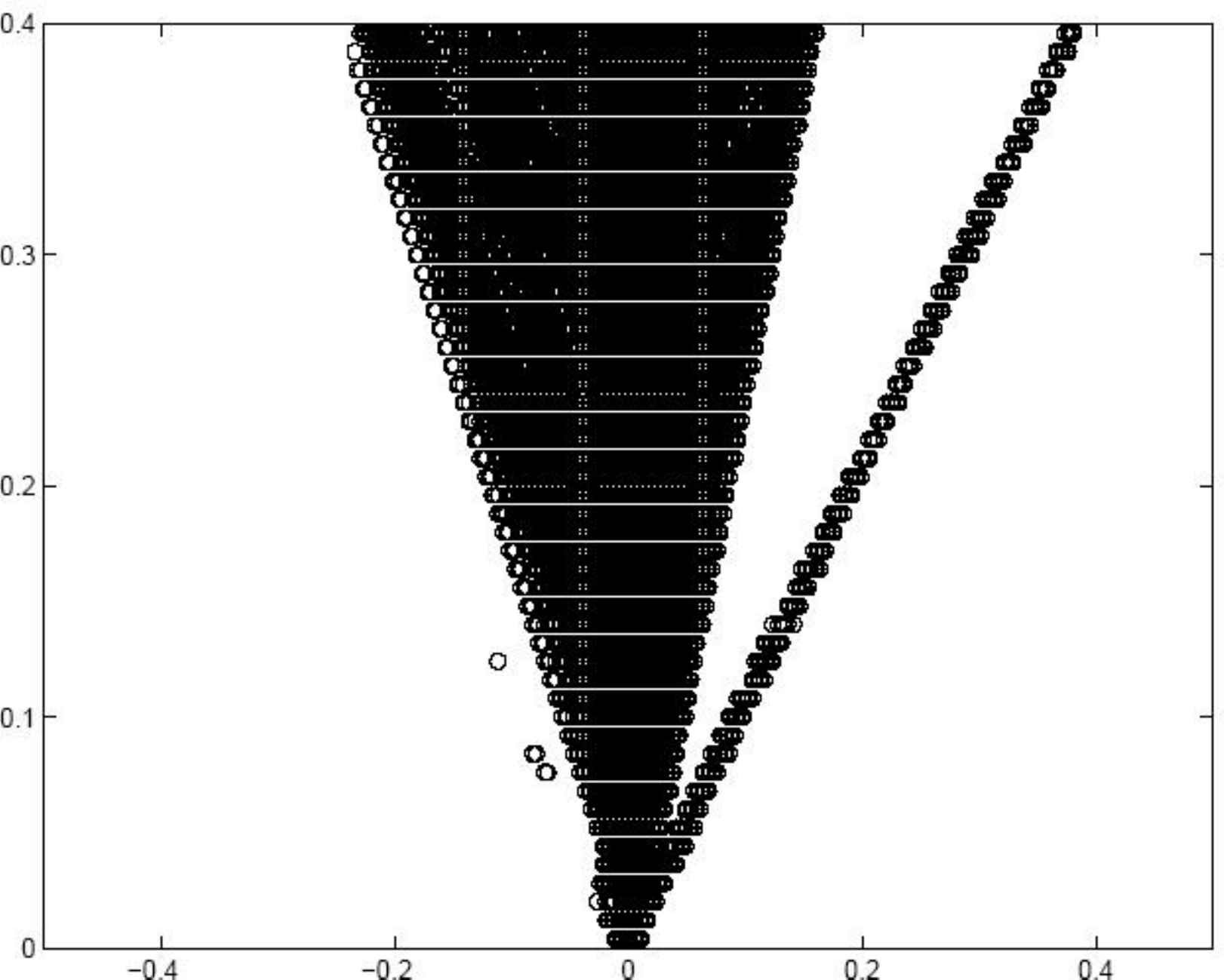}&
\includegraphics[width=0.35\textwidth]{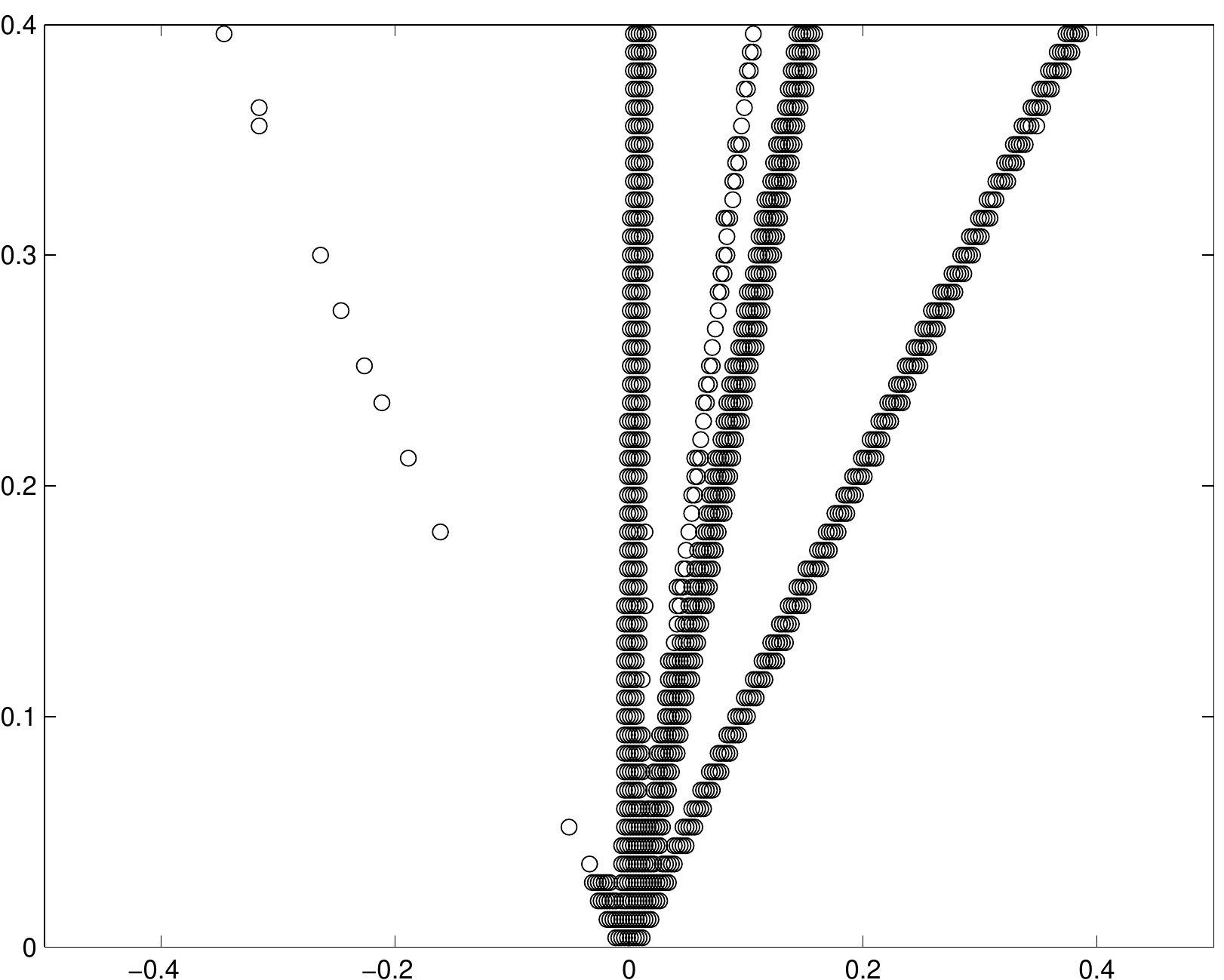}\\
\includegraphics[width=0.35\textwidth]{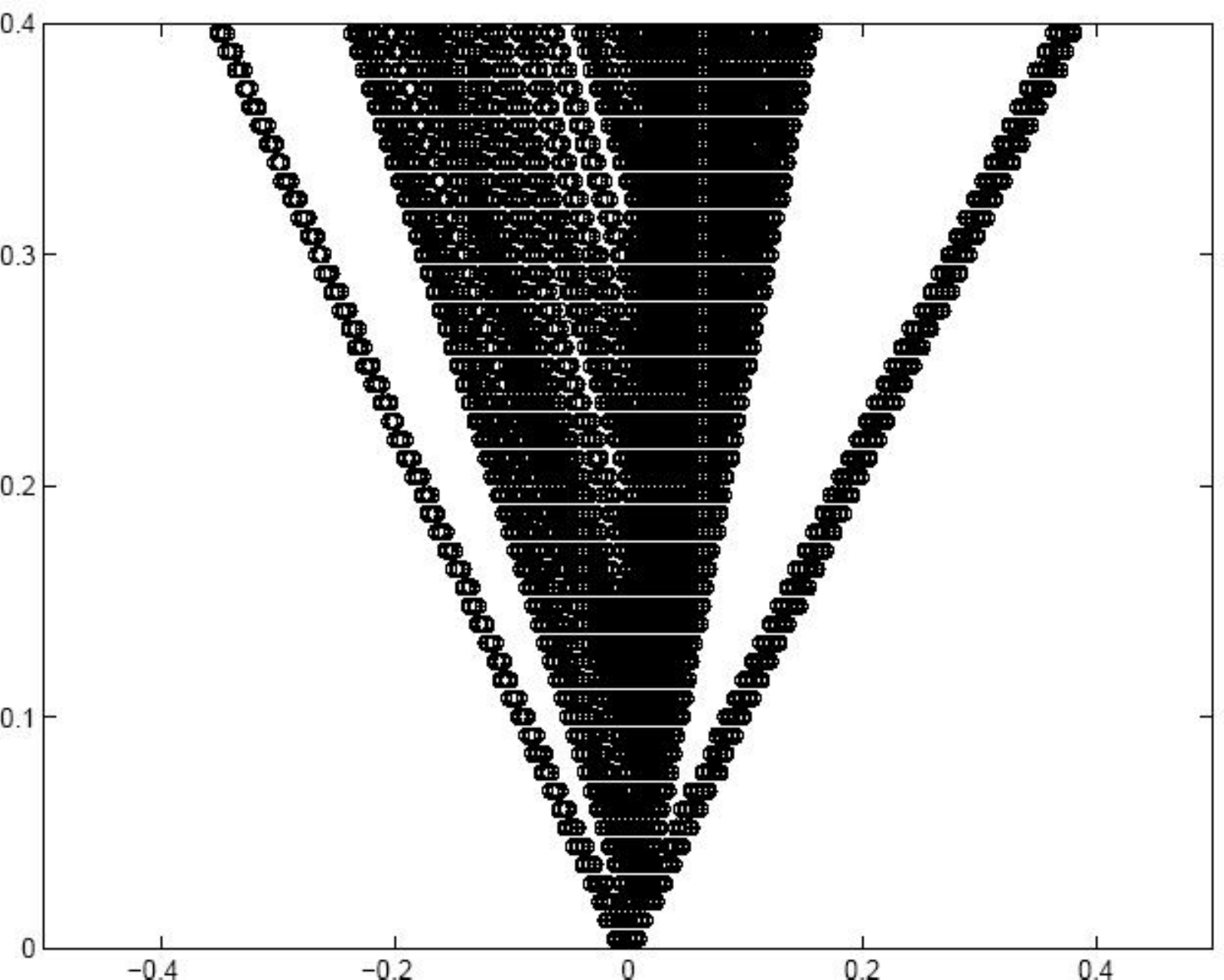}&
\includegraphics[width=0.35\textwidth]{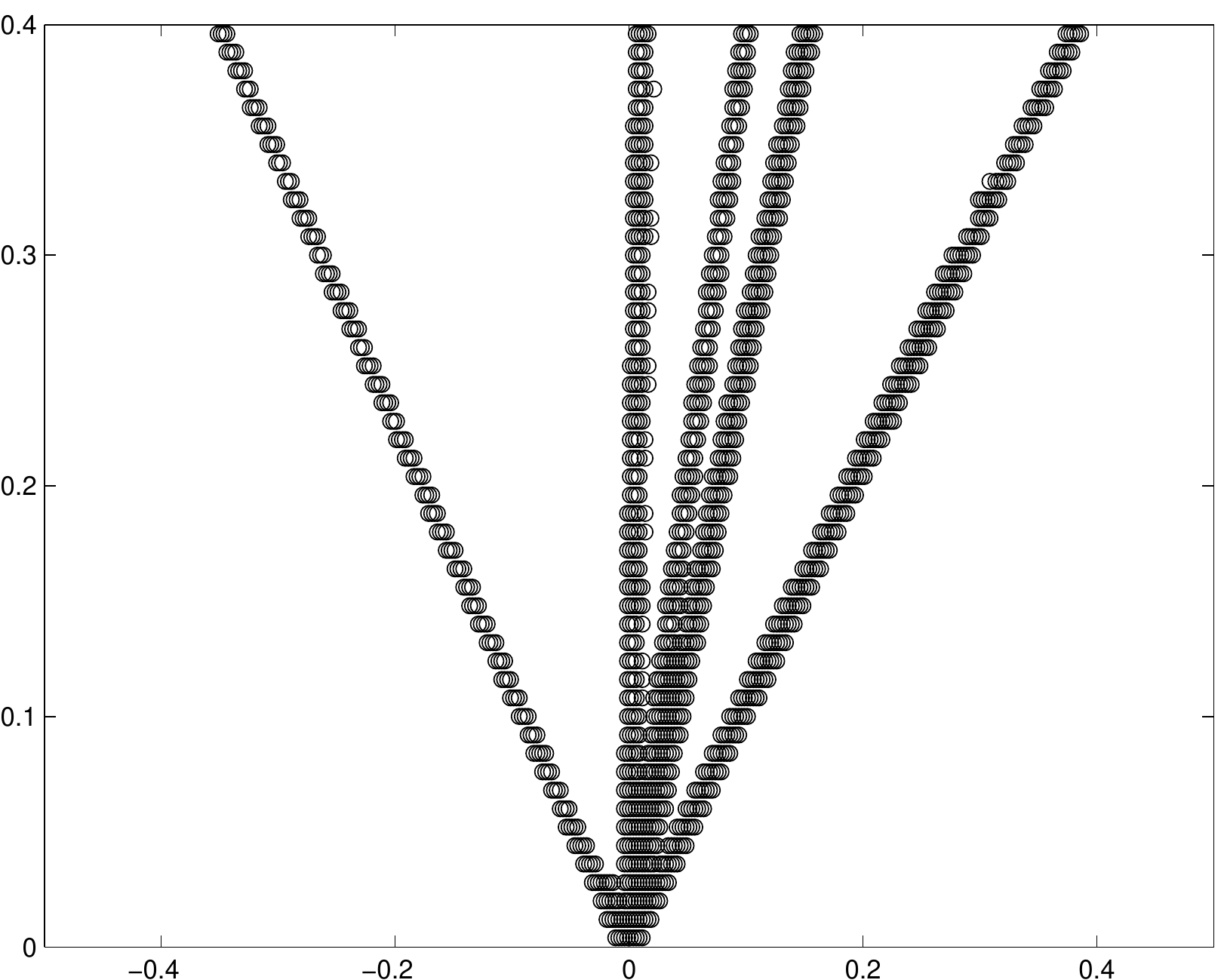}\\
    \end{tabular}
    \caption{Example ~\ref{exRMHDRMT1}: The time evolution of   ``troubled'' cells.
     Left: non-central \DG{}; right: \CDG{}.
      From top to bottom: $K=1,~2,~3$. The cell number   is $800$.}
    \label{fig:RMHDRMT1cell}
  \end{figure}

Figs.~\ref{fig:RMHDRMT1rho},~\ref{fig:RMHDRMT1gam}, and \ref{fig:RMHDRMT1by}
plot the densities $\rho$, Lorentz factors $\gamma$, and magnetic fields $B_y$
at $t=0.4$ obtained by using the non-central and central DG methods.
It is seen that those numerical solutions are   in good agreement with the exact solutions,
 there exist more obvious numerical oscillations at the left-hand side of the compound wave
in the solutions obtained by the  $P^3$-based non-central DG methods, but the
numerical oscillation is not too obvious in the solution of the central DG methods. Moreover,
the non-central DG methods identify  more ``troubled'' cells than the central DG, see Fig.~\ref{fig:RMHDRMT1cell}.

  \begin{Example}[Riemann problem 2]\label{exRMHDRMT3}\rm
The initial data of second Riemann problem are
$$(\rho,v_x,v_y,v_z,B_x,B_y,B_z,p)=\begin{cases}
  (1,0,0,0,5,6,6,30),&x<0,\\
  (1,0,0,0,5,0.7,0.7,1),&x>0,
  \end{cases}
  $$
 with the adiabatic index $\Gamma=5/3$. As time increases,
  the initial discontinuity will decompose into two left-moving rarefaction waves, and a contact discontinuity
  and two right-moving shock waves. 
   \end{Example}

\begin{figure}[!htbp]
    \centering{}
  \begin{tabular}{cc}

    \includegraphics[width=0.35\textwidth]{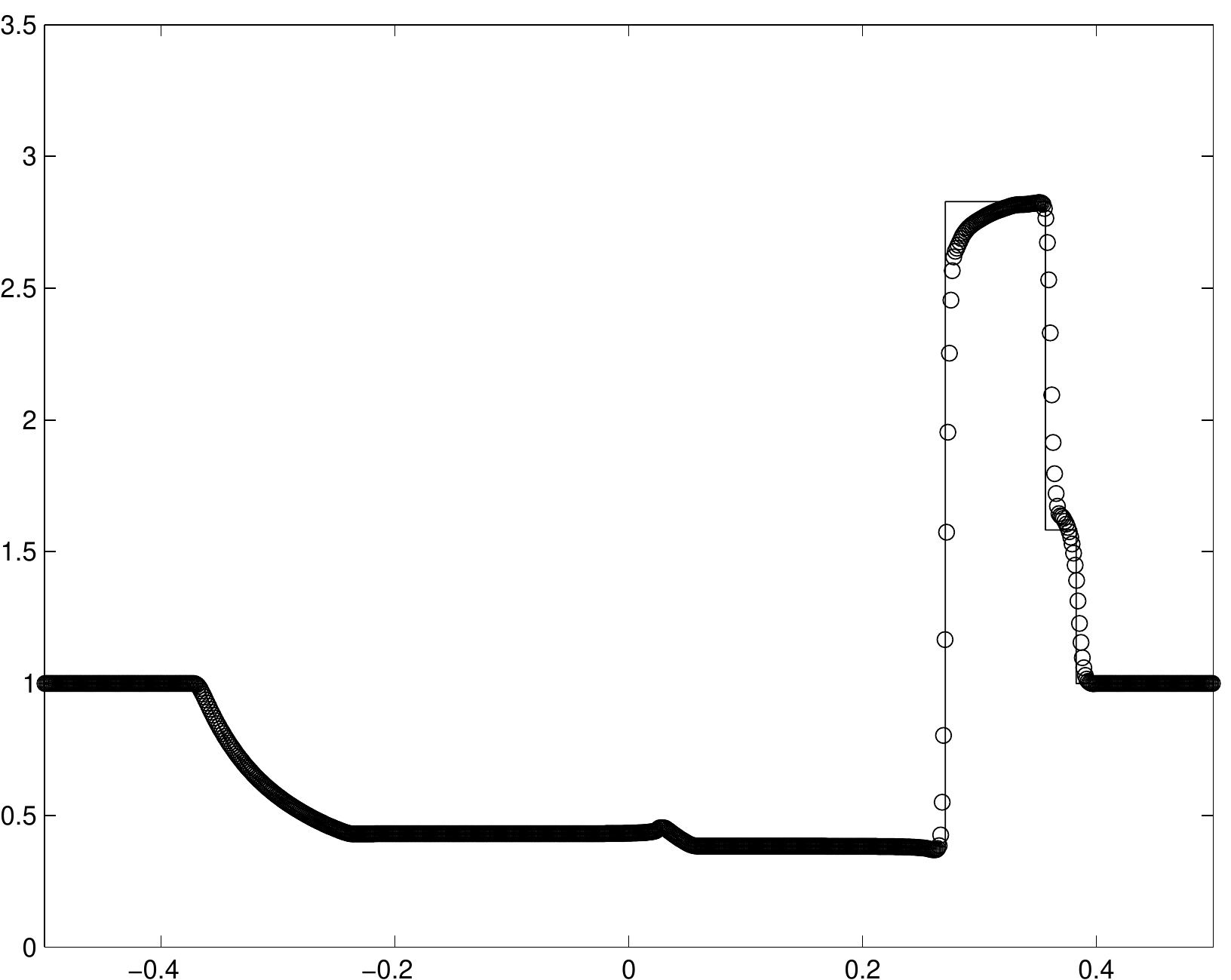}&
\includegraphics[width=0.35\textwidth]{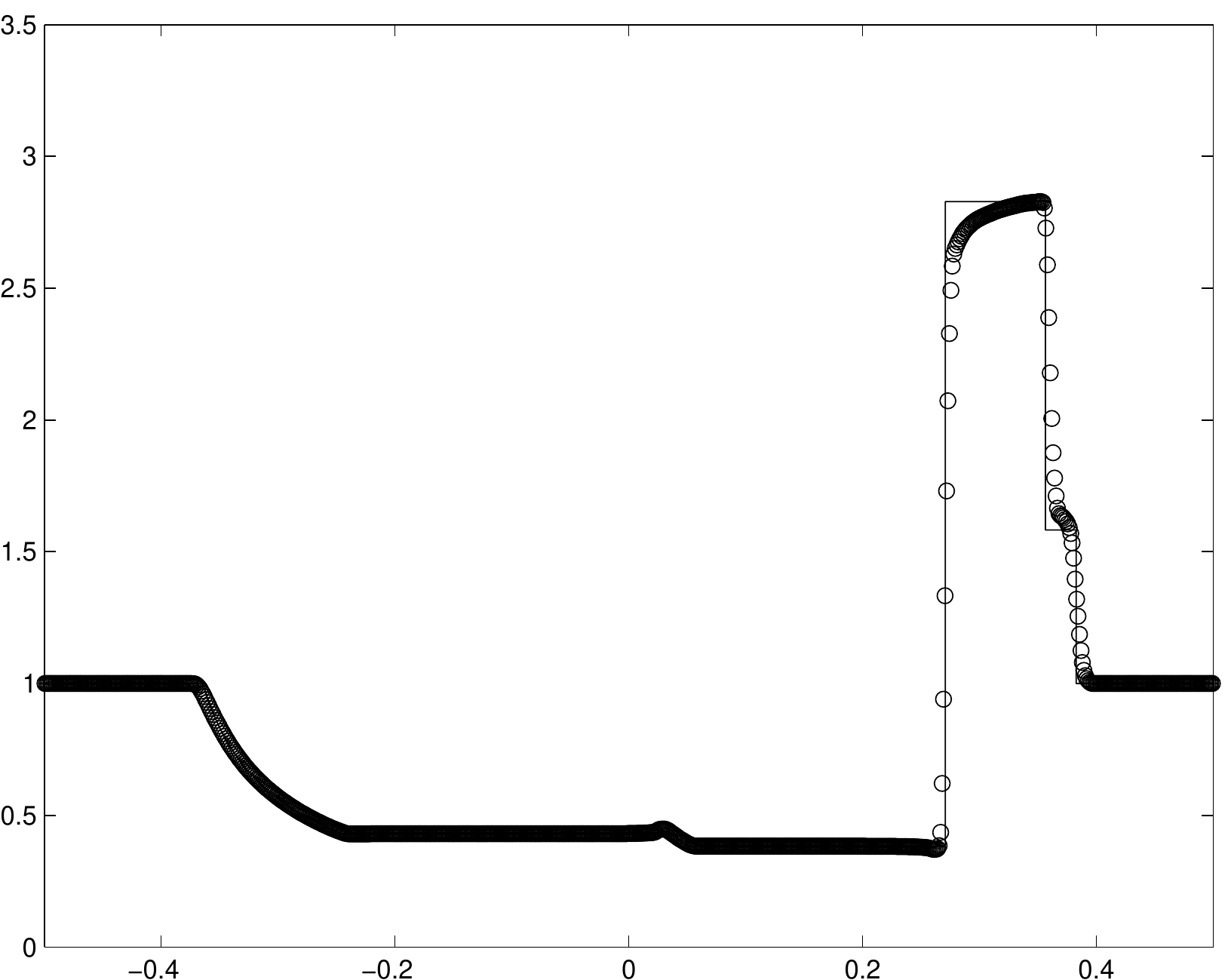}\\
\includegraphics[width=0.35\textwidth]{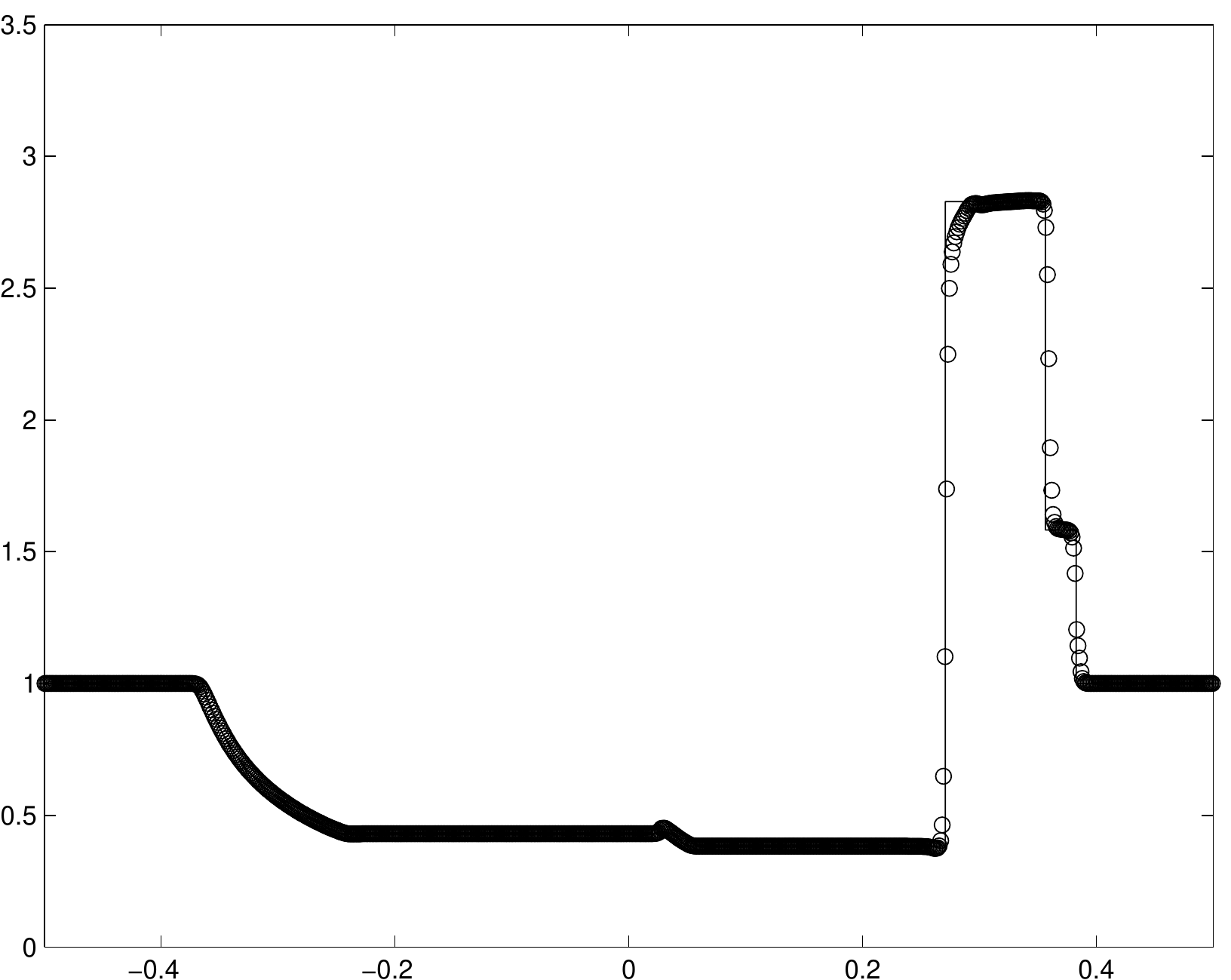}&
\includegraphics[width=0.35\textwidth]{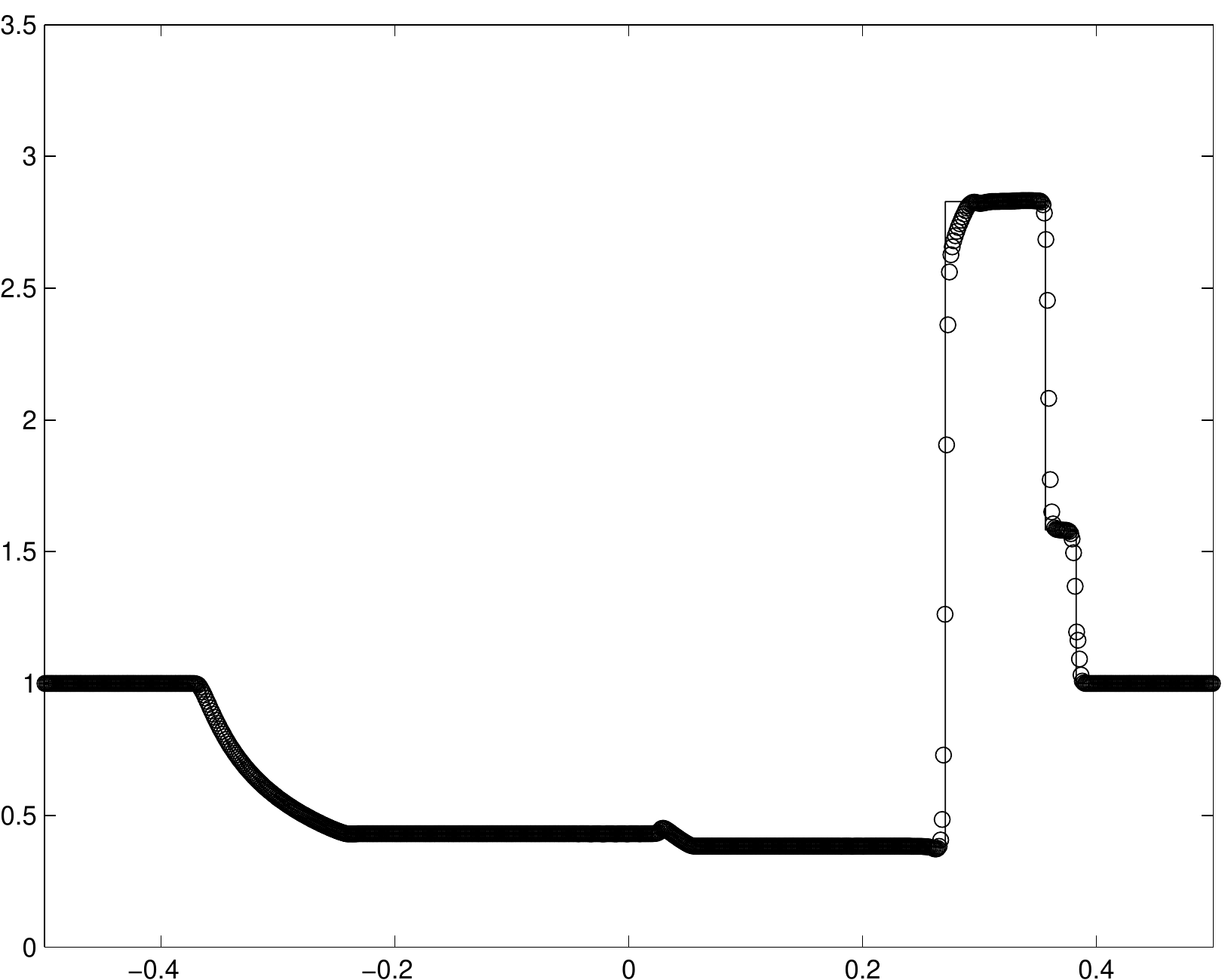}\\
\includegraphics[width=0.35\textwidth]{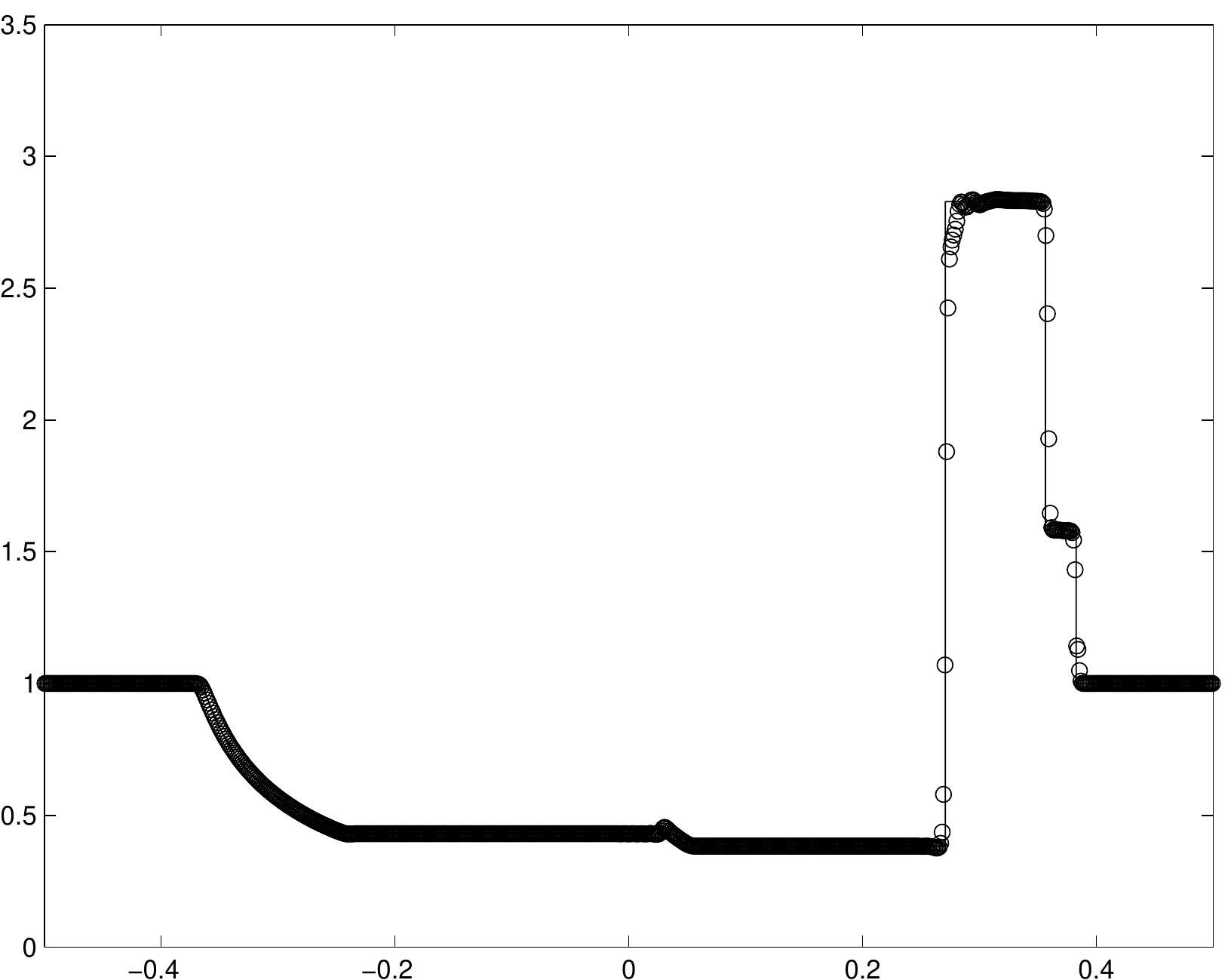}&
\includegraphics[width=0.35\textwidth]{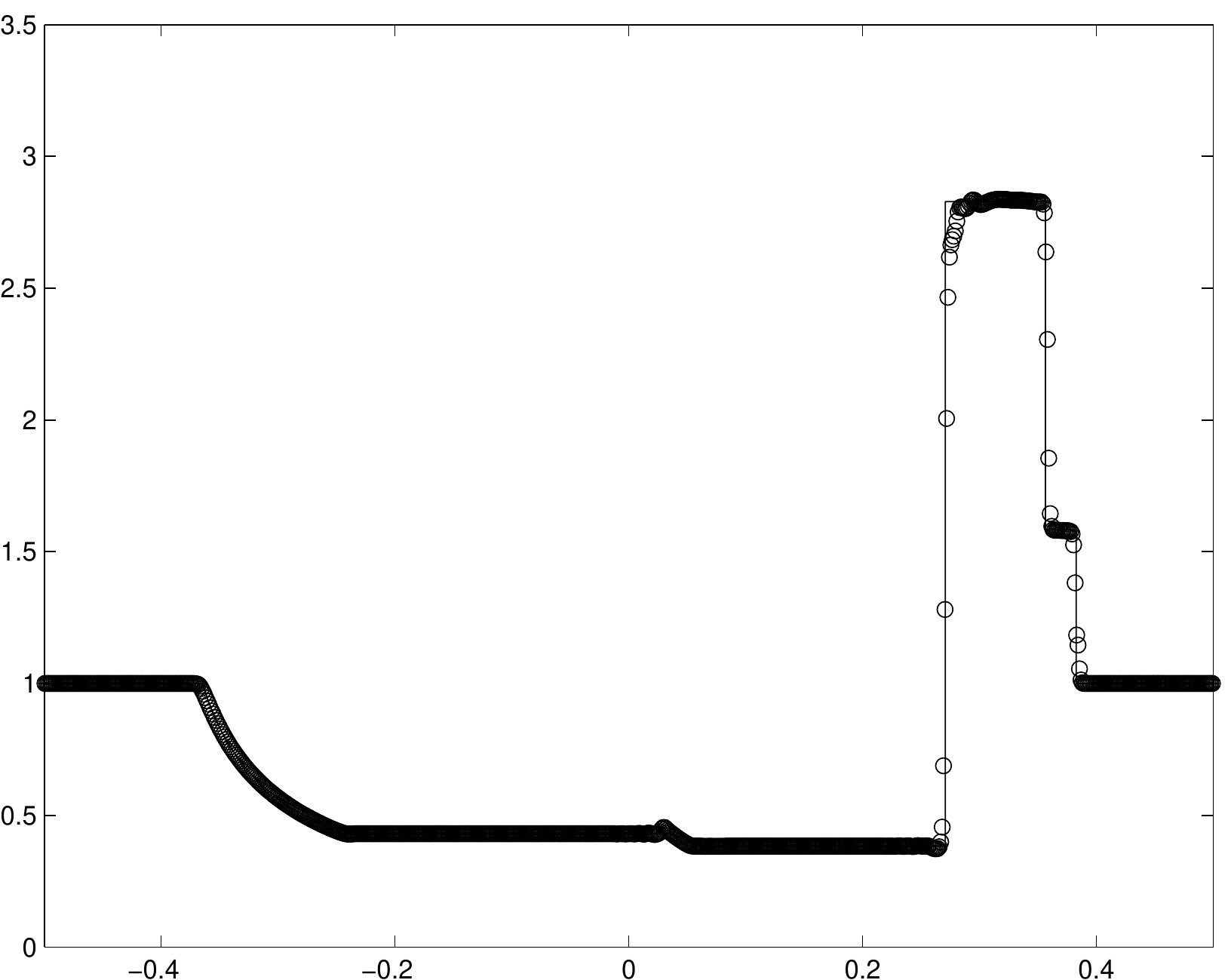}\\
    \end{tabular}
    \caption{Example ~\ref{exRMHDRMT3}: The densities $\rho$ at $t=0.4$.
    The solid line denotes the exact solution, while the symbol ``$\circ$'' is numerical solution
   obtained with $800$ cells. Left: $P^K$-based non-central \DG{}; right: $P^K$-based \CDG{}.
      From top to bottom: $K=1,~2,~3$.
   }
    \label{fig:RMHDRMT3rho}
  \end{figure}

\begin{figure}[!htbp]
    \centering{}
  \begin{tabular}{cc}
    \includegraphics[width=0.35\textwidth]{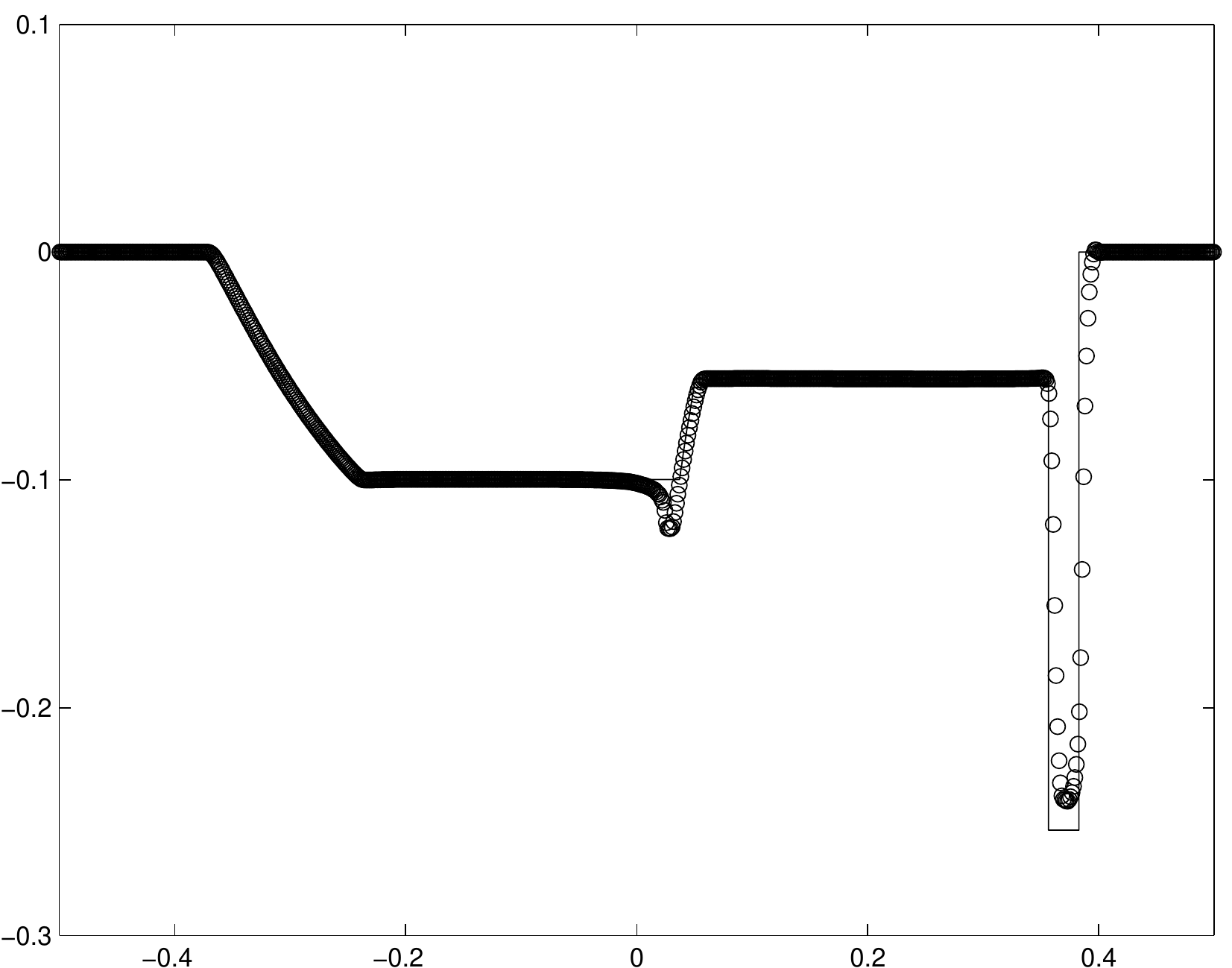}&
\includegraphics[width=0.35\textwidth]{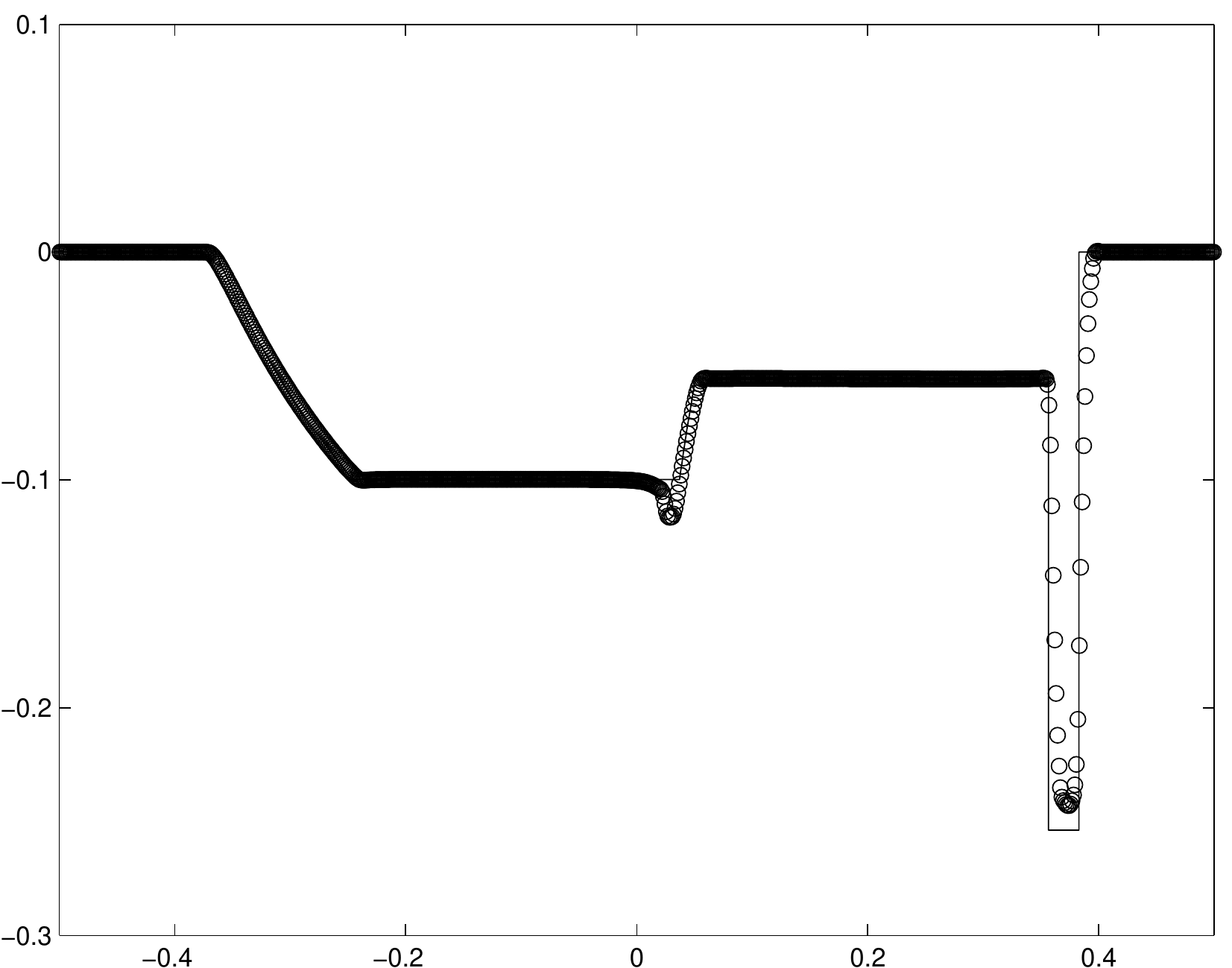}\\
\includegraphics[width=0.35\textwidth]{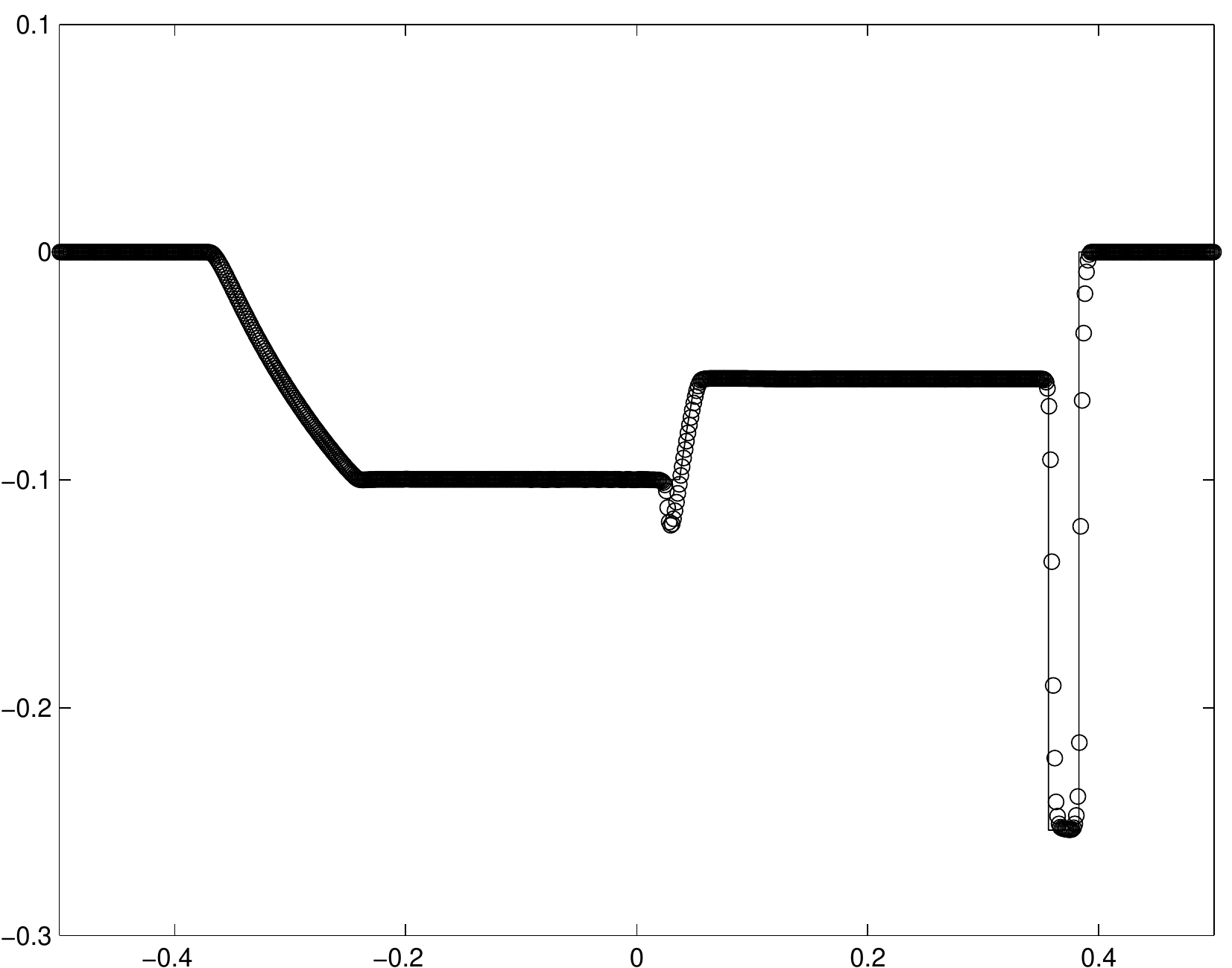}&
\includegraphics[width=0.35\textwidth]{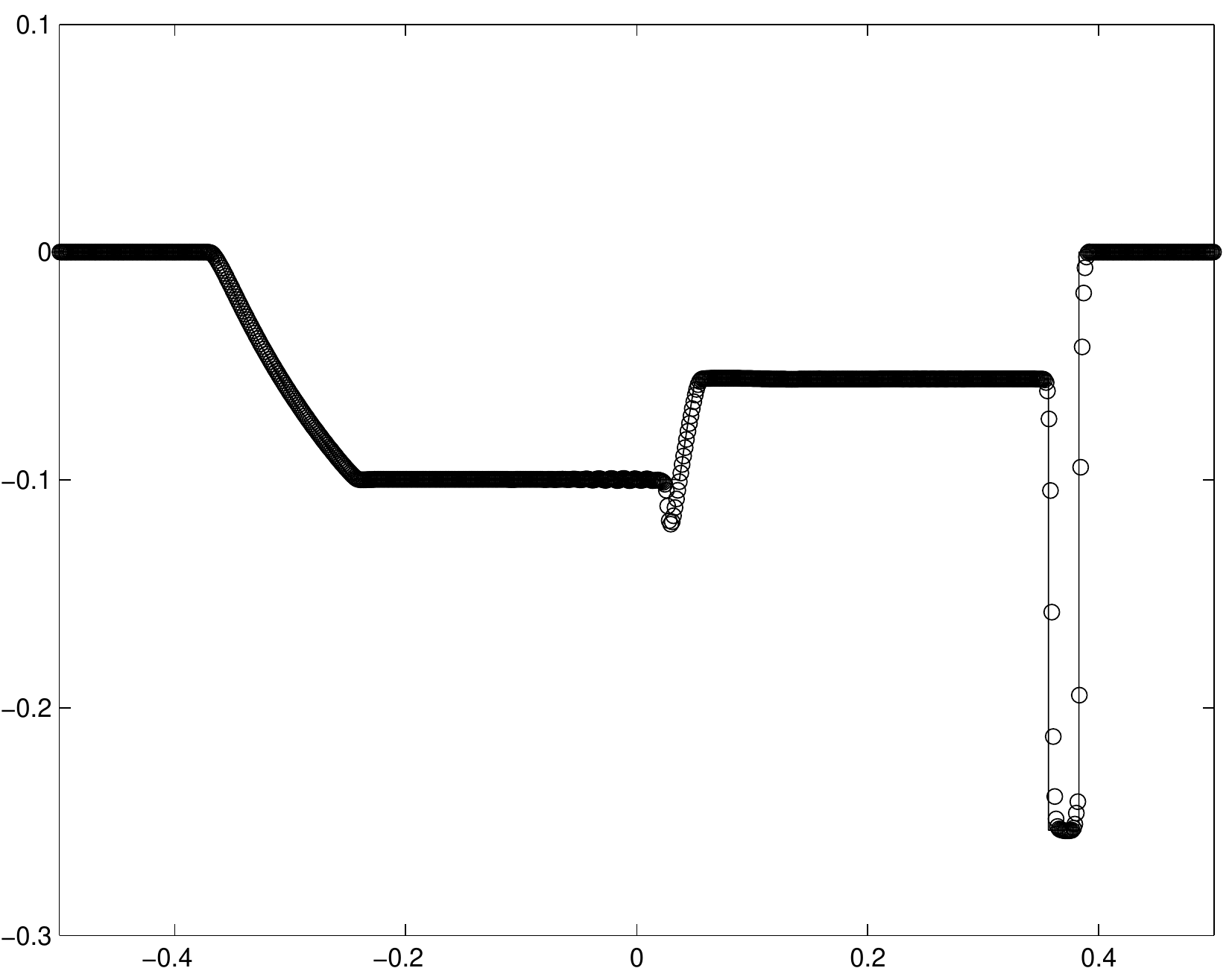}\\
\includegraphics[width=0.35\textwidth]{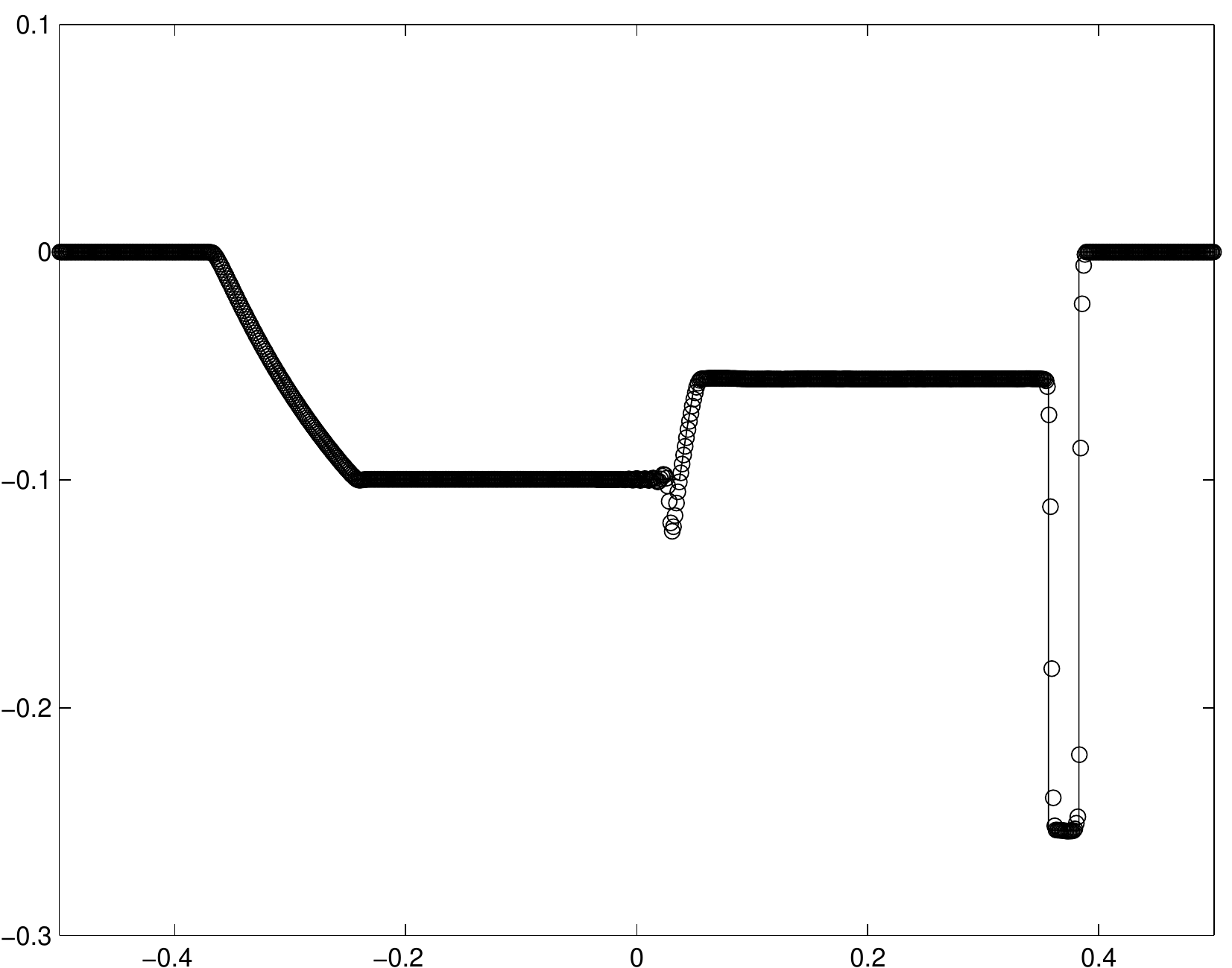}&
\includegraphics[width=0.35\textwidth]{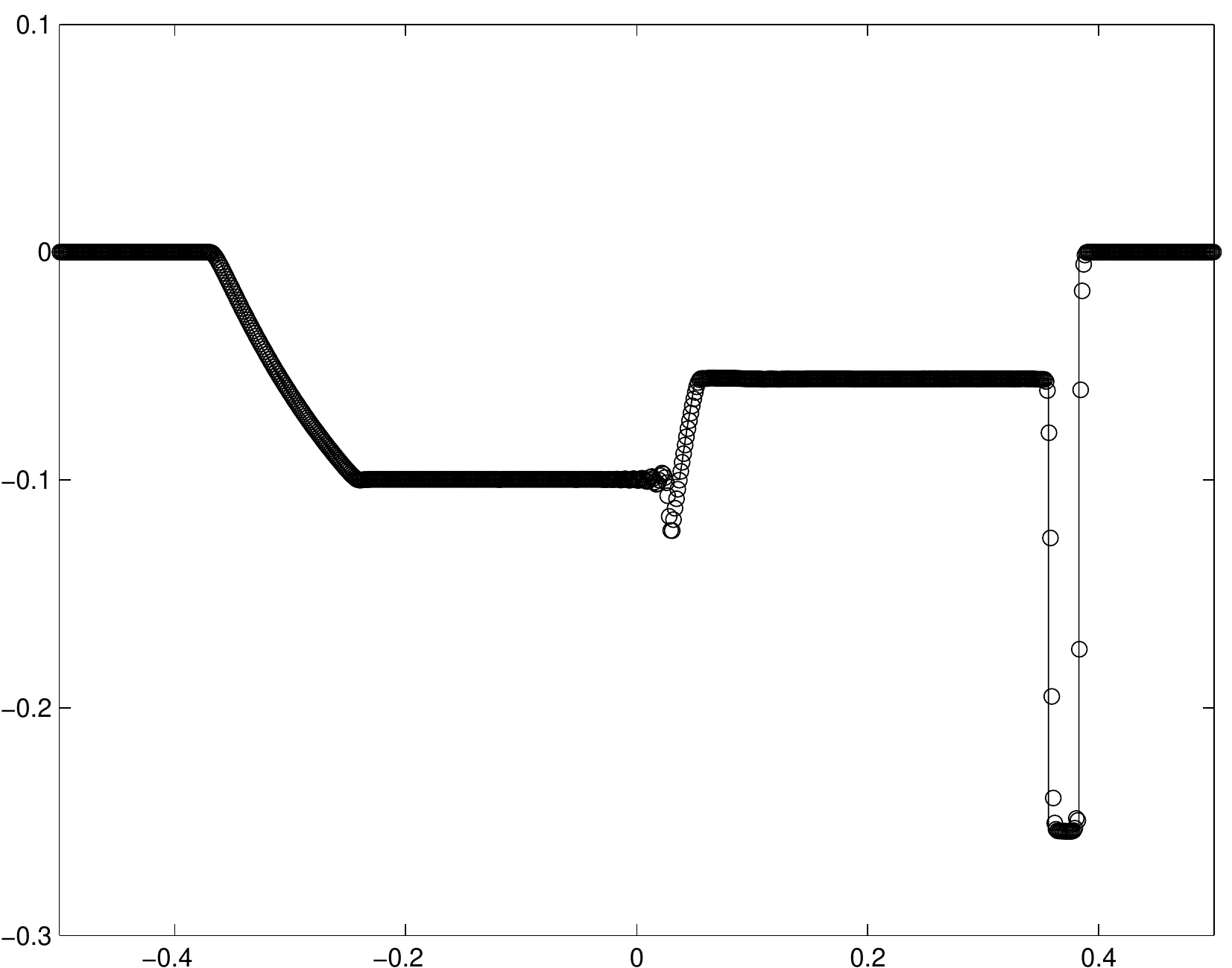}\\
    \end{tabular}
        \caption{Same as Fig.~\ref{fig:RMHDRMT3rho} except for $v_y$.}

    \label{fig:RMHDRMT3vy}
  \end{figure}

\begin{figure}[!htbp]
    \centering{}
  \begin{tabular}{cc}
    \includegraphics[width=0.35\textwidth]{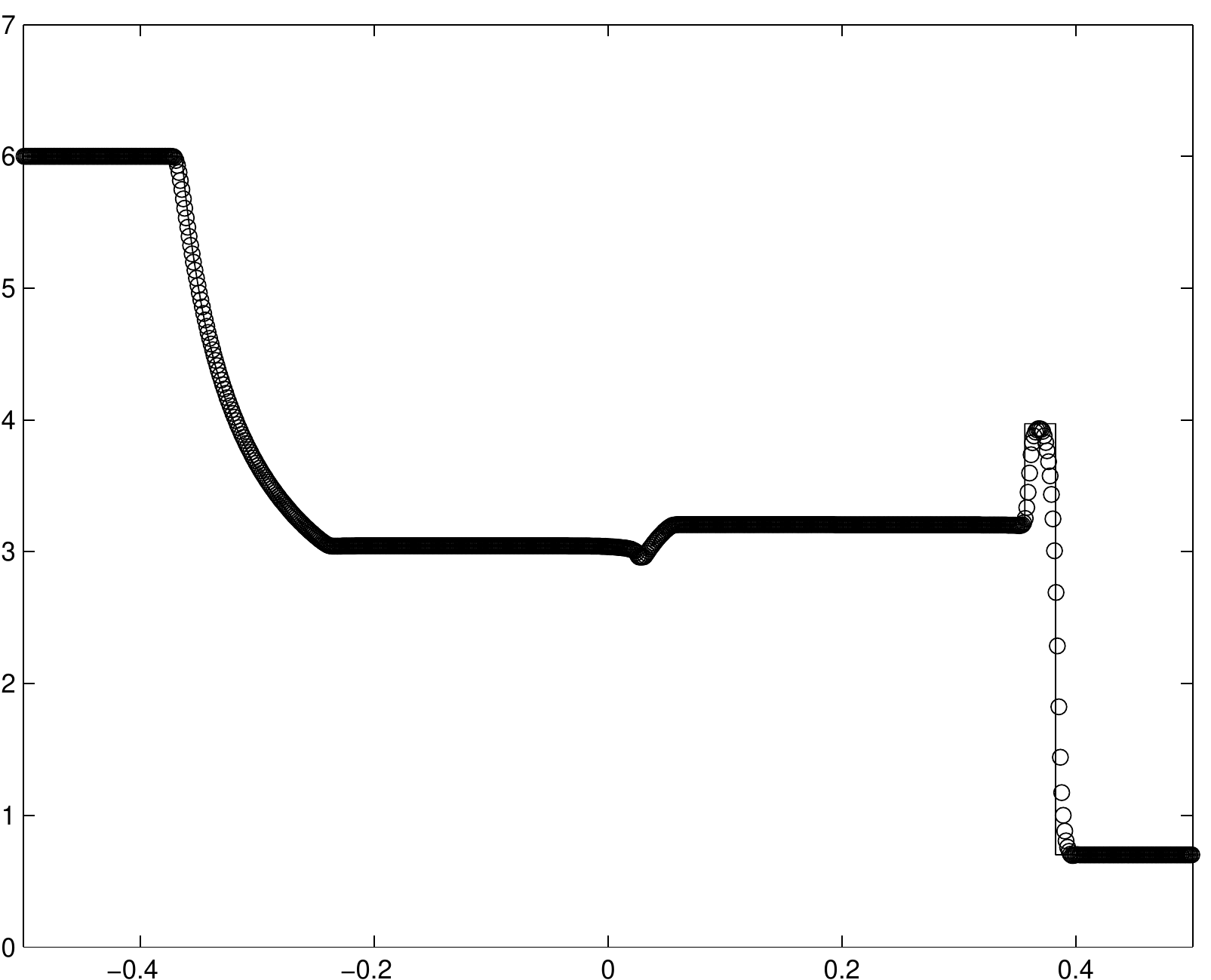}&
\includegraphics[width=0.35\textwidth]{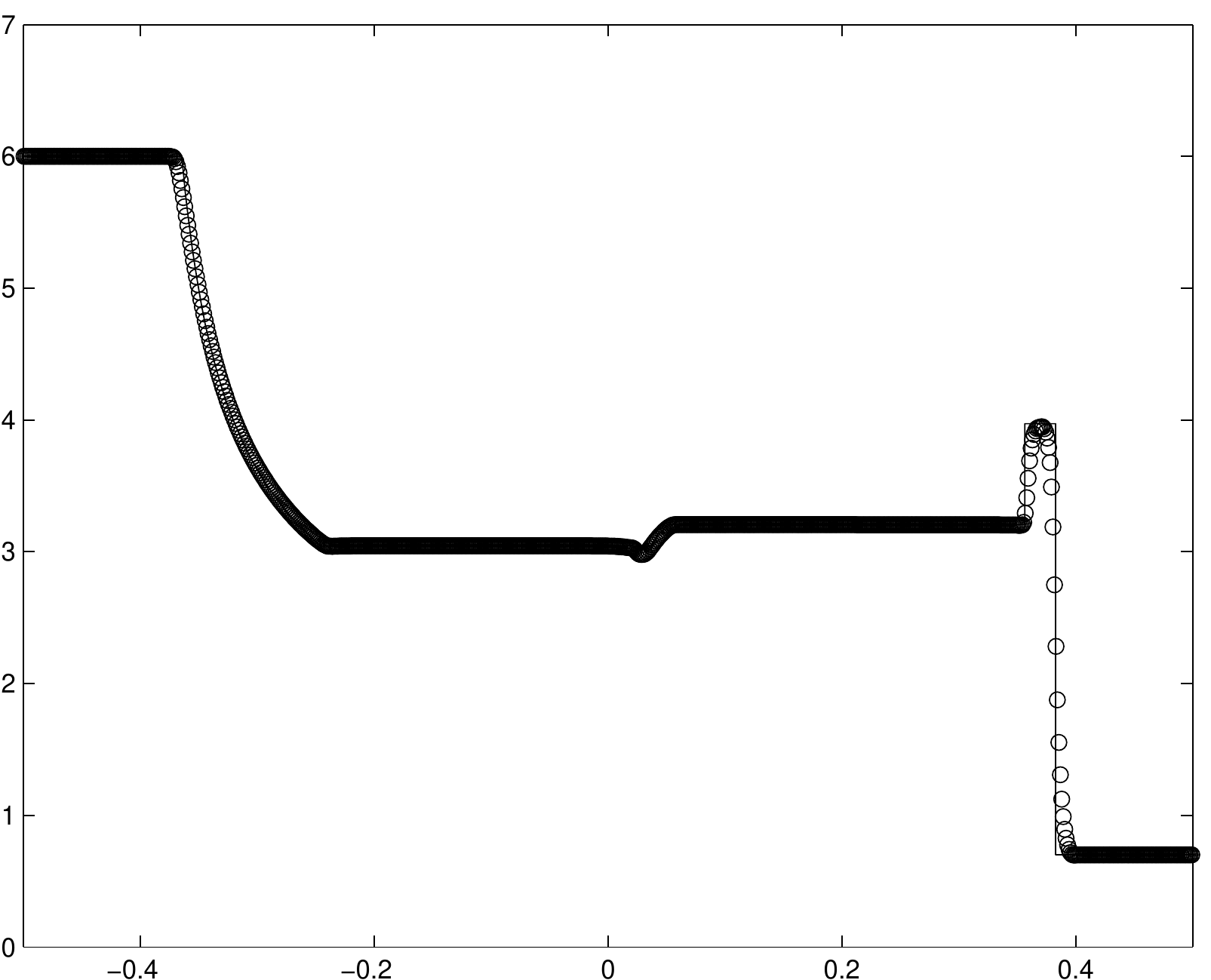}\\
\includegraphics[width=0.35\textwidth]{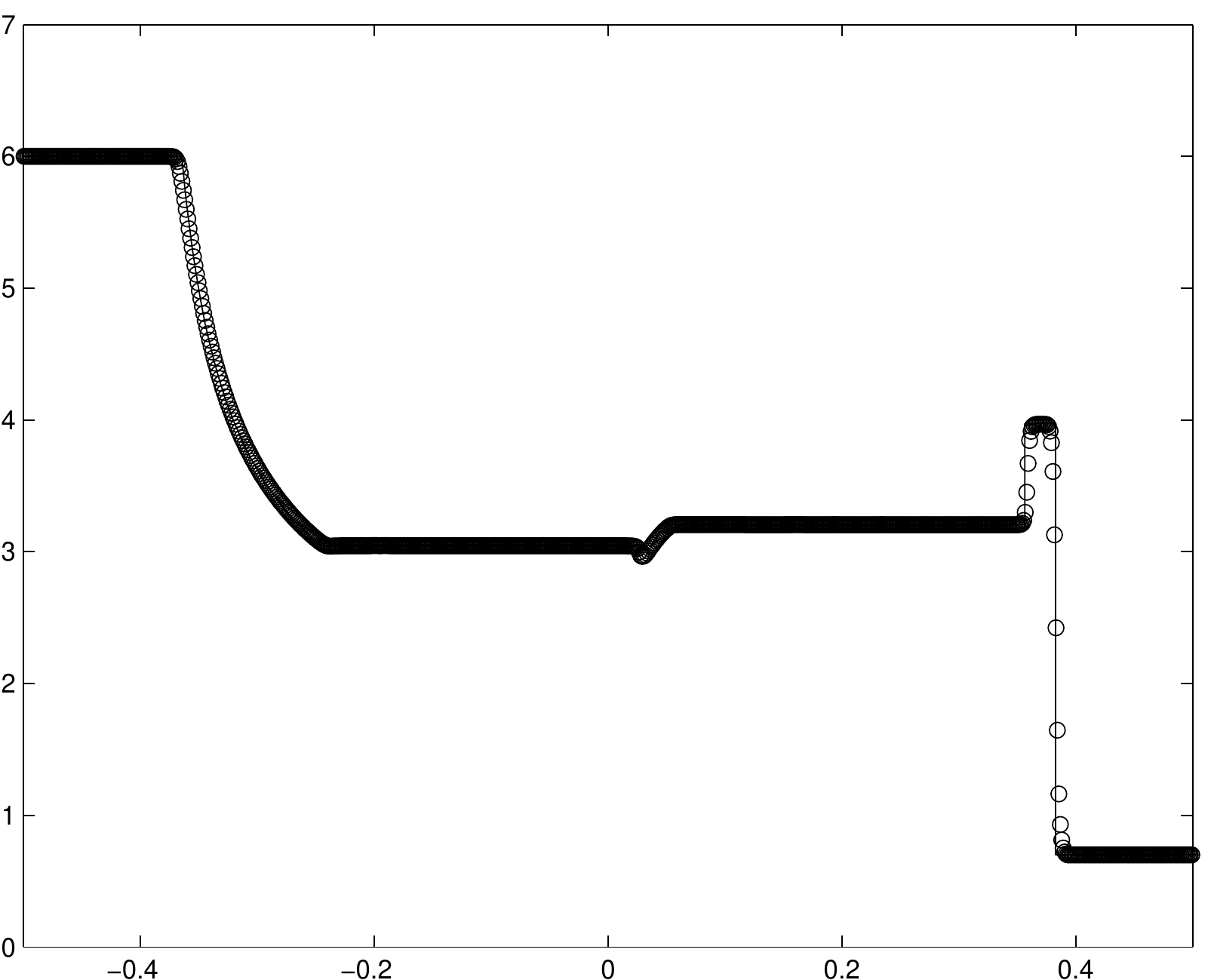}&
\includegraphics[width=0.35\textwidth]{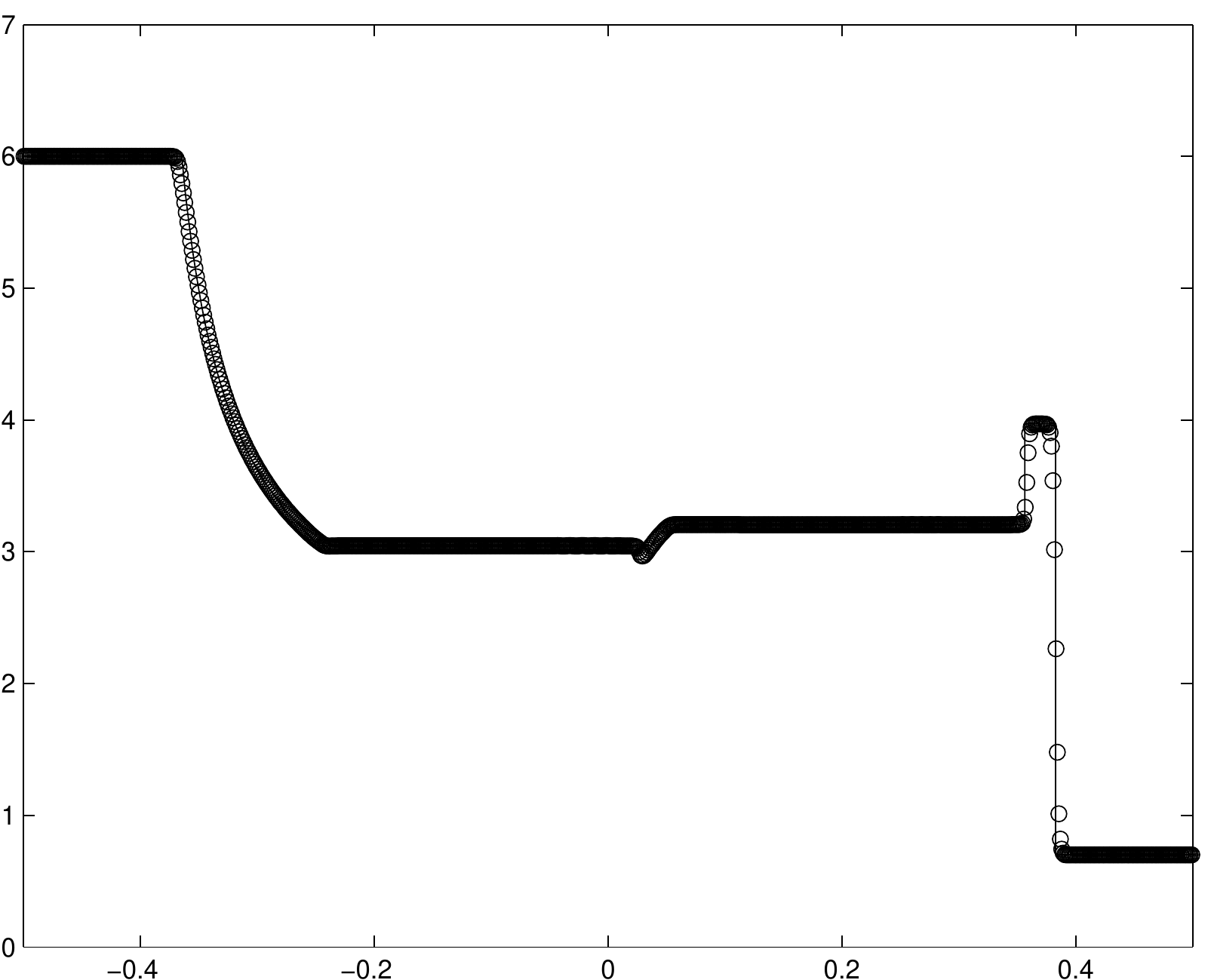}\\
\includegraphics[width=0.35\textwidth]{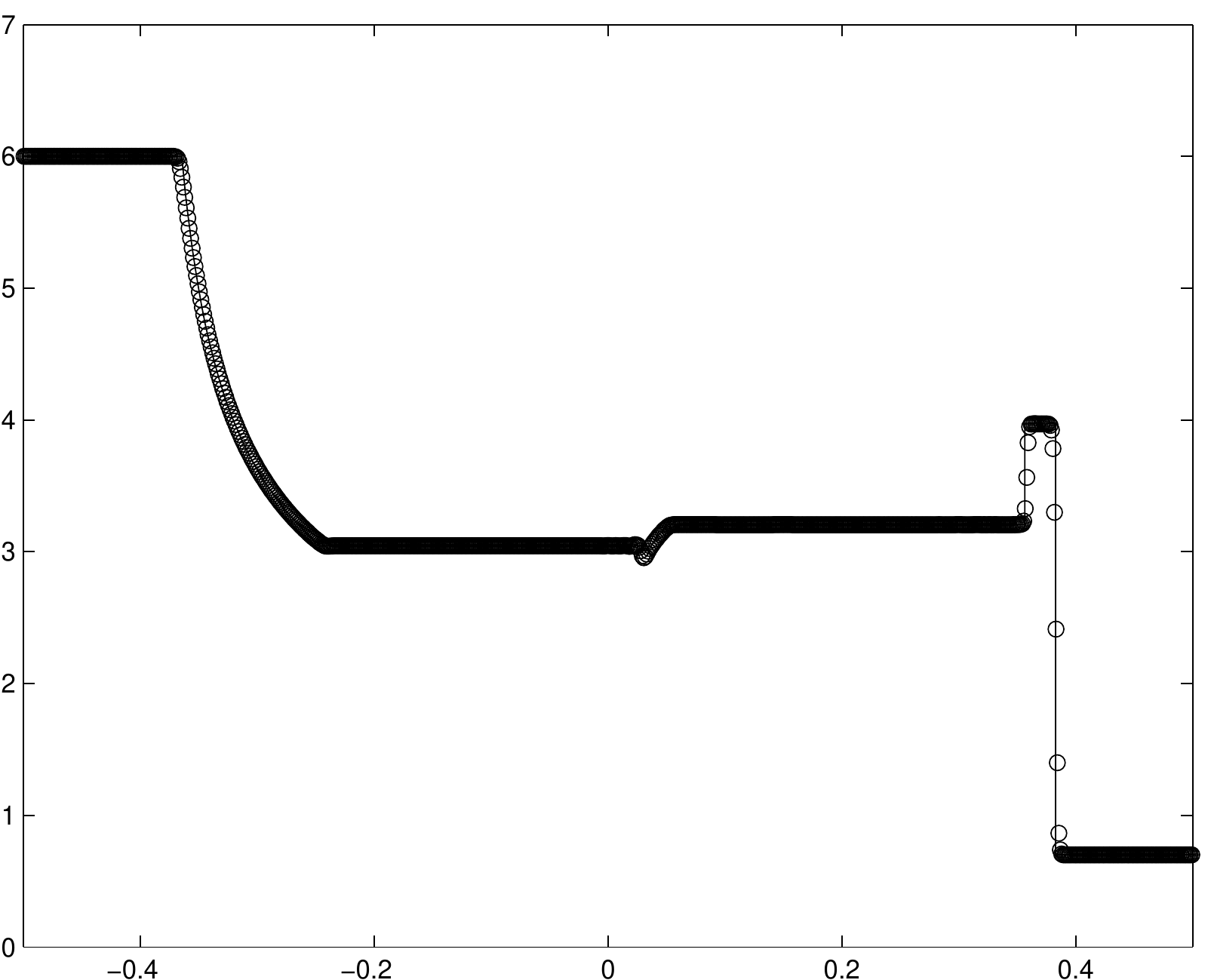}&
\includegraphics[width=0.35\textwidth]{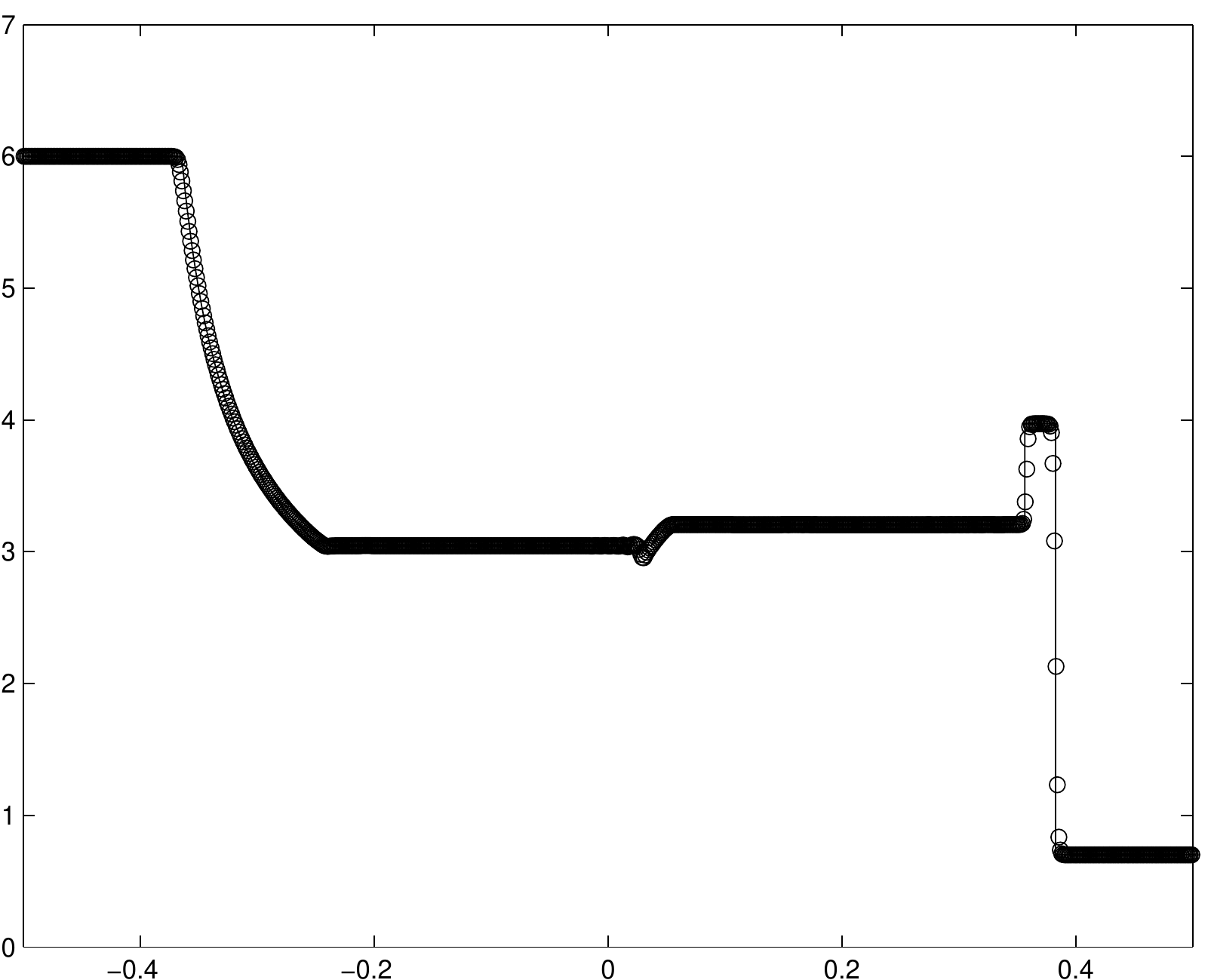}\\
    \end{tabular}
    \caption{Same as Fig.~\ref{fig:RMHDRMT3rho} except for $B_y$.}
    \label{fig:RMHDRMT3by}
  \end{figure}

\begin{figure}[!htbp]
    \centering{}
  \begin{tabular}{cc}
    \includegraphics[width=0.35\textwidth]{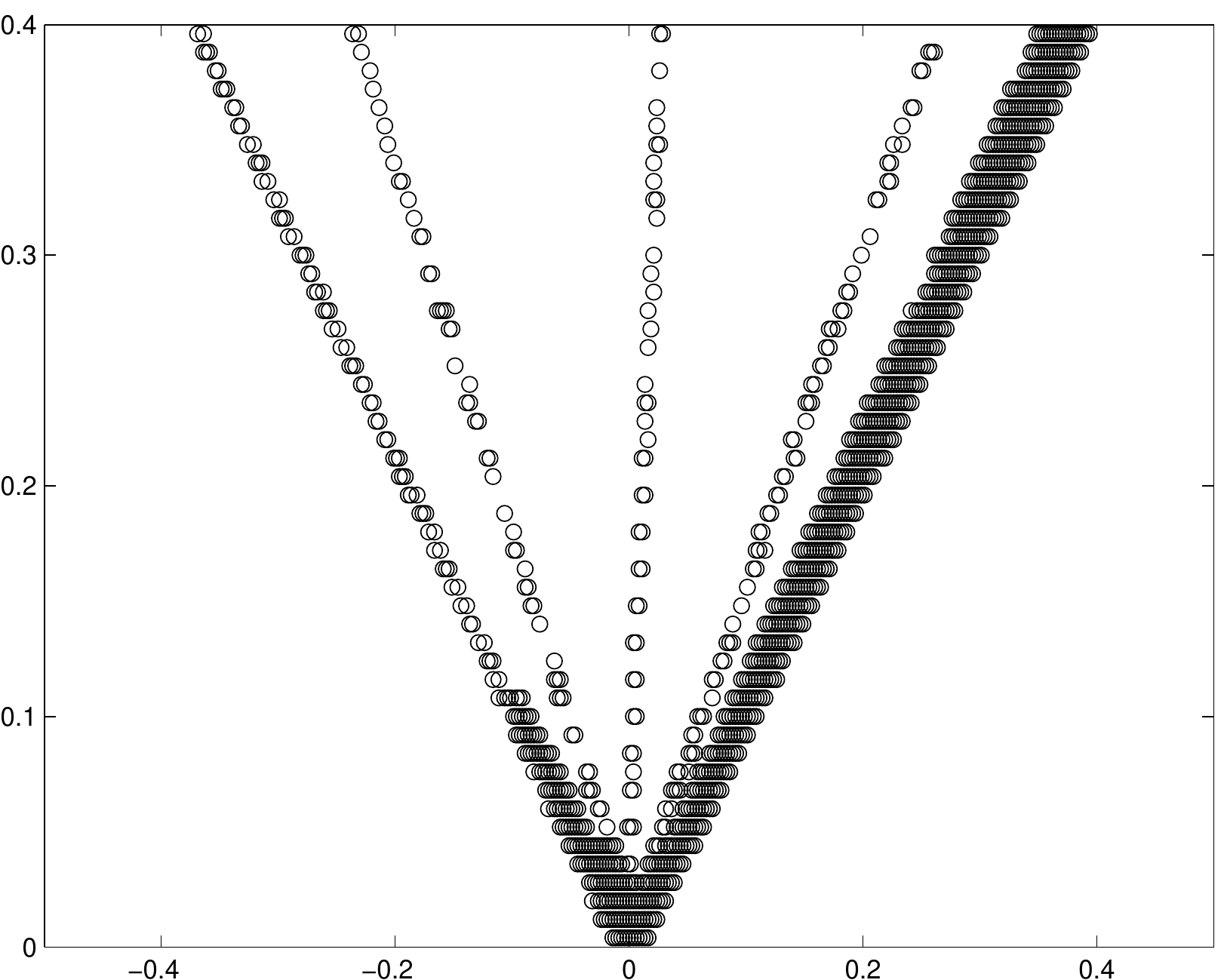}&
\includegraphics[width=0.35\textwidth]{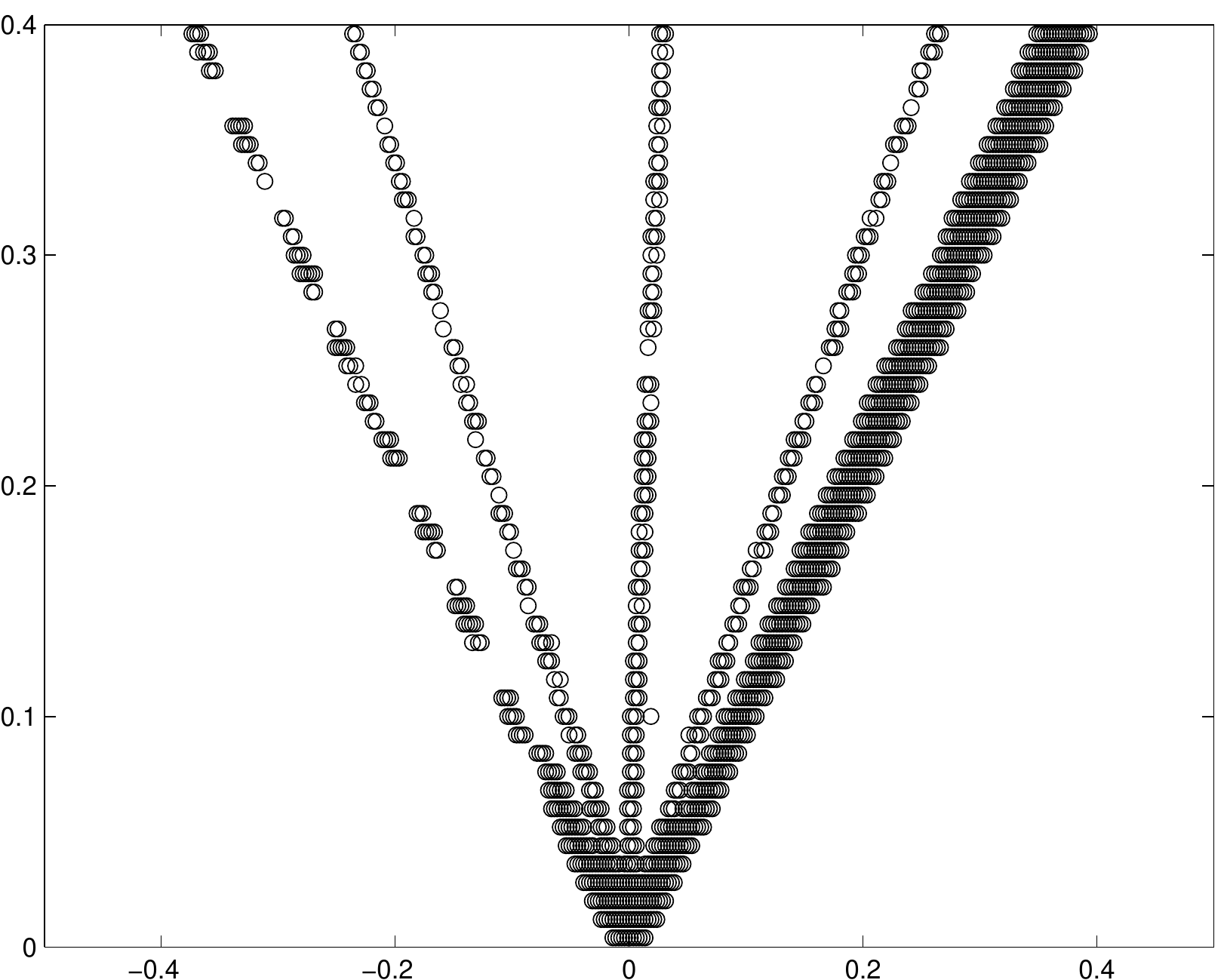}\\
\includegraphics[width=0.35\textwidth]{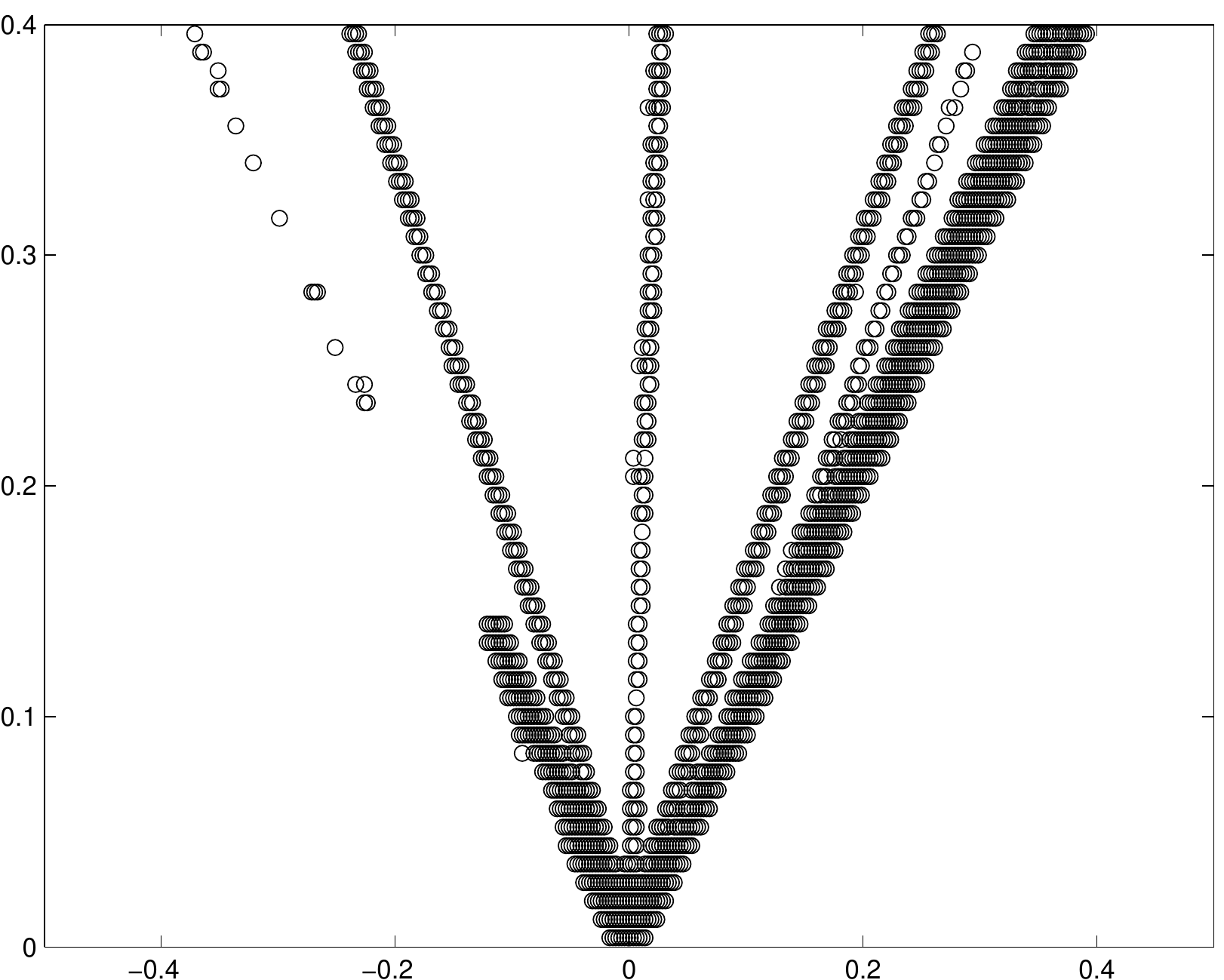}&
\includegraphics[width=0.35\textwidth]{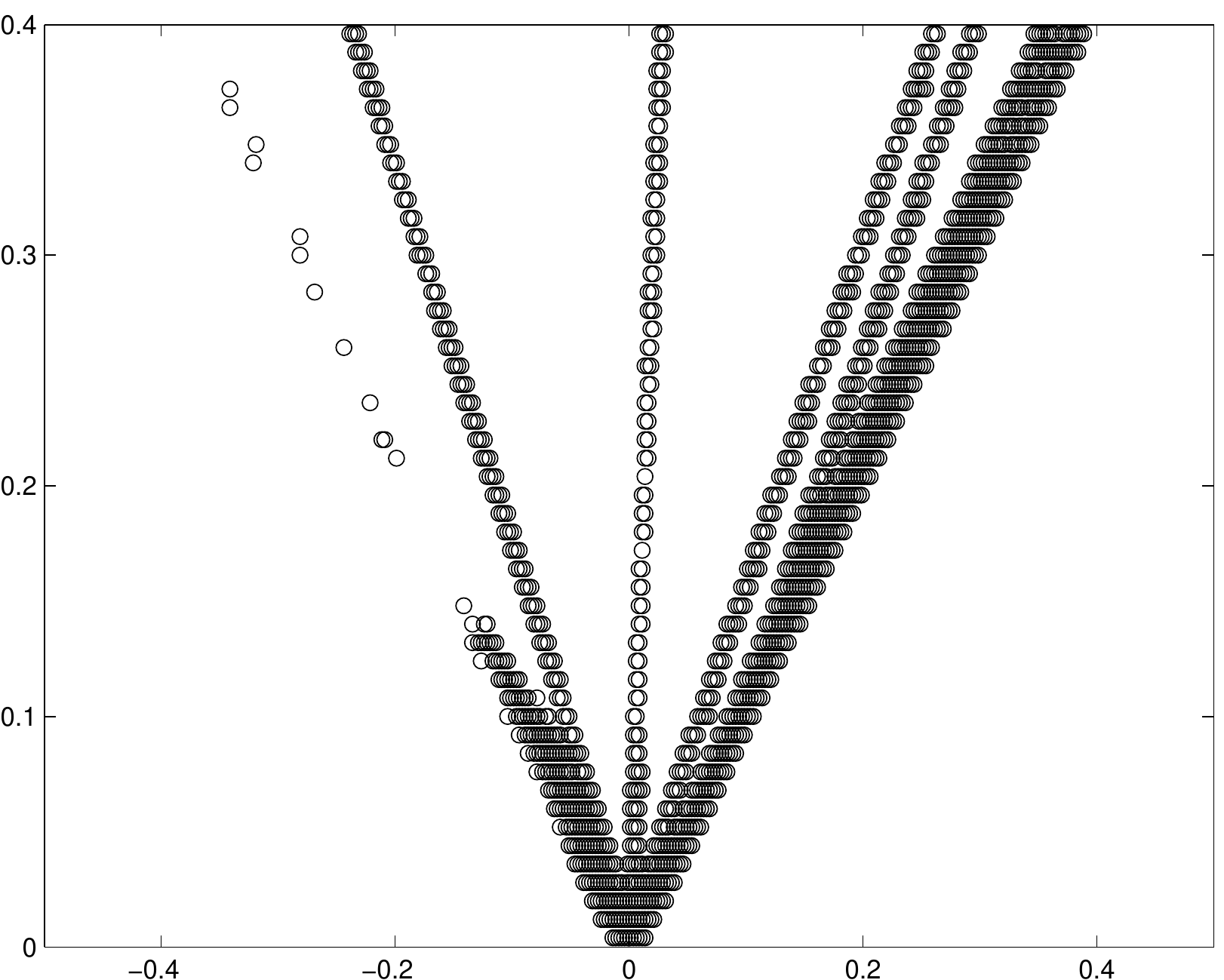}\\
\includegraphics[width=0.35\textwidth]{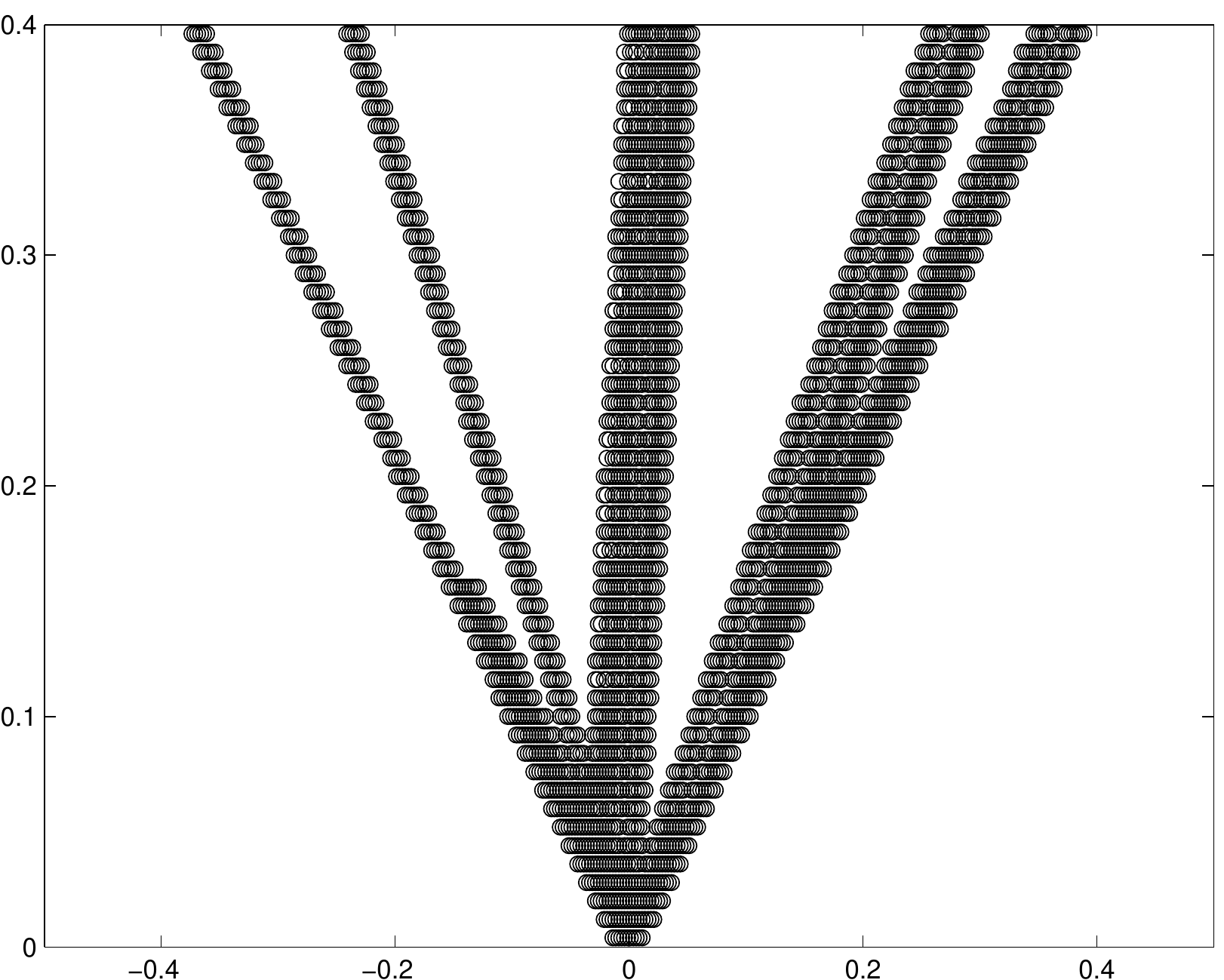}&
\includegraphics[width=0.35\textwidth]{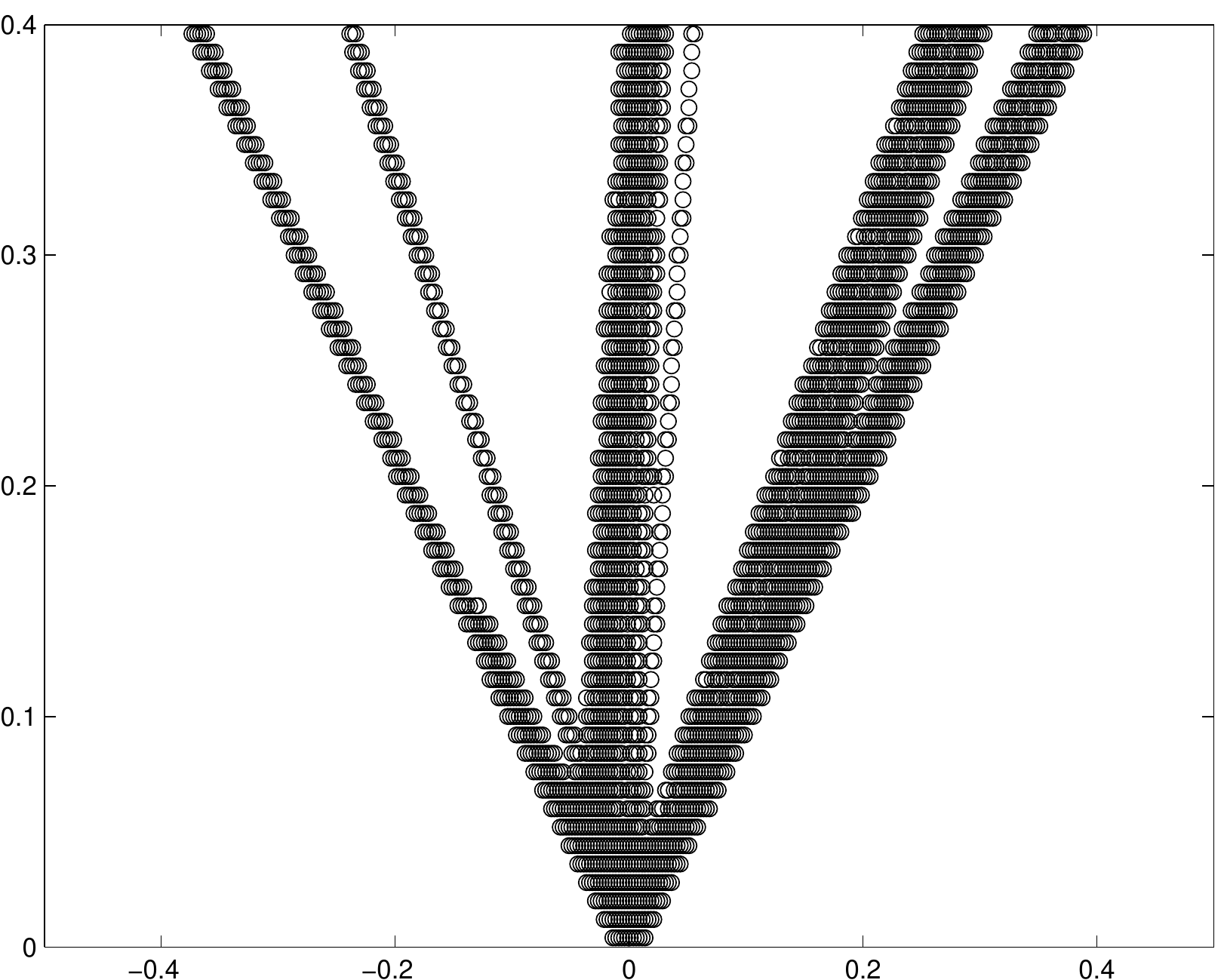}\\
    \end{tabular}
    \caption{Example ~\ref{exRMHDRMT3}: The time evolution of ``troubled'' cells.
       Left: non-central \DG{}; right: \CDG{}.
      From top to bottom: $K=1,~2,~3$. The cell number   is $800$.}
    \label{fig:RMHDRMT3cell}
  \end{figure}

   Figs~\ref{fig:RMHDRMT3rho}, \ref{fig:RMHDRMT3vy}, and \ref{fig:RMHDRMT3by} display
   the densities $\rho$, velocities $v_y$, and magnetic fields $B_y$ at
   $t=0.4$.
It can be seen from the plots that high order methods resolves the contact discontinuity in the density
 and   two right-moving shock waves better than the lower order method,
and the non-central and central DG methods  have essentially the same resolution.
   The ``troubled'' cells identified by the present DG methods are very finite, and
   mainly appear in  the region where the discontinuities in solution are  relatively strong, see
Fig.~\ref{fig:RMHDRMT3cell}.

\begin{figure}[!htbp]
    \centering{}
  \begin{tabular}{cc}

    \includegraphics[width=0.35\textwidth]{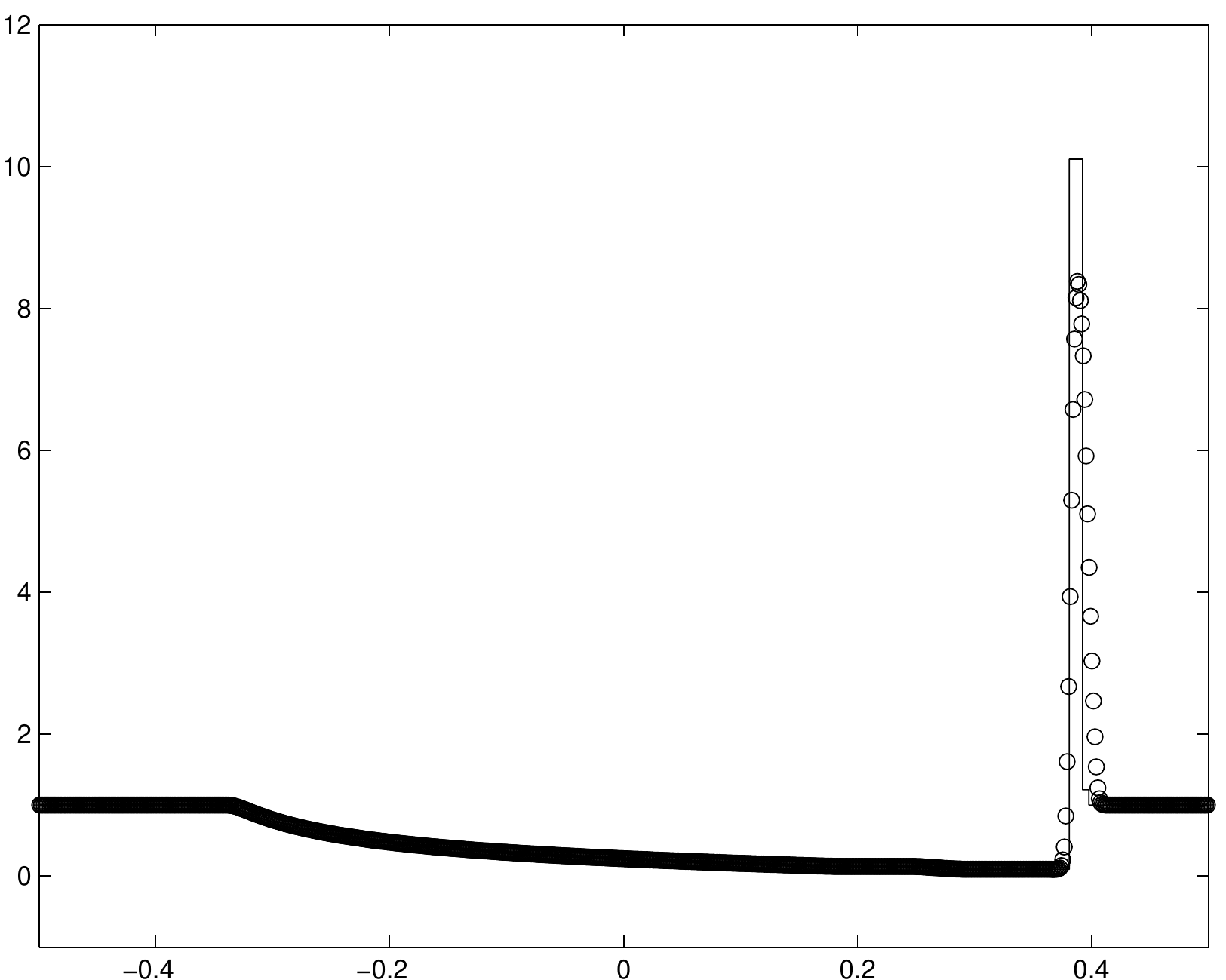}&
\includegraphics[width=0.35\textwidth]{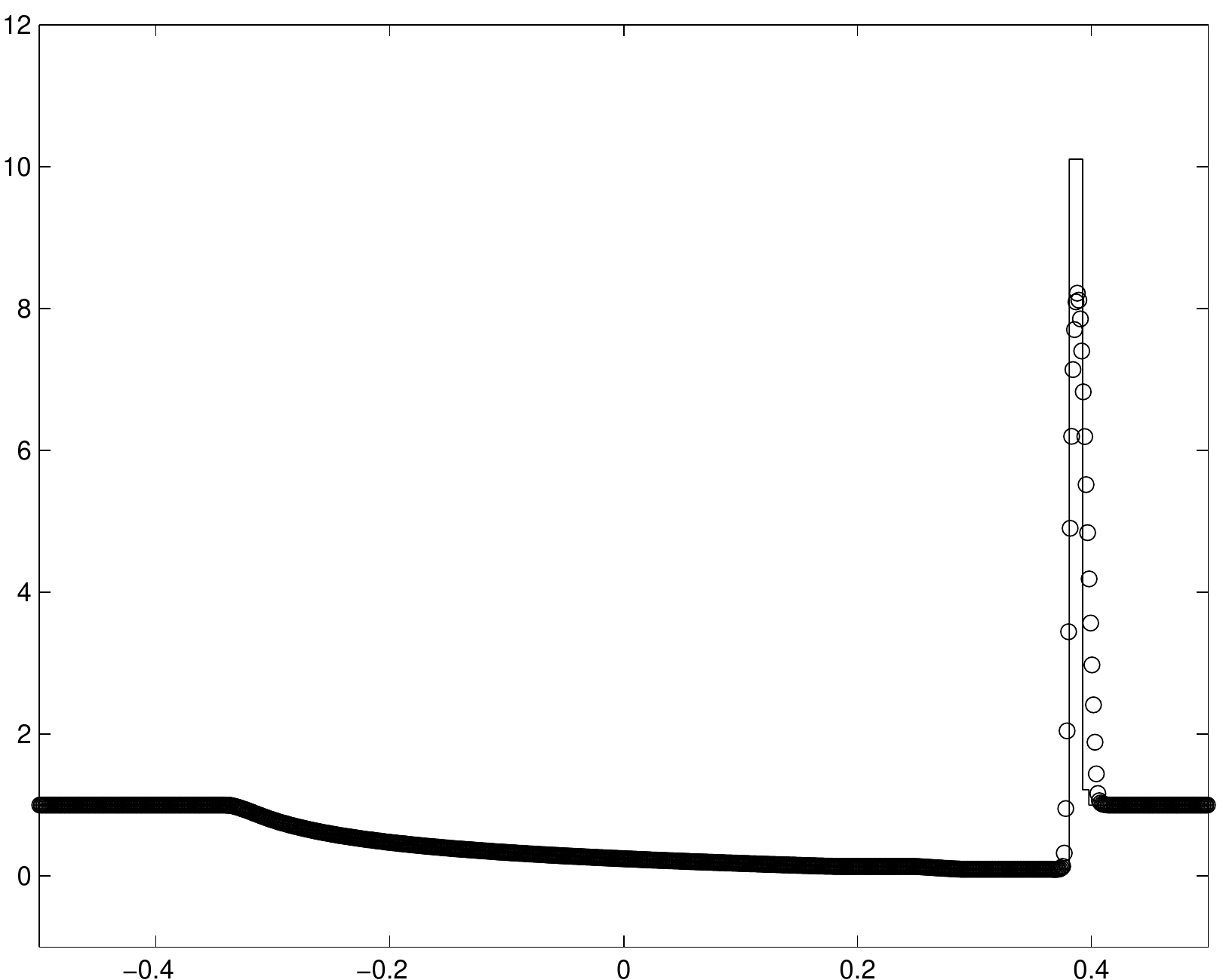}\\
\includegraphics[width=0.35\textwidth]{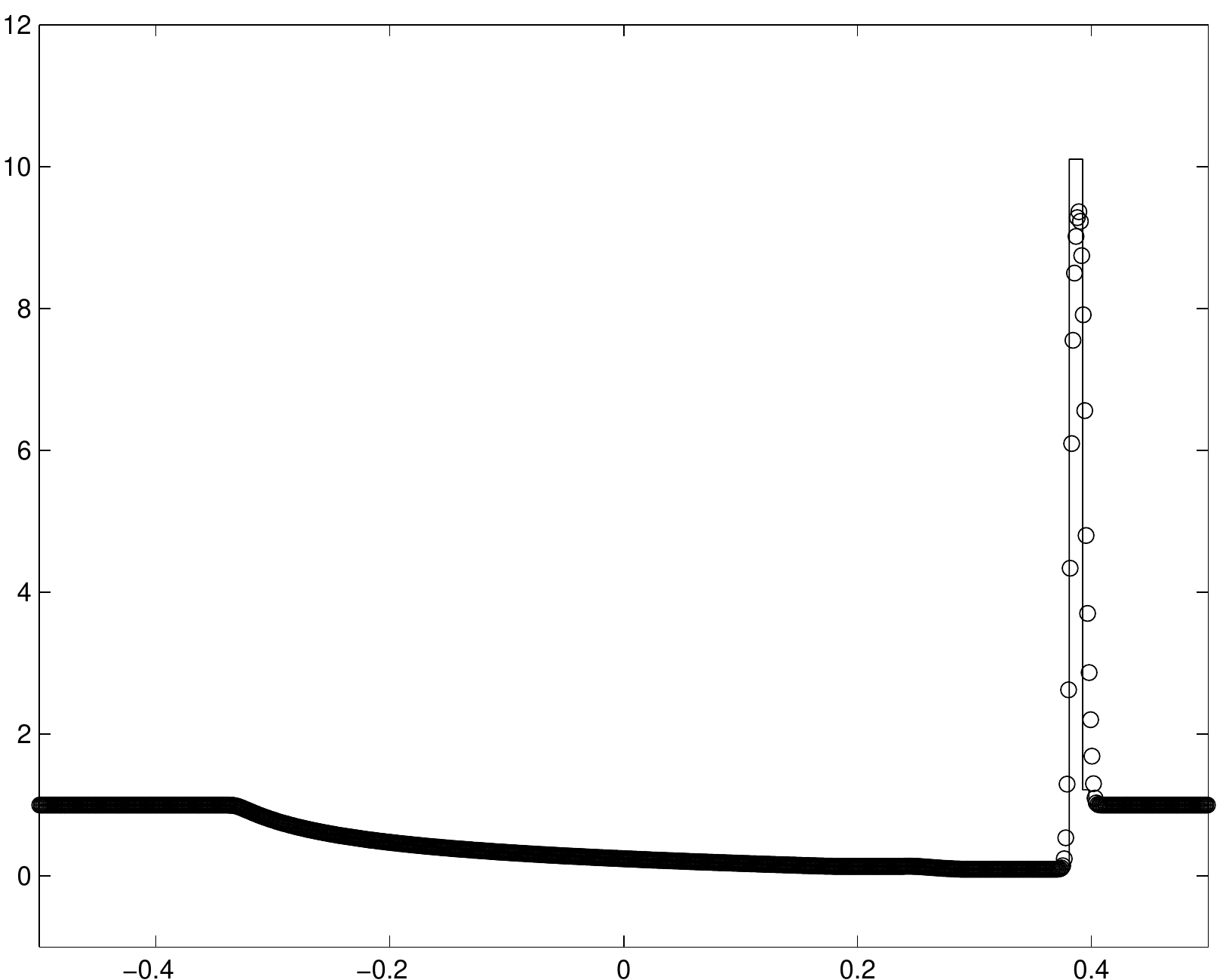}&
\includegraphics[width=0.35\textwidth]{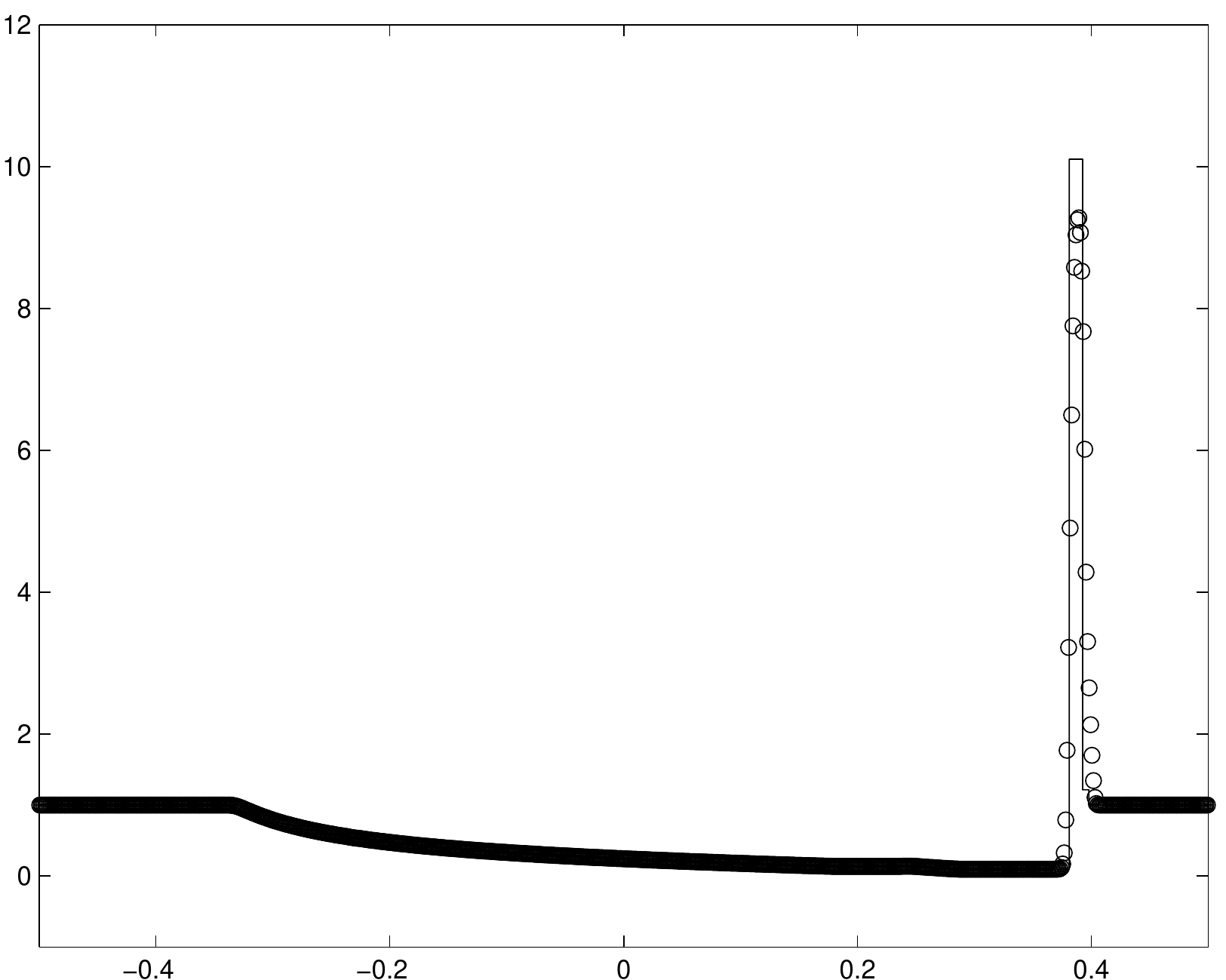}\\
\includegraphics[width=0.35\textwidth]{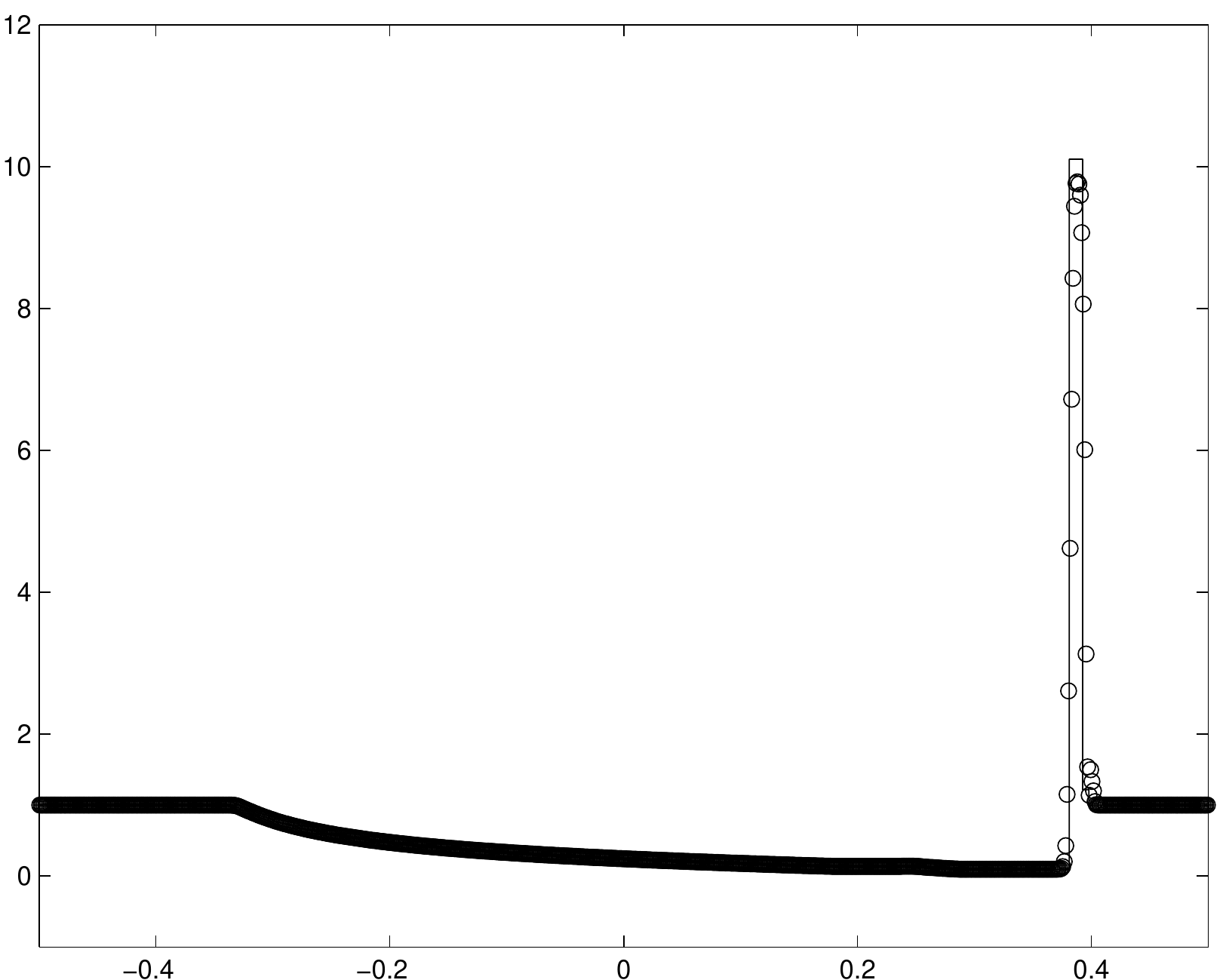}&
\includegraphics[width=0.35\textwidth]{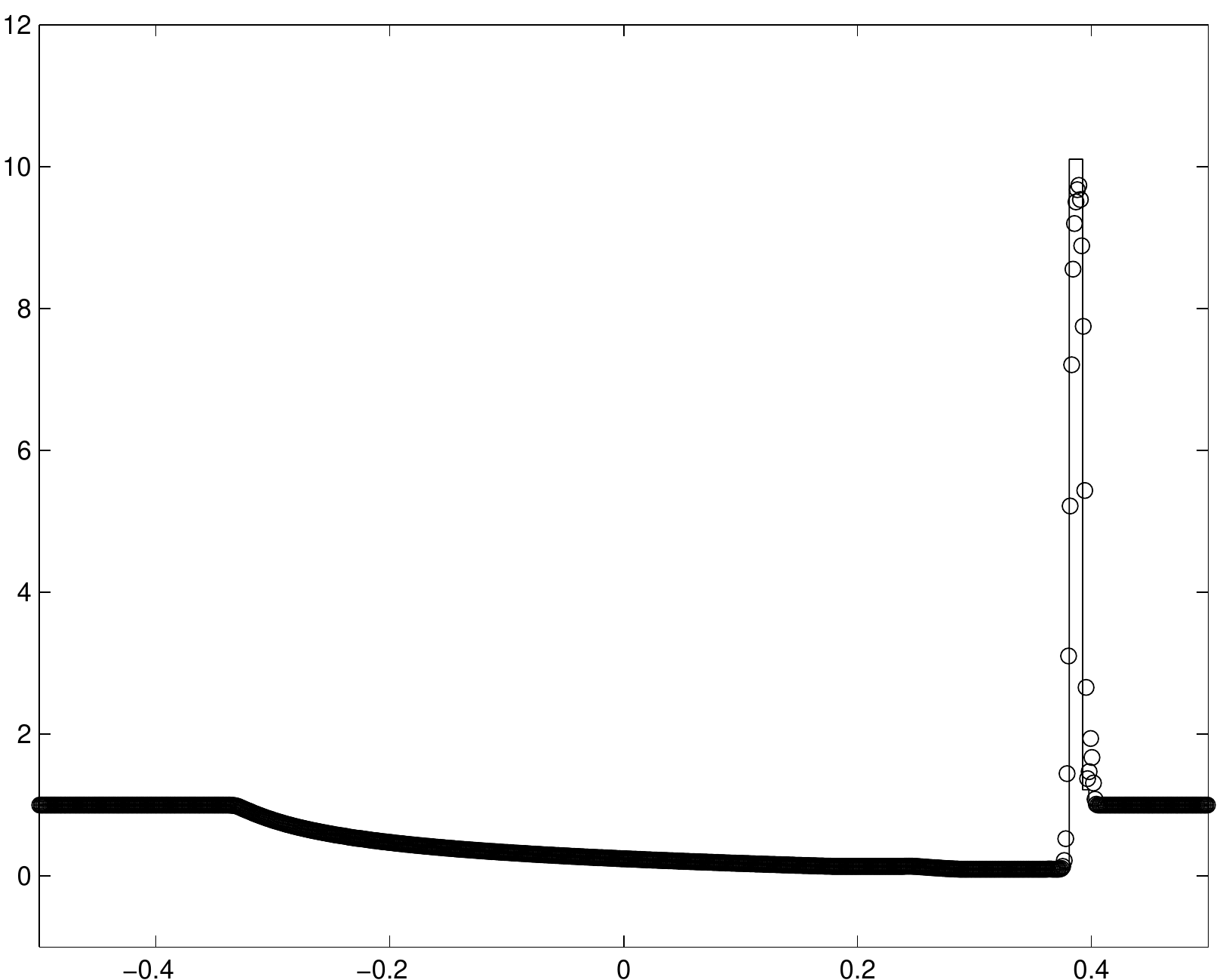}\\
    \end{tabular}
    \caption{Example ~\ref{exRMHDRMT4}: The densities $\rho$ at $t=0.4$.
    The solid line denote the exact solution, while the symbol ``$\circ$'' is numerical solution
   obtained with $800$ cells. Left: $P^K$-based non-central \DG{}; right: $P^K$-based \CDG{}.
      From top to bottom: $K=1,~2,~3$.  }
    \label{fig:RMHDRMT4rho}
  \end{figure}

\begin{figure}[!htbp]
    \centering{}
  \begin{tabular}{cc}
    \includegraphics[width=0.35\textwidth]{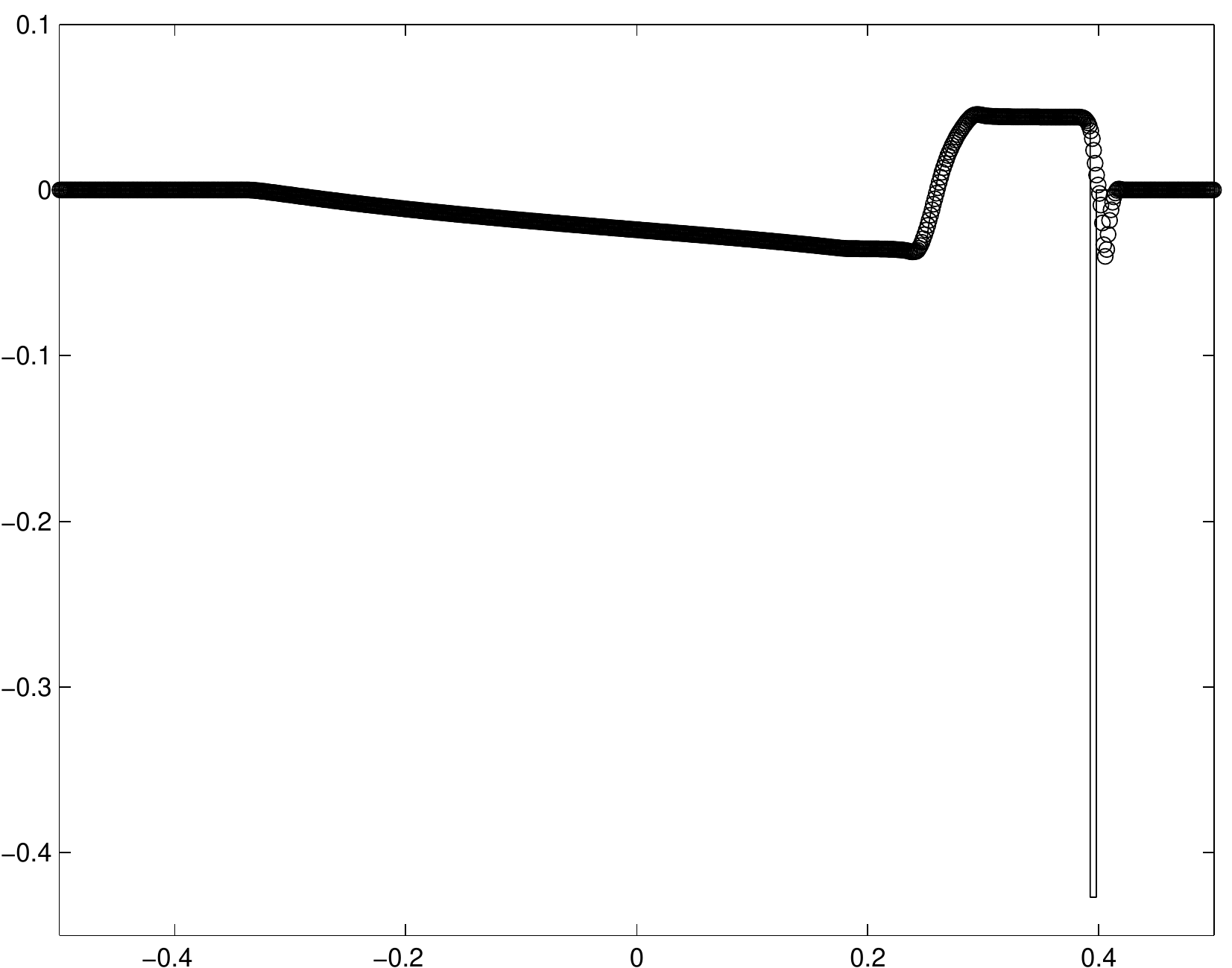}&
\includegraphics[width=0.35\textwidth]{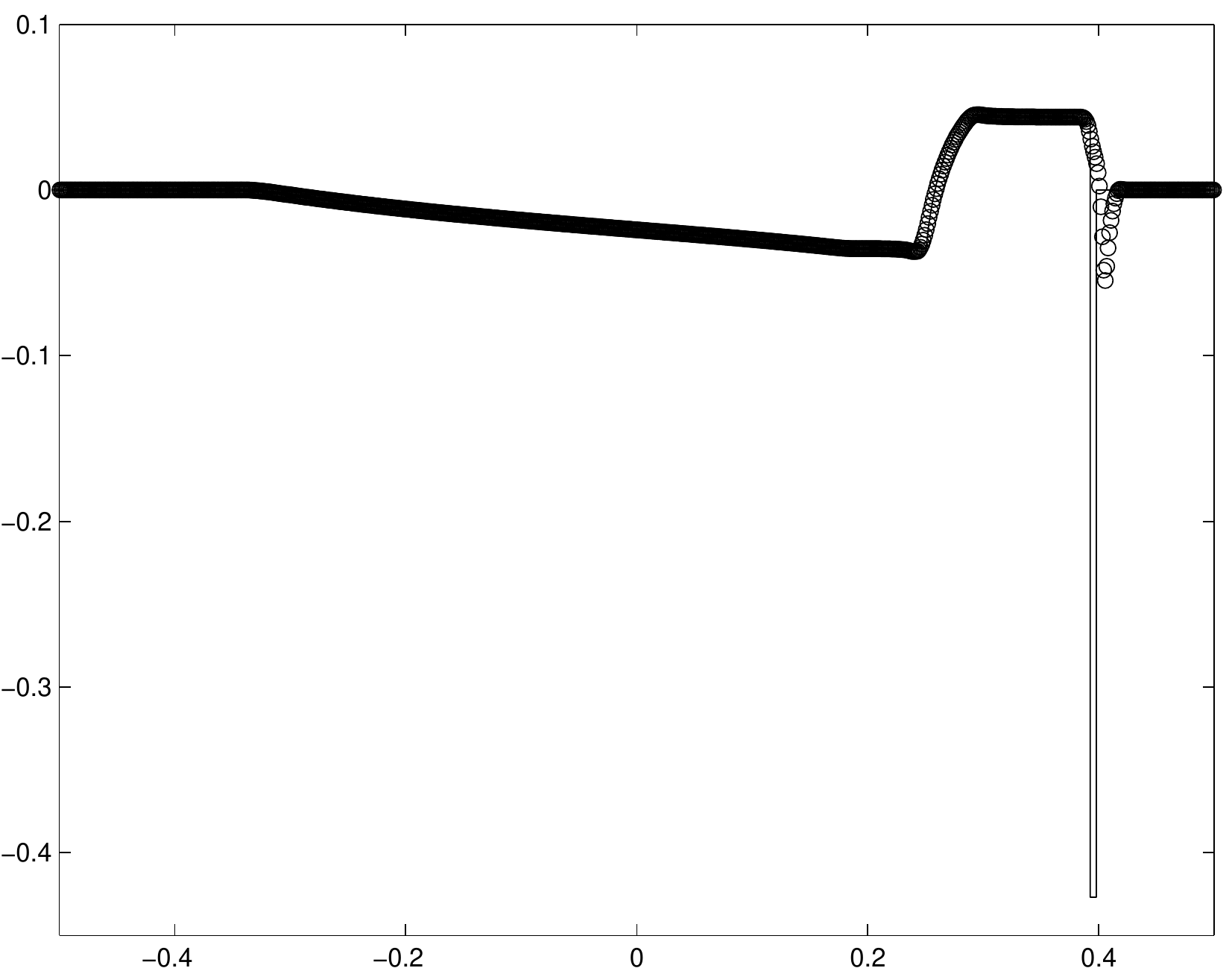}\\
\includegraphics[width=0.35\textwidth]{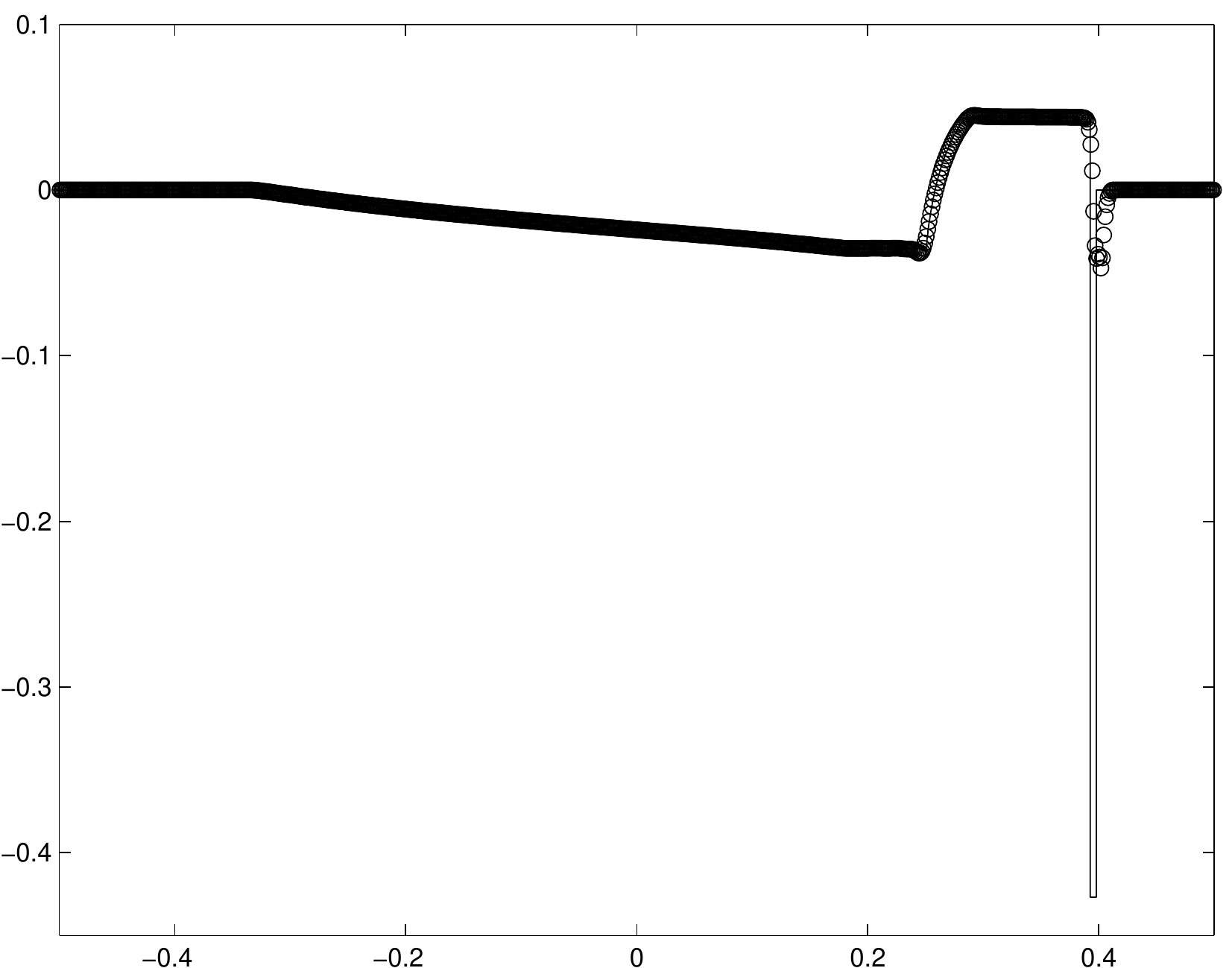}&
\includegraphics[width=0.35\textwidth]{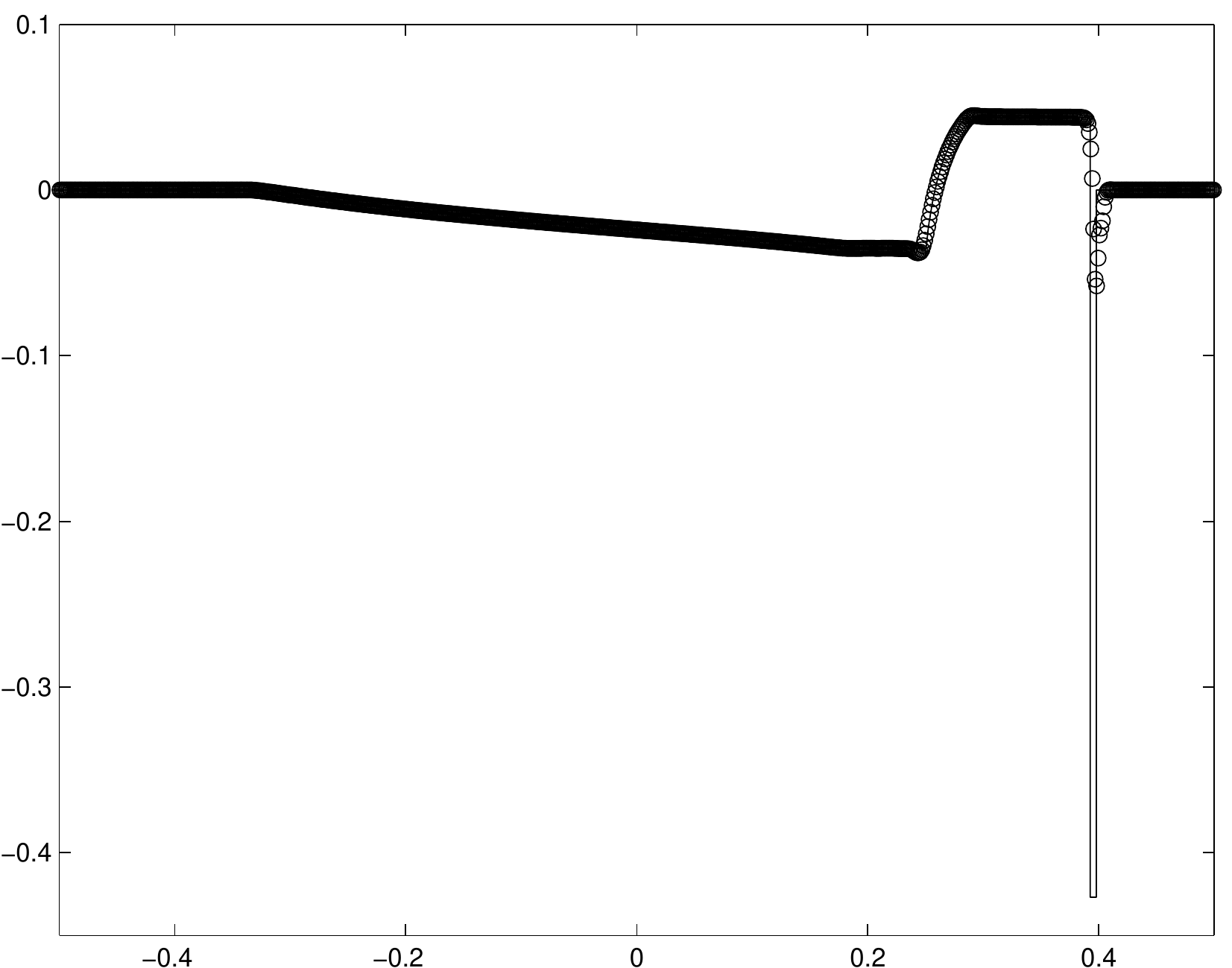}\\
\includegraphics[width=0.35\textwidth]{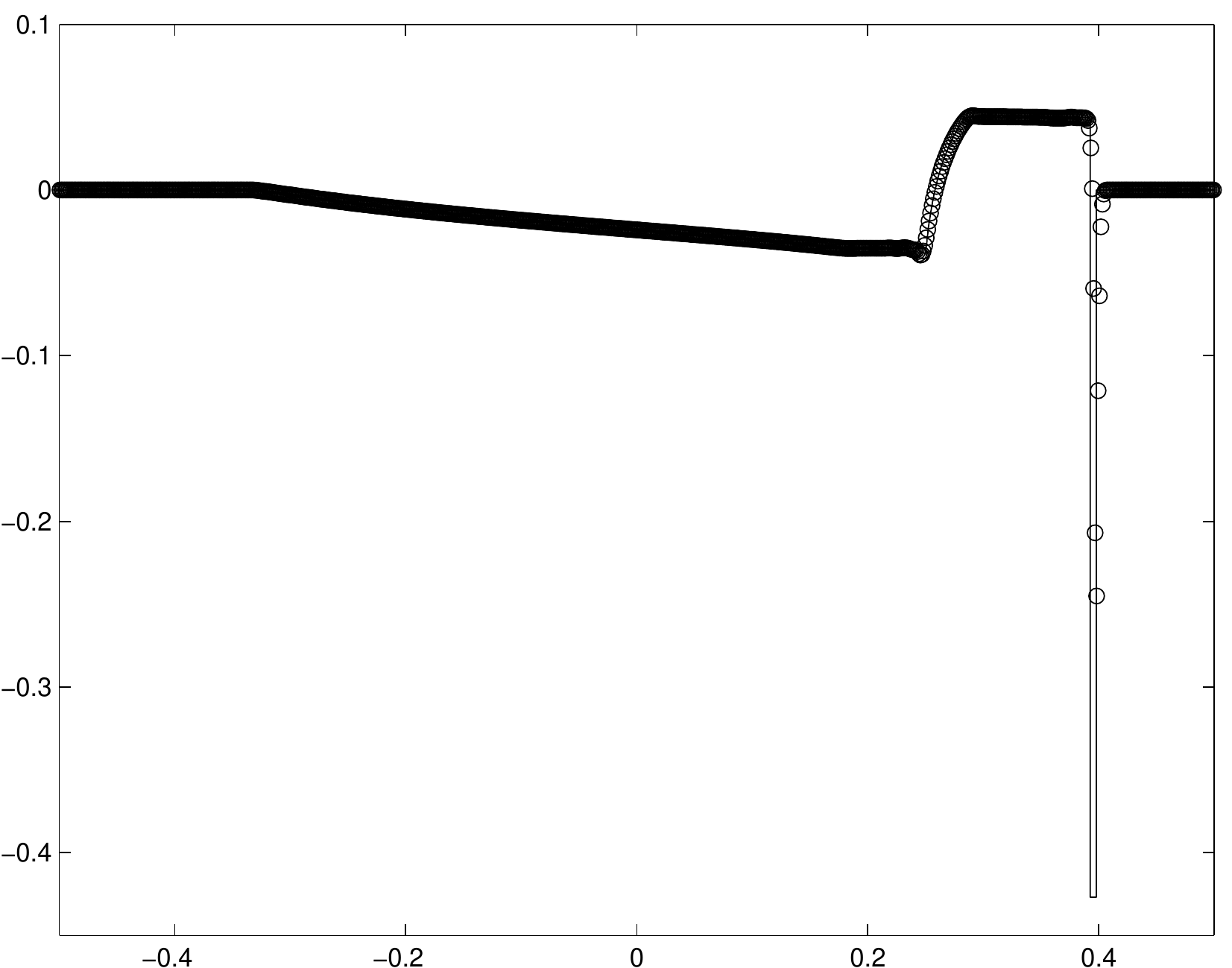}&
\includegraphics[width=0.35\textwidth]{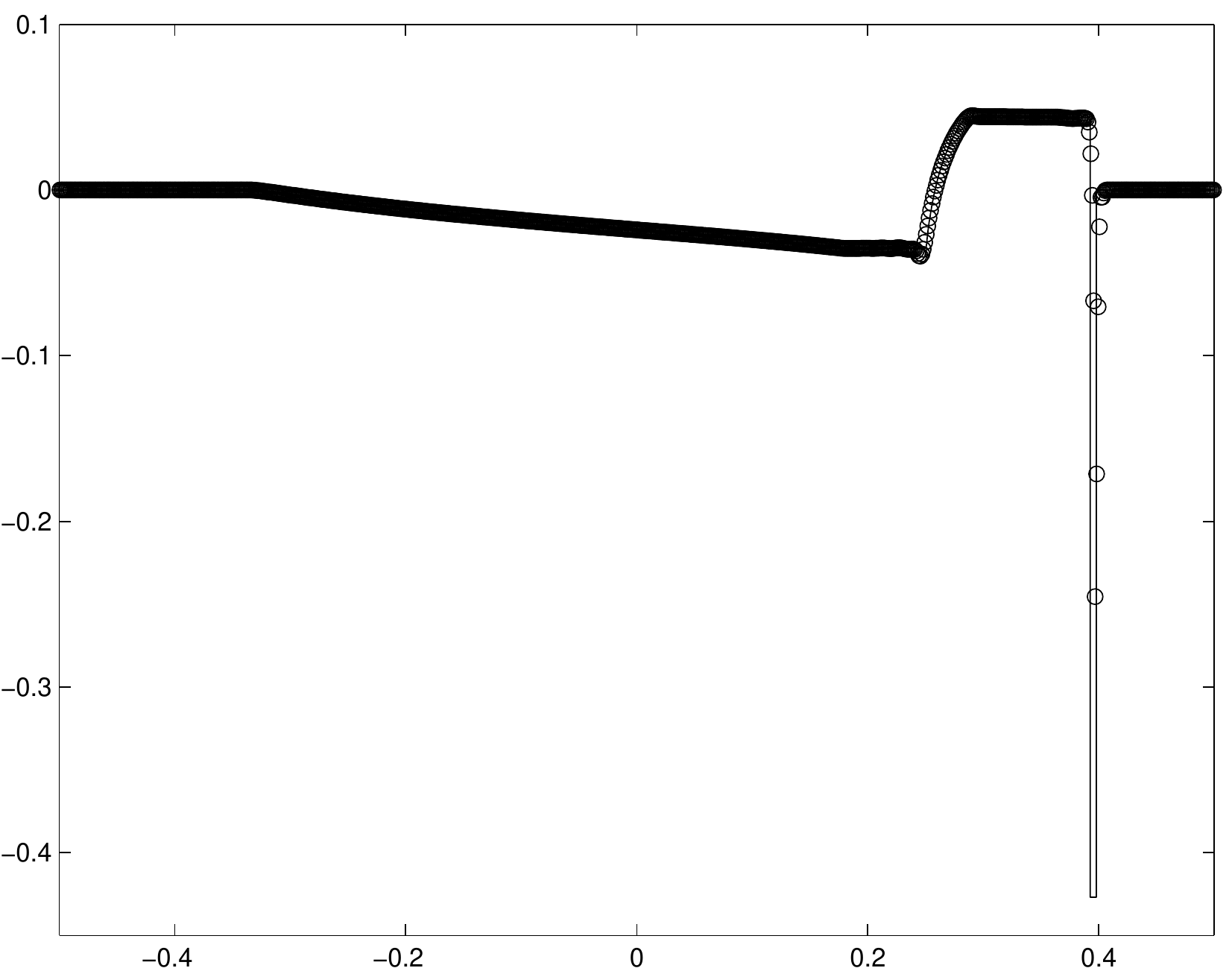}\\
    \end{tabular}
    \caption{Same as Fig. \ref{fig:RMHDRMT4rho} except for the the velocity $v_y$.}
    \label{fig:RMHDRMT4vy}
  \end{figure}

\begin{figure}[!htbp]
    \centering{}
  \begin{tabular}{cc}
    \includegraphics[width=0.35\textwidth]{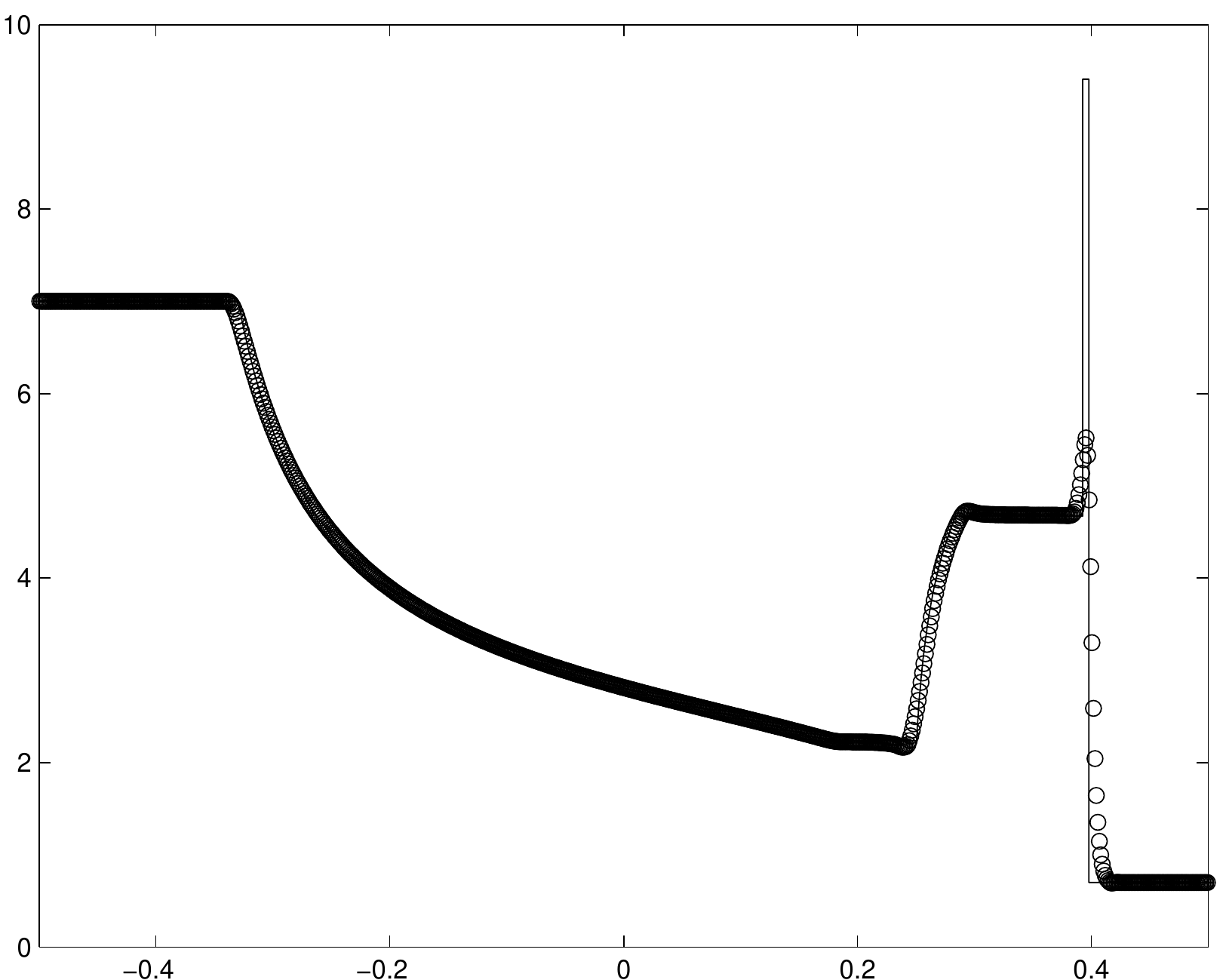}&
\includegraphics[width=0.35\textwidth]{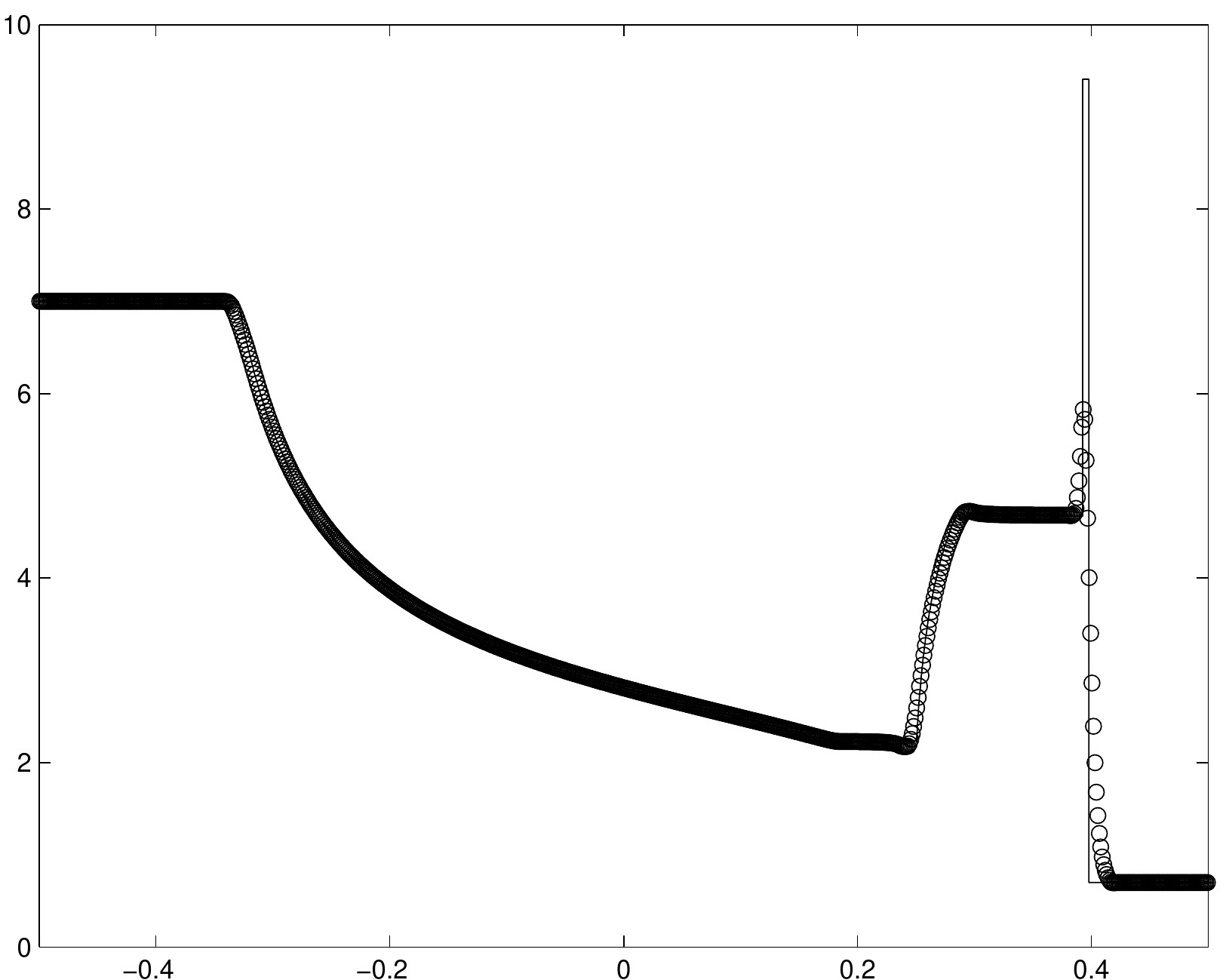}\\
\includegraphics[width=0.35\textwidth]{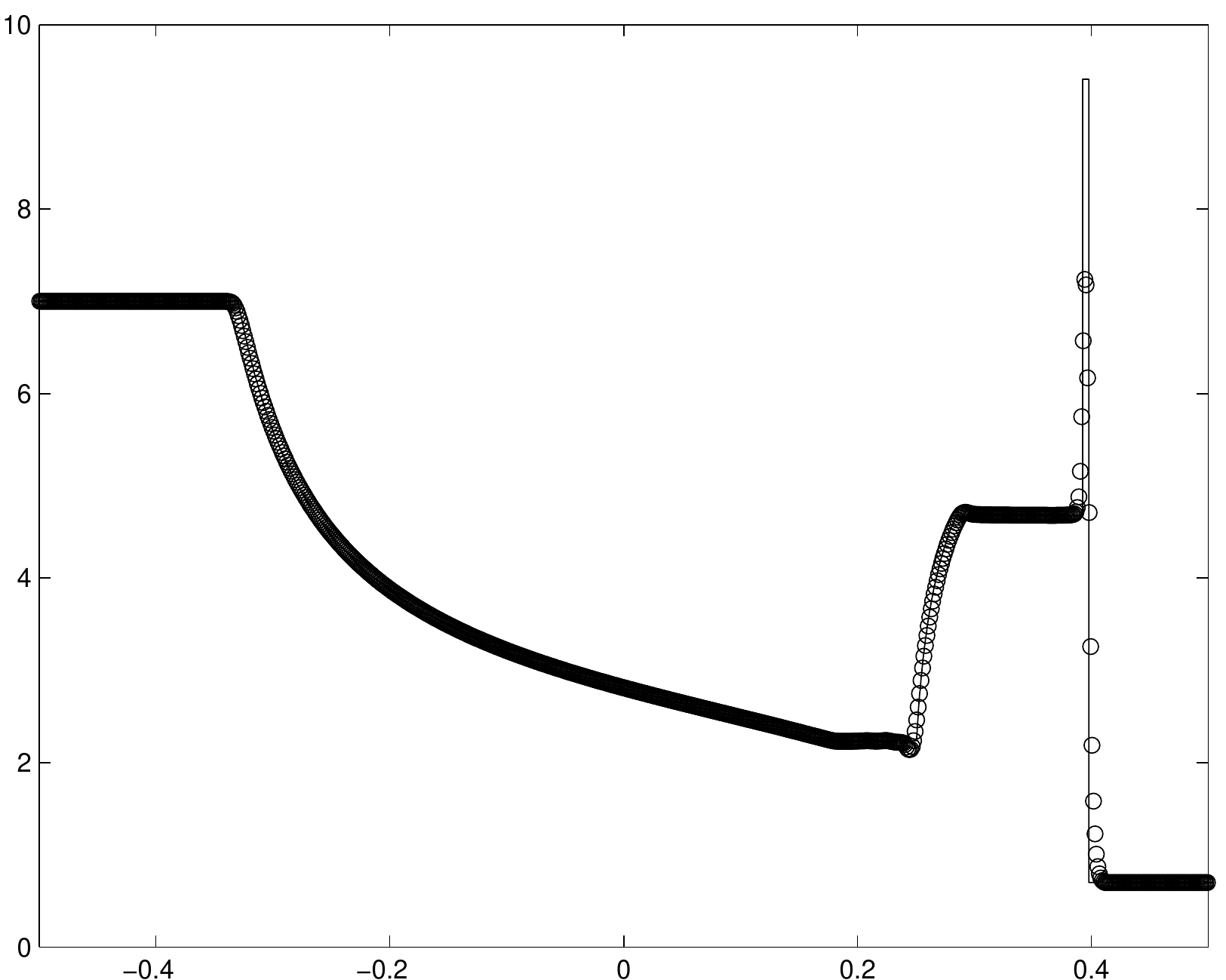}&
\includegraphics[width=0.35\textwidth]{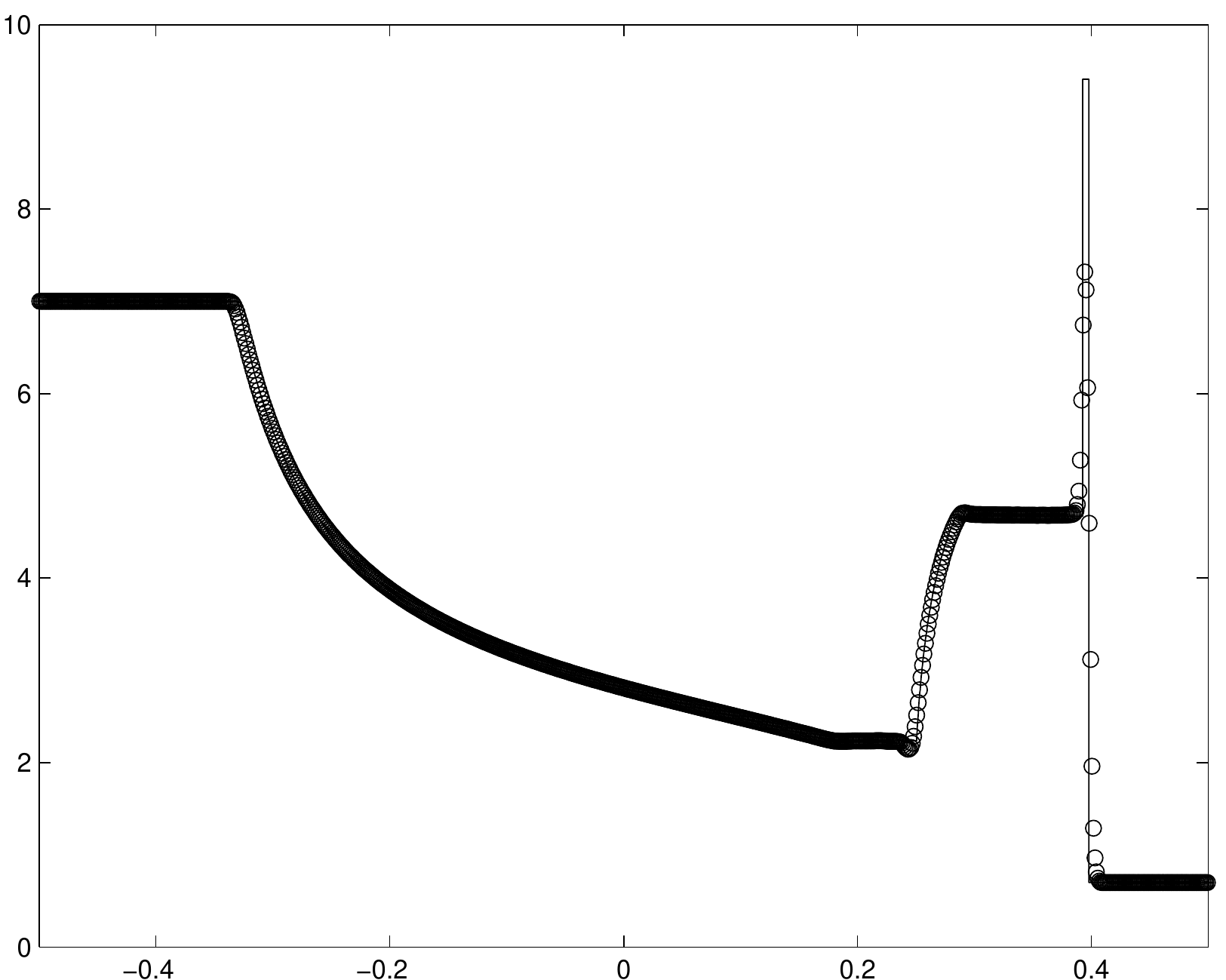}\\
\includegraphics[width=0.35\textwidth]{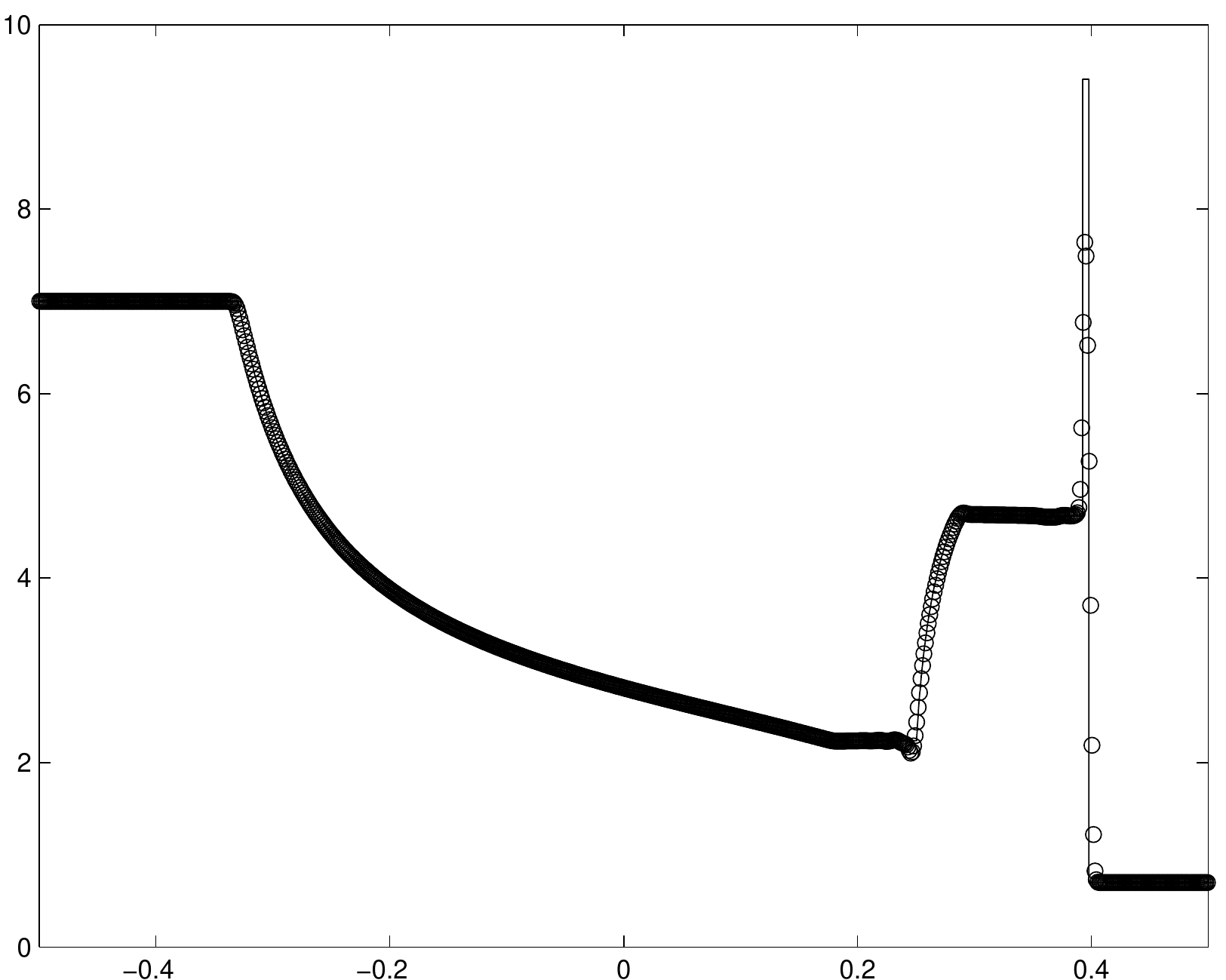}&
\includegraphics[width=0.35\textwidth]{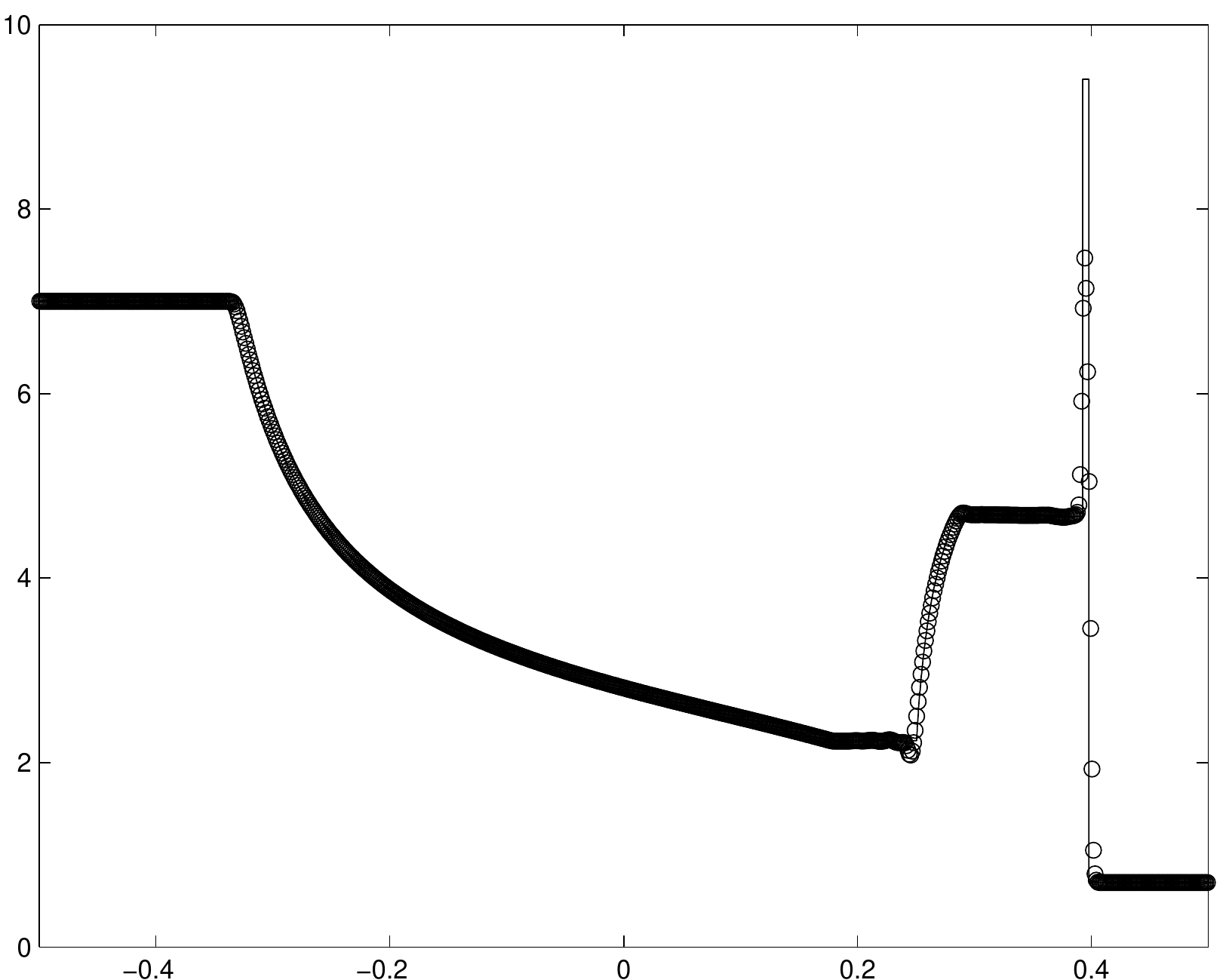}\\
    \end{tabular}
    \caption{Same as Fig. \ref{fig:RMHDRMT4rho} except for the maagnetic field $B_y$.
    }
    \label{fig:RMHDRMT4by}
  \end{figure}

\begin{figure}[!htbp]
    \centering{}
  \begin{tabular}{cc}
    \includegraphics[width=0.35\textwidth]{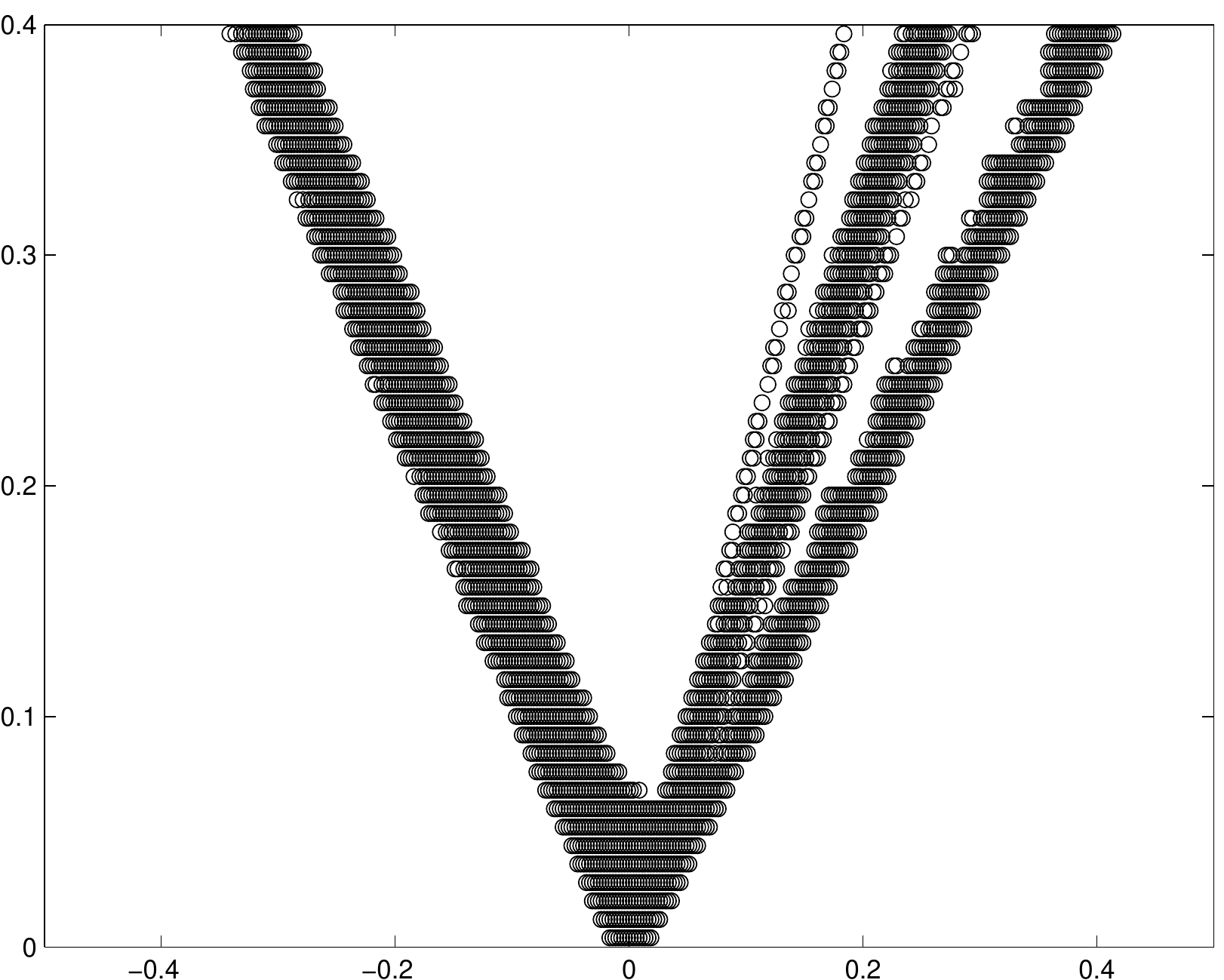}&
\includegraphics[width=0.35\textwidth]{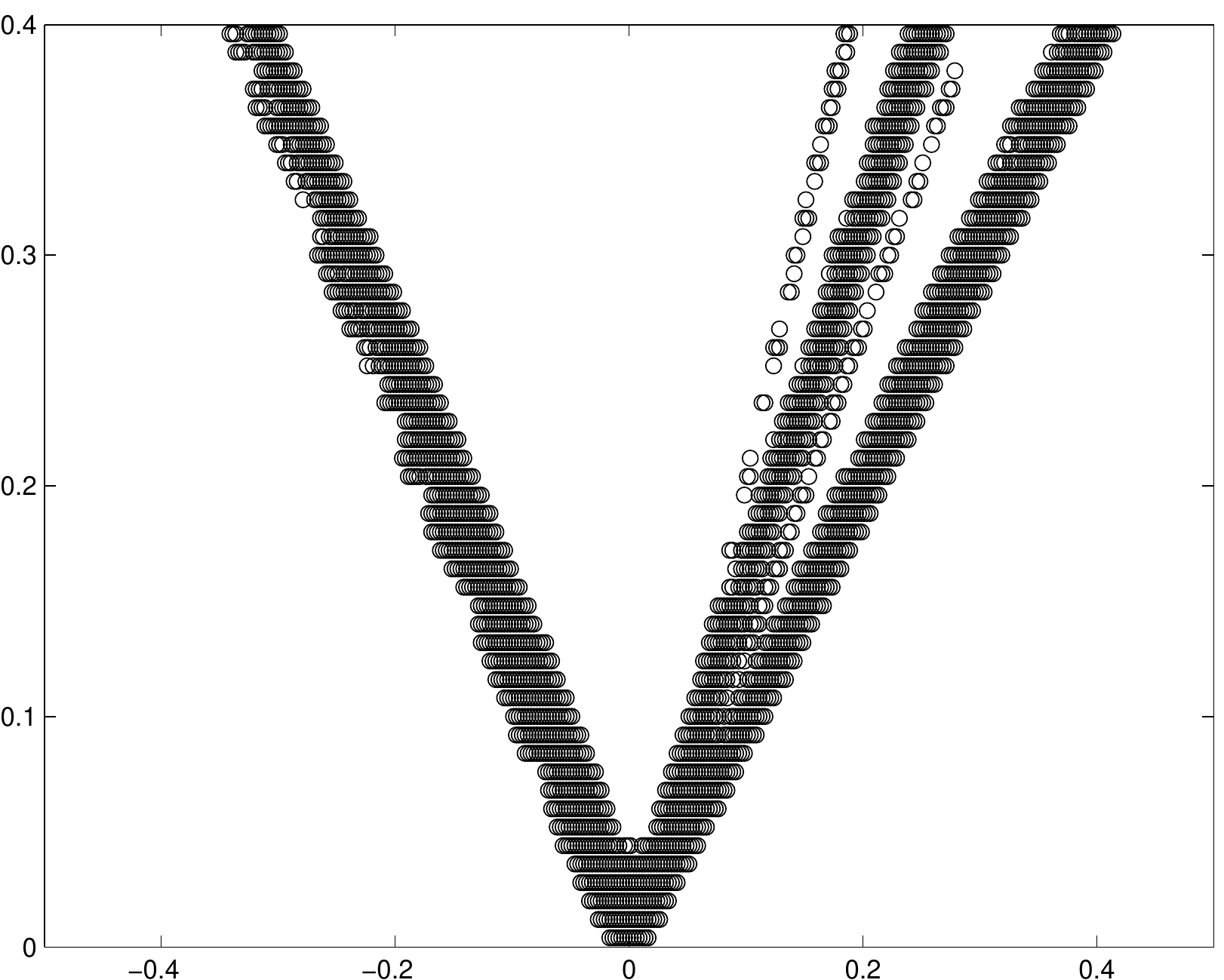}\\
\includegraphics[width=0.35\textwidth]{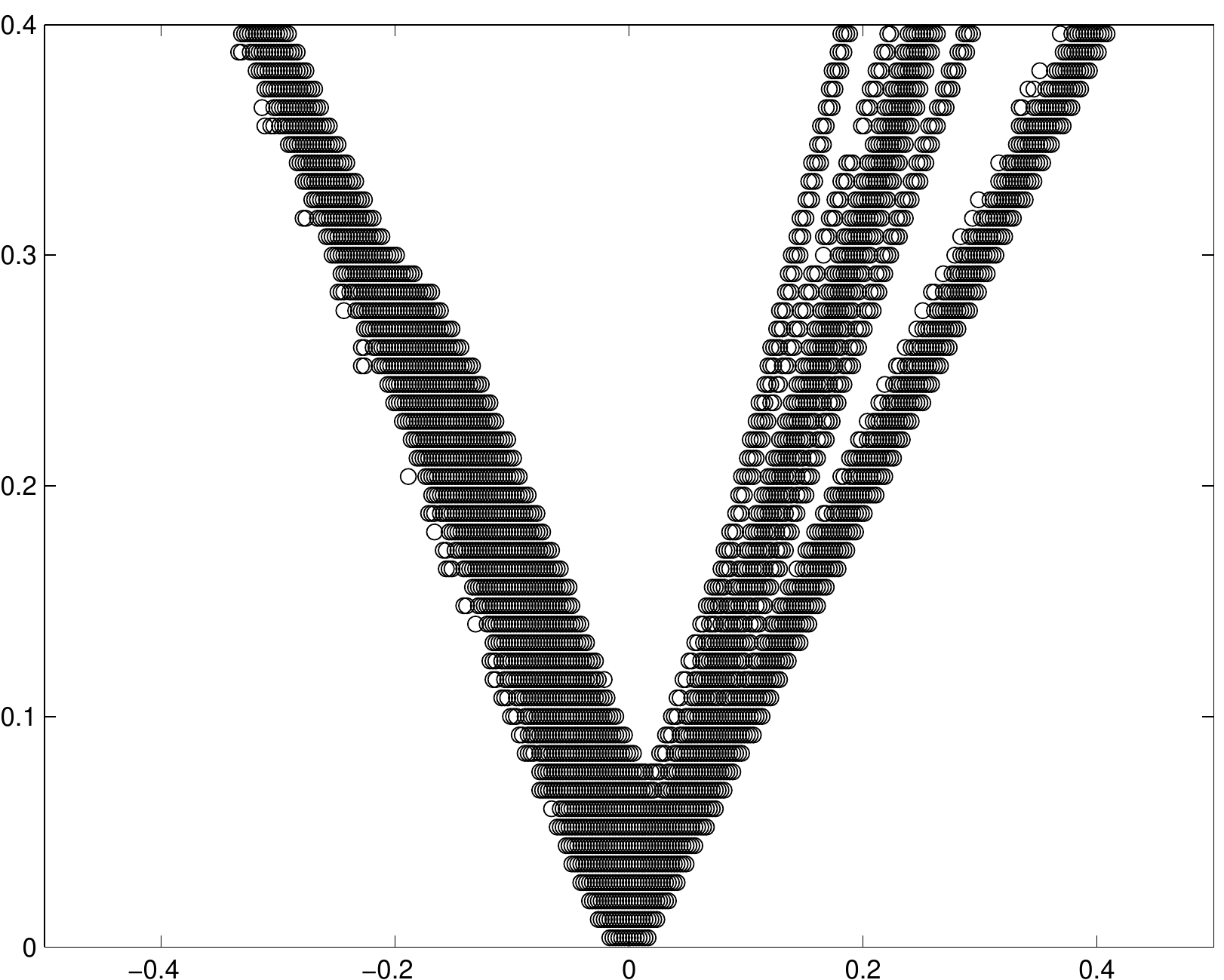}&
\includegraphics[width=0.35\textwidth]{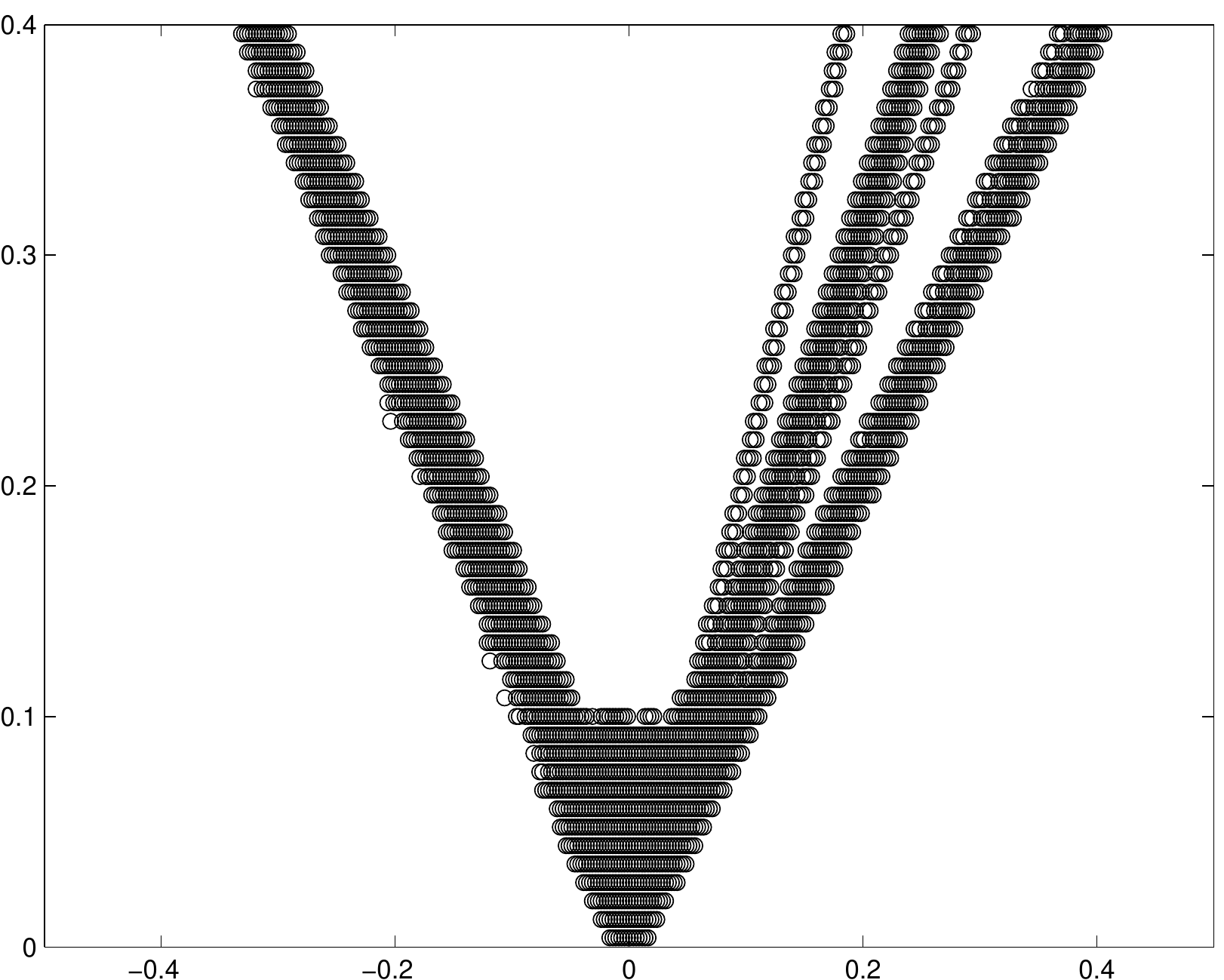}\\
\includegraphics[width=0.35\textwidth]{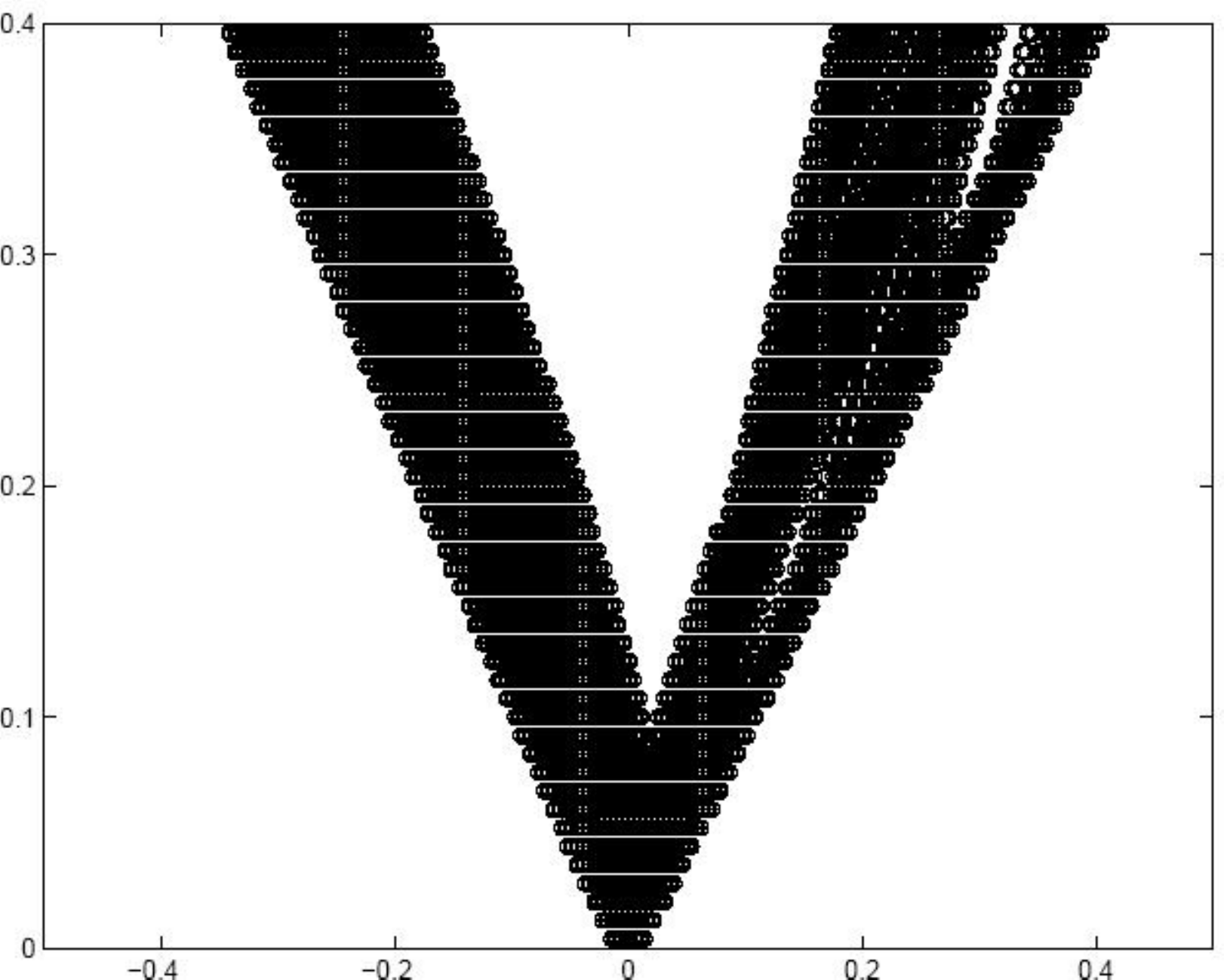}&
\includegraphics[width=0.35\textwidth]{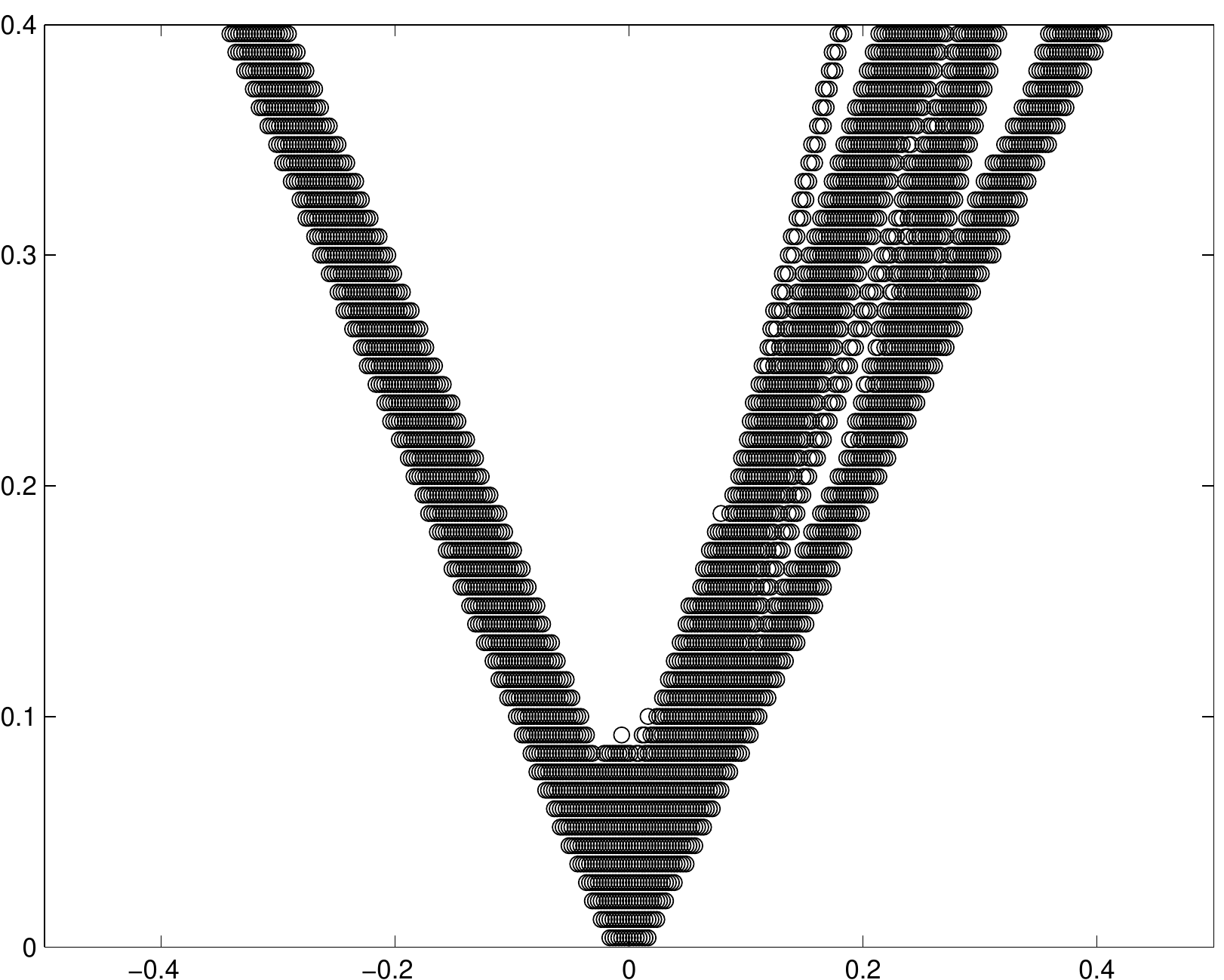}\\
    \end{tabular}
    \caption{Example ~\ref{exRMHDRMT4}:
    The time evolution of ``troubled'' cells.
      Left: non-central \DG{}; right: \CDG{}.
      From top to bottom: $K=1,~2,~3$. The cell number  is $800$. }
    \label{fig:RMHDRMT4cell}
  \end{figure}

 \begin{Example}[Riemann problem 3]\label{exRMHDRMT4}\rm
    The initial data of last 1D Riemann problem are taken as
$$(\rho,v_x,v_y,v_z,B_x,B_y,B_z,p)=\begin{cases}
  (1,0,0,0,10,7,7,1000),&x<0,\\
  (1,0,0,0,10,0.7,0.7,0.1),&x>0,
  \end{cases}
  $$
with the adiabatic index $\Gamma=5/3$.
 It has the same wave structure as Example \ref{exRMHDRMT3},
 but two right-moving shock waves are very strong due to high ratio of initial pressures,
  and their speeds are very close to that of contact discontinuity so that  the difficulty of numerical simulation is seriously increased.

 Figs.~\ref{fig:RMHDRMT4rho}, \ref{fig:RMHDRMT4vy}, and \ref{fig:RMHDRMT4by}
give the densities $\rho$, the velocities  $v_y$, and magnetic fields $B_y$
at $t=0.4$ obtained by using the non-central and central DG methods.
The results show that the densities $\rho$ and
velocities $v_y$ of $P^3$-based DG methods
are obviously better than those of the $P^1$- and $P^2$-based  DG methods.
There is no obvious difference between the non-central and central methods of the same order, and
the numbers of `` troubled'' cells is also basically the same, see Fig.~\ref{fig:RMHDRMT4cell}.
  \end{Example}

\subsection{2D examples}

This section  will solve a smooth problem,
Orszag-Tang problem, blast   problem,    Rotor problem,
and the interaction between the shock wave and cloud by the proposed DG methods.
Unless otherwise stated, the determination of time stepsize is the same as the above and the parameter
$M$ in the modified TVB minmod function is taken as 50.

  \begin{Example}[Smooth problem]\label{exRMHDSmooth2D}\rm
It describes the periodic propagation of a sine wave within the domain $[0,2/\sqrt{3}]\times [0,2]$
and is used to test the accuracy of proposed DG methods.
  The angle between the direction of propagation and the $x$-axis is $\alpha=30^\circ$. The detailed
  initial conditions are     \begin{align*}
  \rho(x,y,0)&=1,\ \ v_x(x,y,0)=-0.1\sin(2\pi \zeta)\sin \alpha,\\
v_y(x,y,0)&=0.1\sin(2\pi \zeta)\cos \alpha,\ \ v_z(x,y,0)=0.1\cos(2\pi \zeta),\\
B_x(x,y,0)&=\cos \alpha +\kappa v_x (x,y,0),\ \ B_y(x,y,0)=\sin
\alpha+\kappa v_y(x,y,0),\\
B_z(x,y,0)&=\kappa v_z(x,y,0),\ \ p(x,y,0)=0.1,
\end{align*}
     where $ \zeta=x\cos\alpha +y\sin\alpha$, $\kappa=\sqrt{1+\rho
       h\gamma^2}$, and corresponding exact solutions are
     \begin{align*}
 \rho(x,y,t)&=1,\quad v_x(x,y,t)=-0.1\sin\big(2\pi (\zeta+t/\kappa)\big)\sin \alpha,\\
v_y(x,y,t)&=0.1\sin\big(2\pi (\zeta+t/\kappa)\big)\cos
\alpha,\quad v_z(x,y,t)=0.1\cos\big(2\pi (\zeta+t/\kappa)\big),\\
B_x(x,y,t)&=\cos \alpha +\kappa v_x (x,y,t),\quad B_y(x,y,t)=\sin
\alpha+\kappa v_y(x,y,t),\\
B_z(x,y,t)&=\kappa v_z(x,y,t),\quad p(x,y,t)=0.1.
\end{align*}

In the present computations,   the fourth-order Runge-Kutta method mentioned in [58] is used
in order to ensure the accuracy in time.    The CFL numbers of
     $P^1$-, $P^2$-, and $P^3$-based non-central DG methods are $0.2$, $0.15$, $0.1$, respectively,
      and the CFL numbers for corresponding  central DG methods are $0.3,~0.25,~0.2$, respectively,
      and $\theta=\Delta t_n/\tau_n=1$.
    Table ~\ref{tab:RMHDsmooth2D} lists  $l^1$ errors and orders  in $B_x$ obtained by the non-central
    and central DG methods.
    As can be seen from the table, two kinds of DG methods have reached
    the expected convergence orders.

   \end{Example}

\begin{table}[!htbp]
  \centering
  \caption{Example \ref{exRMHDSmooth2D}£º
 $l^1$ errors in $B_x$ and orders at $t=1$ obtained by
non-central and central \DG{}.   The fourth-order accurate Runge-Kutta matheod and $N \times 2N$ uniform cells are used. }
\begin{tabular}{|c|c|c|c|c|c|c|c|c|c|}
  \hline
 \multirow{3}{20pt}{}
 &\multirow{2}{2pt}{}
 &\multicolumn{4}{|c|}{without limiter}&\multicolumn{4}{|c|}{with limiter in global}\\
 \cline{3-10}
 & & \multicolumn{2}{|c|}{non-central DG}&\multicolumn{2}{|c|}{central DG}&
 \multicolumn{2}{|c|}{non-central DG}&\multicolumn{2}{|c|}{central DG}\\
 \cline{2-10}
 & N &  $l^1$ error& order &  $l^1$ error & order &  $l^1$ error& order &  $l^1$ error & order \\
\hline
\multirow{6}{20pt}{$P^{1}$}
&10& 4.03e-03& --& 3.48e-03& --& 4.05e-02& --& 3.65e-02& --\\
\cline{2-10}
&20&6.99e-04& 2.53&4.51e-04& 2.95&1.27e-02& 1.68&1.18e-02& 1.63\\
\cline{2-10}
&40&1.46e-04& 2.25&5.98e-05& 2.91&4.70e-03& 1.43&4.35e-03& 1.44\\
\cline{2-10}
&80&3.37e-05& 2.12&9.54e-06& 2.65&1.29e-03& 1.86&1.20e-03& 1.86\\
\cline{2-10}
&160&8.19e-06& 2.04&1.72e-06& 2.47&3.27e-04& 1.98&3.03e-04& 1.98\\
\cline{2-10}
&320&2.03e-06& 2.02&4.11e-07& 2.07&8.07e-05& 2.02&7.28e-05& 2.06\\
\hline
\multirow{6}{20pt}{$P^{2}$}
&10& 3.11e-04& --& 8.01e-04& --& 2.93e-03& --& 8.39e-03& --\\
\cline{2-10}
&20&3.87e-05& 3.01&9.81e-05& 3.03&1.22e-04& 4.59&6.69e-04& 3.65\\
\cline{2-10}
&40&4.83e-06& 3.00&1.22e-05& 3.01&1.03e-05& 3.57&5.37e-05& 3.64\\
\cline{2-10}
&80&6.08e-07& 2.99&1.52e-06& 3.00&1.83e-06& 2.49&5.74e-06& 3.23\\
\cline{2-10}
&160&7.74e-08& 2.97&1.89e-07& 3.00&2.48e-07& 2.89&6.90e-07& 3.06\\
\cline{2-10}
&320&9.87e-09& 2.97&2.37e-08& 3.00&3.16e-08& 2.97&8.56e-08& 3.01\\
\hline
\multirow{6}{20pt}{$P^{3}$}
&10& 2.18e-05& --& 3.32e-05& --& 5.51e-04& --& 5.68e-04& --\\
\cline{2-10}
&20&1.77e-06& 3.62&2.23e-06& 3.90&1.06e-05& 5.70&1.57e-05& 5.17\\
\cline{2-10}
&40&1.03e-07& 4.10&1.43e-07& 3.97&2.44e-07& 5.43&6.96e-07& 4.50\\
\cline{2-10}
&80&5.67e-09& 4.19&8.97e-09& 3.99&1.10e-08& 4.47&4.25e-08& 4.04\\
\cline{2-10}
&160&3.11e-10& 4.19&5.62e-10& 4.00&6.44e-10& 4.10&2.69e-09& 3.98\\
\cline{2-10}
&320&1.89e-11& 4.04&3.51e-11& 4.00&4.00e-11& 4.01&1.70e-10& 3.98\\
\hline
\end{tabular}
  \label{tab:RMHDsmooth2D}
  \end{table}

\begin{figure}[!htbp]
    \centering{}
  \begin{tabular}{cc}
    \includegraphics[width=0.35\textwidth]{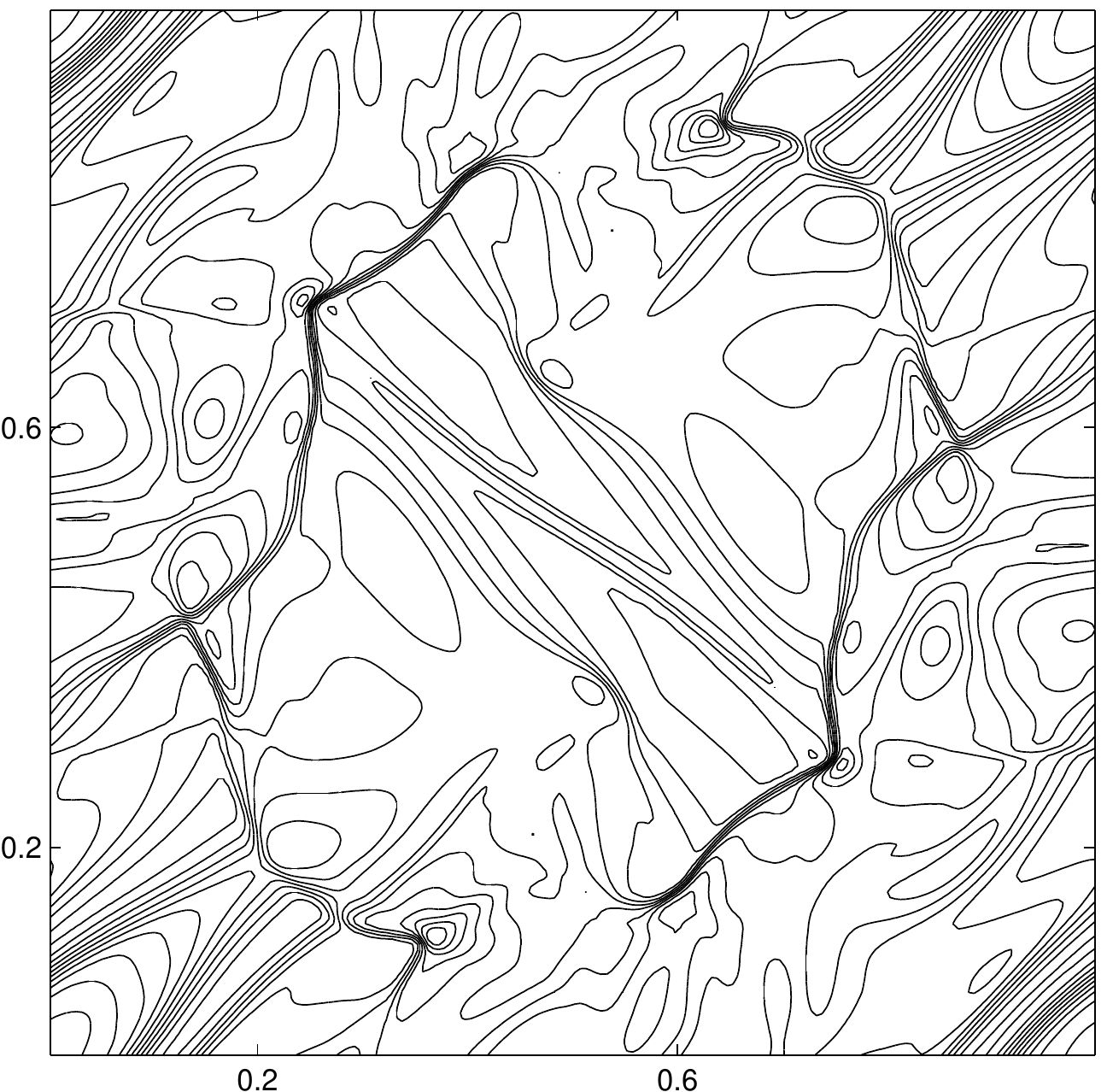}&
 \includegraphics[width=0.35\textwidth]{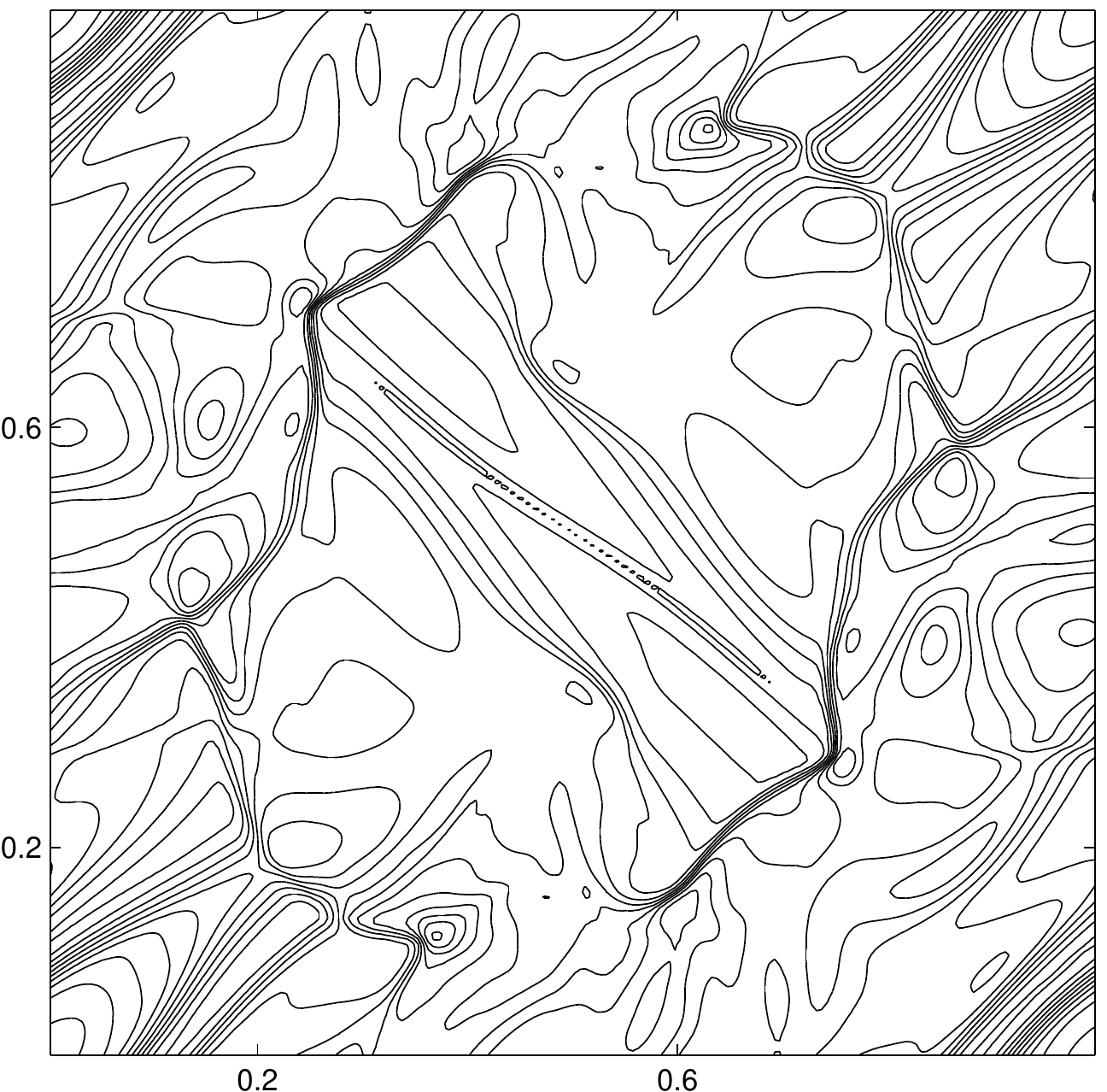}\\
    \includegraphics[width=0.35\textwidth]{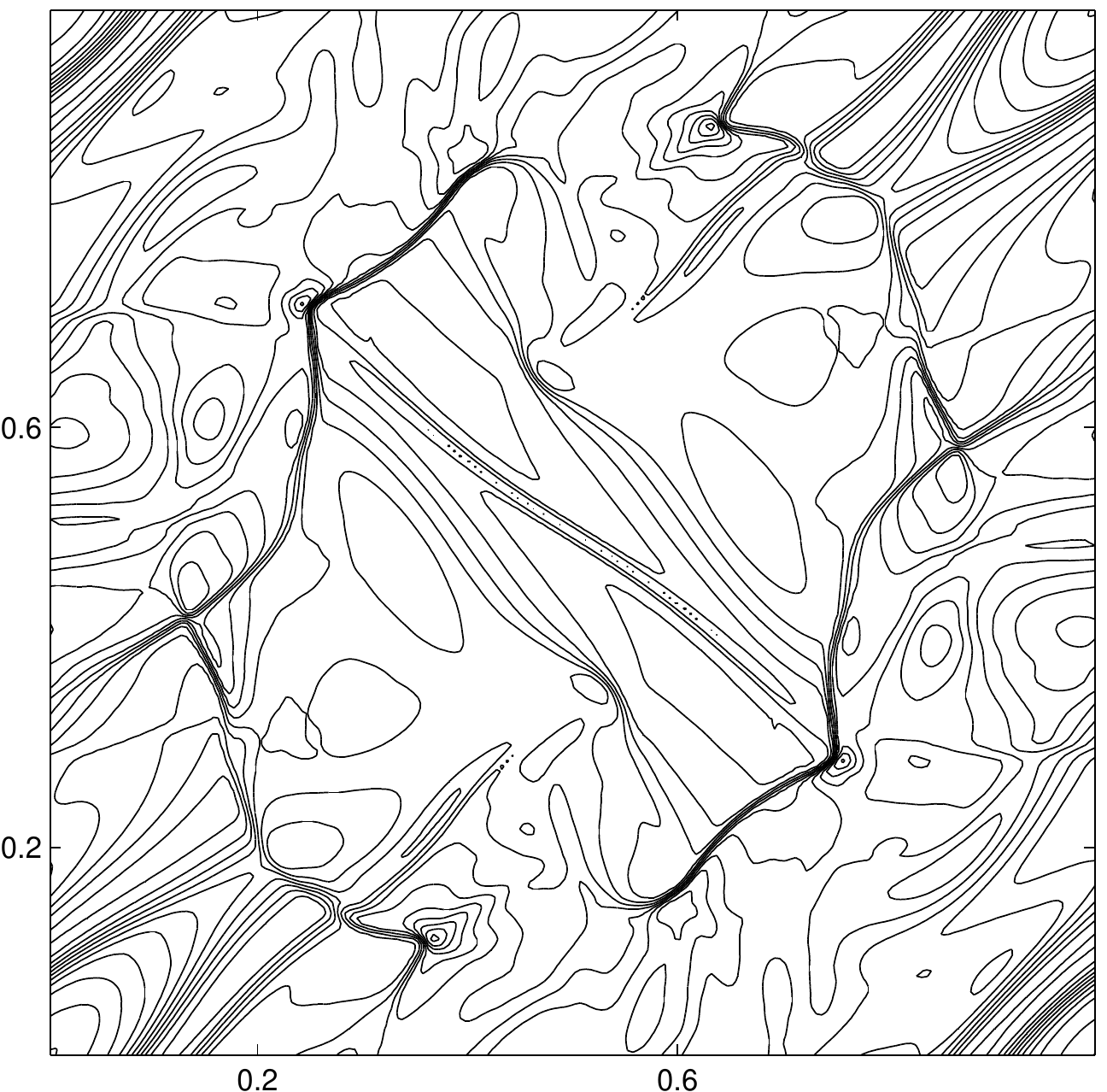}&
 \includegraphics[width=0.35\textwidth]{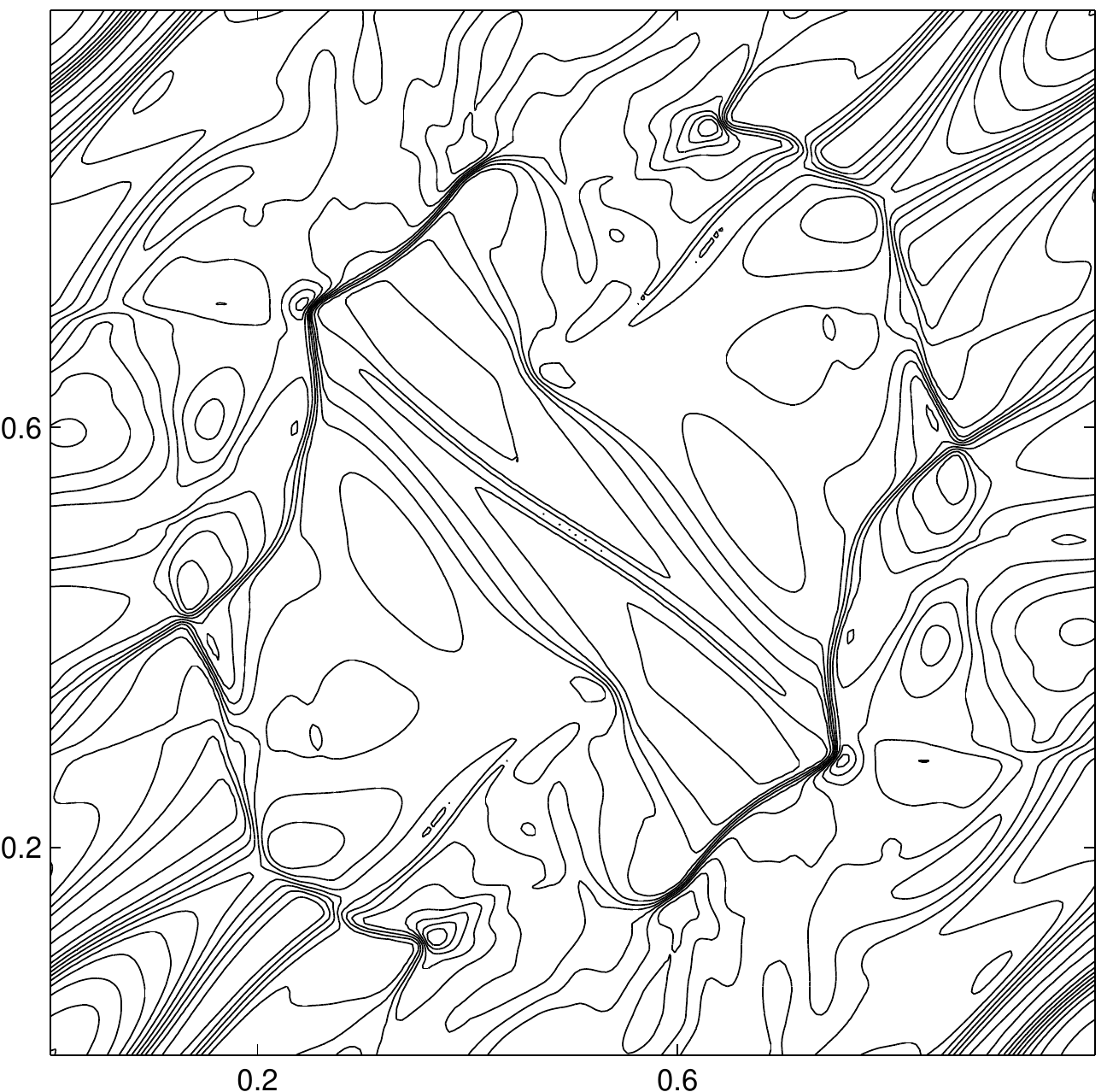}\\
 \includegraphics[width=0.35\textwidth]{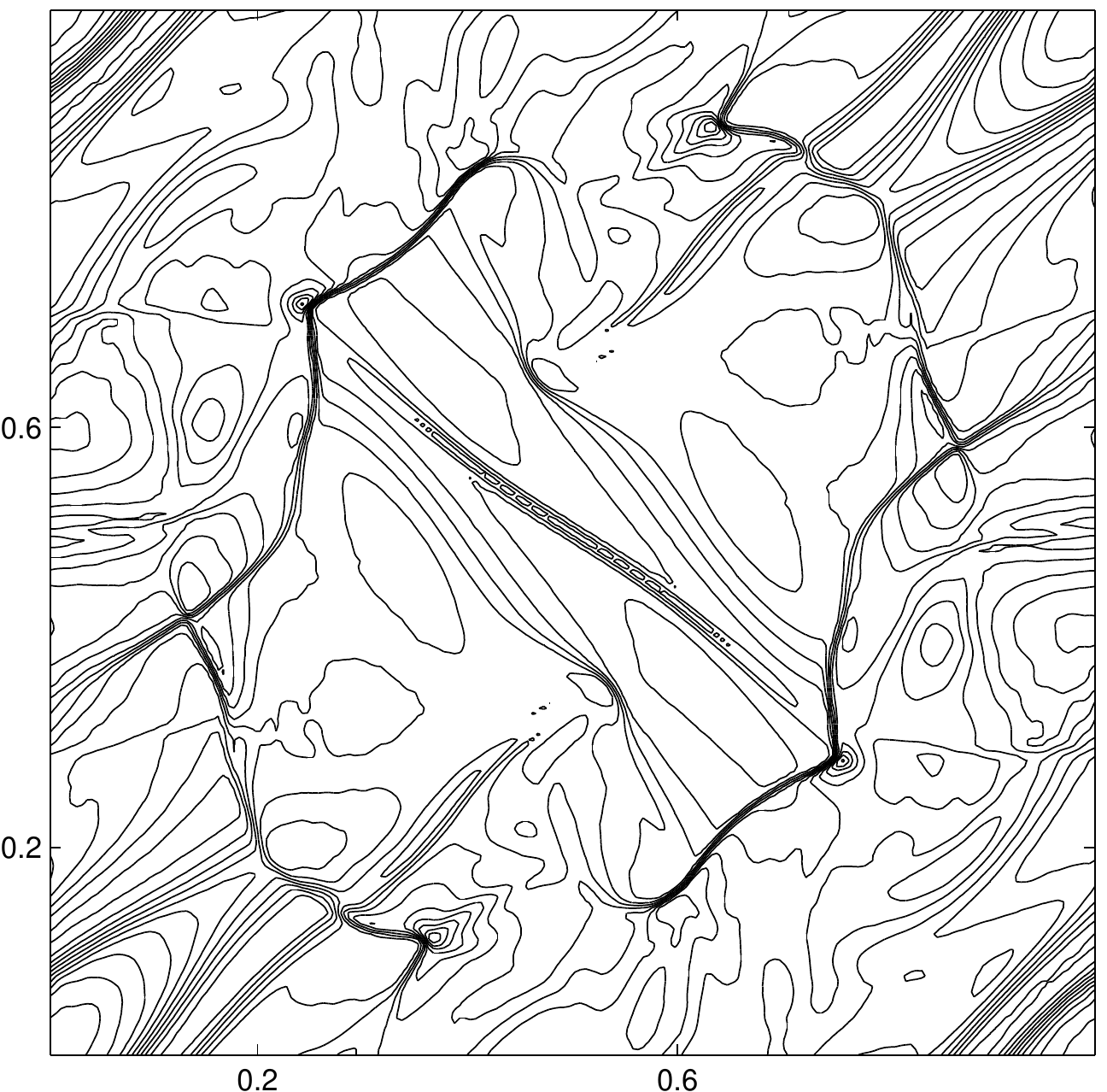}&
  \includegraphics[width=0.35\textwidth]{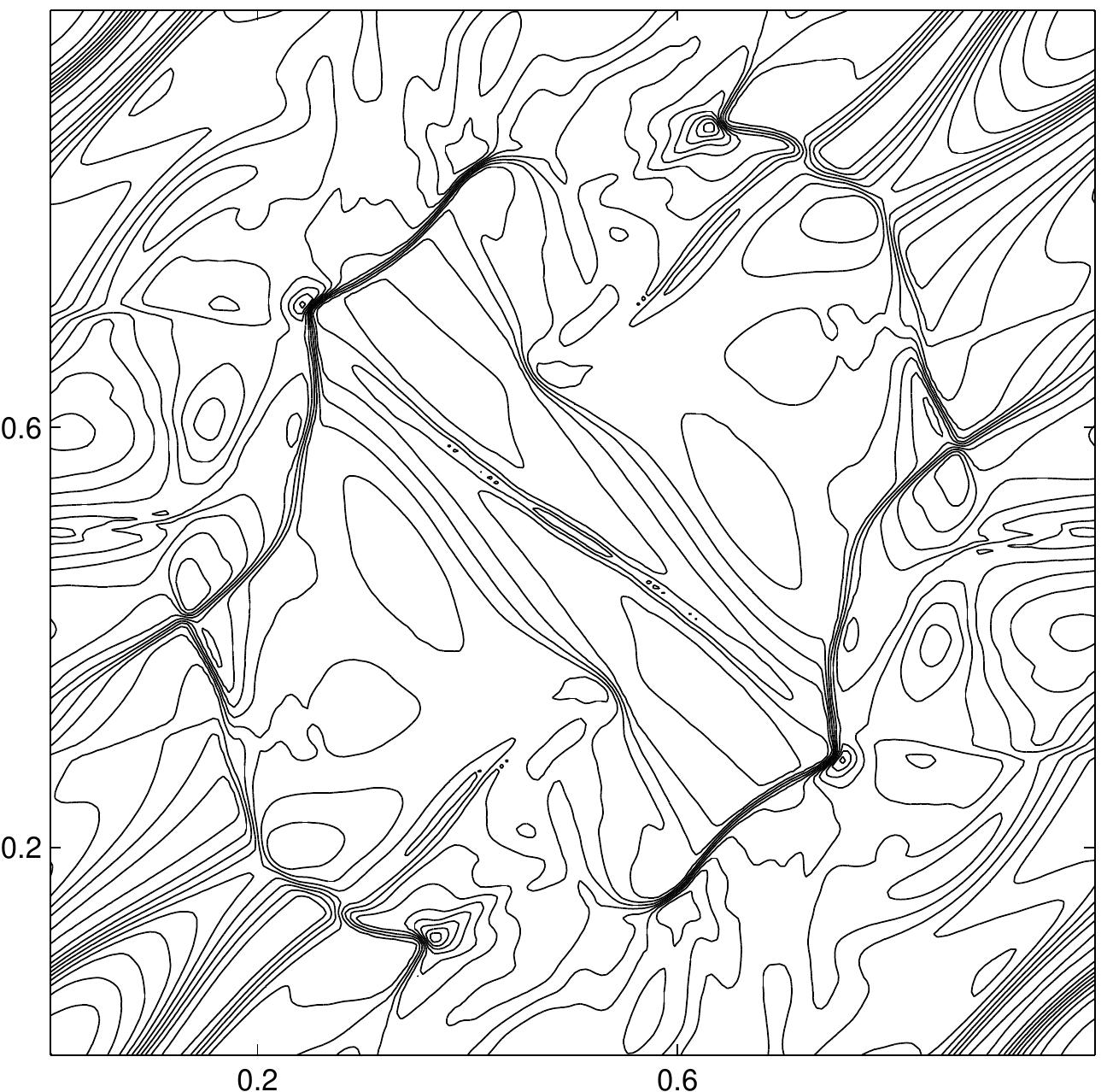}\\
    \end{tabular}
    \caption{Example \ref{exRMHDOT}£º
    The contour plots of density $\rho$ at $t=1$ obtained by using $200\times 200$ cells
     (15 equally spaced contour lines from 0.06 to 0.48).
     Left: $P^K$-based \DG{}; right: $P^K$-based \CDG{}.
     From top to bottom: $K=1,~2,~3$.   }
    \label{fig:RMHDOTrho}
  \end{figure}

\begin{figure}[!htbp]
    \centering{}
  \begin{tabular}{cc}
    \includegraphics[width=0.35\textwidth]{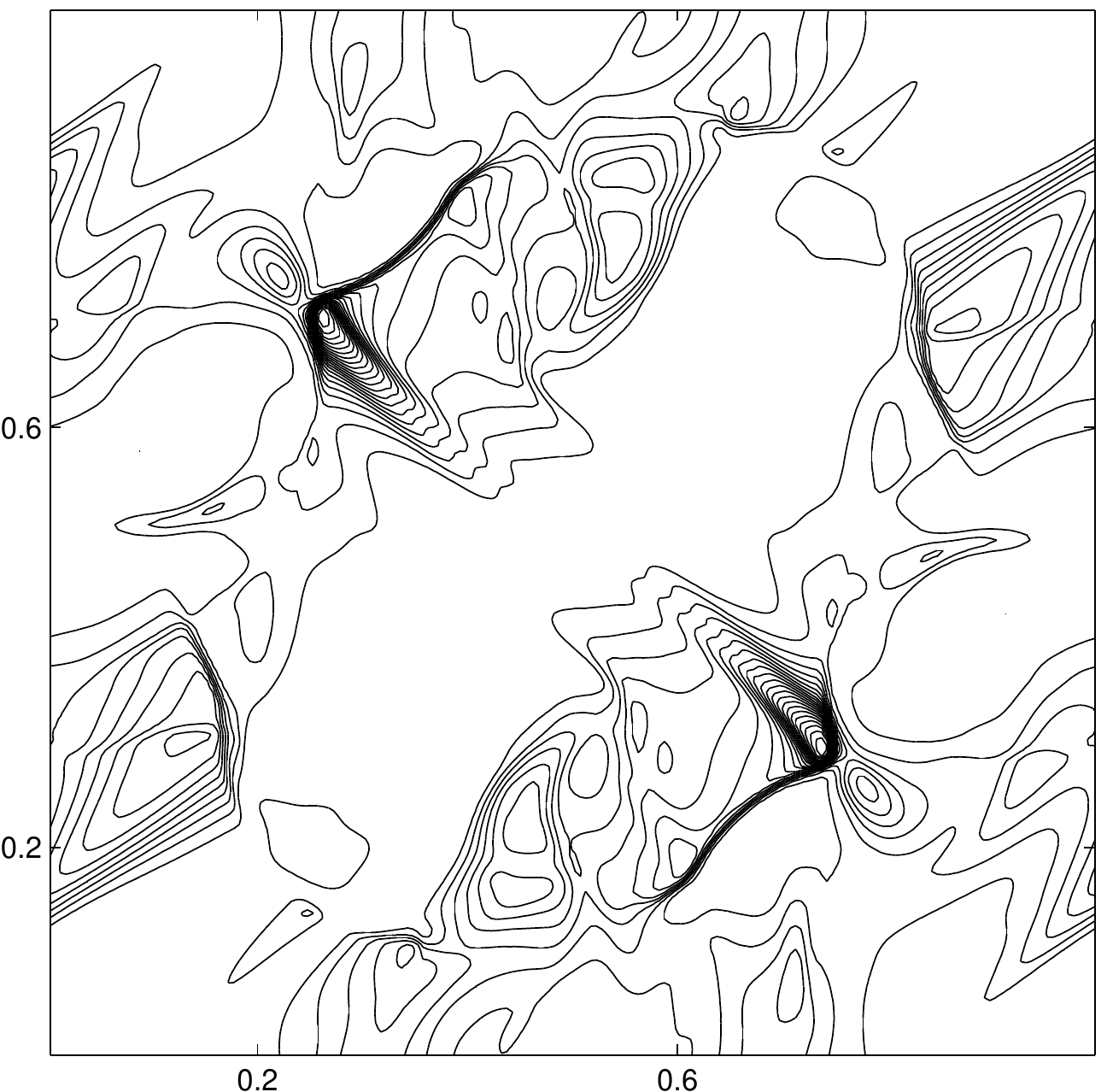}&
 \includegraphics[width=0.35\textwidth]{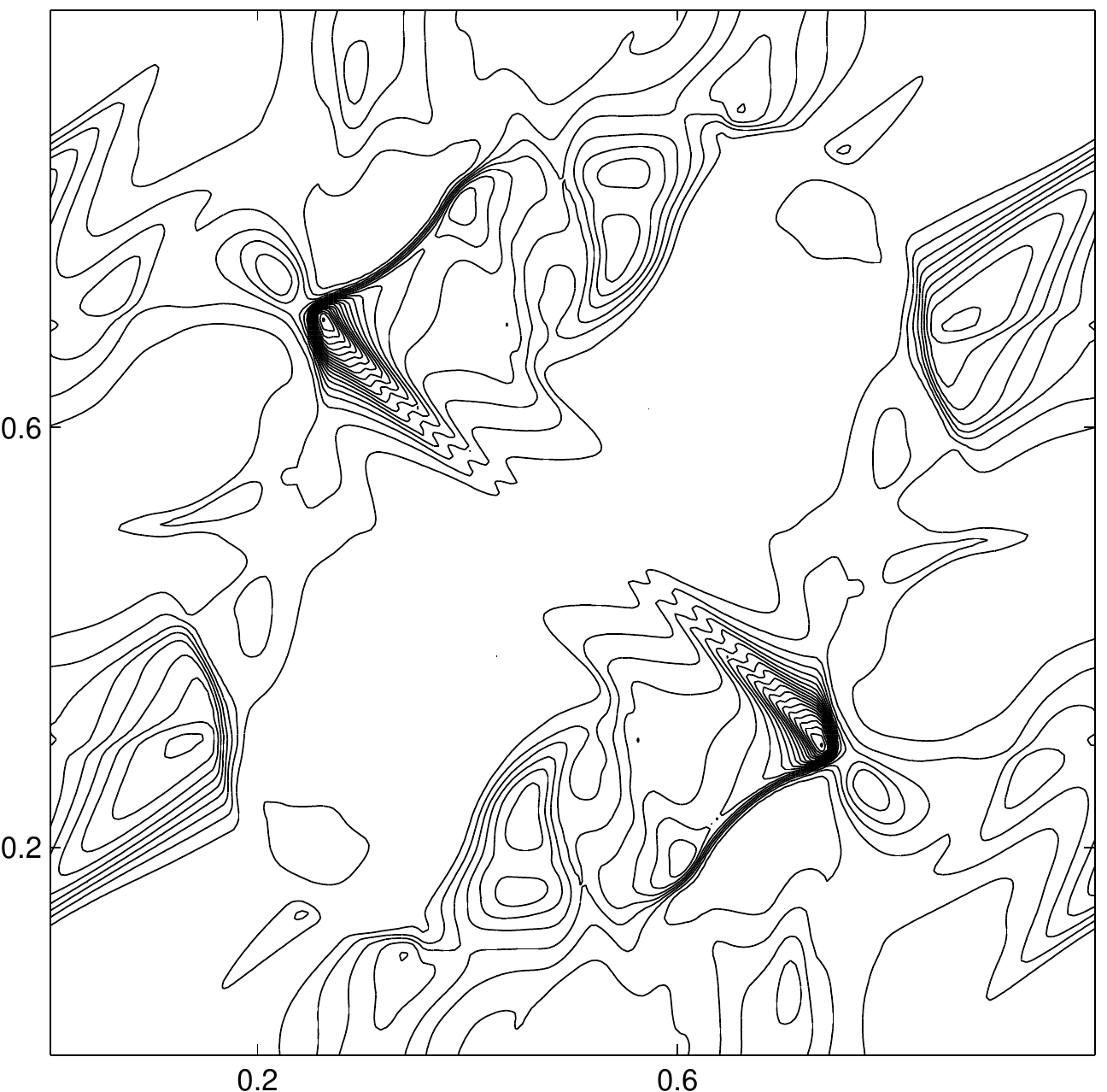}\\
    \includegraphics[width=0.35\textwidth]{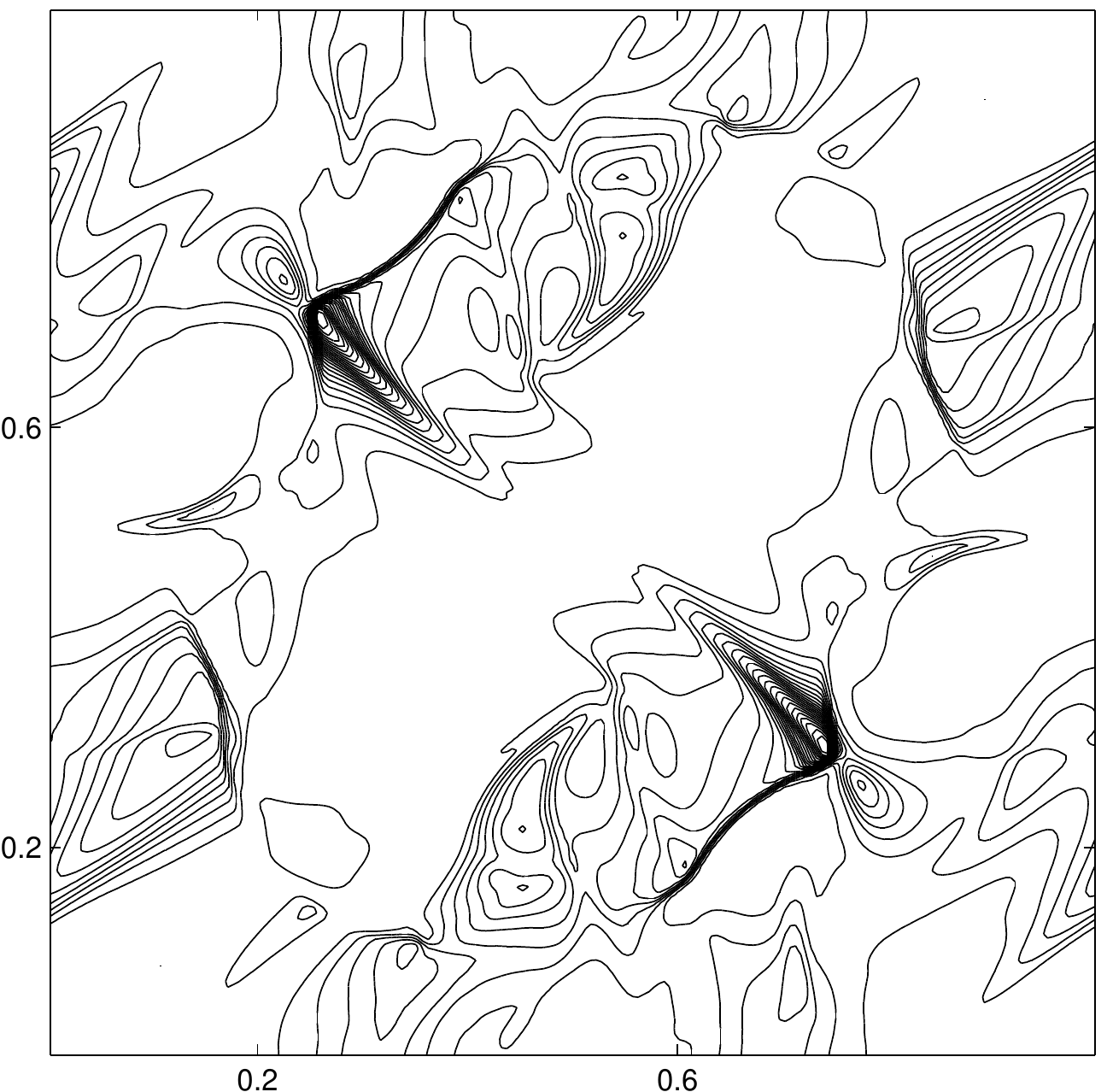}&
 \includegraphics[width=0.35\textwidth]{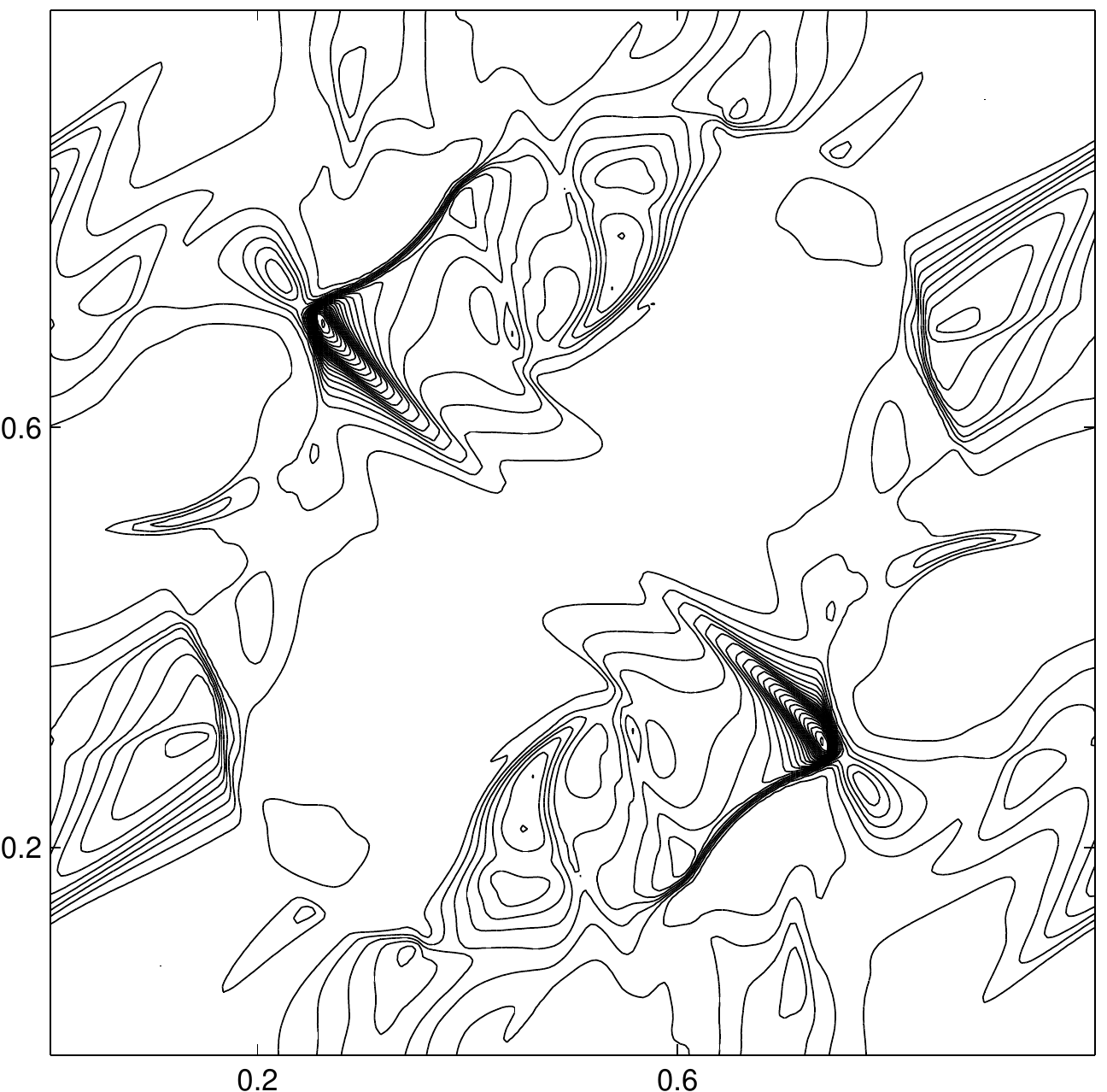}\\
 \includegraphics[width=0.35\textwidth]{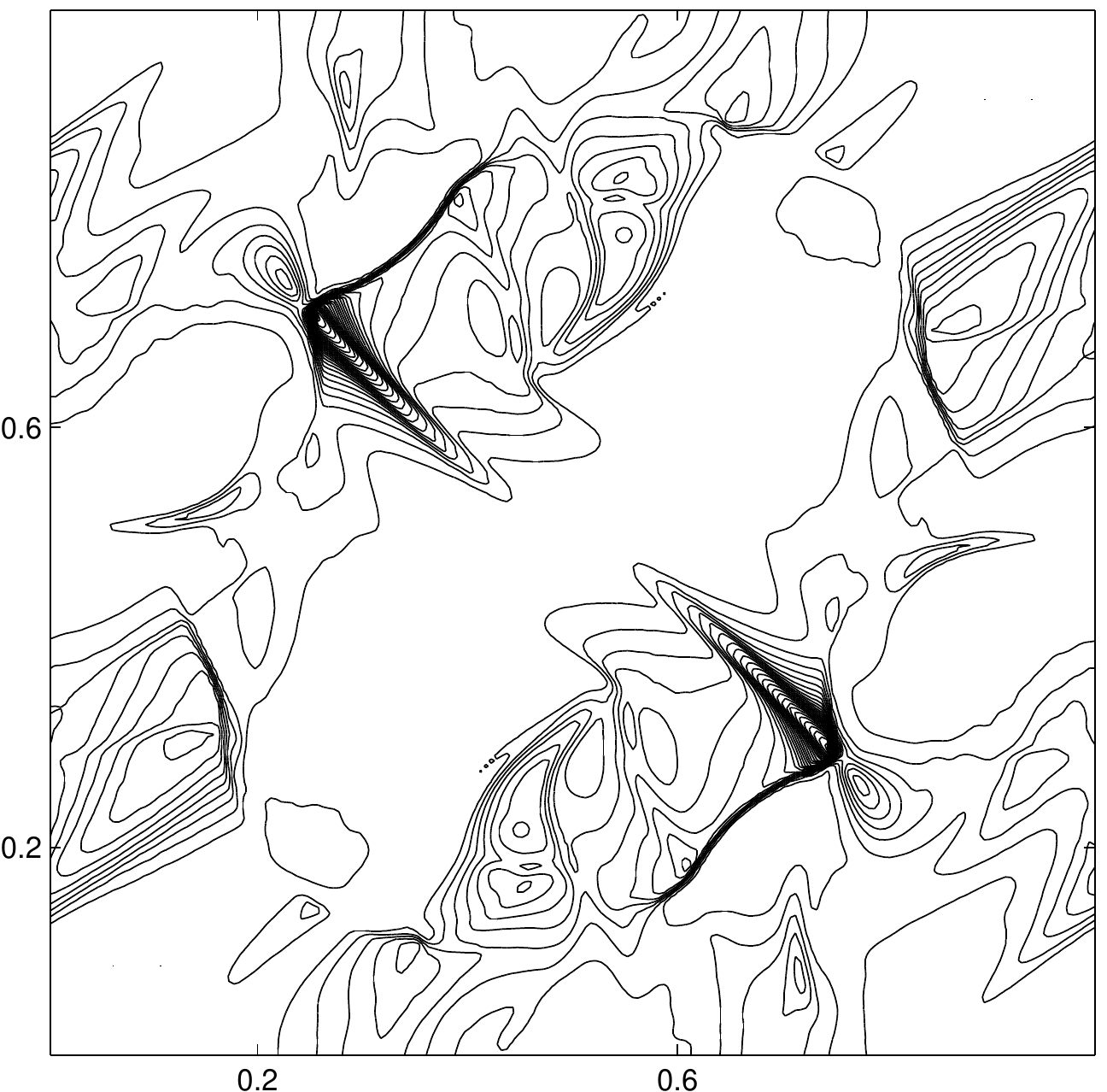}&
 \includegraphics[width=0.35\textwidth]{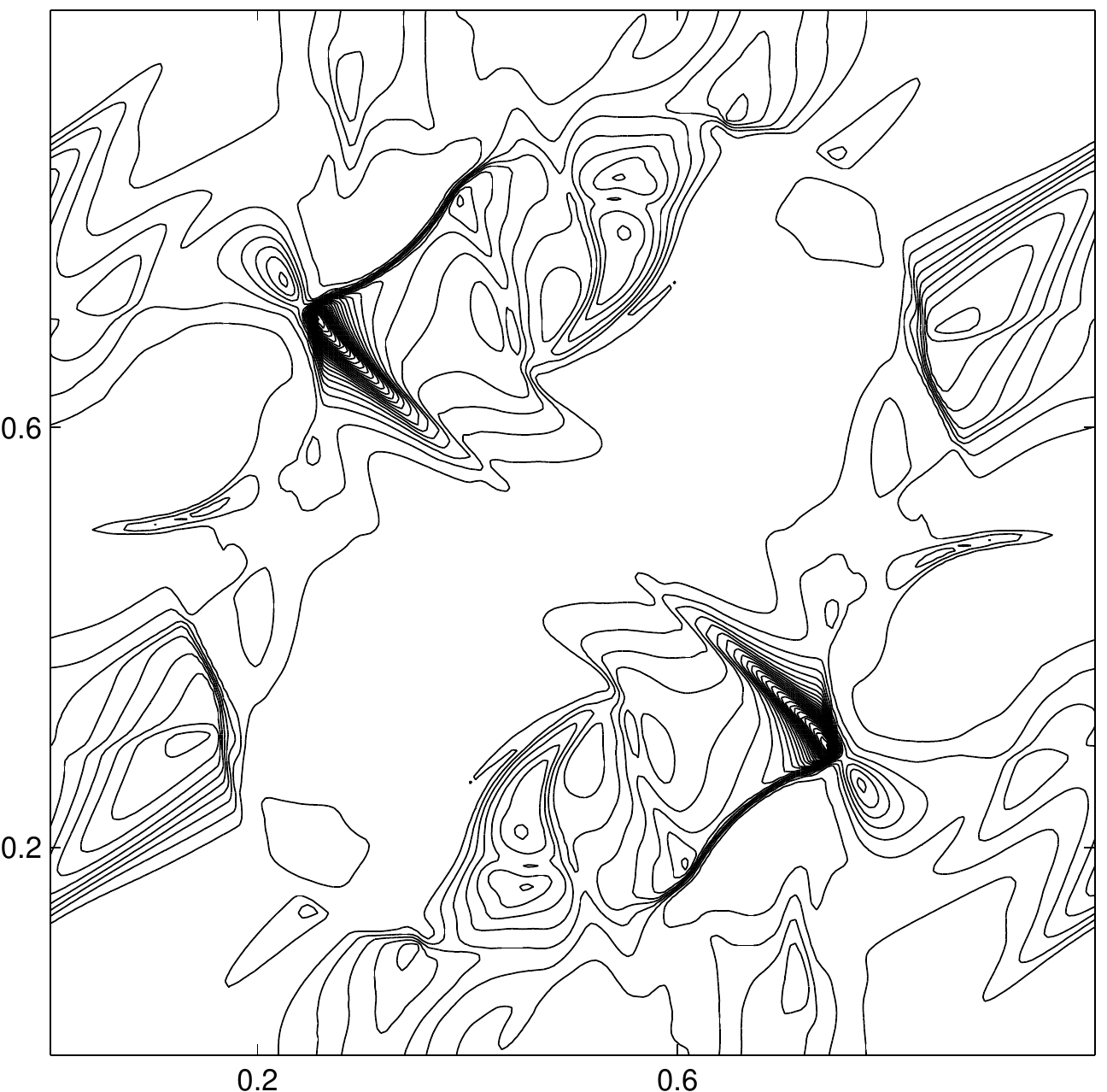}\\
    \end{tabular}
    \caption{Same as Fig.~~\ref{fig:RMHDOTrho} except for the Lorentz factor $\gamma$
    (30 equally spaced contour lines from 1 to 2.2). }
    \label{fig:RMHDOTgam}
  \end{figure}

\begin{figure}[!htbp]
    \centering{}
    \begin{tabular}{cc}
    \includegraphics[width=0.45\textwidth]{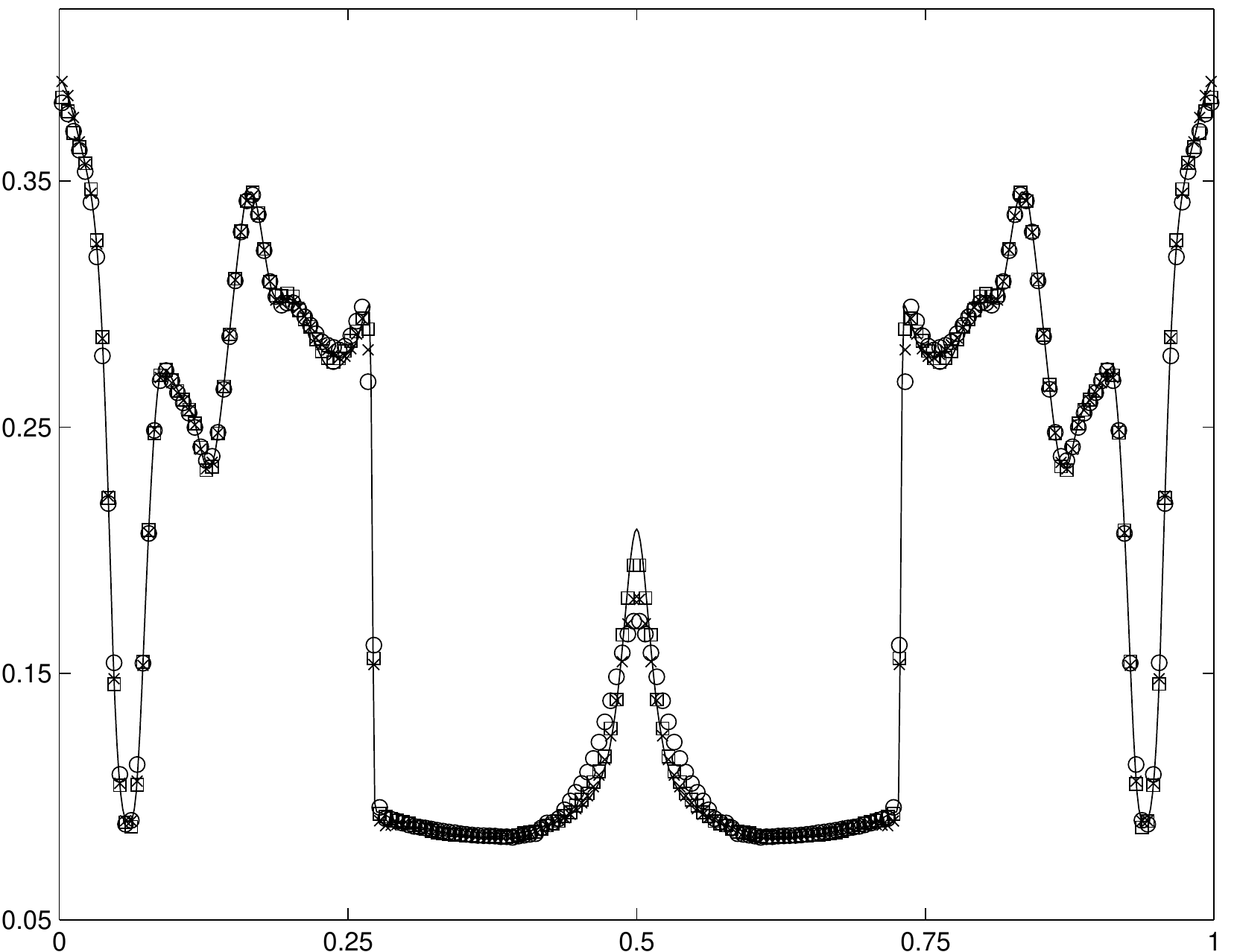}&
        \includegraphics[width=0.45\textwidth]{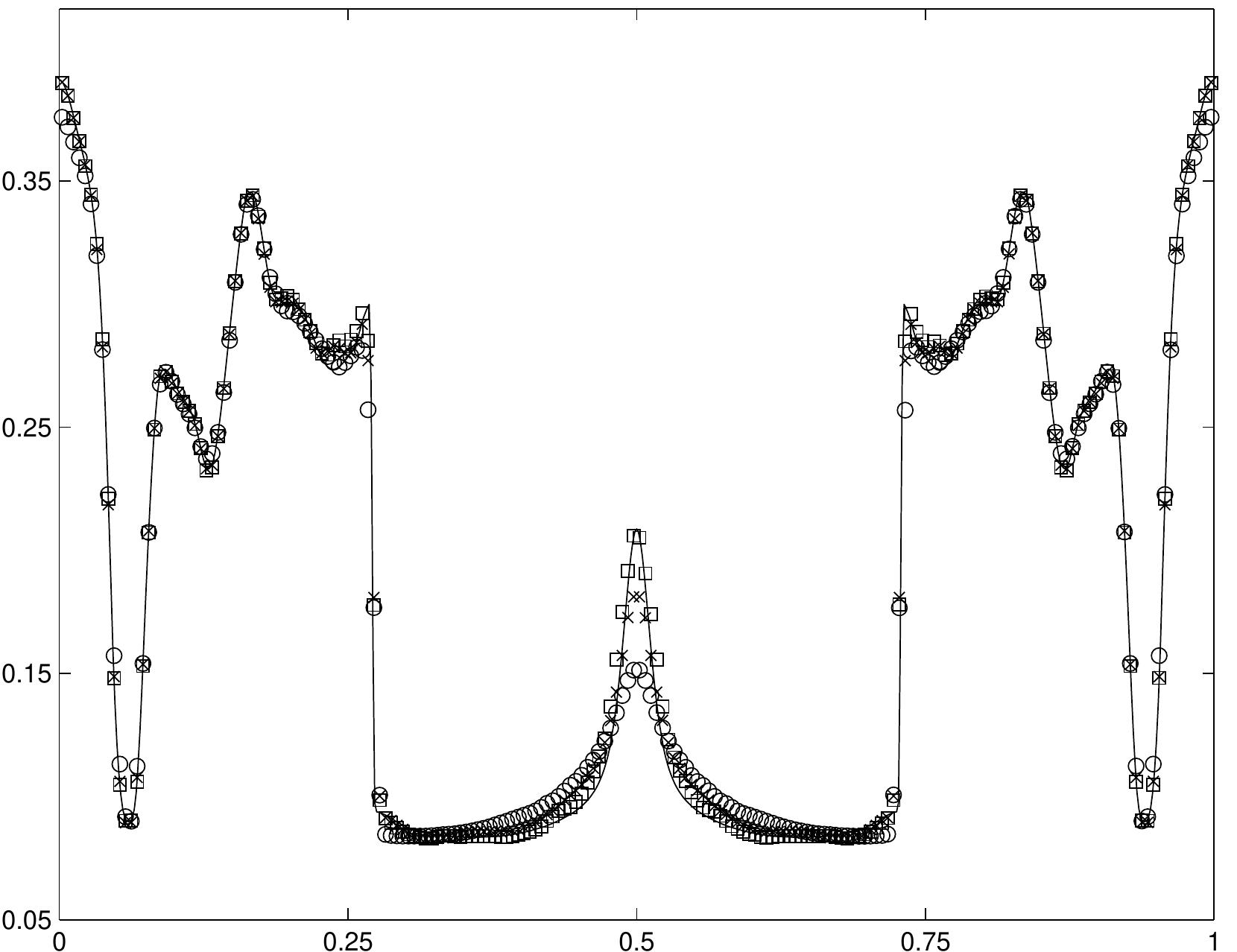}\\
         \includegraphics[width=0.45\textwidth]{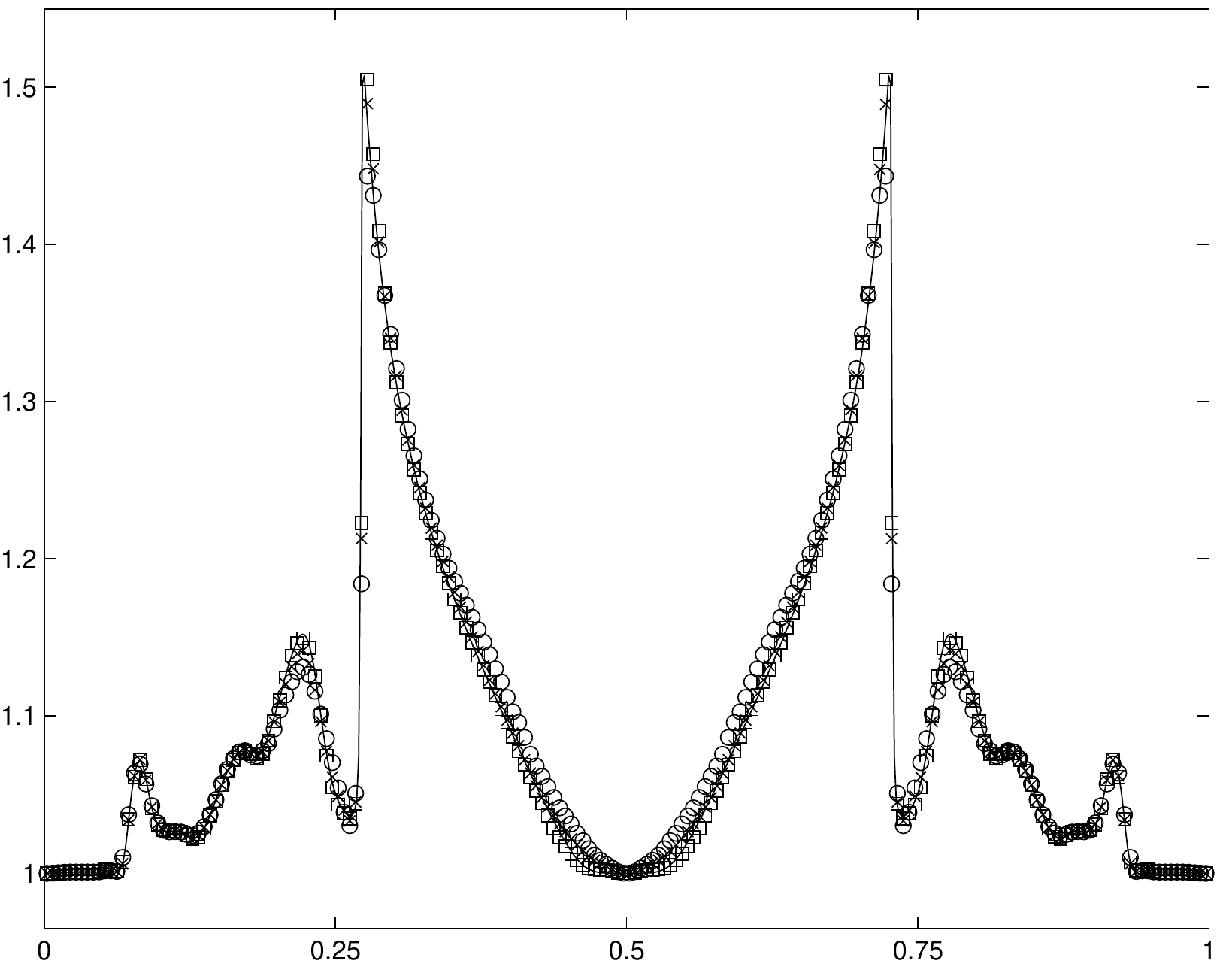}&
        \includegraphics[width=0.45\textwidth]{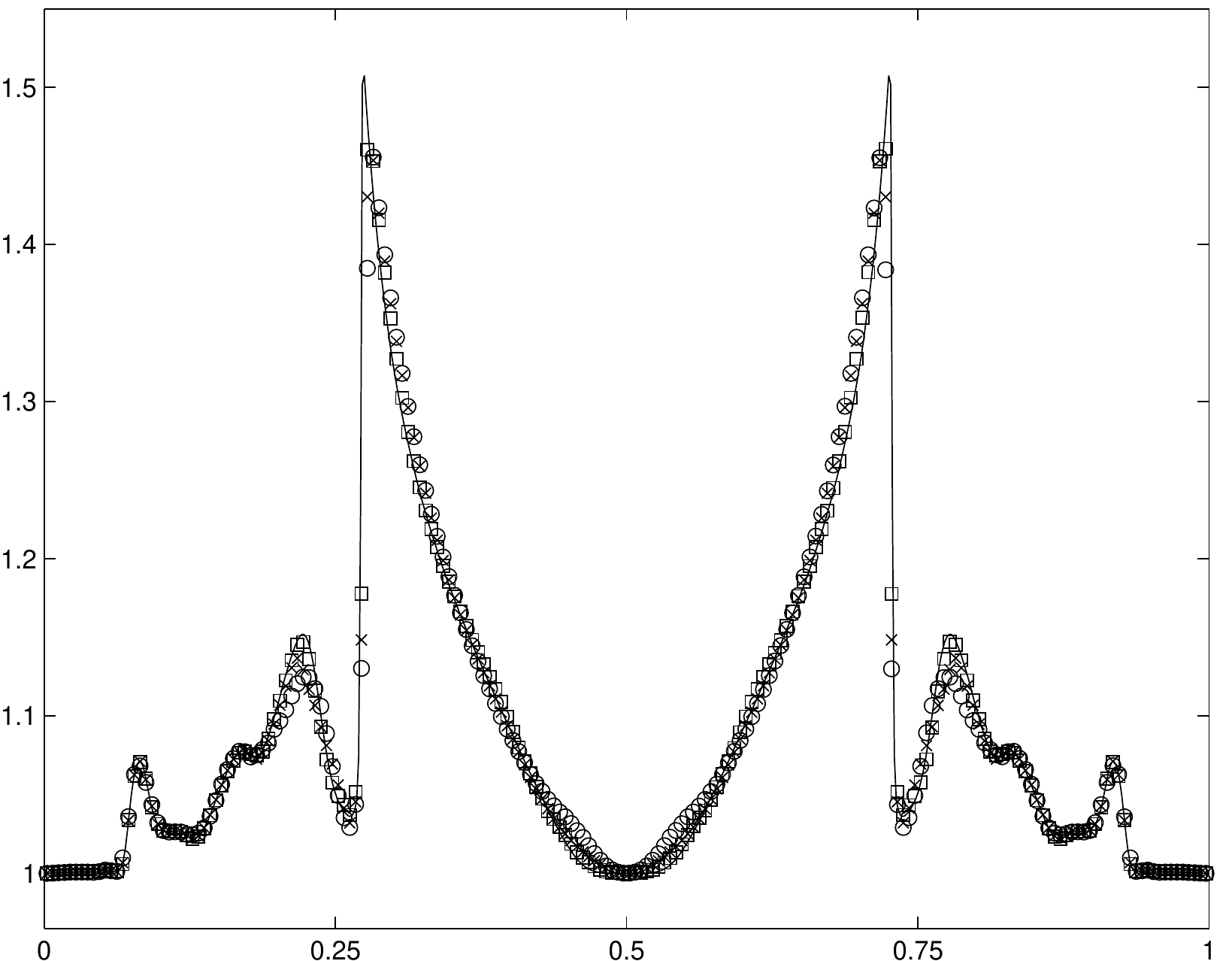}
        \end{tabular}
    \caption{Example \ref{exRMHDOT}:
    The densities $\rho$ (top) and Lorentz factors $\gamma$ (bottom) at $t=1$ along the line $y=1-x$.
    The solid line denotes the reference solution obtained by using the MUSCL scheme with $600\times 600$ cells,
    while the symbols ``$\circ$'', ``$\times$'', and ``$\square$''
denote the numerical solutions obtained by using the $P^1$-, $P^2$-, and $P^3$ methods with $200\times200$
cells. Left: \DG{}; right: \CDG{}. }
    \label{fig:OTyxrhogam}
  \end{figure}

\begin{figure}[!htbp]
    \centering{}
  \begin{tabular}{cc}
    \includegraphics[width=0.35\textwidth]{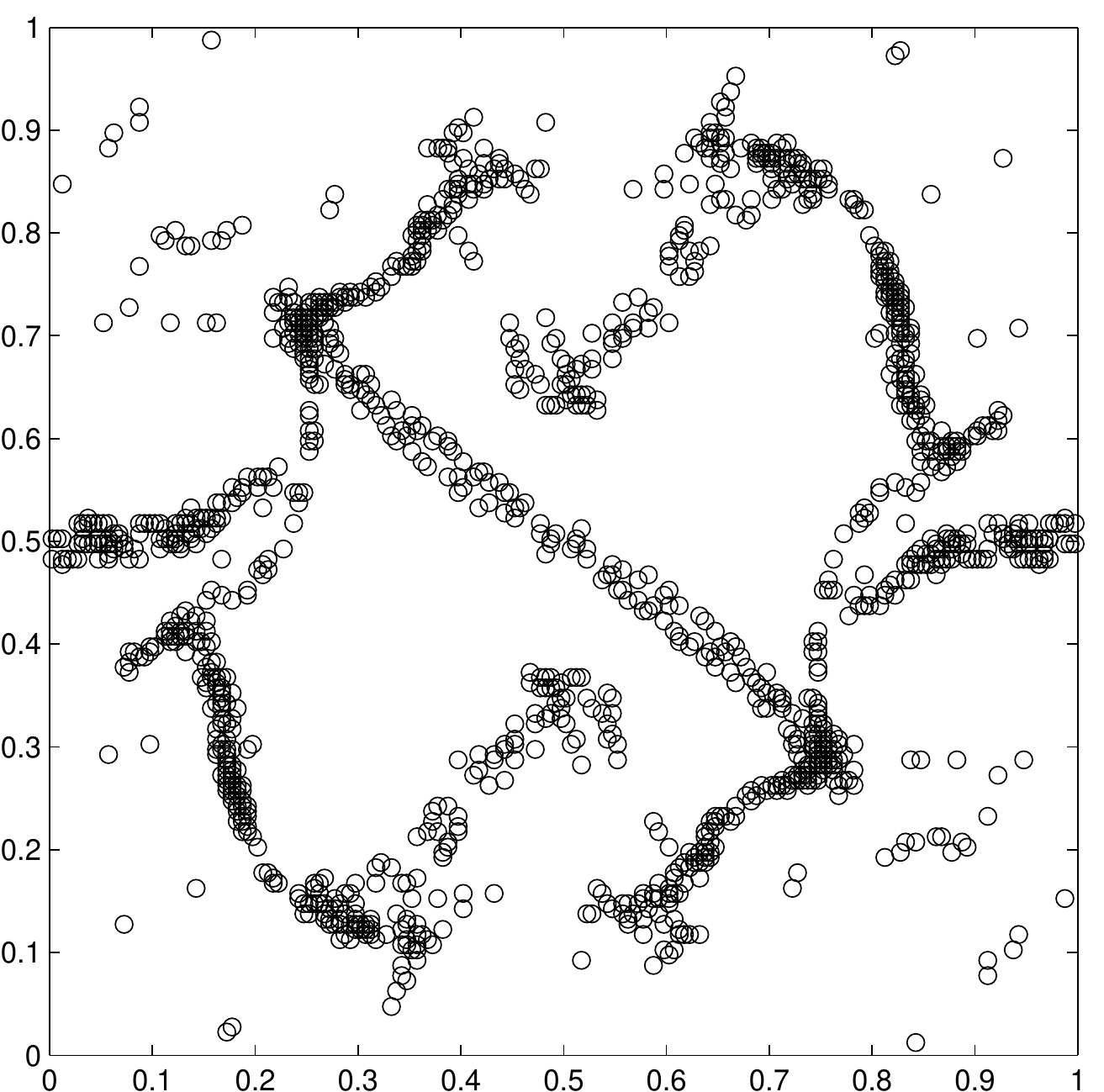}&
 \includegraphics[width=0.35\textwidth]{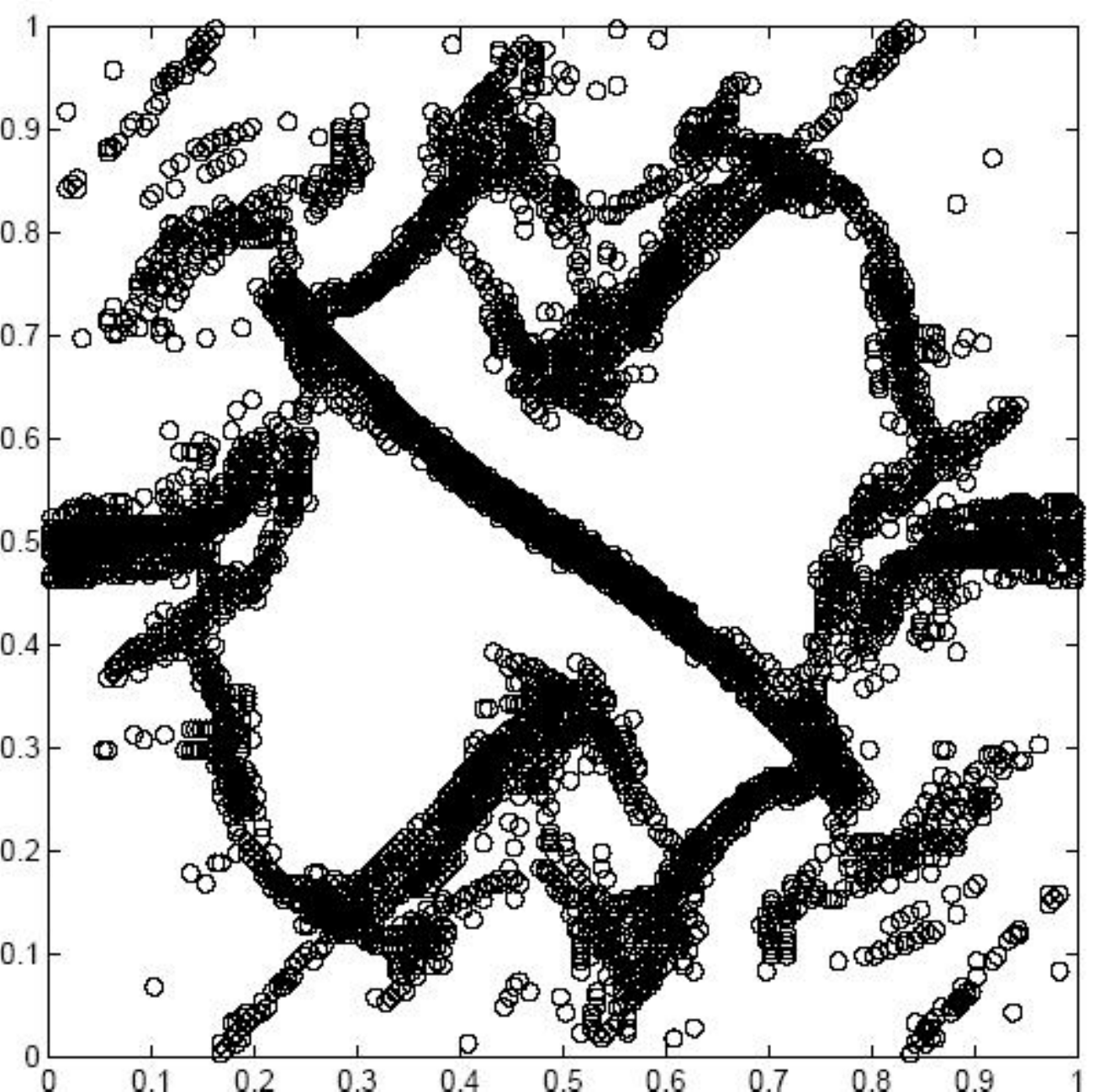}\\
    \includegraphics[width=0.35\textwidth]{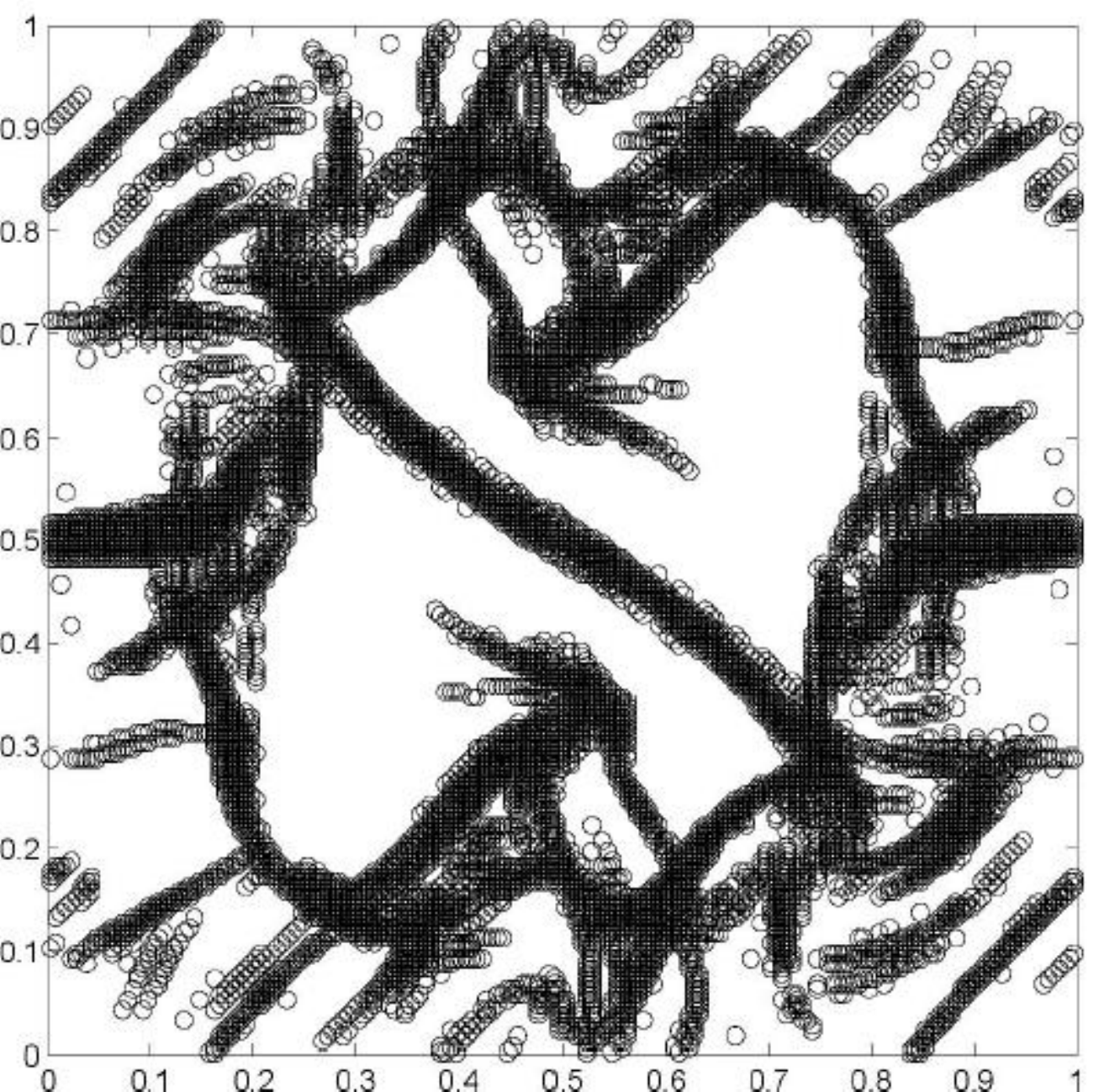}&
 \includegraphics[width=0.35\textwidth]{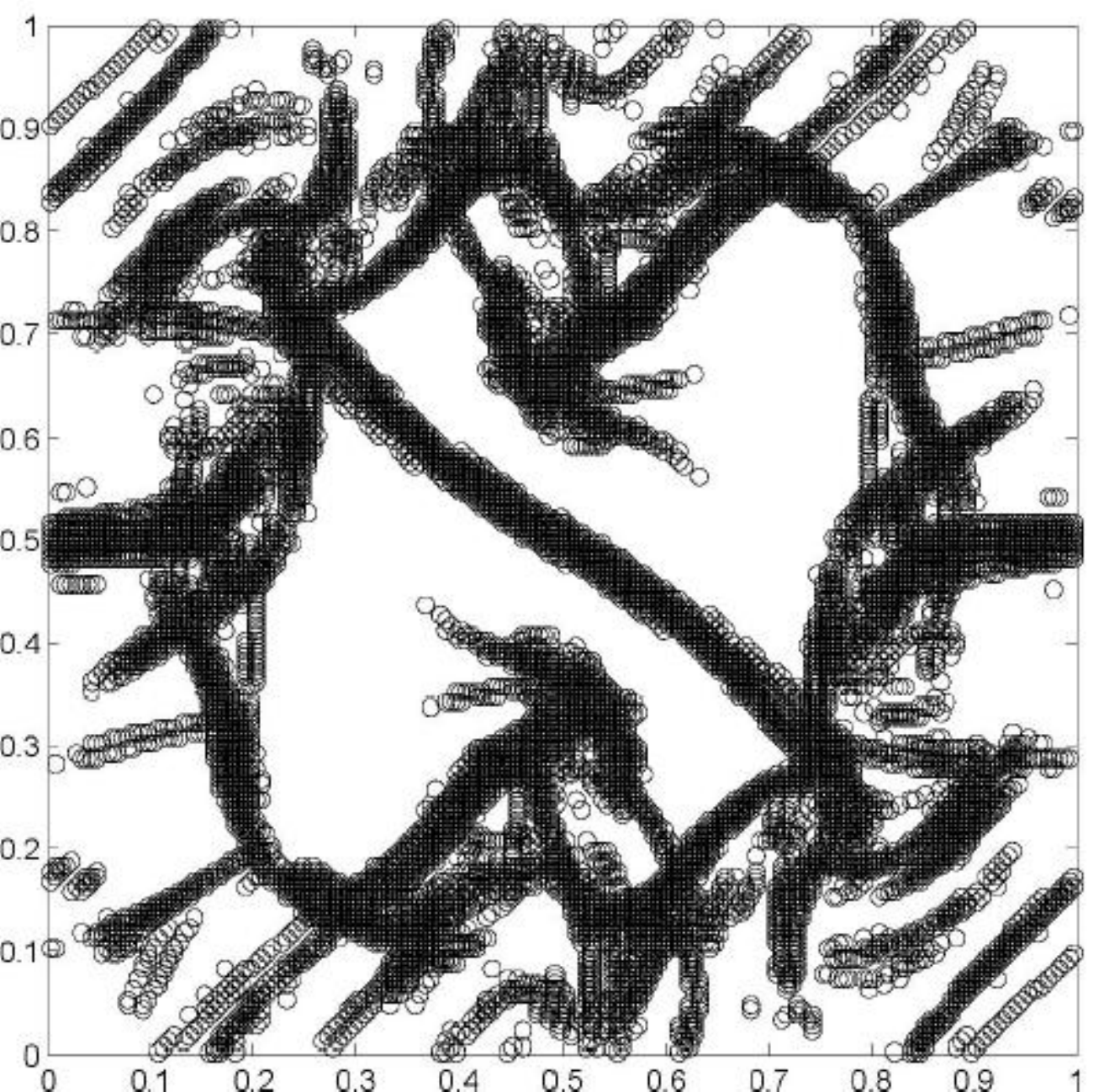}\\
 \includegraphics[width=0.35\textwidth]{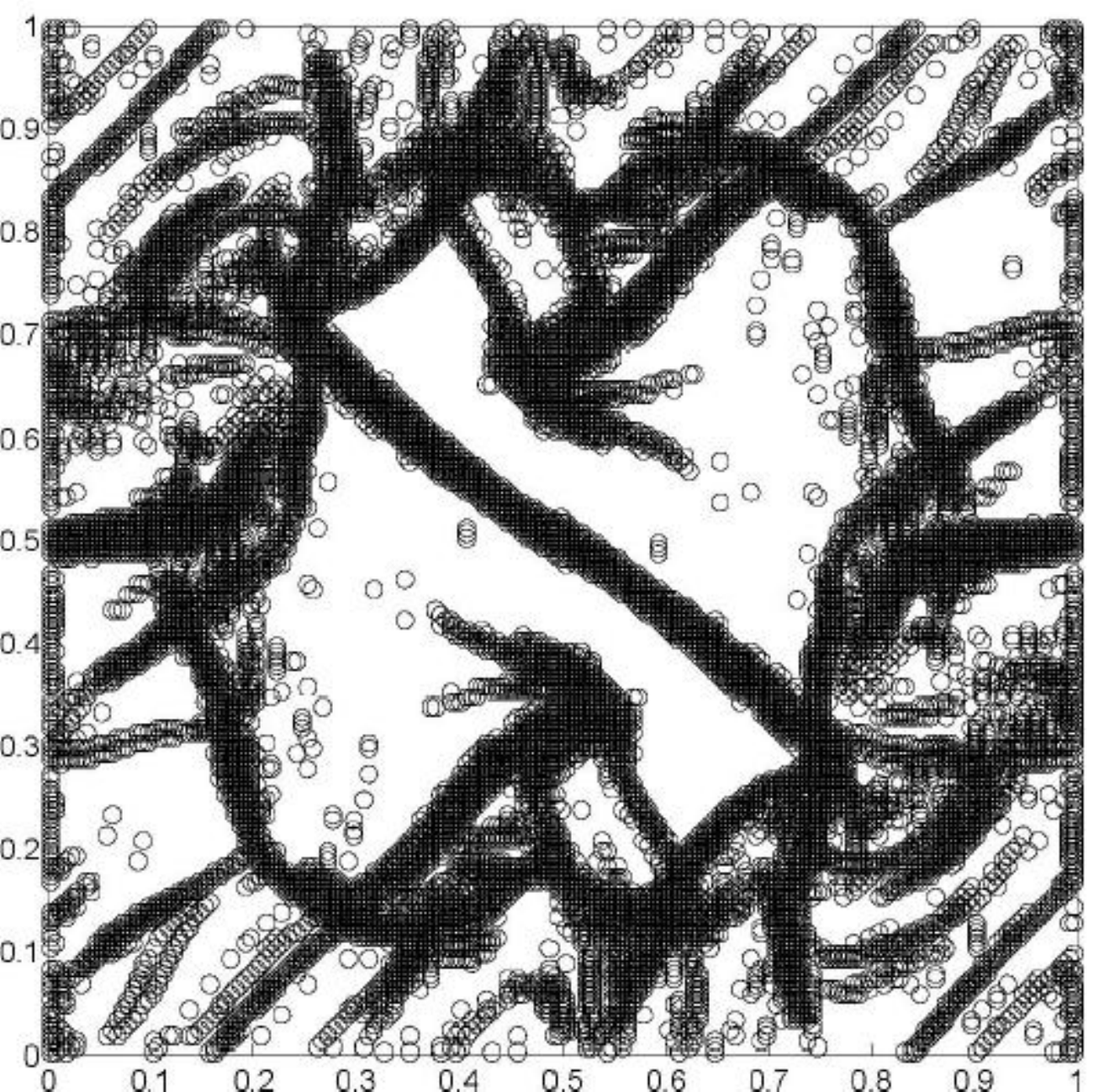}&
  \includegraphics[width=0.35\textwidth]{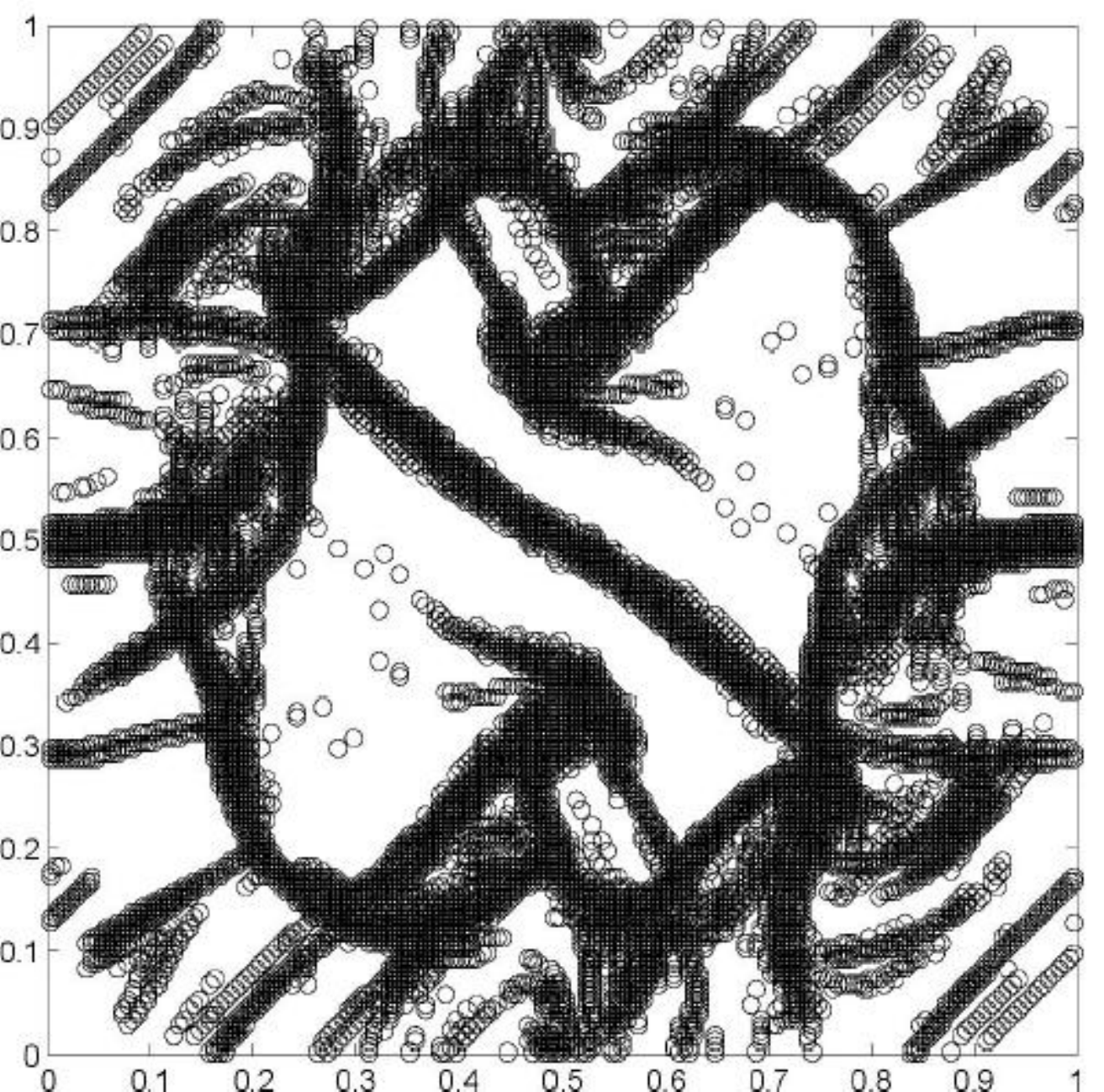}\\
    \end{tabular}
    \caption{Same as Fig.~\ref{fig:RMHDOTrho} except for the ``troubled'' cells. 
    }
    \label{fig:RMHDOTcell}
  \end{figure}

\begin{Example}[Orszag-Tang problem]\label{exRMHDOT}\rm
This RMHD Orszag-Tang problem is an extension of the non-relativistic version.
The computational domain is chosen as $\Omega=[0,1]\times[0,1]$, and the initial data are
\begin{align*}
&\rho(x,y,0)=\frac{25}{36\pi},\ \ v_x(x,y,0)=0.5\sin(2\pi y),\ \
v_y(x,y,0)=0.5\sin(2\pi x),\\
&v_z(x,y,0)=0,\ \ B_x(x,y,0)=-\frac{1}{\sqrt{4\pi}}\sin(2\pi y),\ \
B_y(x,y,0)=\frac{1}{\sqrt{4\pi}}\sin(4\pi x),\\
&B_z(x,y,0)=0,\ \ p(x,y,0)=\frac{5}{12\pi},
\end{align*}
with the adiabatic index $\Gamma=5/3$.
The   solution of problem is smooth initially, but the complicated wave
structure is formed as the time increases, and it has the turbulence behavior.

  The CFL numbers of $P^1$-, $P^2$-, and $P^3$-based non-central DG methods
   are $0.2,~0.15,~0.1$, respectively, while those of  corresponding
 central DG methods are chosen as $0.3,~0.25,~0.2$, respectively, and
 $\theta=\Delta t_n/\tau_n=1$.

    Figs.~\ref{fig:RMHDOTrho} and \ref{fig:RMHDOTgam} show
    the contours of
      densities $\rho$ and Lorentz factors $\gamma$
    at $t=1$ obtained by using the proposed DG methods with $200\times 200$ cells.
    Fig.~\ref{fig:OTyxrhogam} plots
    them and the reference solutions along the line $y=1-x$,
     where the reference solutions are obtained by using the MUSCL scheme with $600\times 600$ uniform
     cells.
     It is seen that the solutions obtained by higher-order DG methods are in good
     agreement with the reference solutions. The distribution of identified ``troubled'' cells
     is also consistent with the solutions, see Fig.~\ref{fig:RMHDOTcell}.
  \end{Example}


\begin{figure}[!htbp]
    \centering{}
  \begin{tabular}{cc}
    \includegraphics[width=0.3\textwidth]{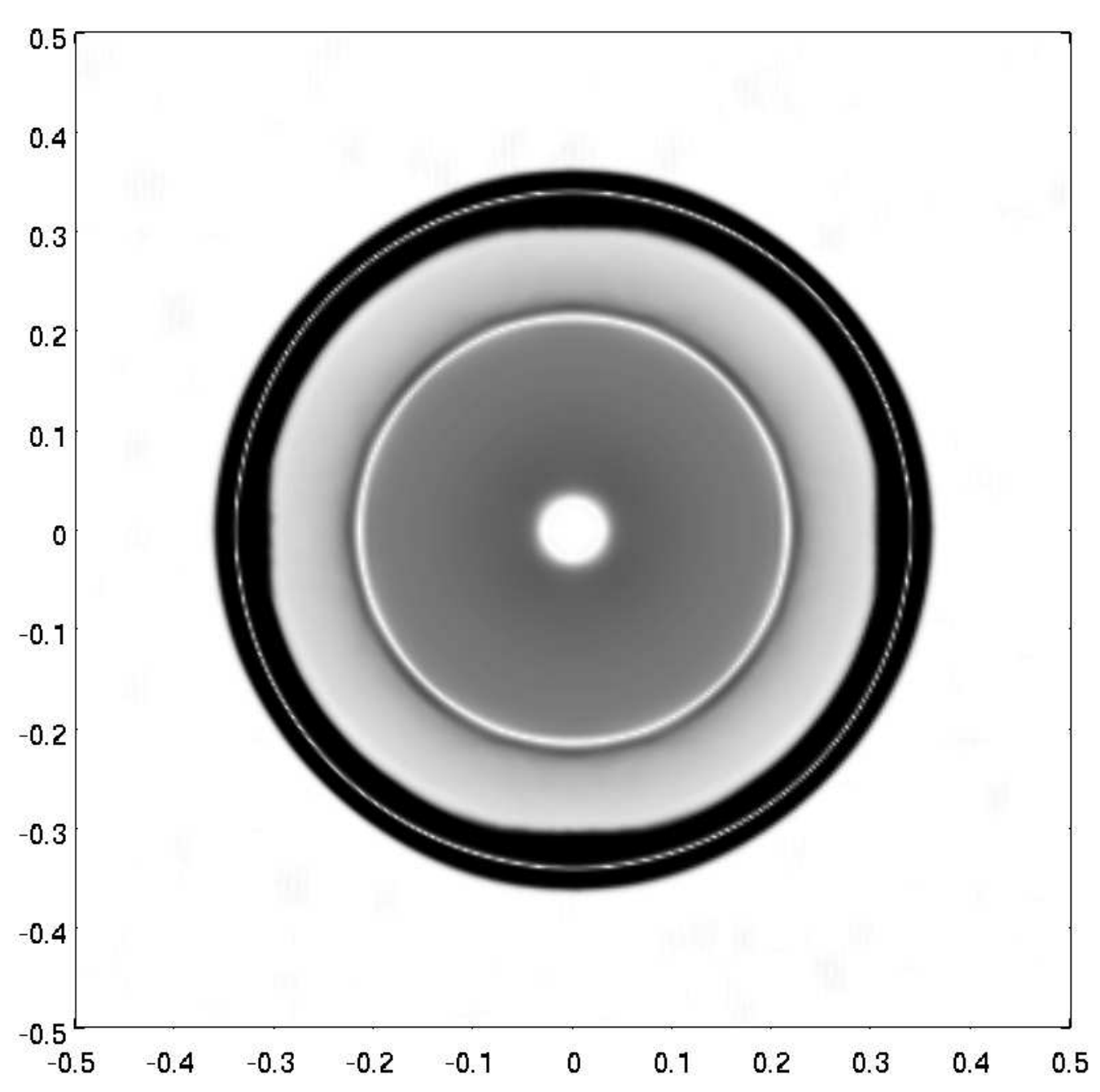}&
        \includegraphics[width=0.3\textwidth]{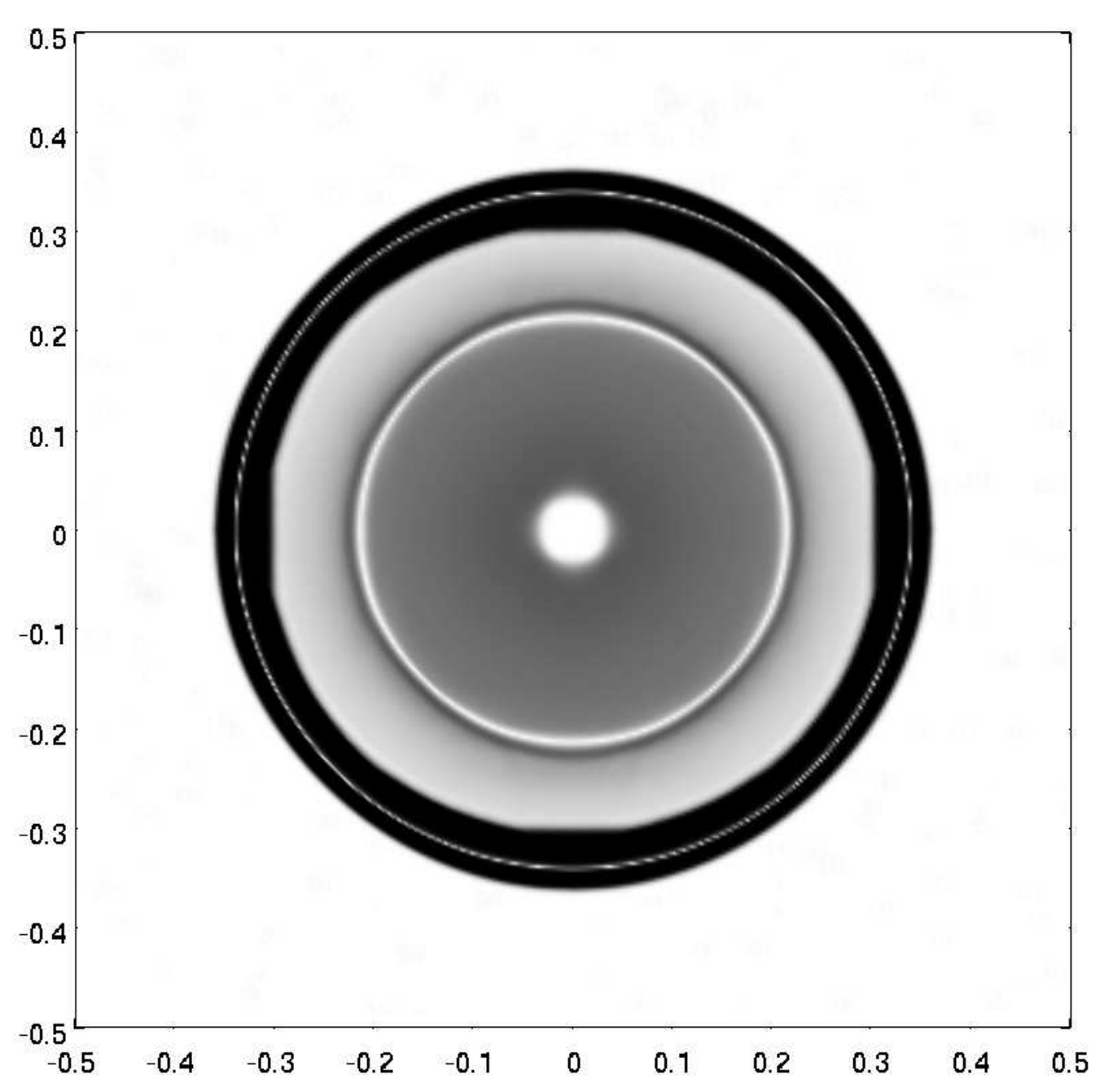}\\
        \includegraphics[width=0.3\textwidth]{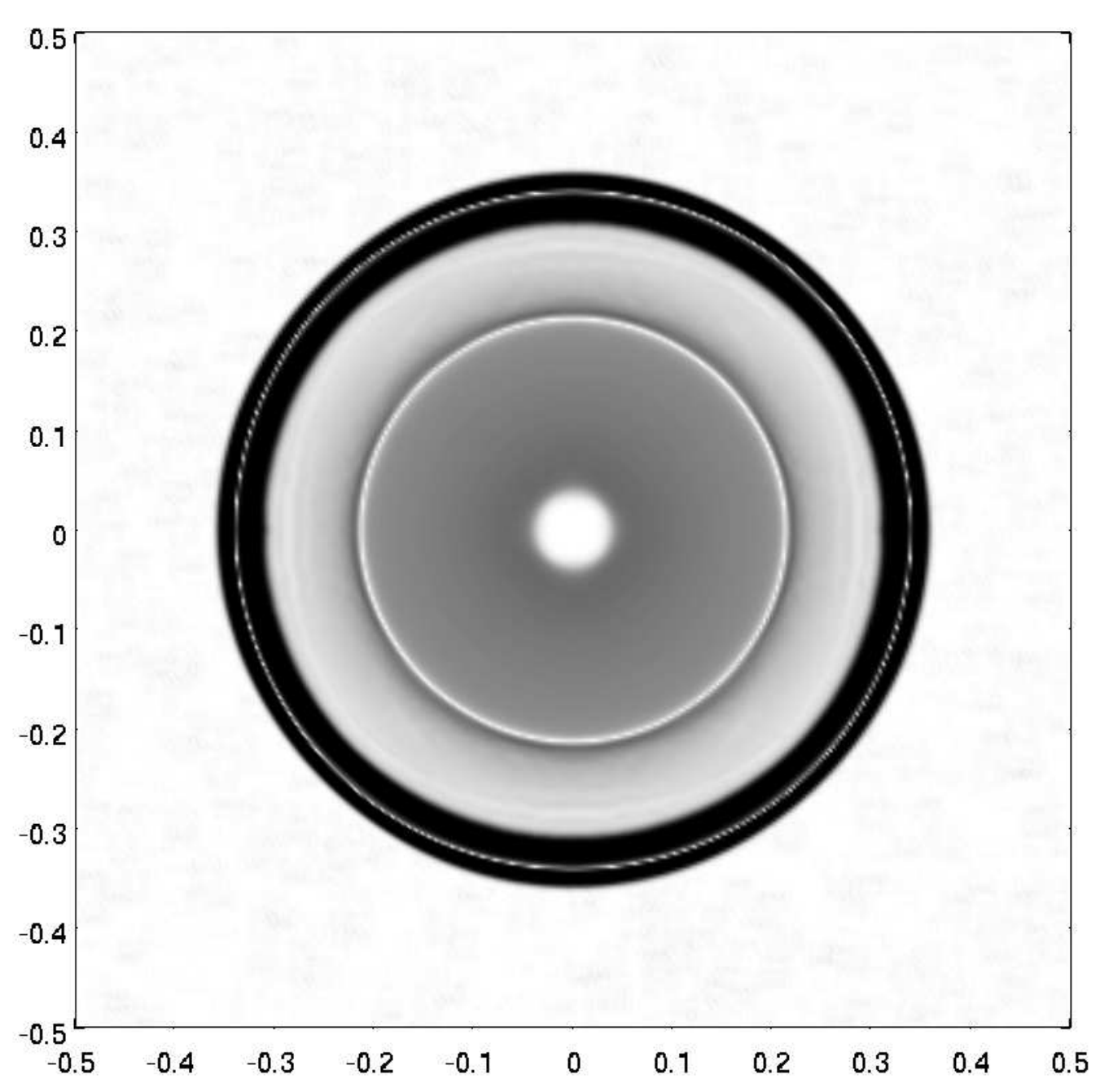}&
        \includegraphics[width=0.3\textwidth]{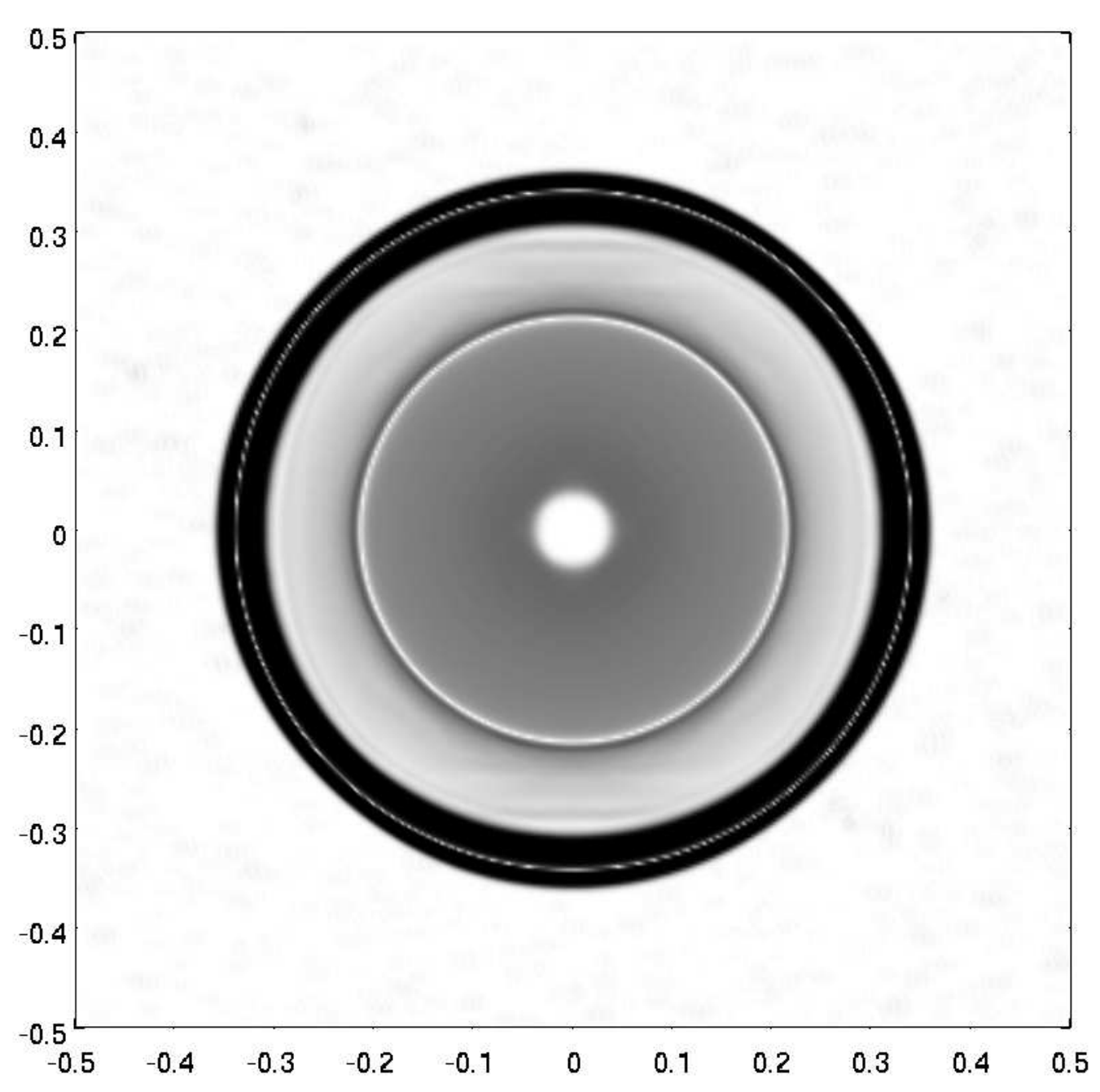}\\
        \includegraphics[width=0.3\textwidth]{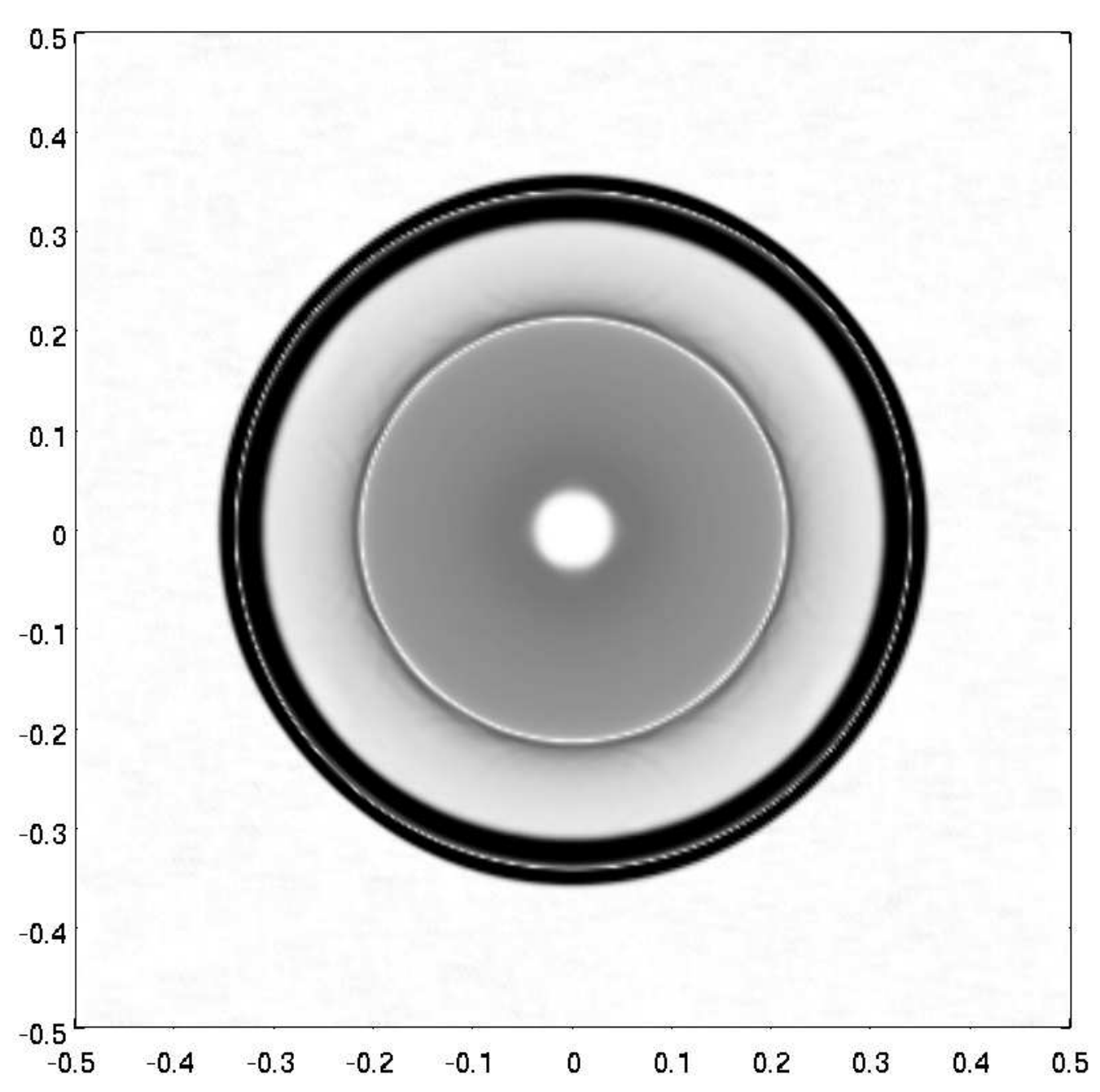}&
        \includegraphics[width=0.3\textwidth]{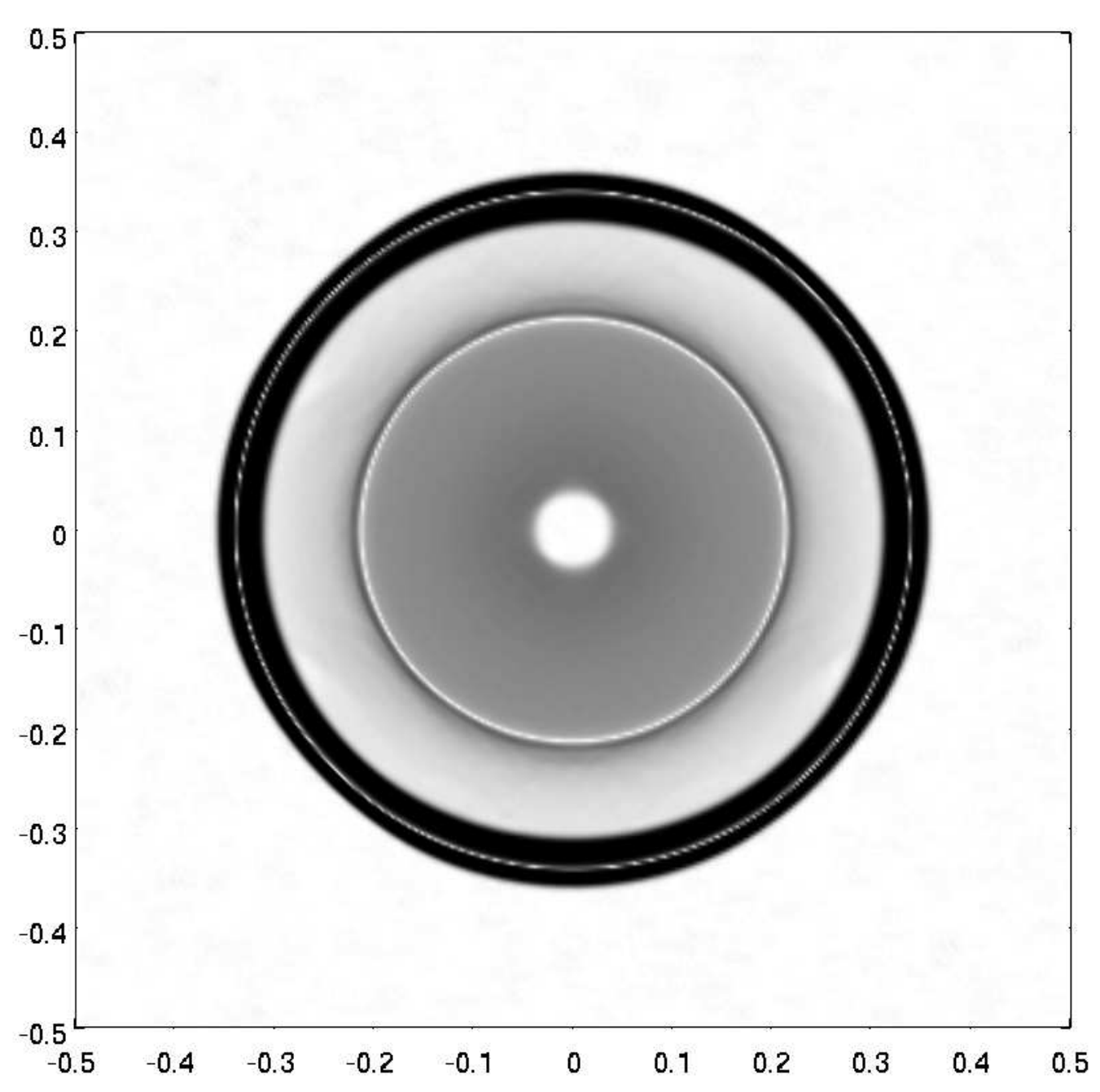}\\
    \end{tabular}
    \caption{Example \ref{exRMHDBlast}:
         Schlieren images of density $\rho$ at $t=0.3$ obtained with $300\times 300$
cells. Left: $P^K$-based RKDG methods; right: $P^K$-based Runge-Kutta CDG methods. From
top to bottom: $K = 1, 2, 3$.
 }
    \label{fig:RMHDBLrho}
  \end{figure}

\begin{figure}[!htbp]
    \centering{}
  \begin{tabular}{cc}
    \includegraphics[width=0.3\textwidth]{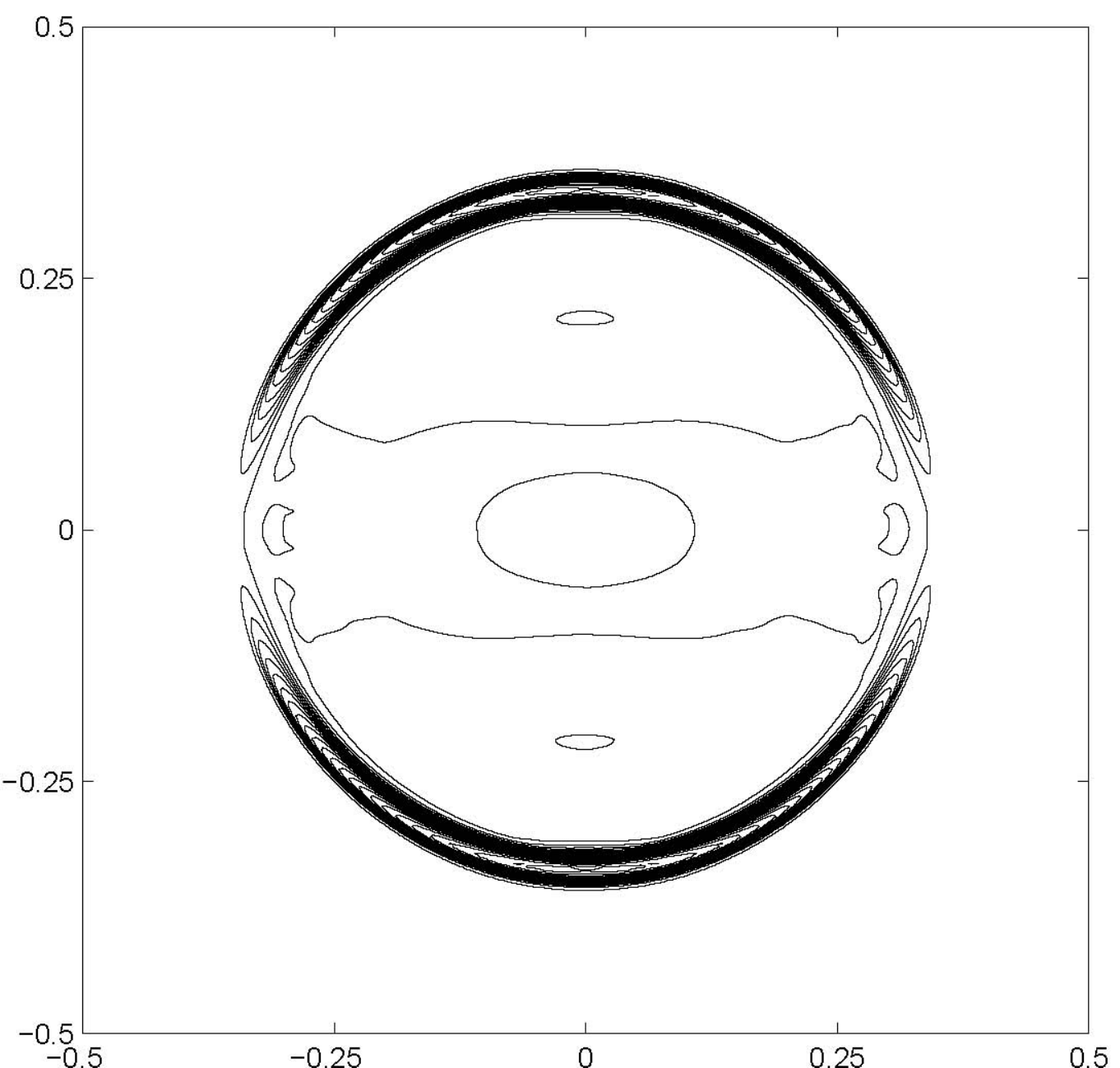}&
        \includegraphics[width=0.3\textwidth]{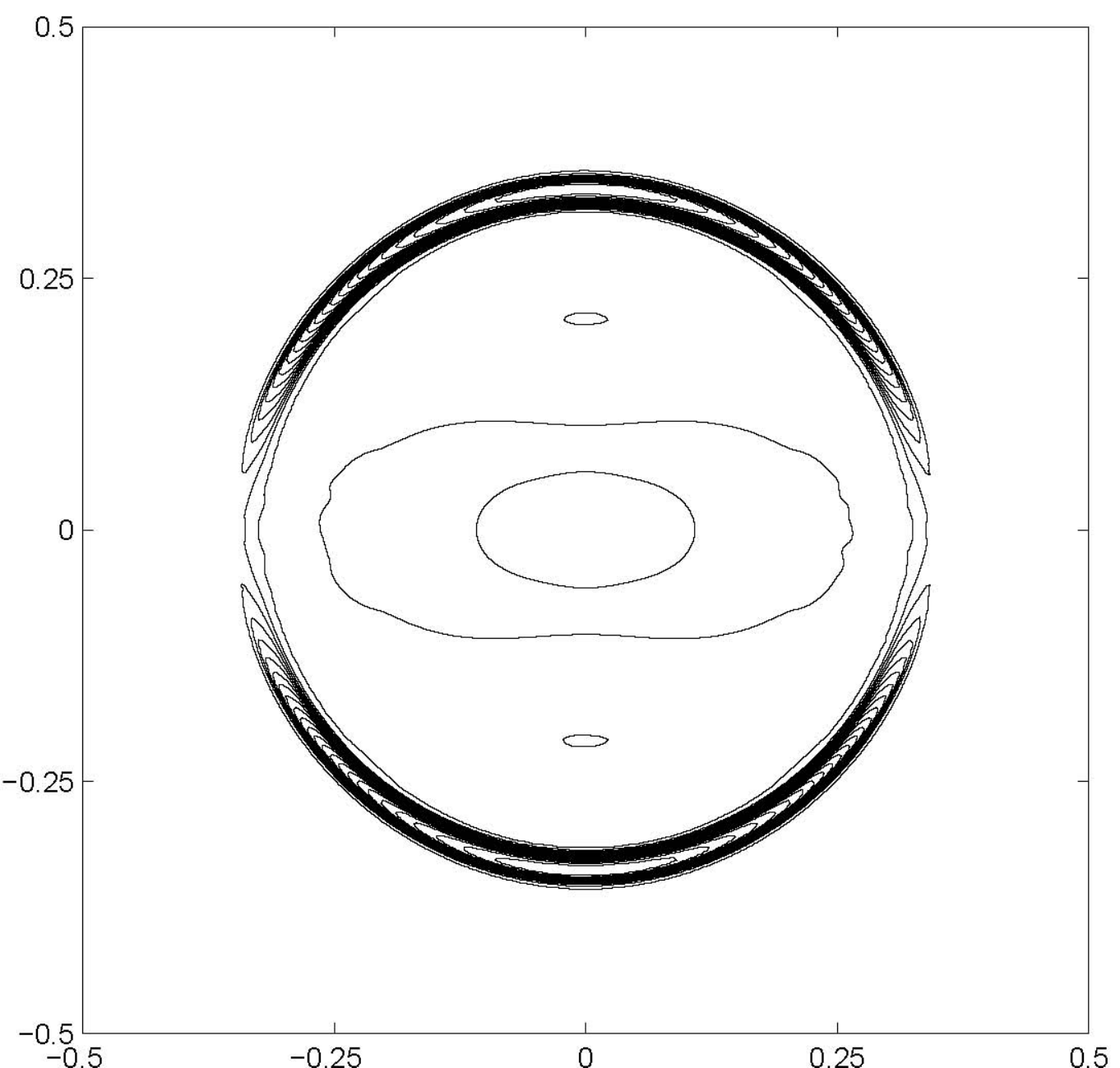}\\
        \includegraphics[width=0.3\textwidth]{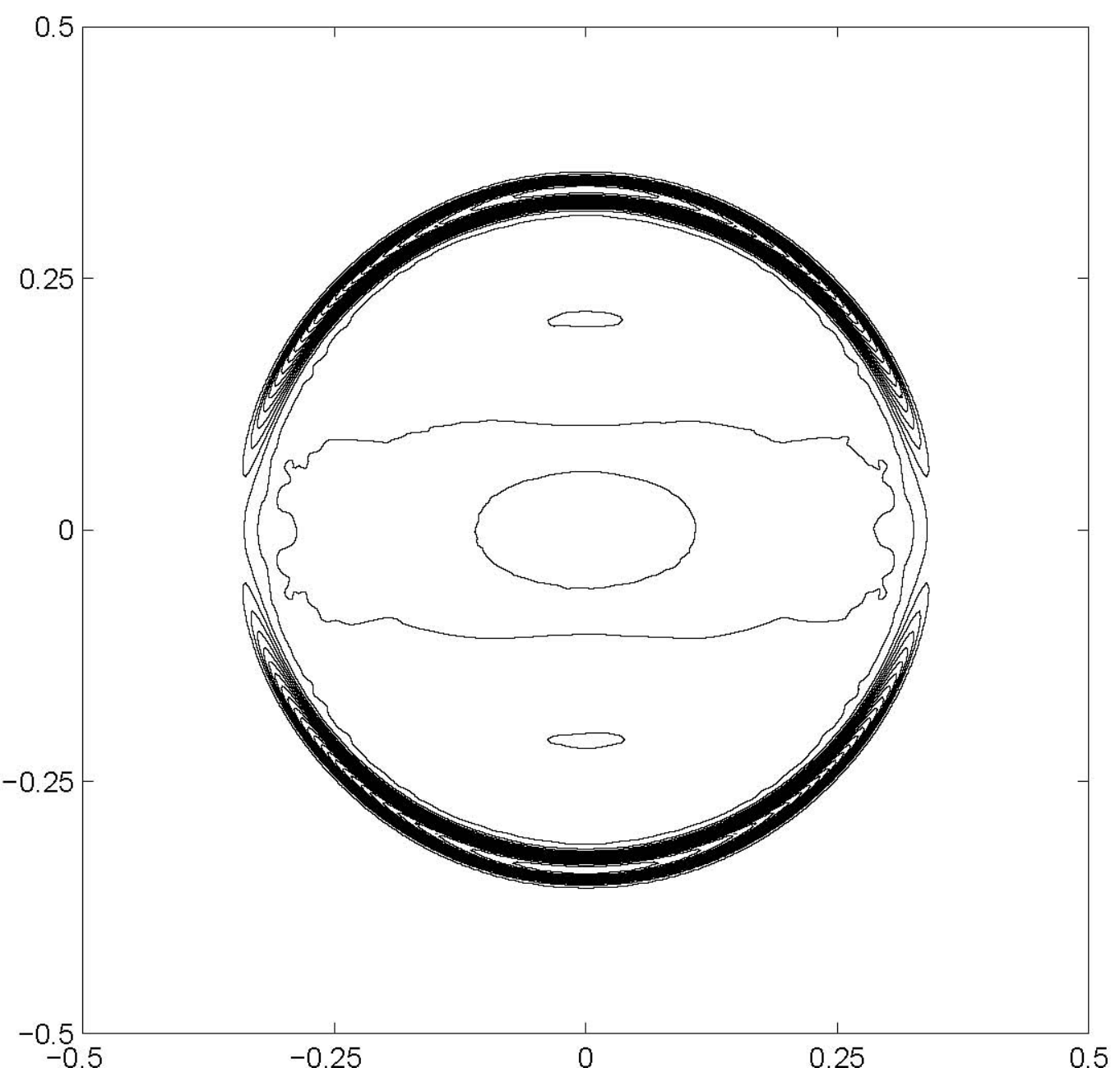}&
        \includegraphics[width=0.3\textwidth]{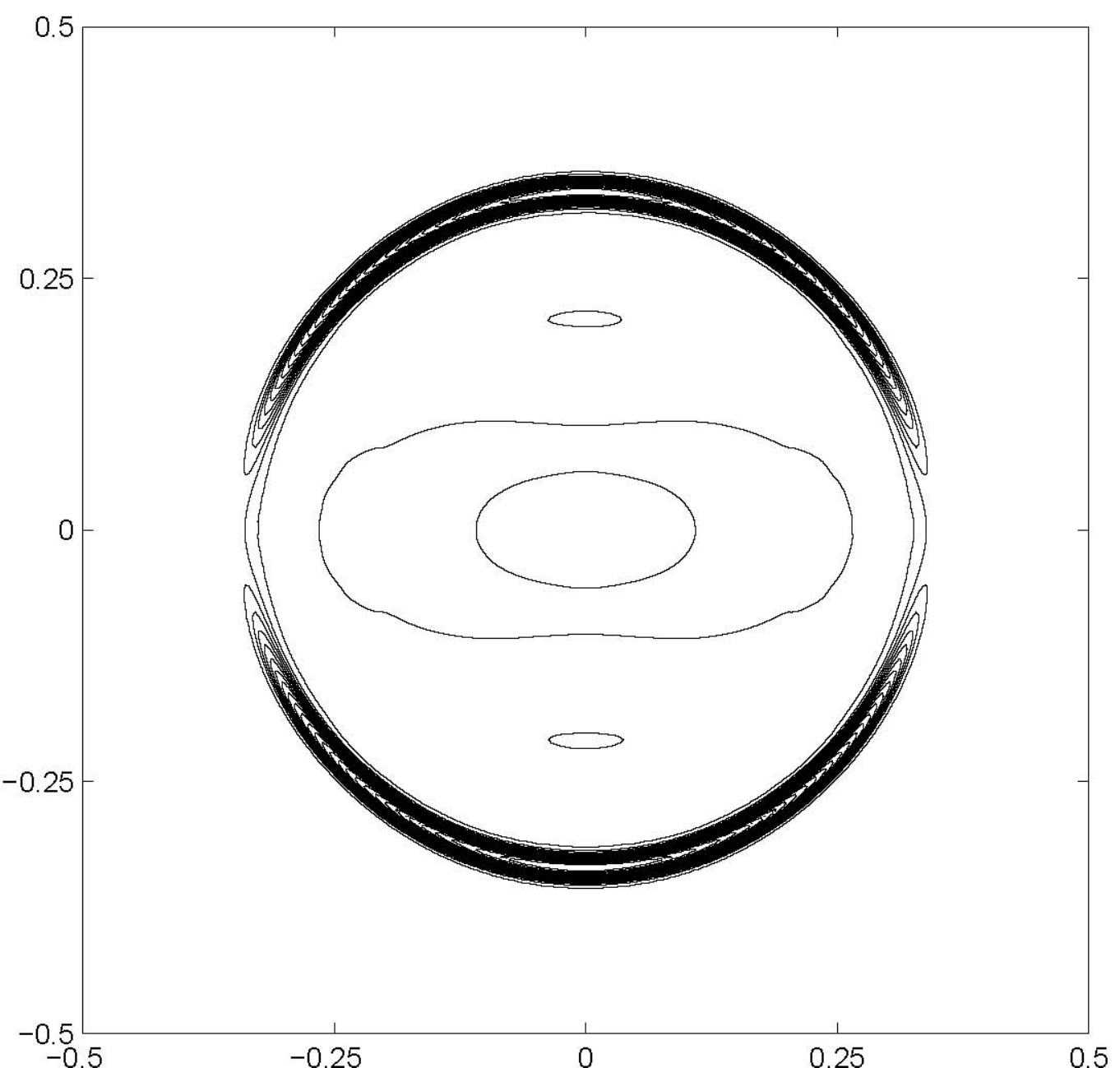}\\
        \includegraphics[width=0.3\textwidth]{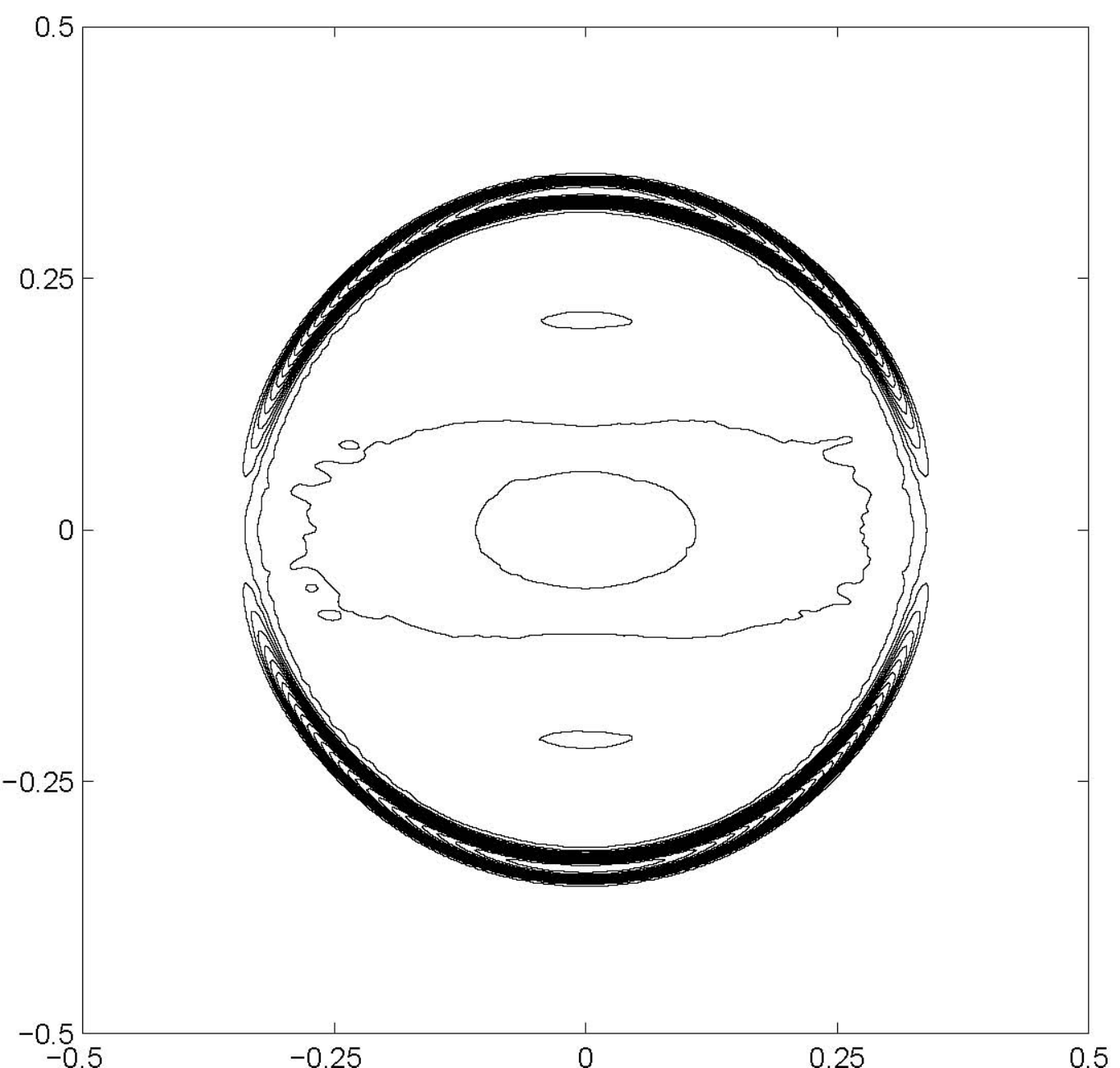}&
        \includegraphics[width=0.3\textwidth]{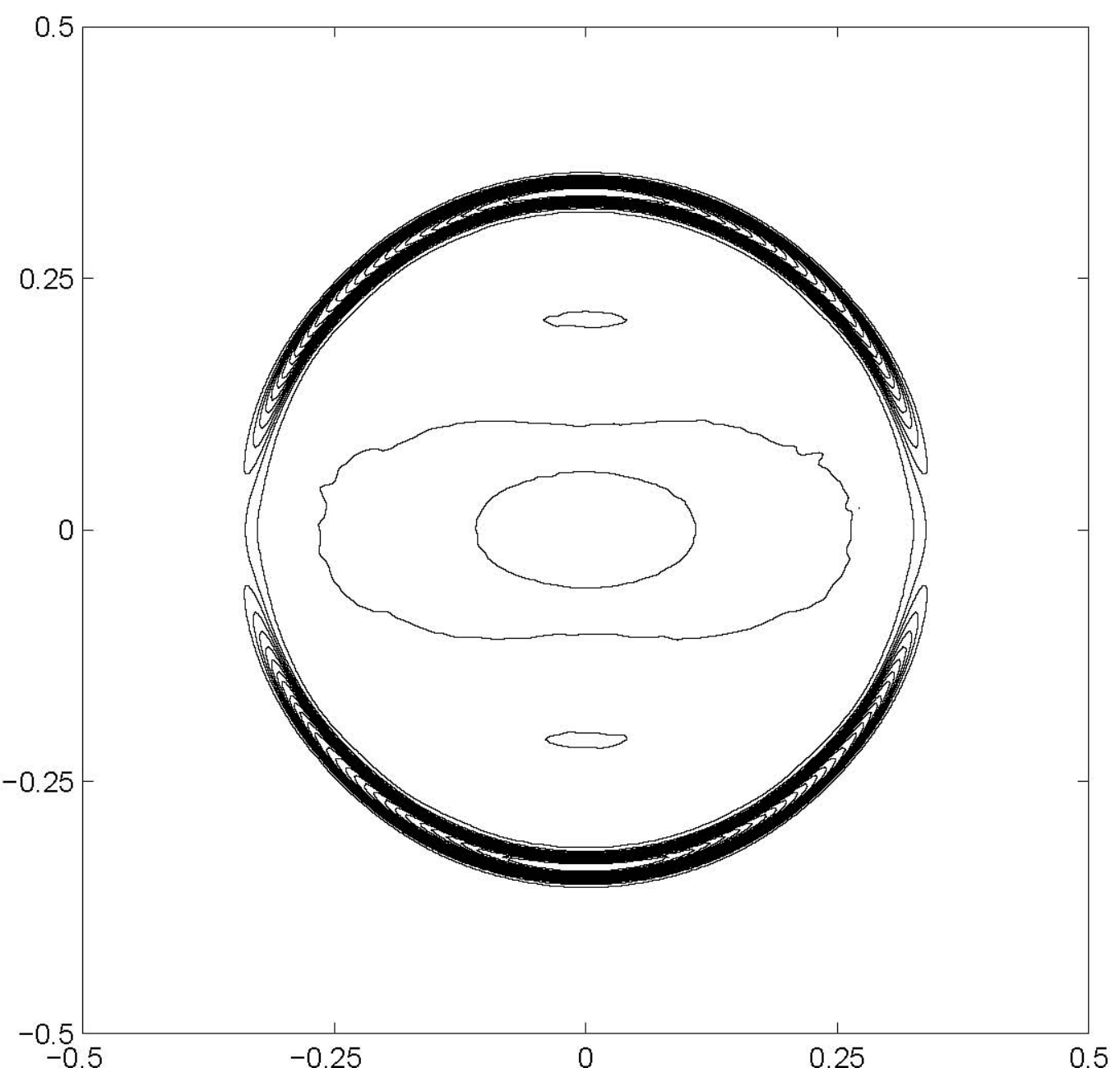}\\
    \end{tabular}
    \caption{Same as Fig.~\ref{fig:RMHDBLrho} except for the
    magnetic field $B_x$  (15 equally spaced contour lines from 0.02 to 0.35).}
    \label{fig:RMHDBLbx}
  \end{figure}

\begin{figure}[!htbp]
    \centering{}
    \begin{tabular}{cc}
    \includegraphics[width=0.3\textwidth]{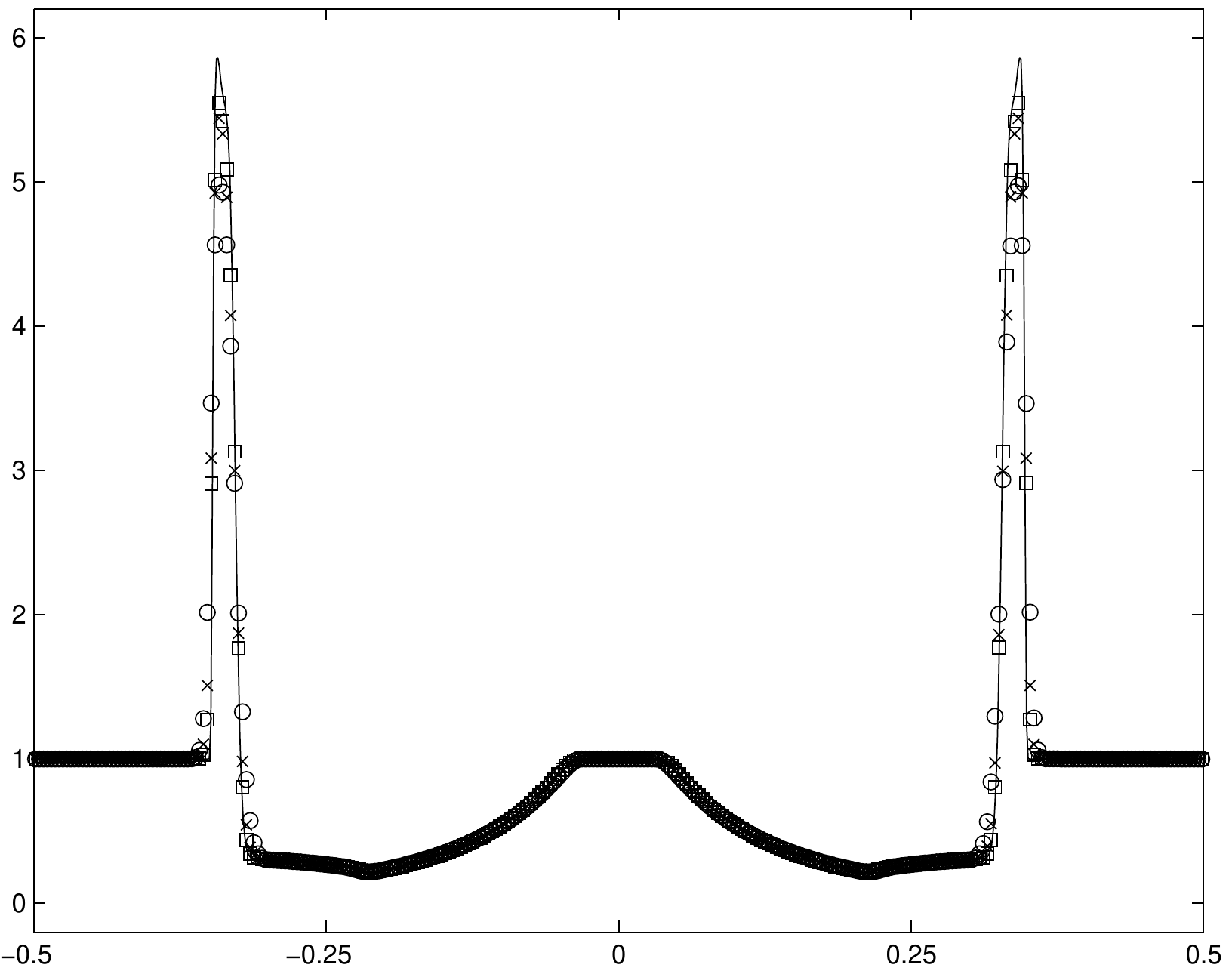}&
        \includegraphics[width=0.3\textwidth]{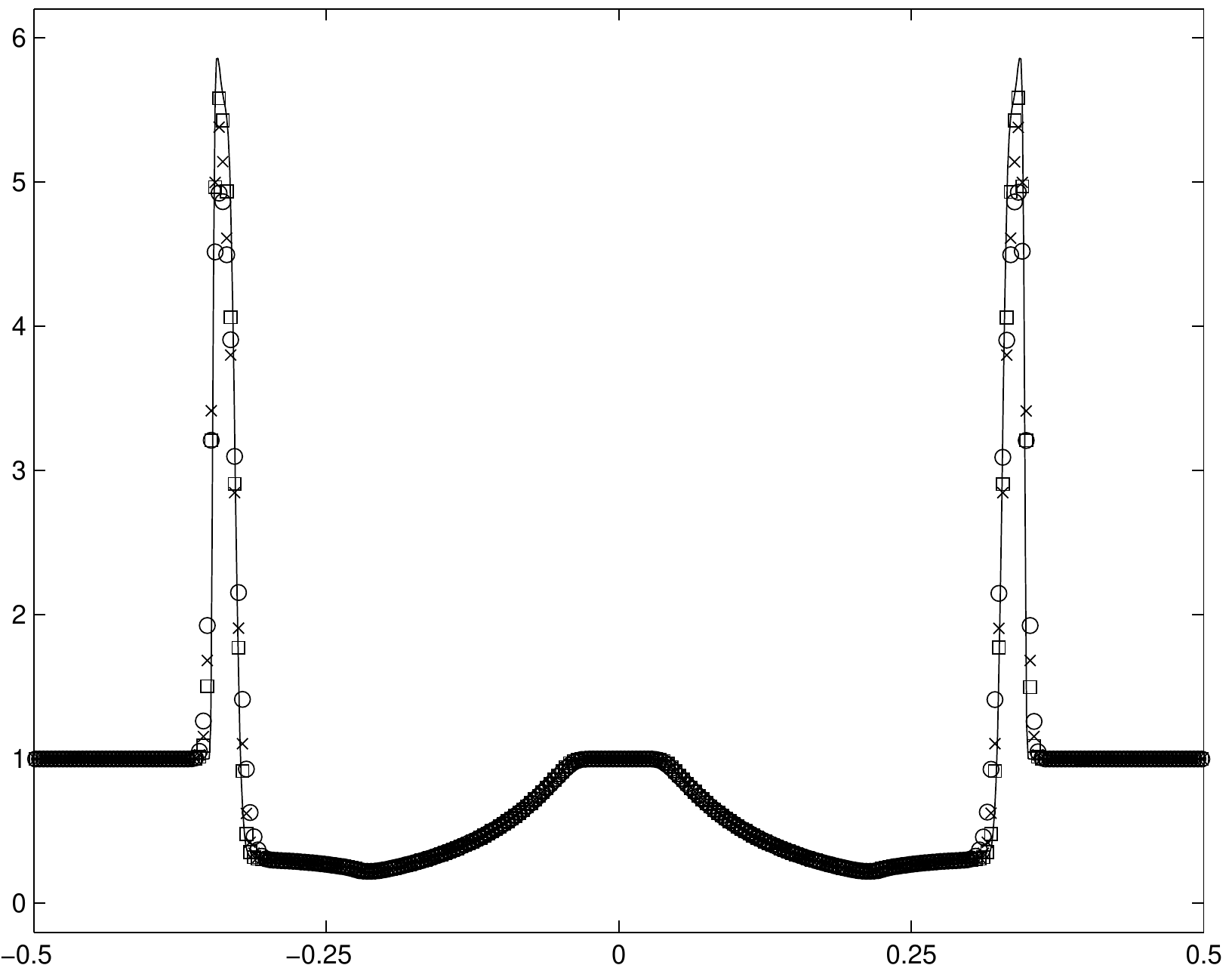}
        \end{tabular}
    \caption{Example \ref{exRMHDBlast}:
The density   $\rho$ at
    $t=0.3$ along the line $x=0$. The solid line denote the reference solution obtained by using the MUSCL scheme
    with $800\times 800$ cells, while the symbol ``$\circ$'', ``$\times$'', and ``$\square$'' denote
    the solutions obtained by using the $P^1$-, $P^2$-, and $P^3$-based methods with $300\times 300$ cells
    respectively. Left: \DG{}; right:\CDG{}. }
    \label{fig:BLcmprho}
  \end{figure}

\begin{figure}[!htbp]
    \centering{}
    \begin{tabular}{cc}
    \includegraphics[width=0.35\textwidth]{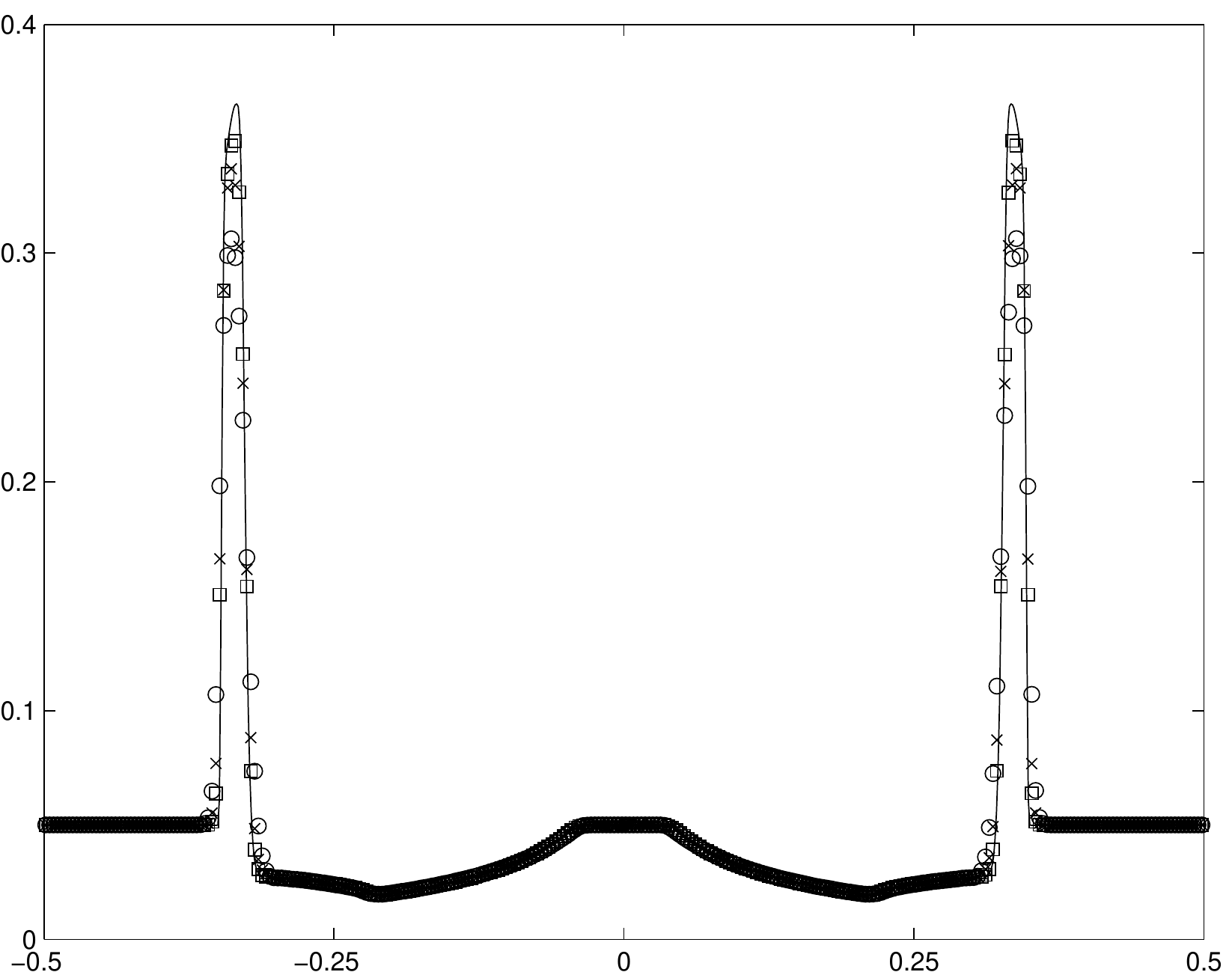}&
        \includegraphics[width=0.35\textwidth]{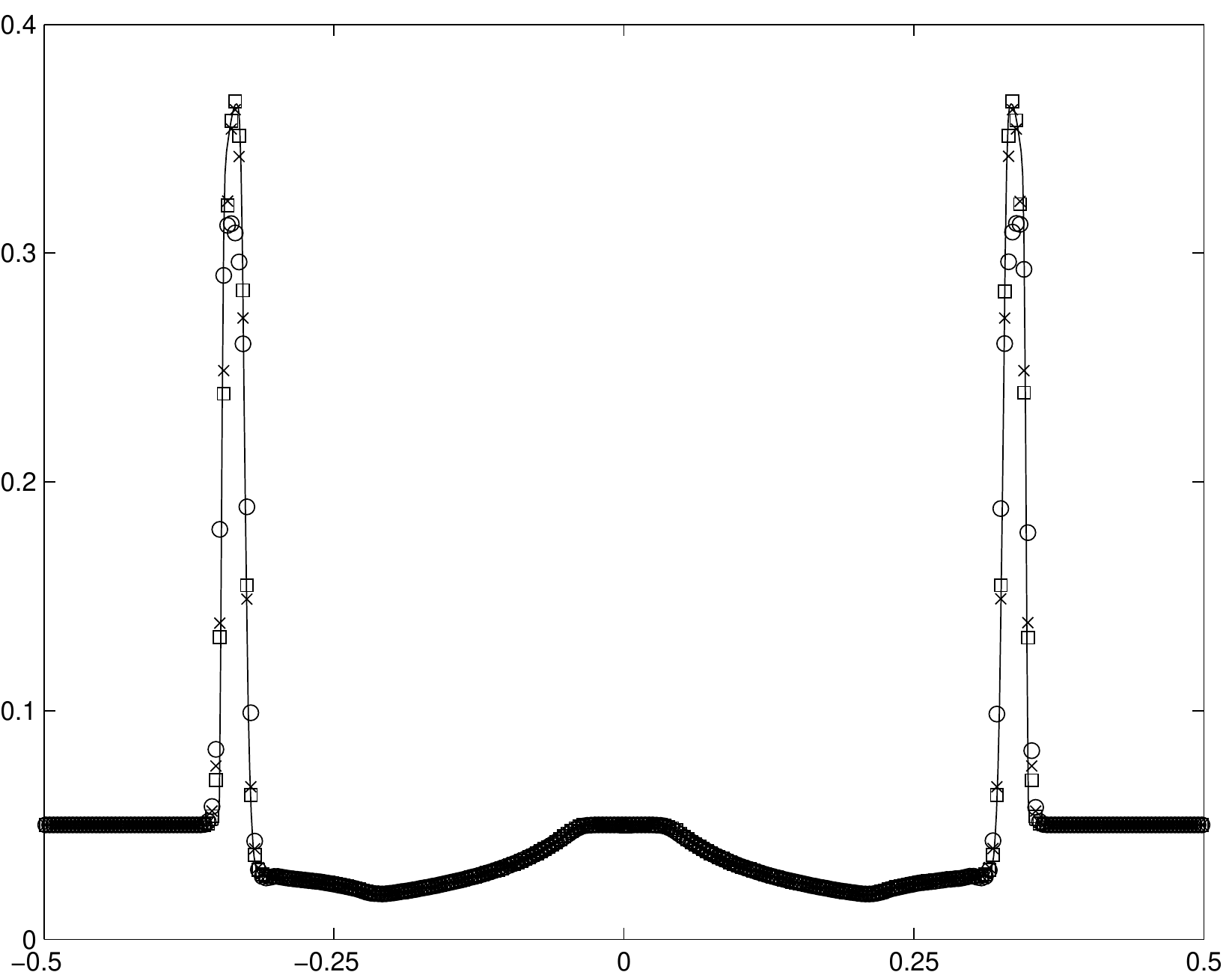}\\
        \end{tabular}
    \caption{Same as Fig.~\ref{fig:BLcmprho} except for the
    magnetic field $B_x$ along the line $x=0$.}
    \label{fig:BLcmpbx}
  \end{figure}

\begin{Example}[Blast problem]\label{exRMHDBlast}\rm
It has become  a very useful test for the multidimensional numerical schemes.
Our setup is the same as that in \cite{AderMood}.
The adiabatic index $\Gamma$ and  computational domain are
 taken as $4/3$ and $\Omega=[-0.5,0.5]\times[-0.5,0.5]$
with four outflow conditions, respectively. Initially, the magnetic field is constant,
 the fluid velocity is zero everywhere,  and
 there is an explosion region with radius one at the center of $\Omega$.
The detailed description of initial data is as follows
$$(\rho,v_x,v_y,v_z,B_x,B_y,B_z,p)=\begin{cases} (1,0,0,0,0.05,0,0,1), & r<0.2,\\
      (1,0,0,0,0.05,0,0,10^{-3}), & r\ge 0.2,\end{cases}$$
where  $r=\sqrt{x^2+y^2}$.

 Figs. \ref{fig:RMHDBLrho} and \ref{fig:RMHDBLbx}
 give the densities $\rho$ and magnetic fields $B_x$
at $t=0.3$ obtained by the proposed DG methods with $300\times 300$ cells,
while FigS.~\ref{fig:BLcmprho} and
\ref{fig:BLcmpbx} show them along the line $x=0$
and the reference solutions, which are obtained by using the MUSCL scheme with
$800\times 800$ uniform cells.
The CFL numbers of $P^1$-, $P^2$-, and $P^3$-based non-central DG methods haven been chosen as $0.2,~0.15,~0.1$, respectively, while those of corresponding central DG methods are $0.3,~0.25,~0.2$, and $\theta=\Delta t_n/\tau_n=1$.
As can be seen from those plots,
the higher order methods can better resolve the discontinuities,
and the numerical solutions   are in good agreement with the reference solutions.
 Table \ref{tab:cellperRMHDSc} records
  the percentage of ``troubled'' cells, and shows that the
  ``troubled'' cells of relatively limited are identified.
      \end{Example}

   \begin{table}[!htbp]
  \centering
  \caption{The percentage of ``troubled'' cells identified by non-central and central \DG{}. }
\begin{tabular}{|c|c|c|c|c|c|c|}
 \hline
    \multirow{2}{10pt}{}& \multicolumn{2}{|c|}{Example ~\ref{exRMHDBlast}} & \multicolumn{2}{|c|}{Example ~\ref{exRMHDRotor}} & \multicolumn{2}{|c|}{Example ~\ref{exRMHDSC}}
  \\
  \cline{2-7}
 & non-central & central    & non-central & central & non-central & central\\
  \hline
  $P^1$& 3.52 & 5.42 & 0.61  &  8.18  &    8.81  &  8.58\\
  \hline
  $P^2$& 9.15 & 10.2 & 12.5 &   14.5 &   20.5 & 17.1 \\
  \hline
 $P^3$& 11.8 & 13.0 & 16.7  & 18.8 & 26.3 &  27.4\\
 \hline
  \end{tabular}
  \label{tab:cellperRMHDSc}
  \end{table}

    \begin{table}[!htbp]
  \centering
  \caption{CPU times of non-central and central \DG{} (second). }
\begin{tabular}{|c|c|c|c|c|c|c|}
   \hline
    \multirow{2}{10pt}{}  & \multicolumn{2}{|c|}{Example ~\ref{exRMHDBlast}}
    & \multicolumn{2}{|c|}{Example ~\ref{exRMHDRotor}} & \multicolumn{2}{|c|}{Example ~\ref{exRMHDSC}}
  \\
  \cline{2-7}
 & non-central & central         & non-central & central  & non-central & central\\
  \hline
  $P^1$& 1.12e5 & 2.11e5  & 2.81e5 & 6.25e5 & 1.76e6 & 3.44e6 \\
  \hline
  $P^2$& 5.52e5 & 8.00e5 & 6.55e5 & 1.84e6 & 3.56e6  & 1.10e7 \\
  \hline
 $P^3$& 1.16e6 & 2.59e6 & 1.73e6 & 5.98e6 & 1.10e7 &   2.25e7\\
 \hline
  \end{tabular}
  \label{tab:RMHDSccmpRC}
  \end{table}

\begin{figure}[!htbp]
    \centering{}
  \begin{tabular}{cc}
    \includegraphics[width=0.35\textwidth]{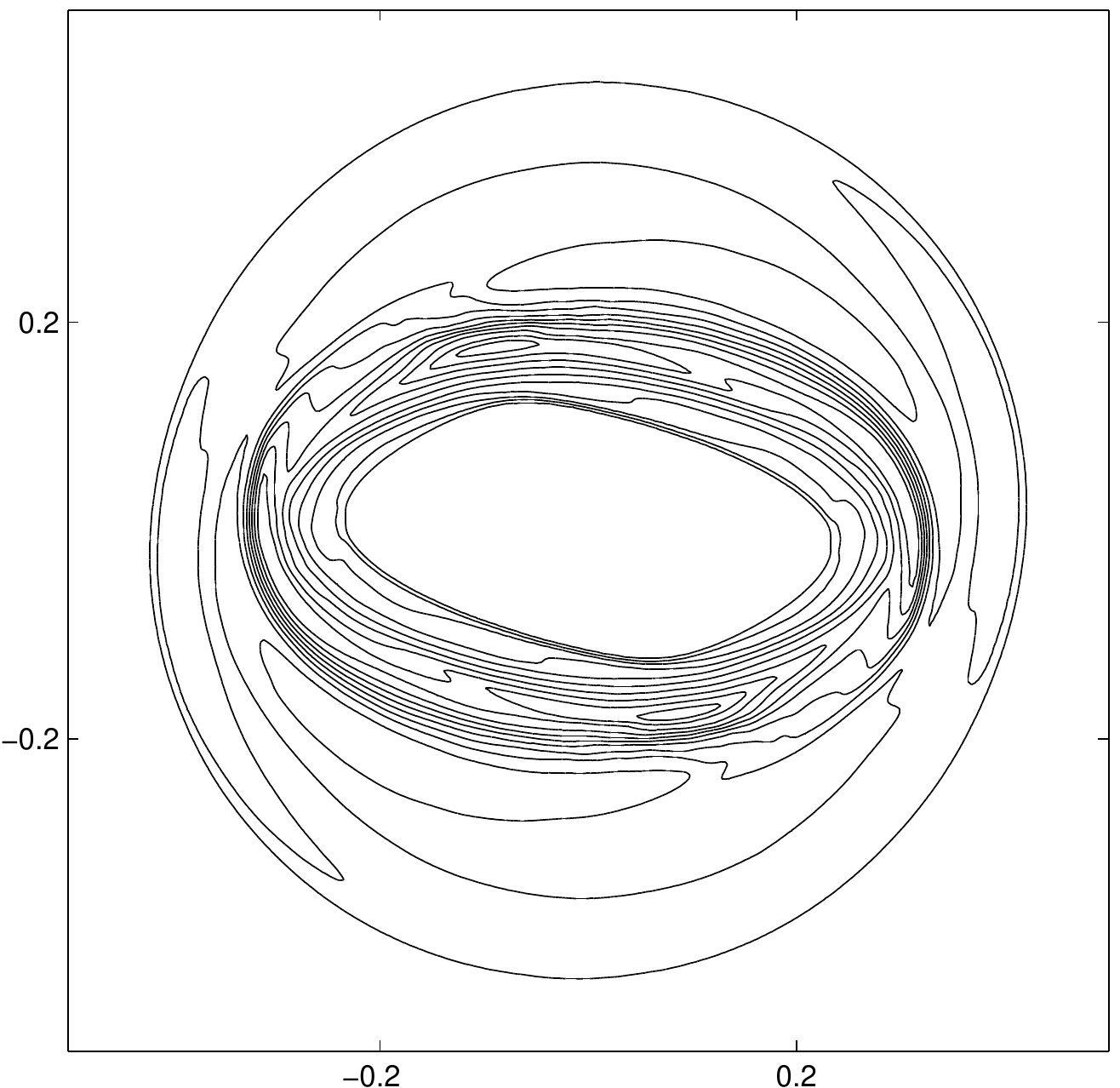}&
\includegraphics[width=0.35\textwidth]{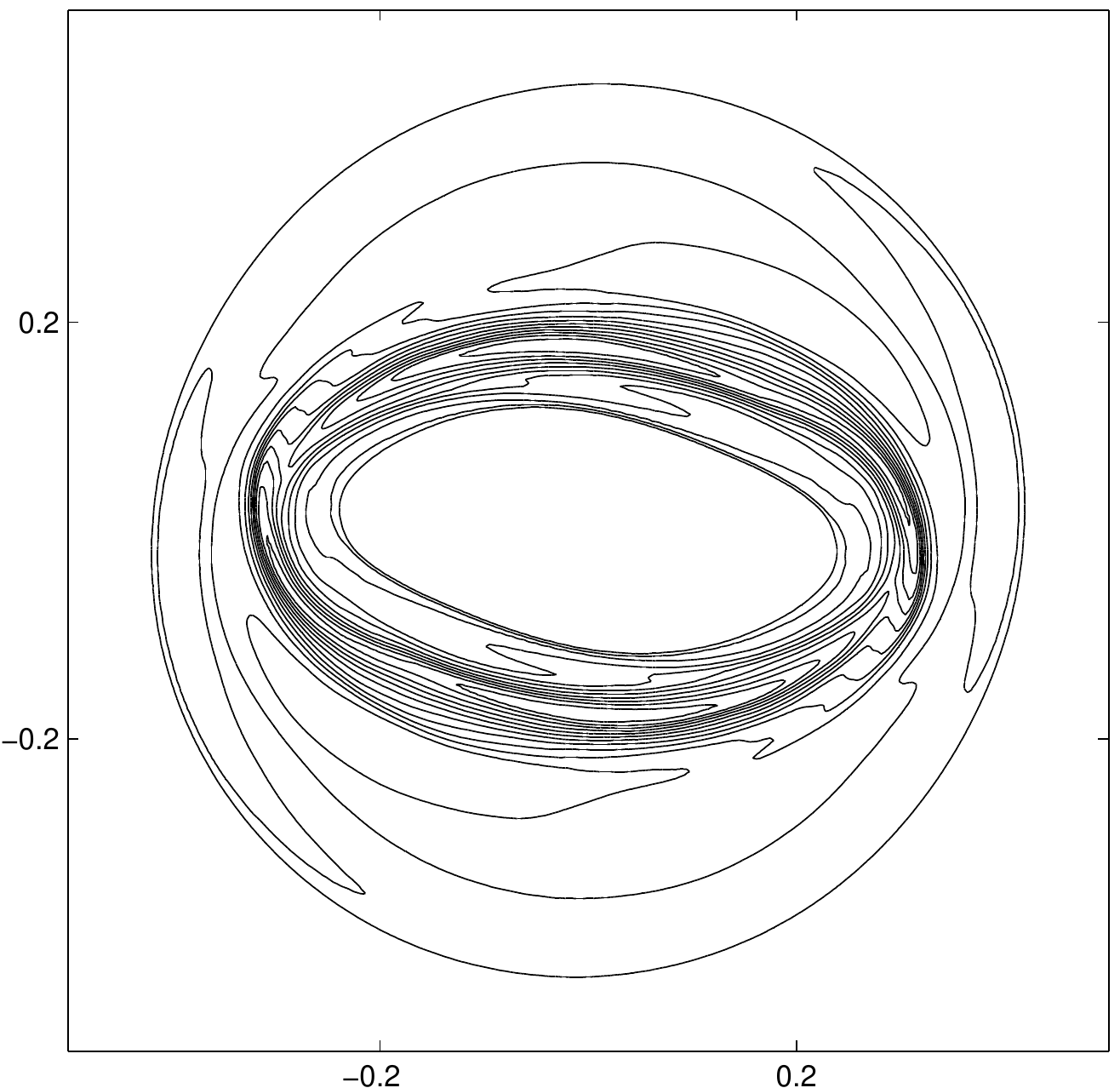}\\
    \includegraphics[width=0.35\textwidth]{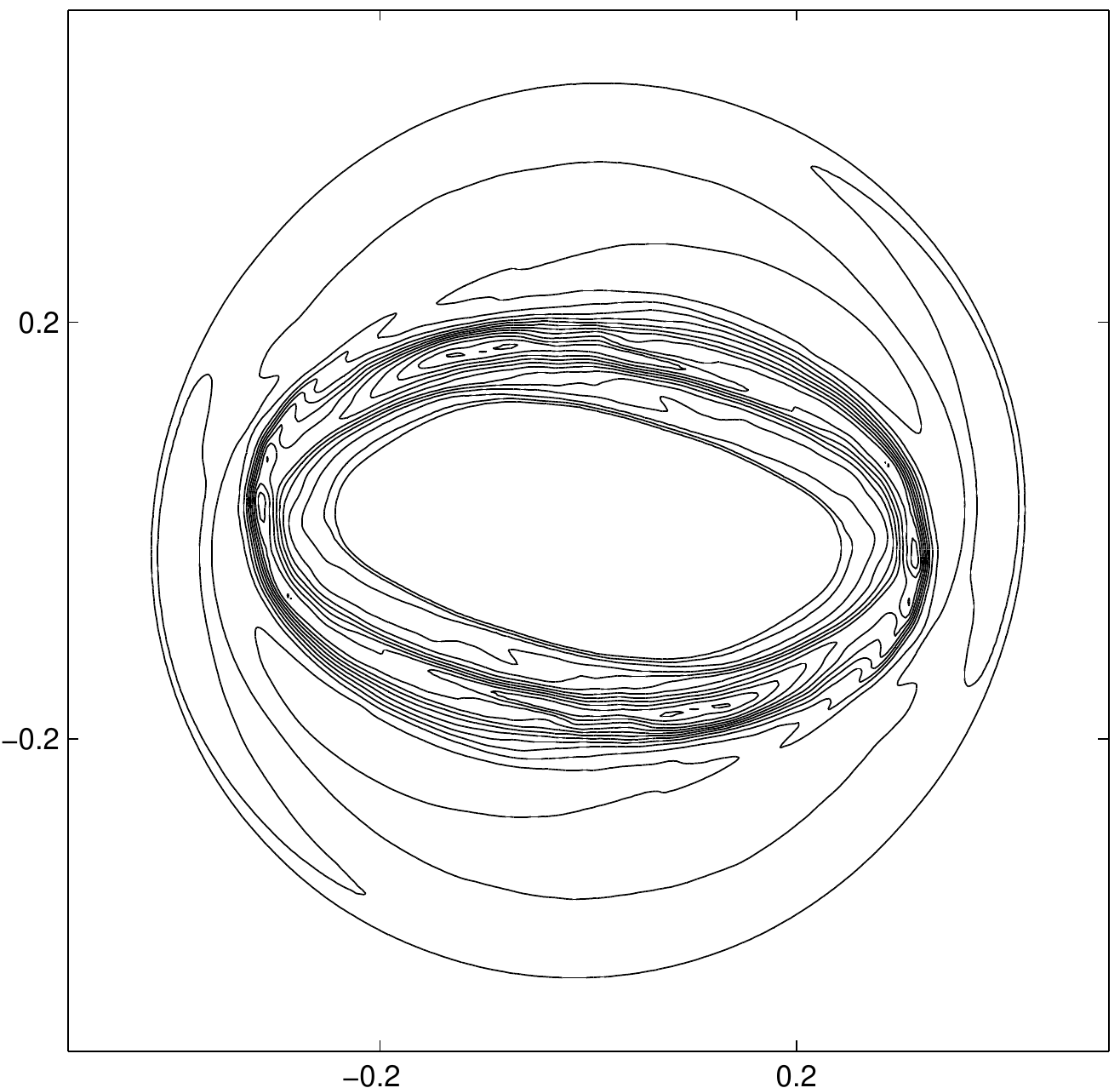}&
\includegraphics[width=0.35\textwidth]{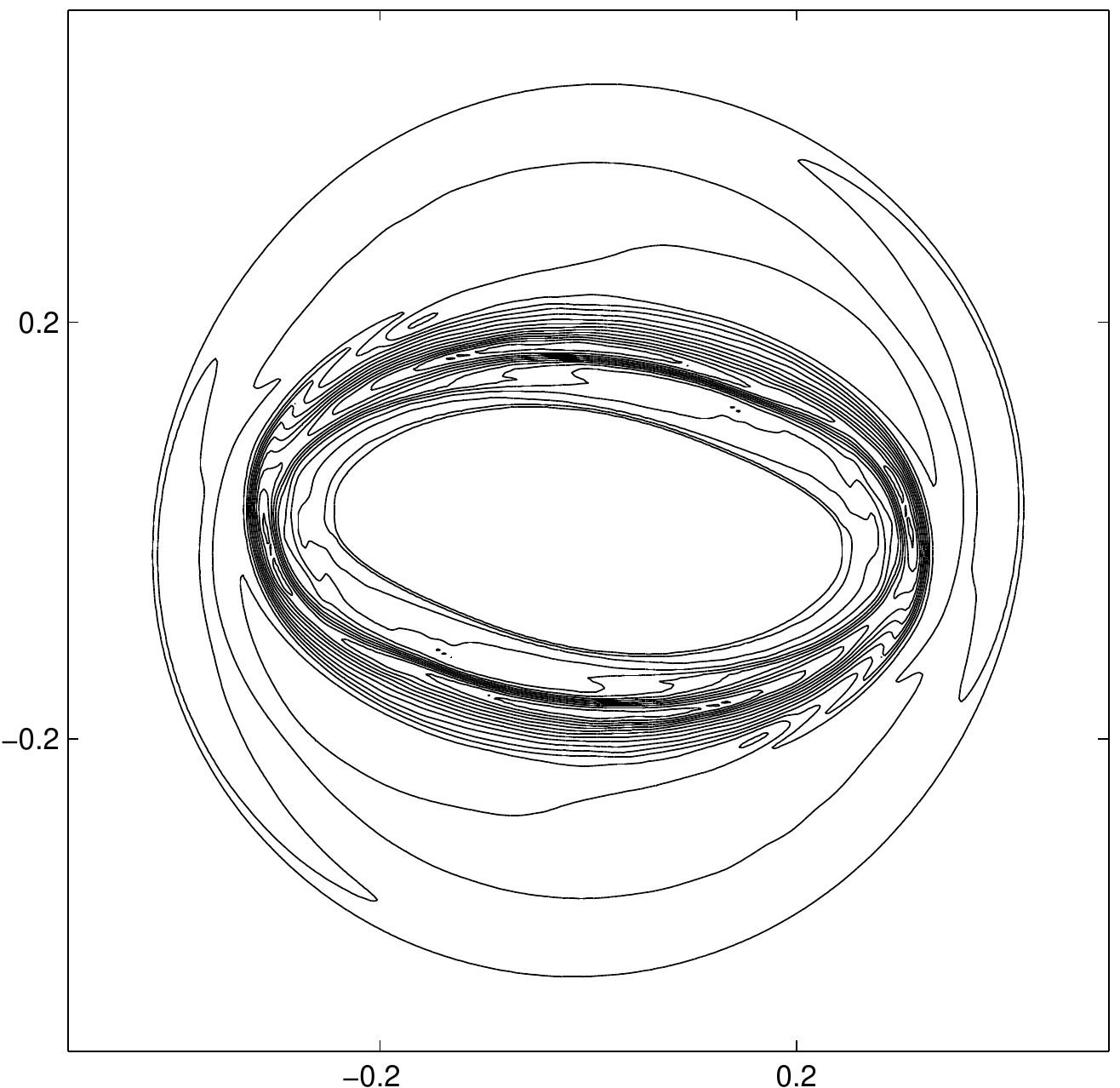}\\
 \includegraphics[width=0.35\textwidth]{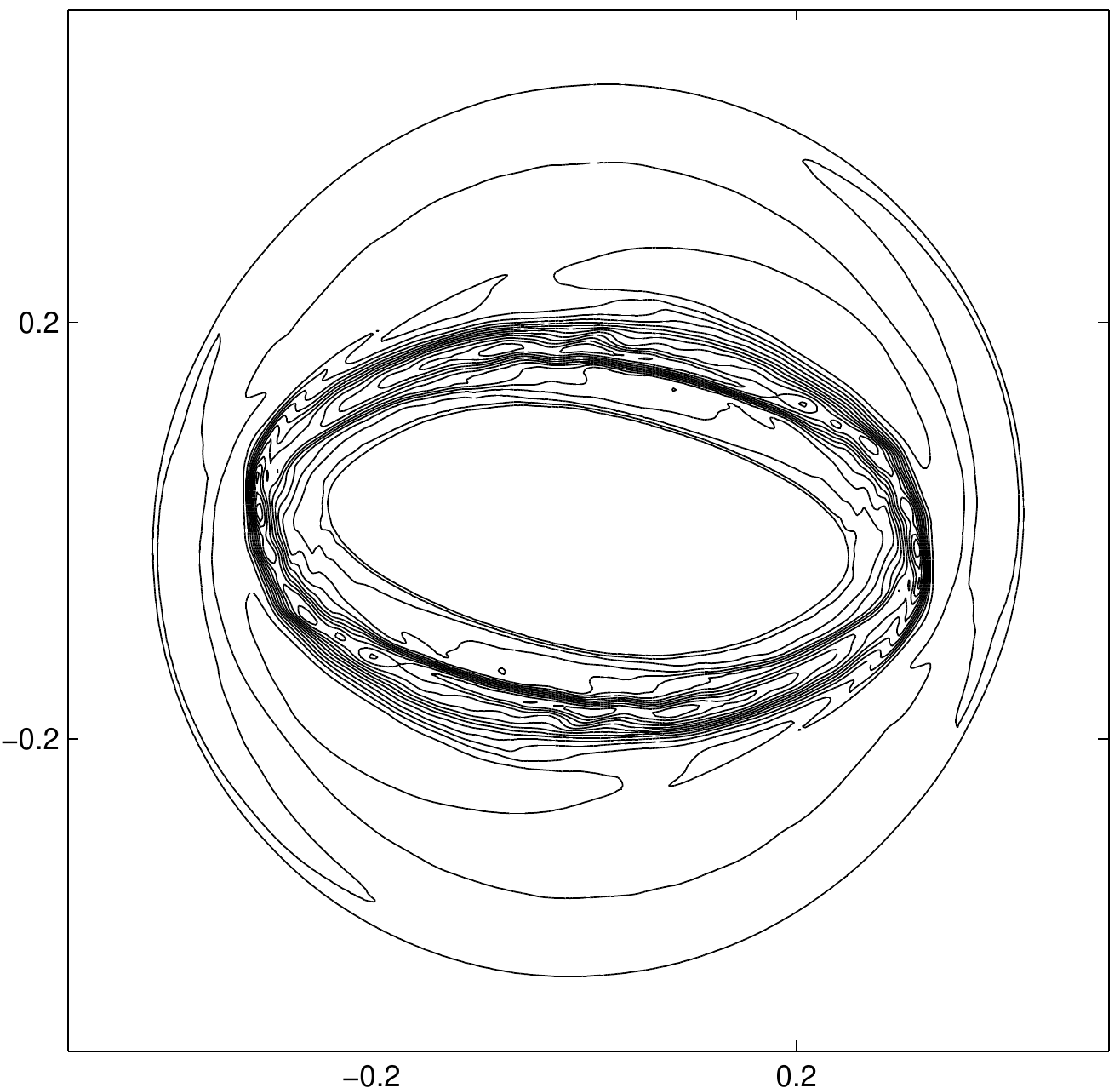}&
\includegraphics[width=0.35\textwidth]{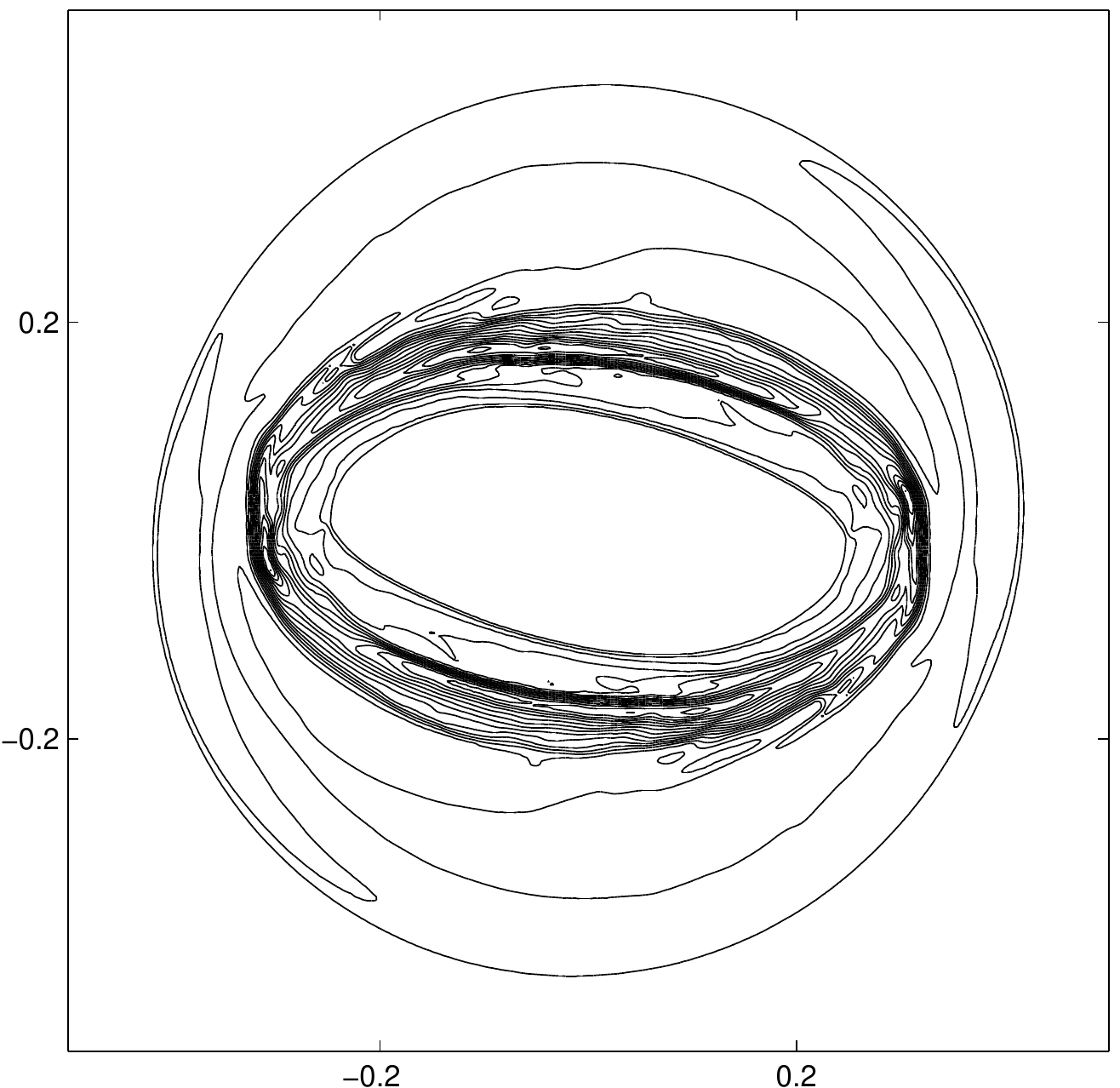}
    \end{tabular}
    \caption{Example~\ref{exRMHDRotor}:
    The contour plots of density $\rho$ at $t = 0.4$ obtained with $300\times 300$ cells
    (15 equally spaced contour lines from 0.3 to 8.2).
     Left: $P^K$-based \DG{}; right: $P^K$-based \CDG{}.  From top to bottom: $K=1,~2,~3$. }
    \label{fig:RMHDRotorrho}
  \end{figure}

\begin{figure}[!htbp]
    \centering{}
  \begin{tabular}{cc}
    \includegraphics[width=0.35\textwidth]{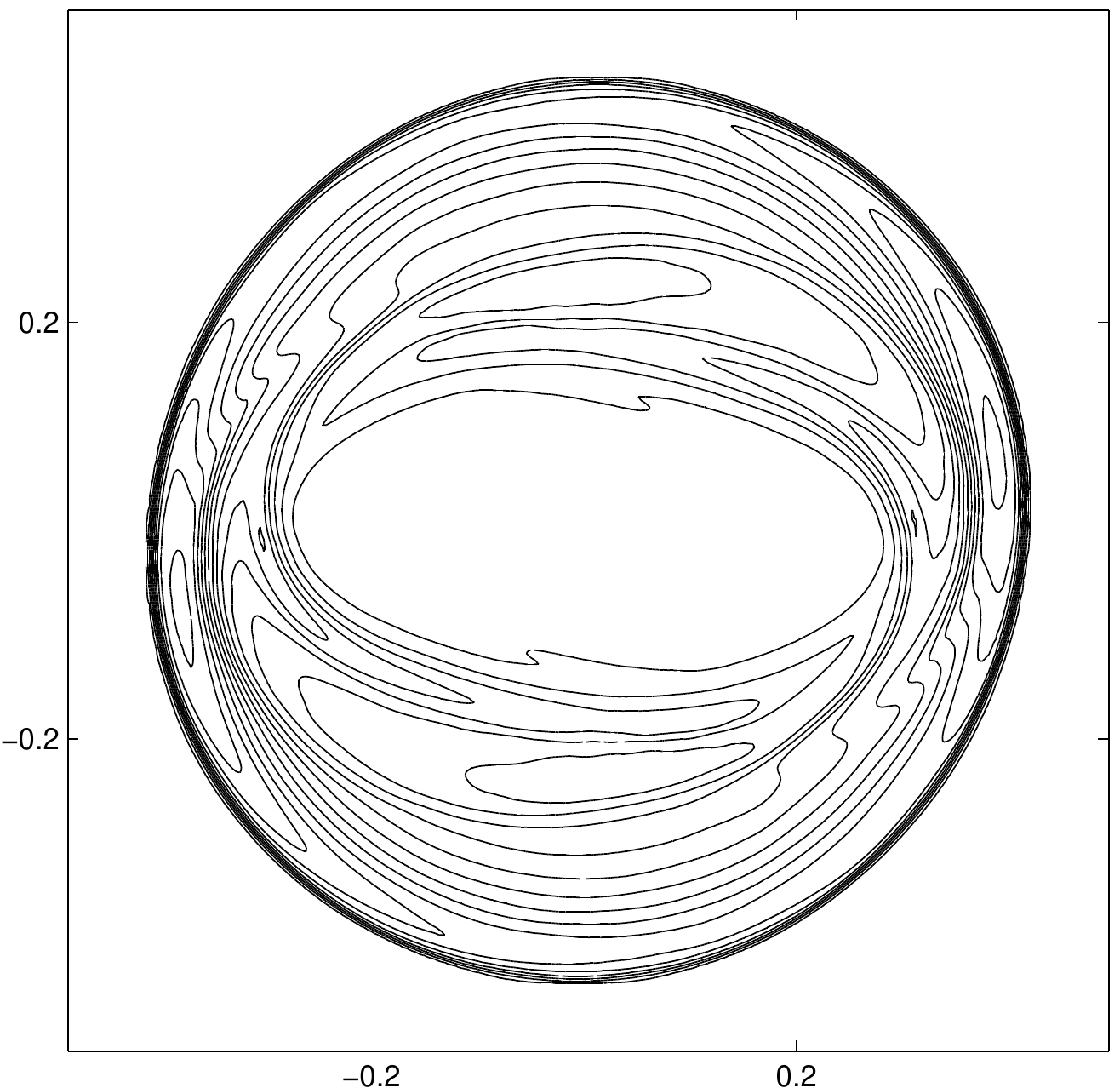}&
\includegraphics[width=0.35\textwidth]{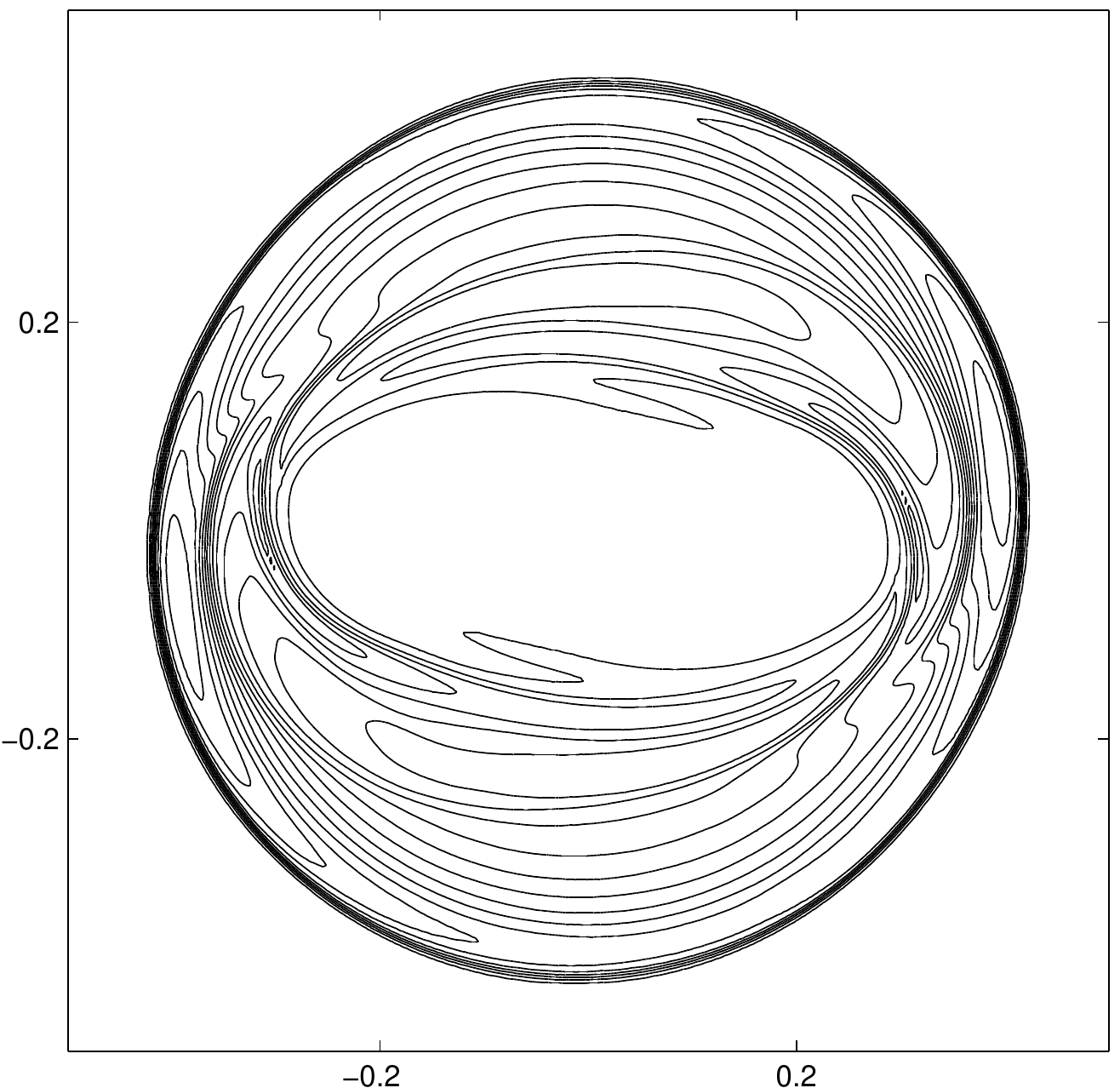}\\
    \includegraphics[width=0.35\textwidth]{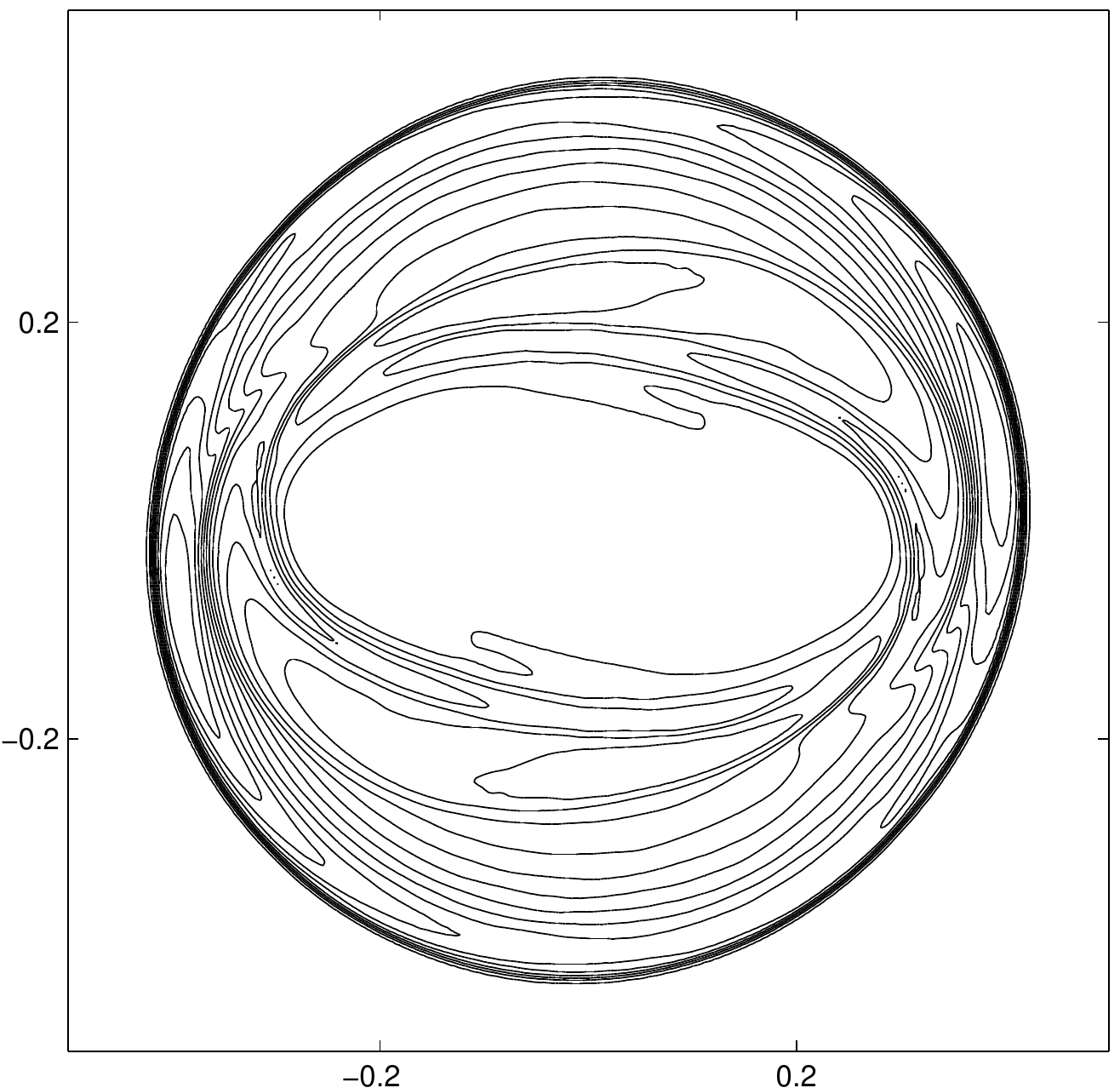}&
\includegraphics[width=0.35\textwidth]{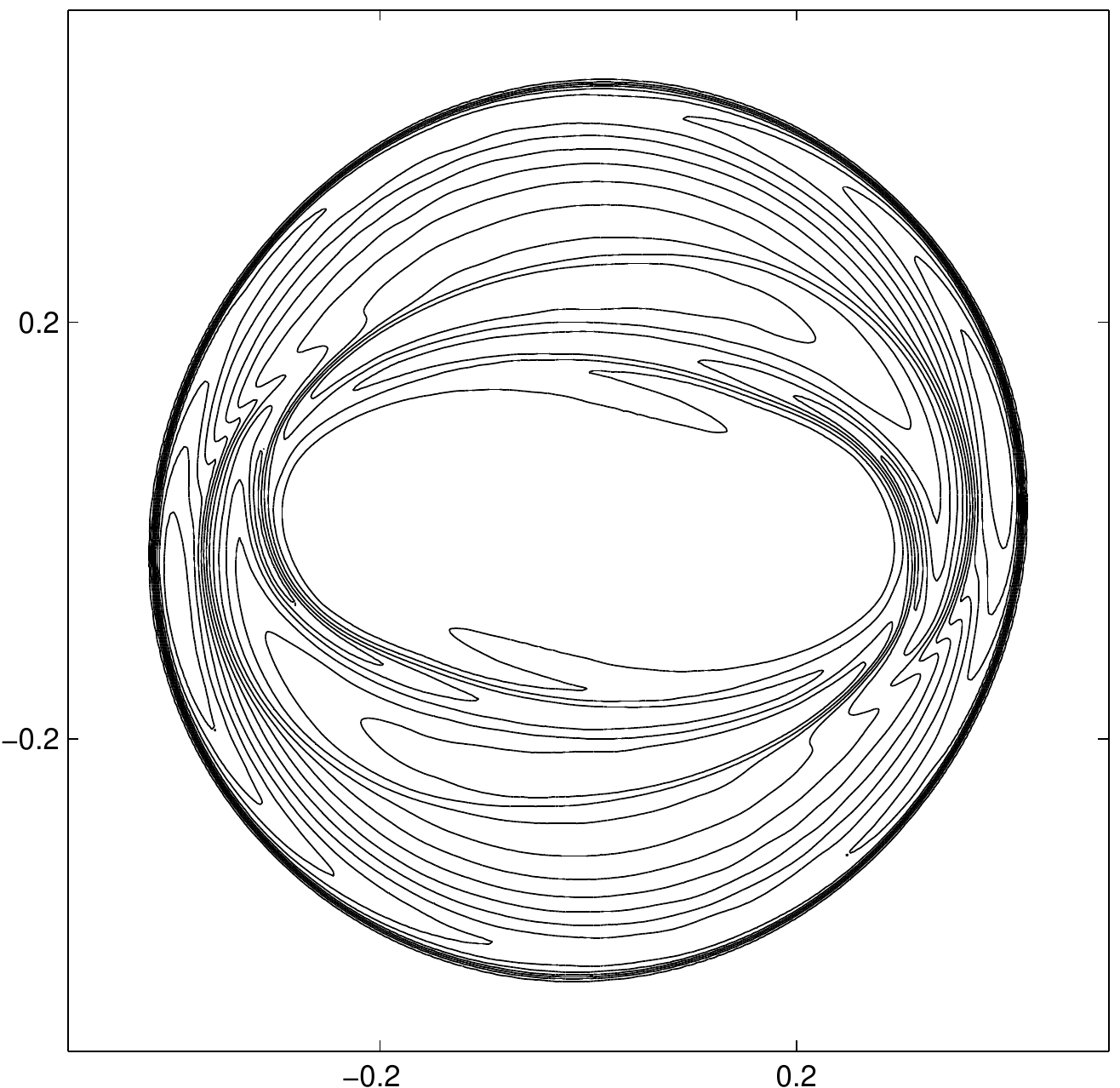}\\
 \includegraphics[width=0.35\textwidth]{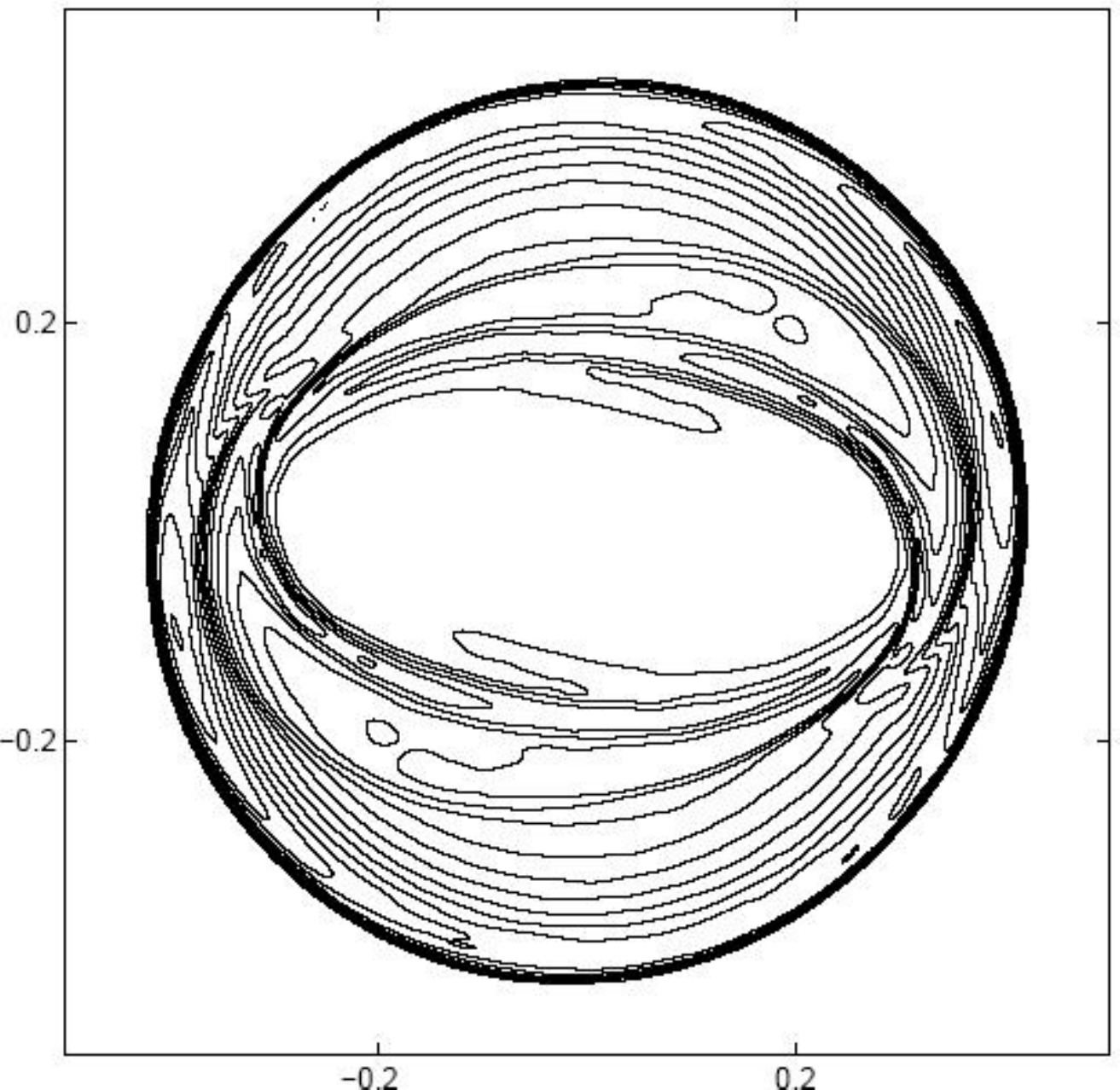}&
\includegraphics[width=0.35\textwidth]{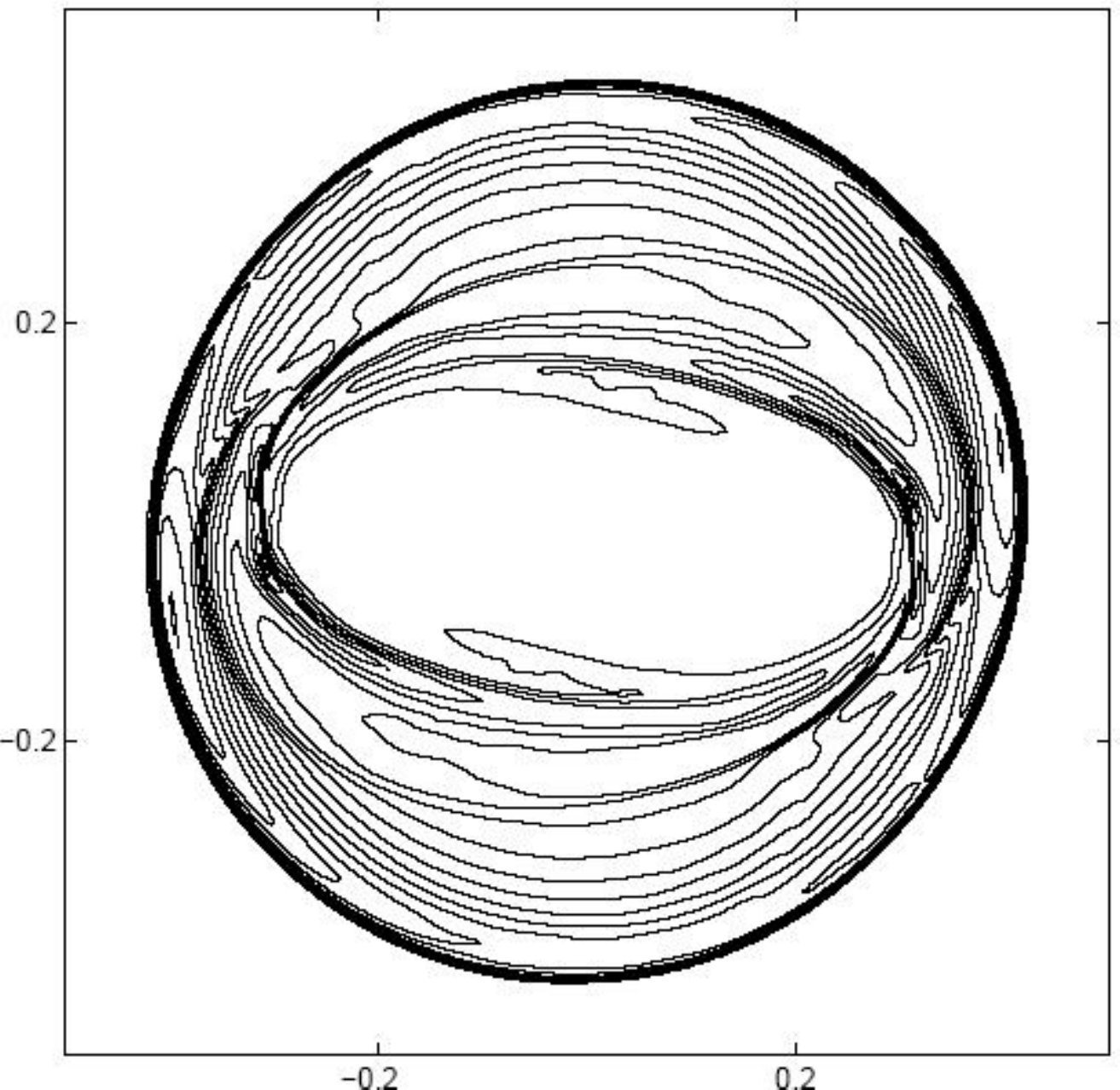}
    \end{tabular}
    \caption{
    Same as Fig.~\ref{fig:RMHDRotorrho} except for
    the gas pressure $p$ (15 equally spaced contour lines from 0 to 4.2).
     }
    \label{fig:RMHDRotorpre}
  \end{figure}

\begin{figure}[!htbp]
    \centering{}
  \begin{tabular}{cc}
    \includegraphics[width=0.35\textwidth]{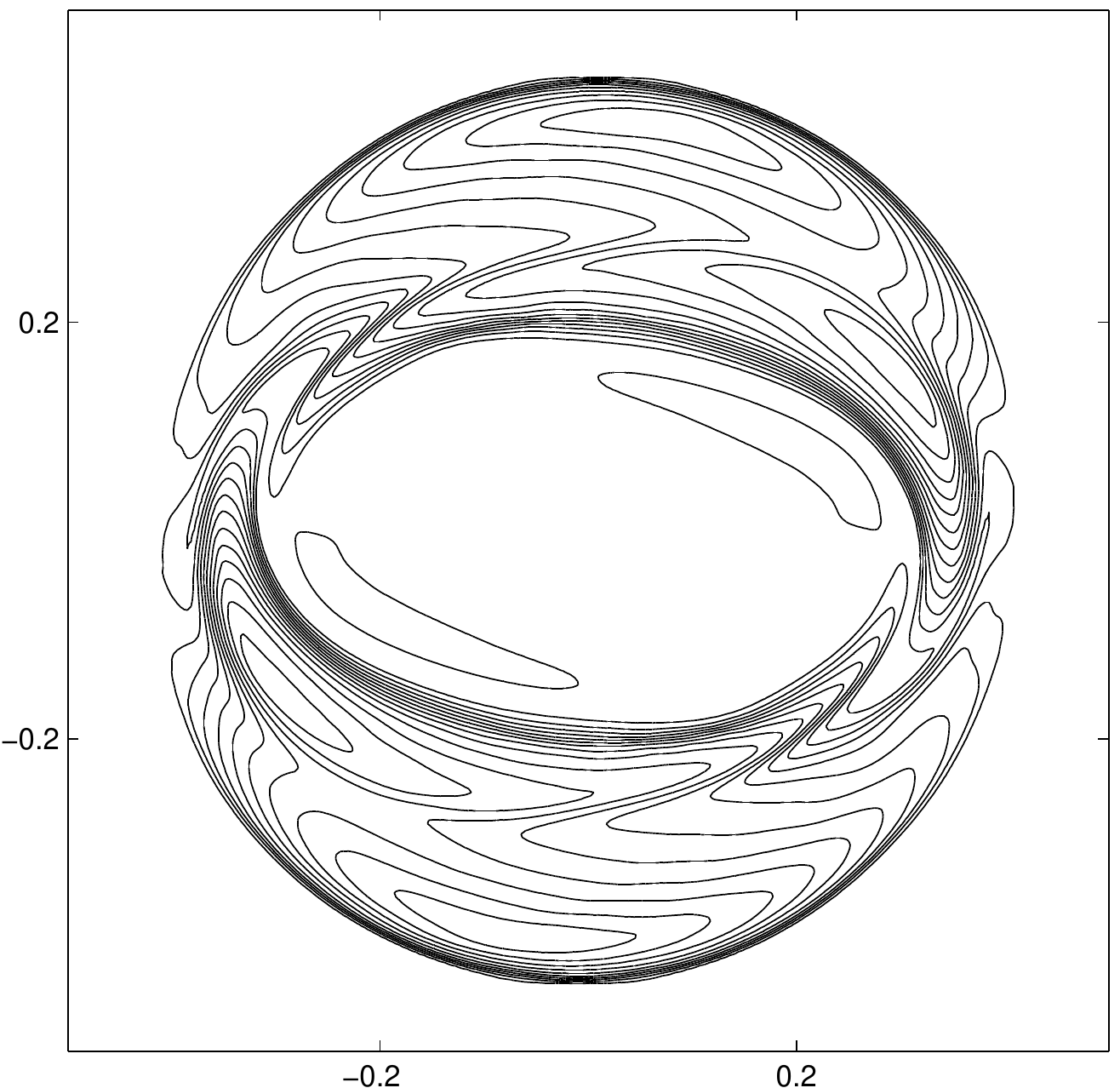}&
\includegraphics[width=0.35\textwidth]{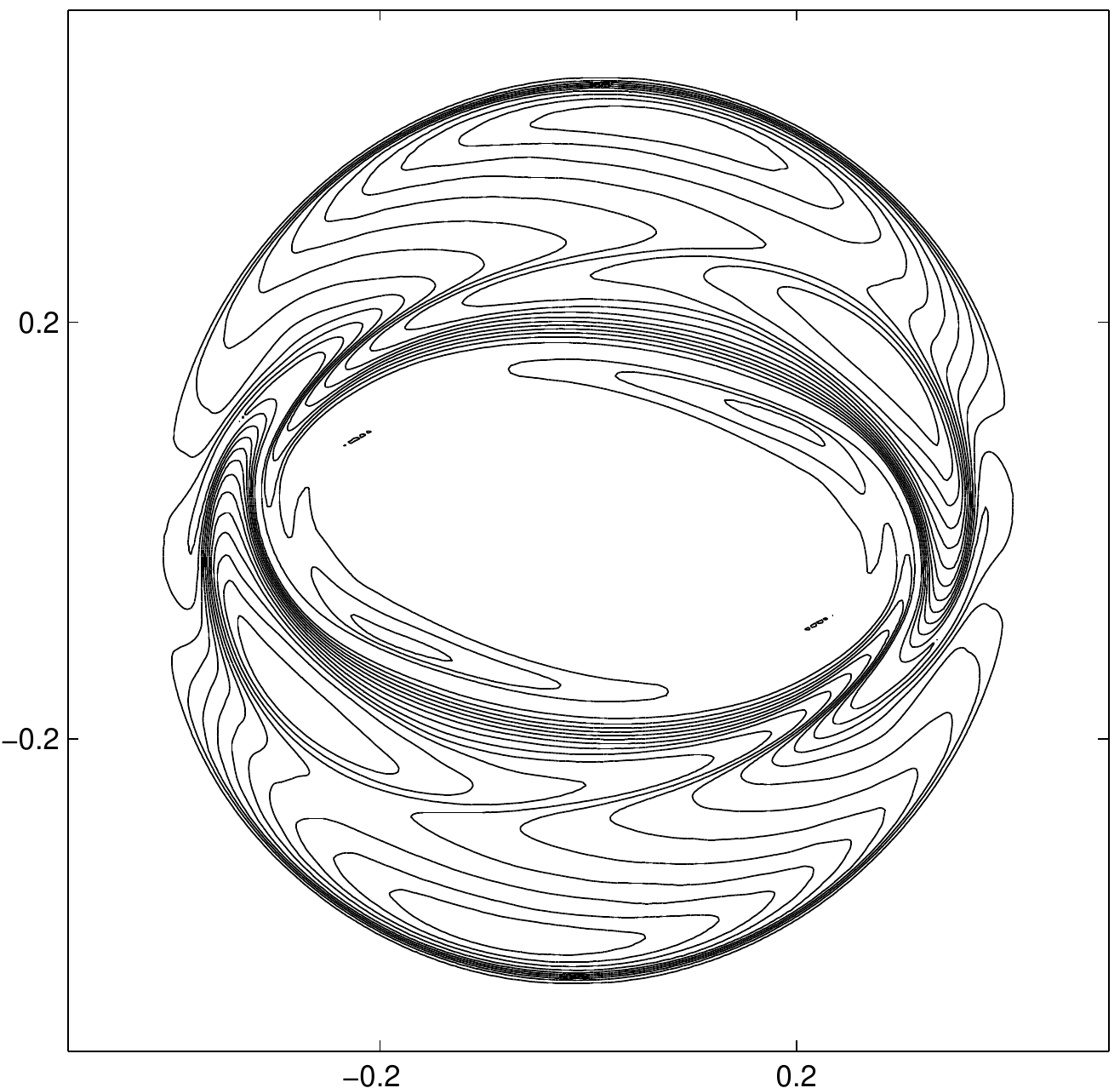}\\
   \includegraphics[width=0.35\textwidth]{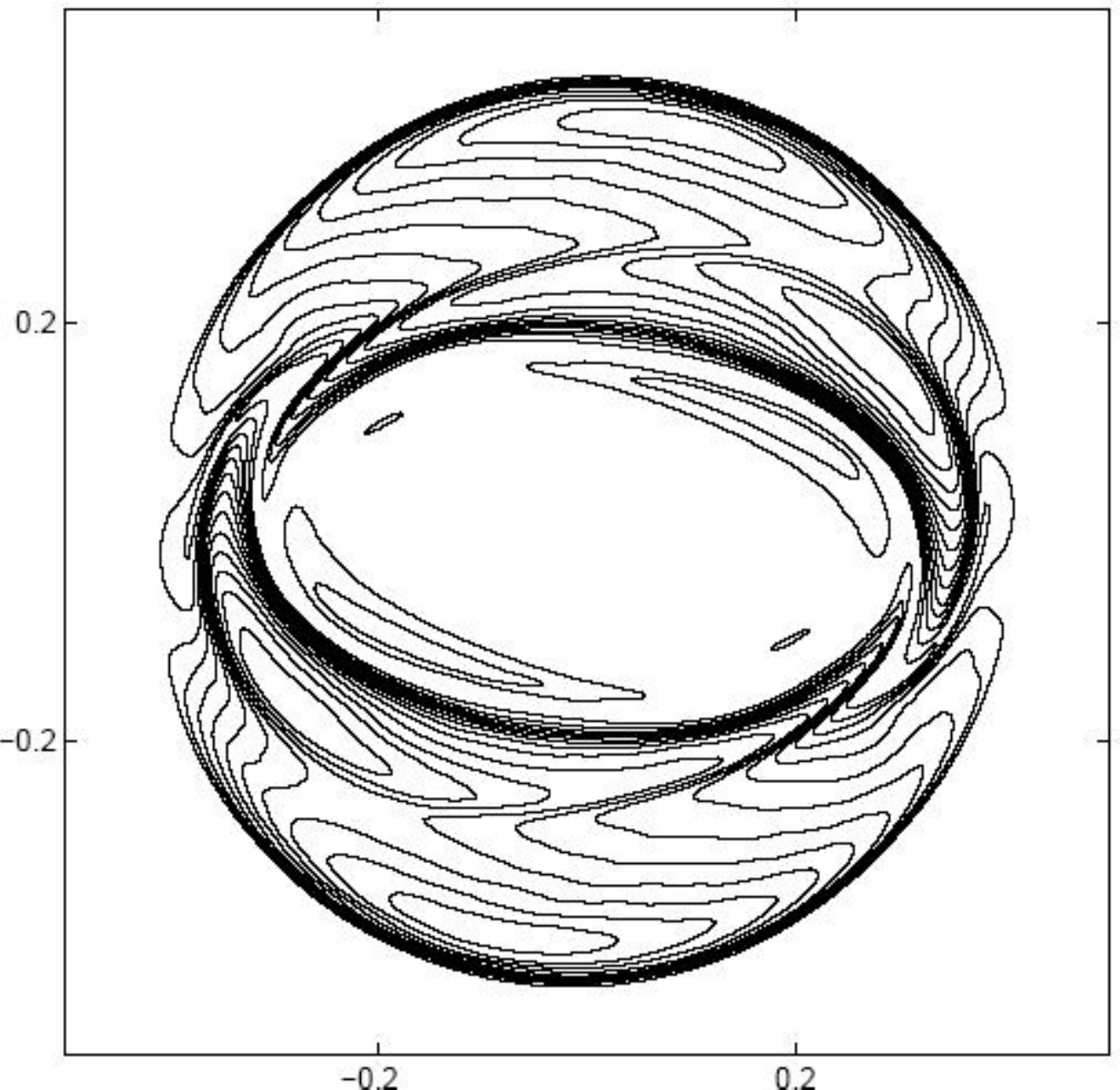}&
\includegraphics[width=0.35\textwidth]{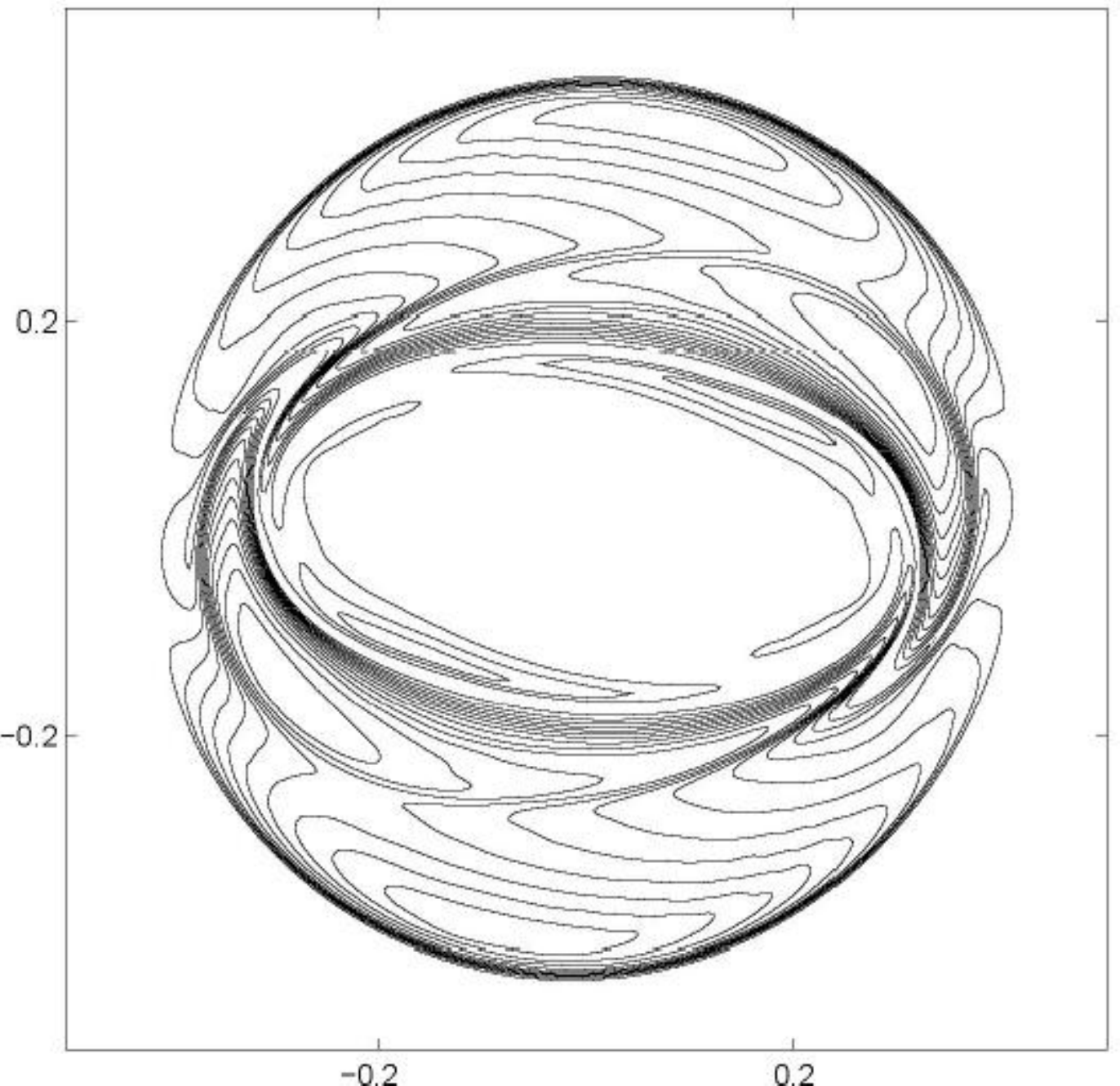}\\
\includegraphics[width=0.35\textwidth]{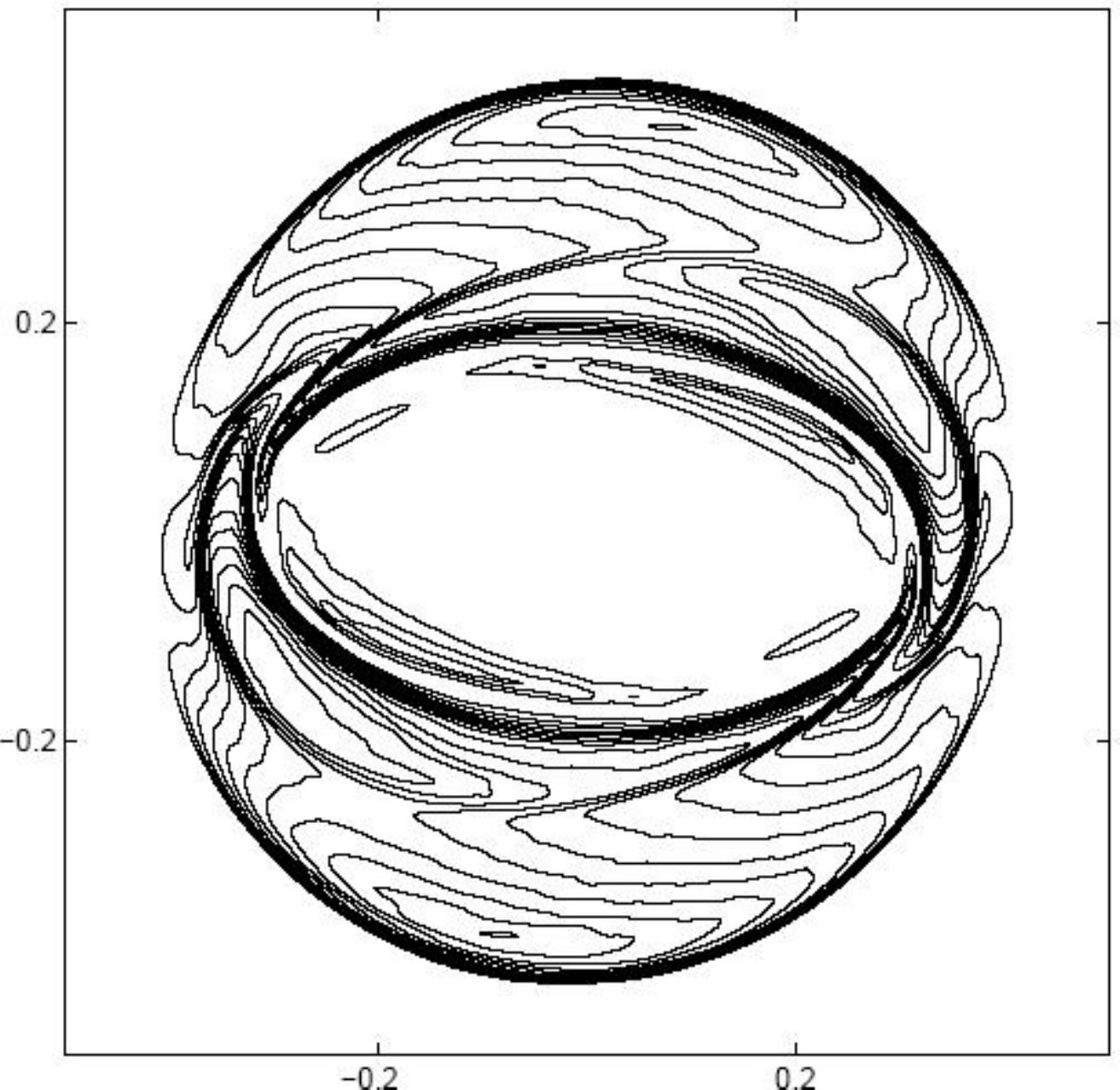}&
\includegraphics[width=0.35\textwidth]{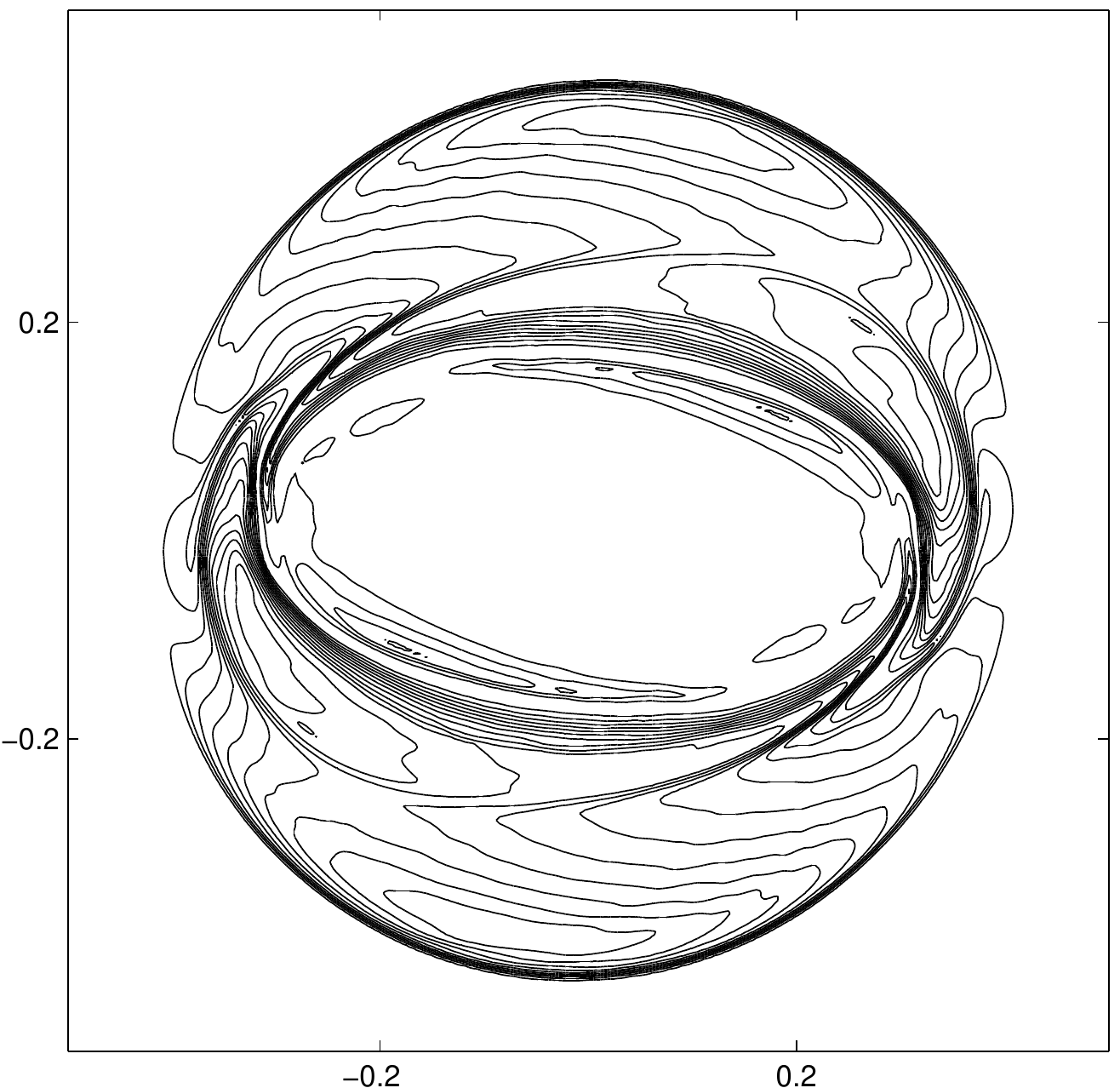}
    \end{tabular}
    \caption{Same as Fig.~\ref{fig:RMHDRotorrho} except for
    the magnetic pressure $p_m$ (15 equally spaced contour lines from 0 to 2.1).
 }
    \label{fig:RMHDRotorpm}
  \end{figure}

\begin{figure}[!htbp]
    \centering{}
  \begin{tabular}{cc}
    \includegraphics[width=0.35\textwidth]{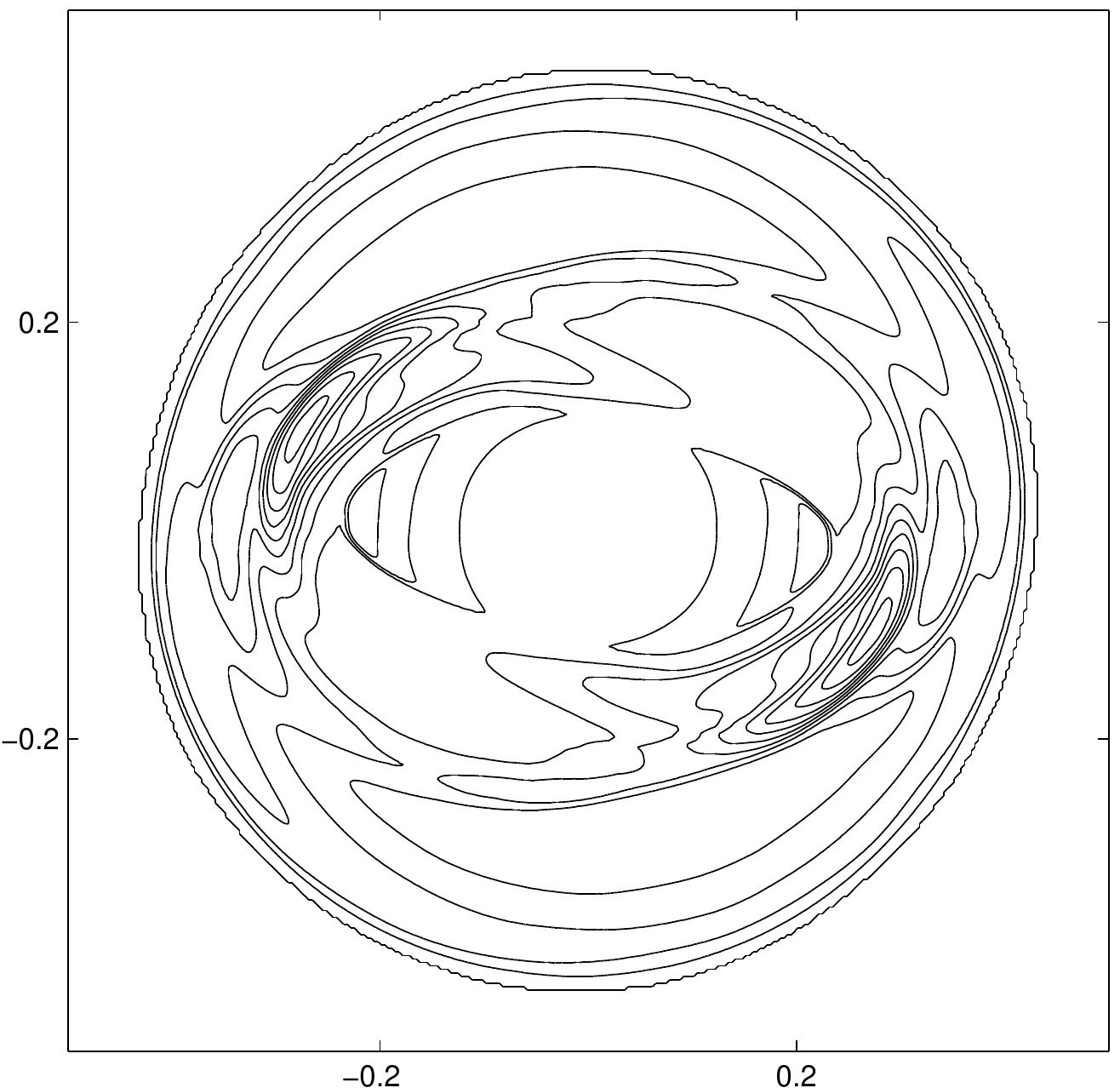}&
\includegraphics[width=0.35\textwidth]{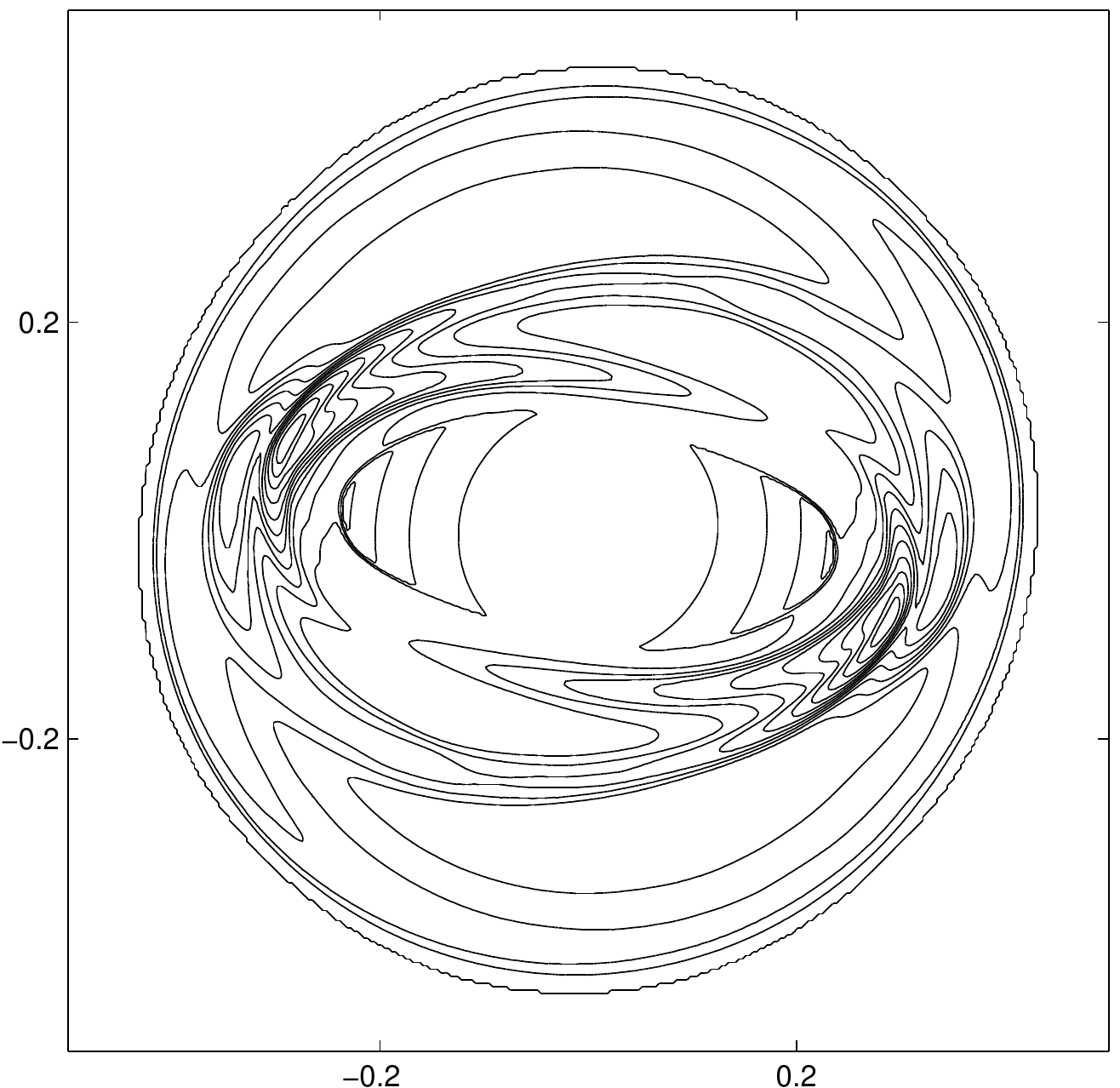}\\
\includegraphics[width=0.35\textwidth]{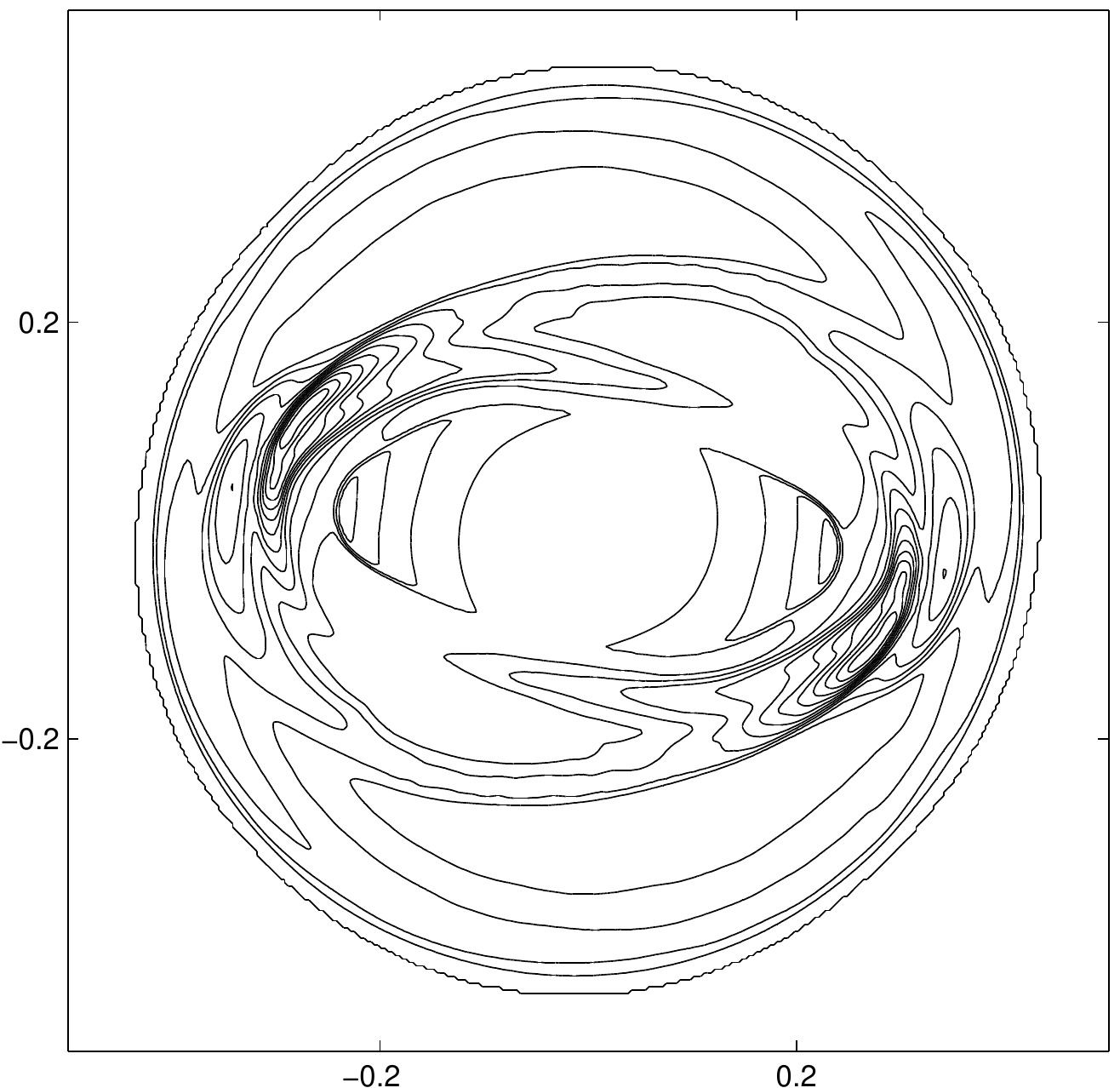}&
\includegraphics[width=0.35\textwidth]{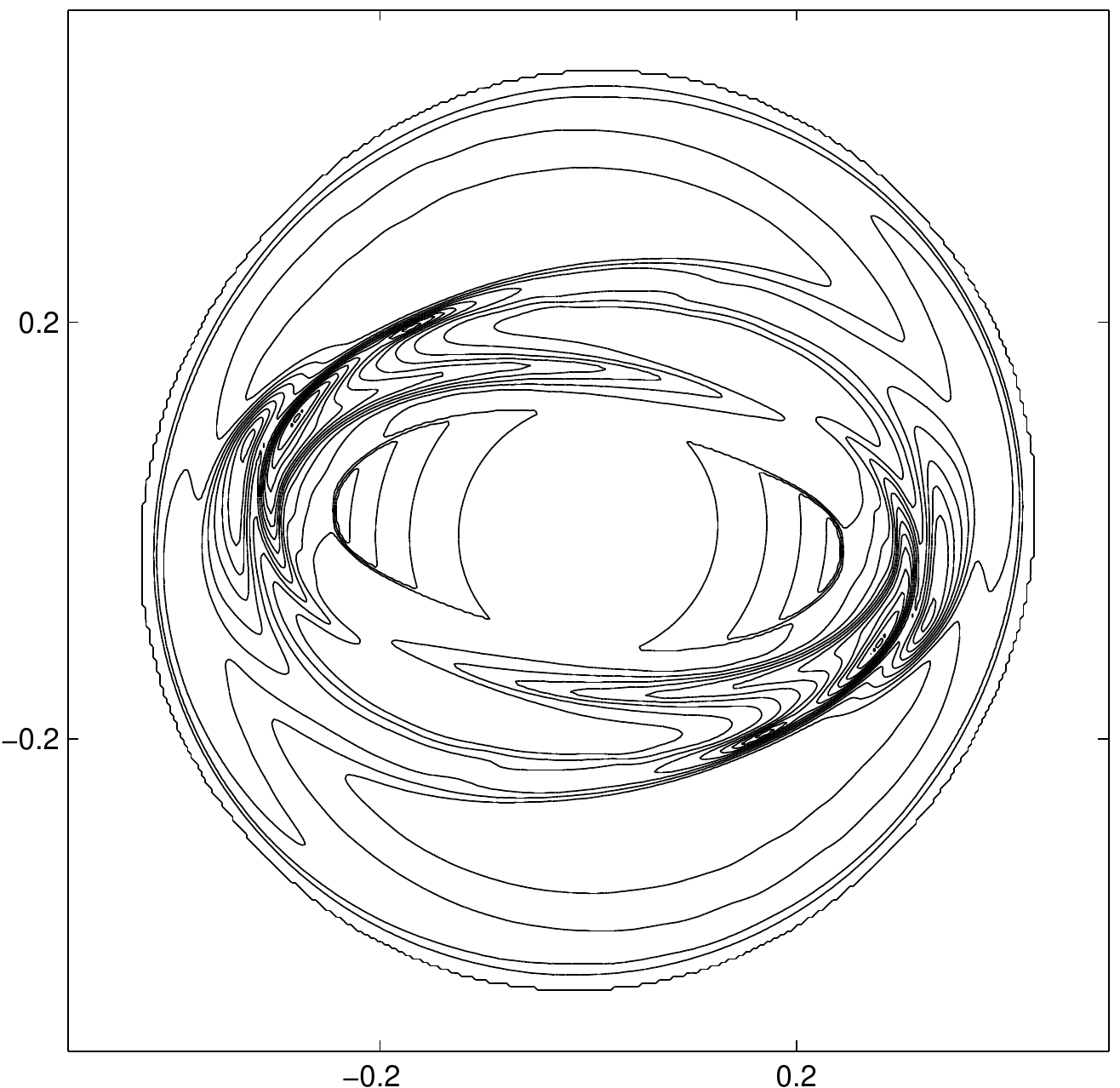}\\
\includegraphics[width=0.35\textwidth]{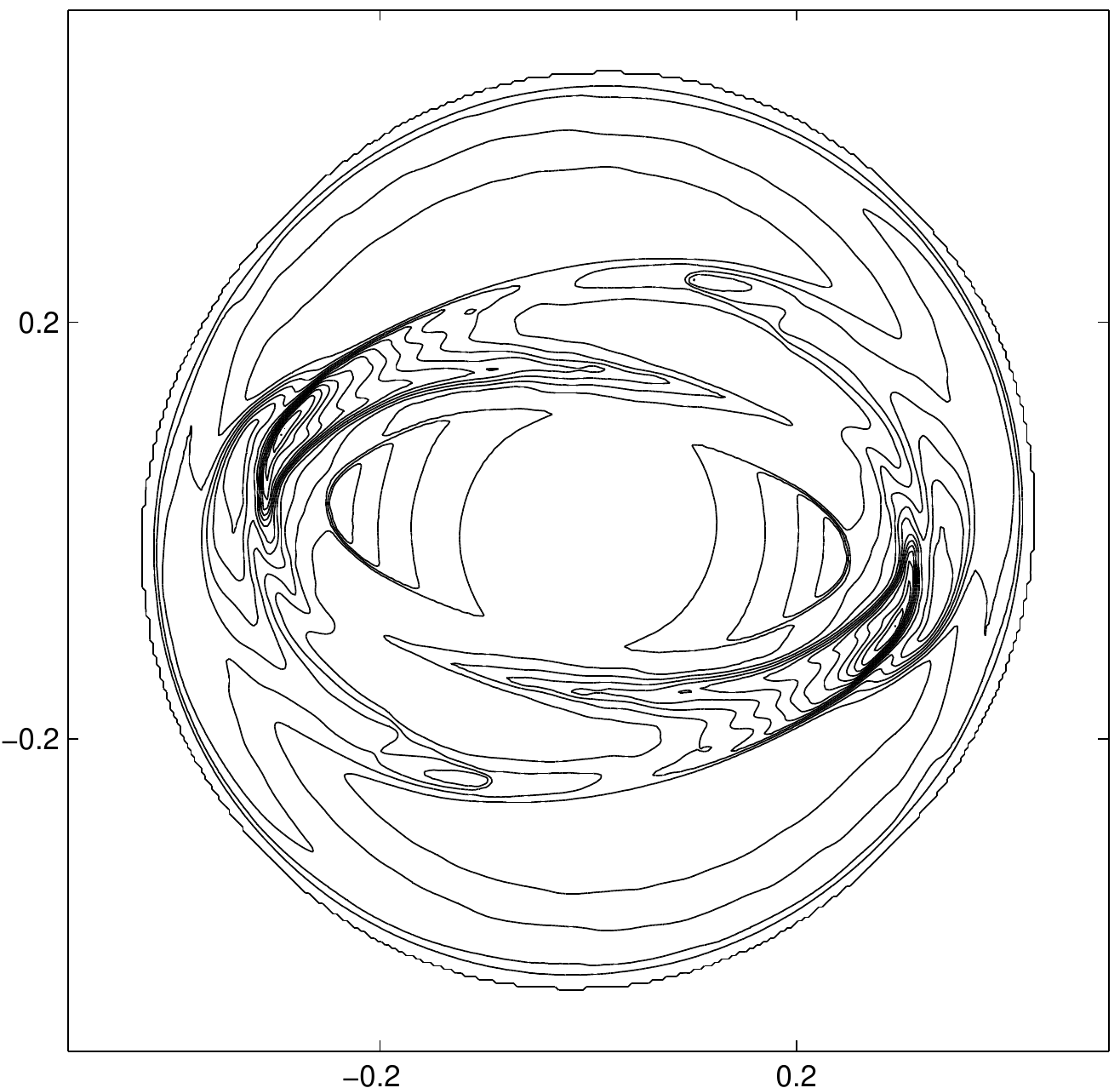}&
\includegraphics[width=0.35\textwidth]{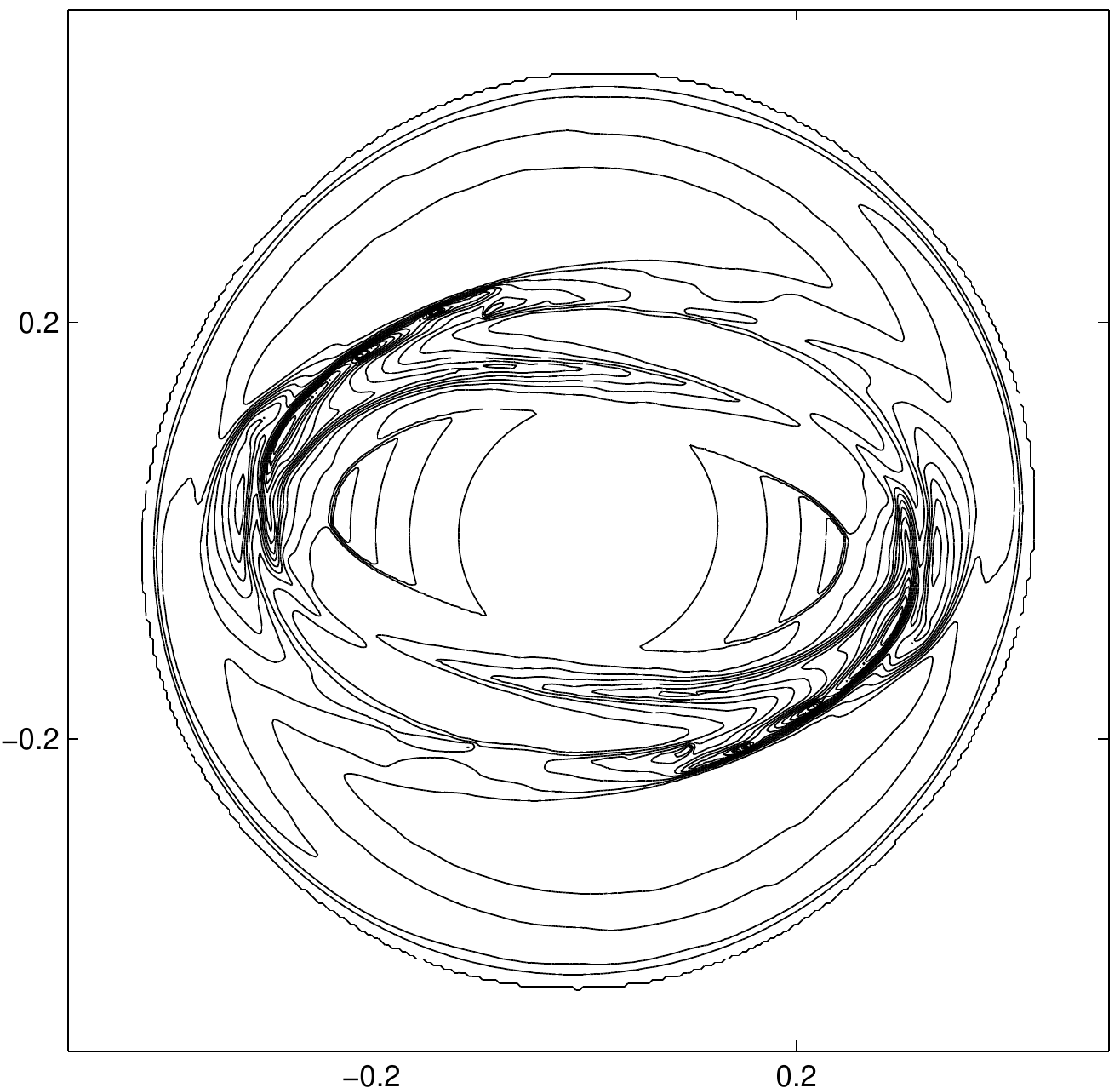}
    \end{tabular}
    \caption{Same as Fig.~\ref{fig:RMHDRotorrho} except for
    the Lorentz factor (15 equally spaced contour lines from 1 to 1.9).
      }
    \label{fig:RMHDRotorgam}
  \end{figure}

\begin{figure}[!htbp]
    \centering{}
    \begin{tabular}{cc}
    \includegraphics[width=0.35\textwidth]{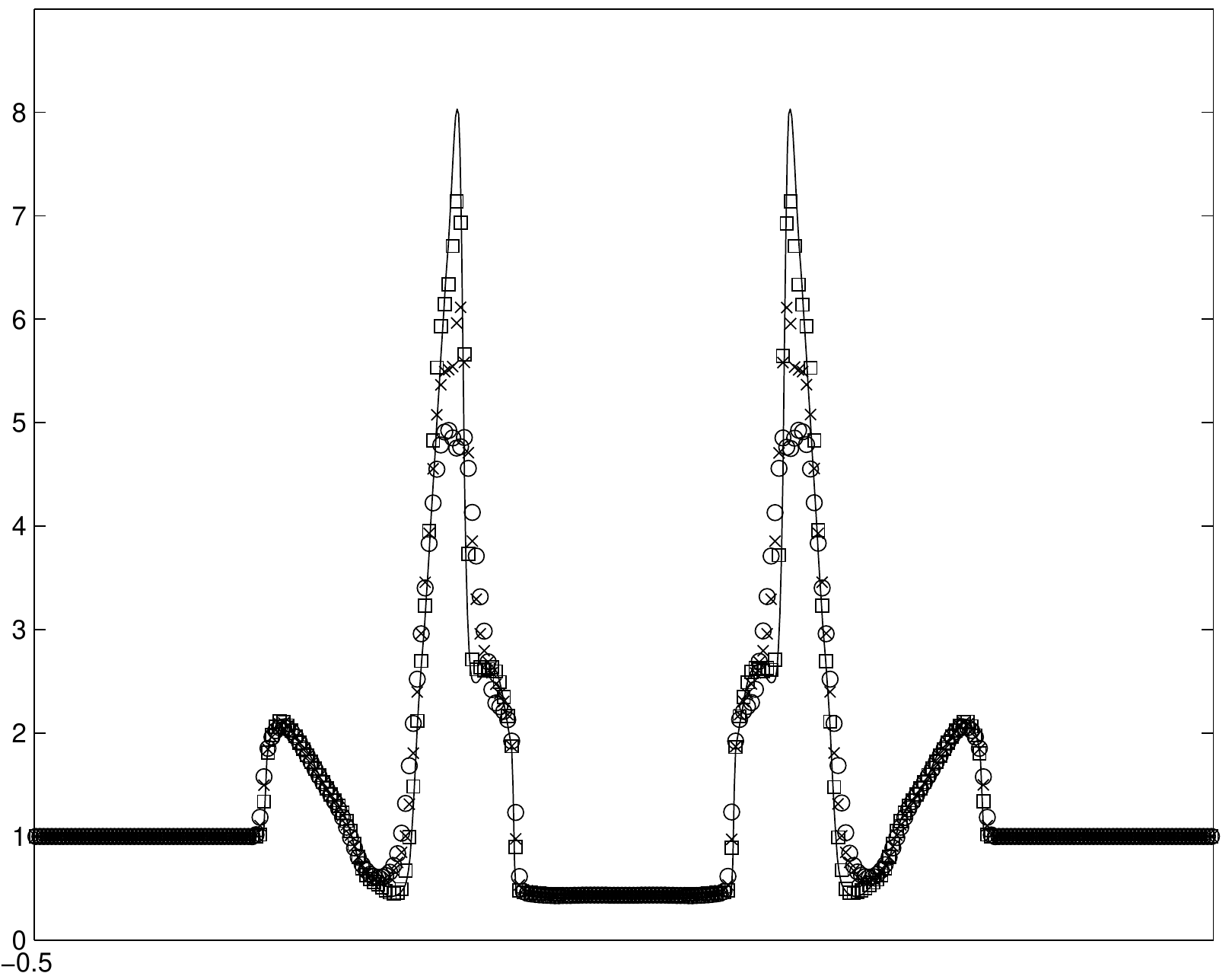}&
        \includegraphics[width=0.35\textwidth]{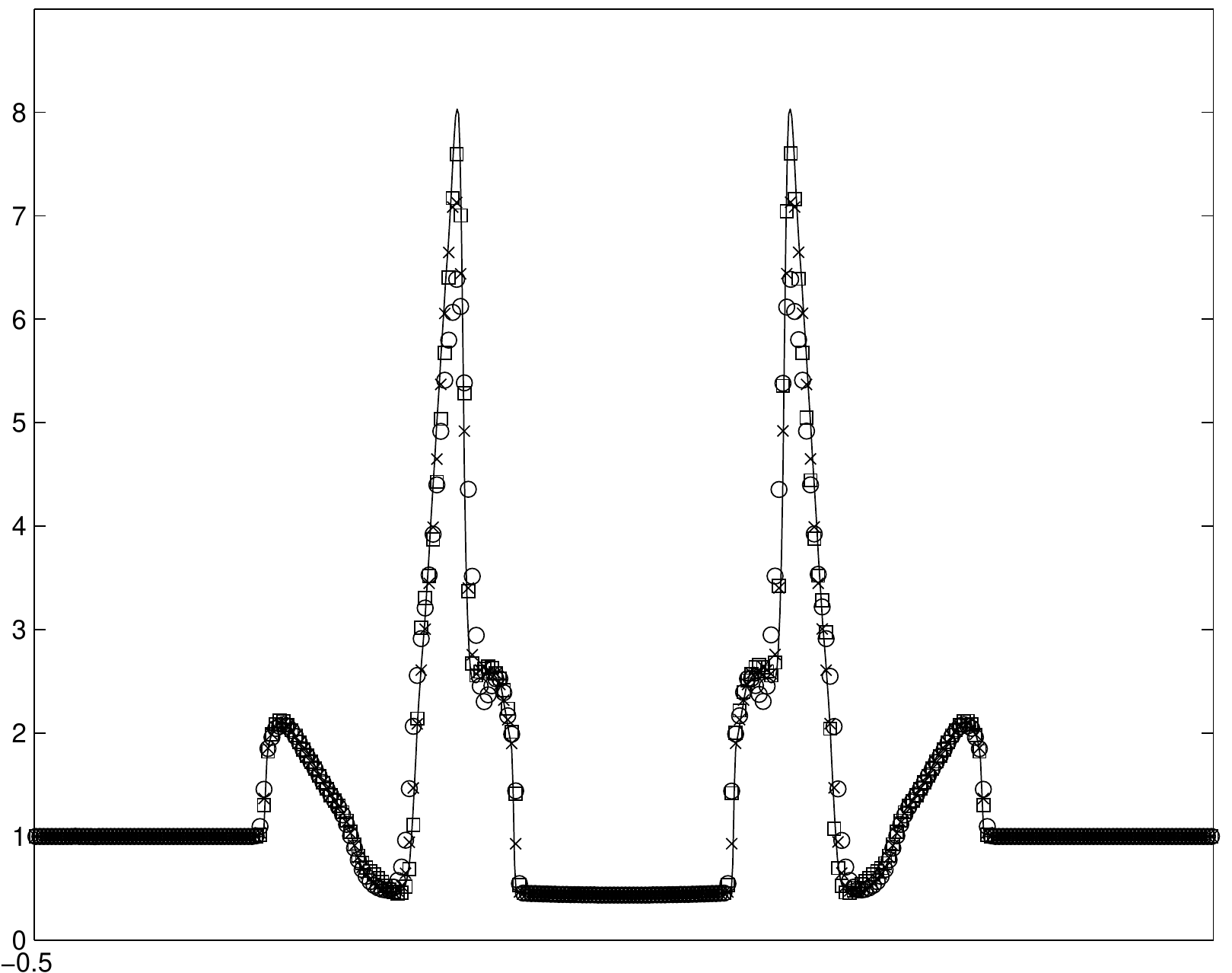}\\
         \includegraphics[width=0.35\textwidth]{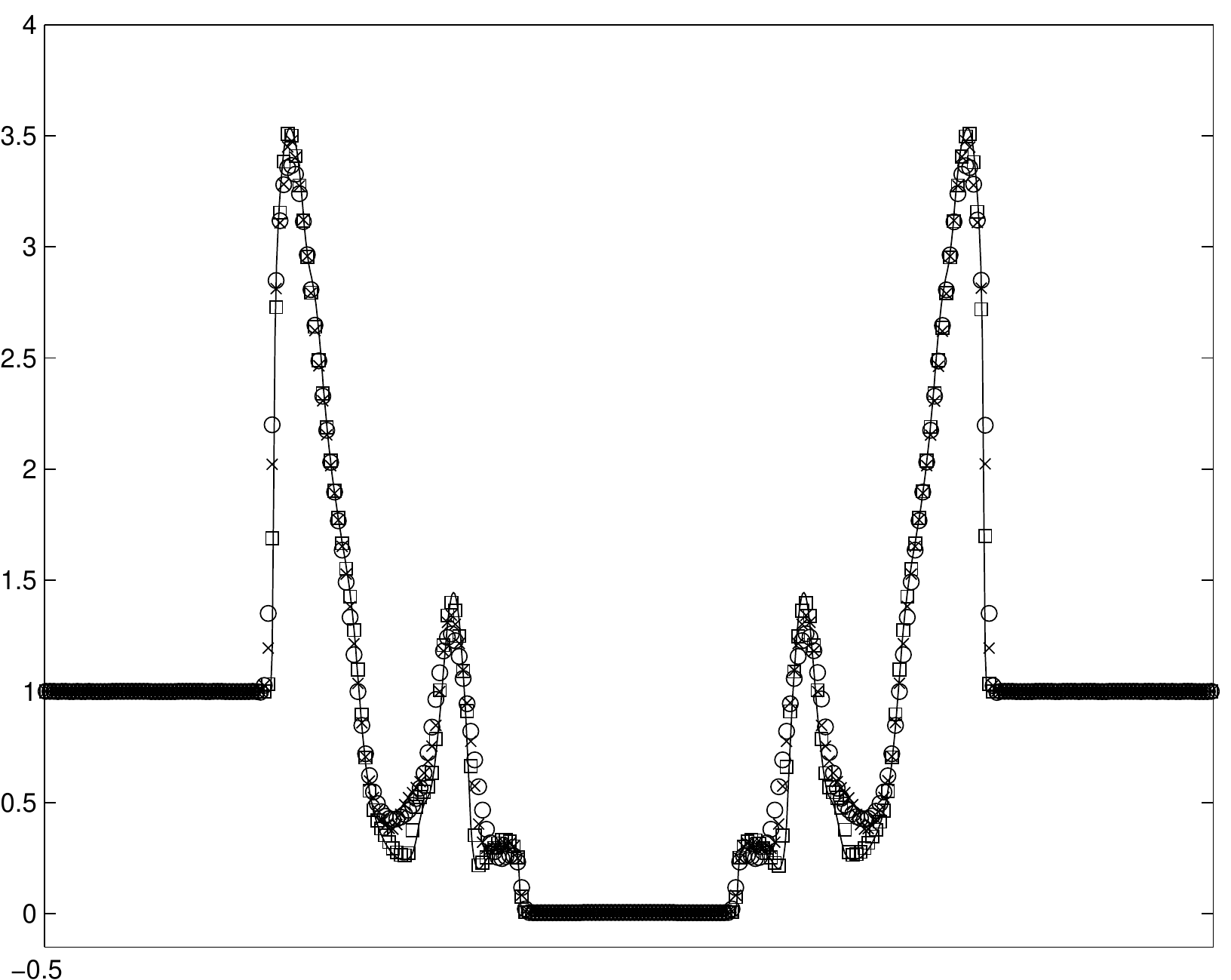}&
        \includegraphics[width=0.35\textwidth]{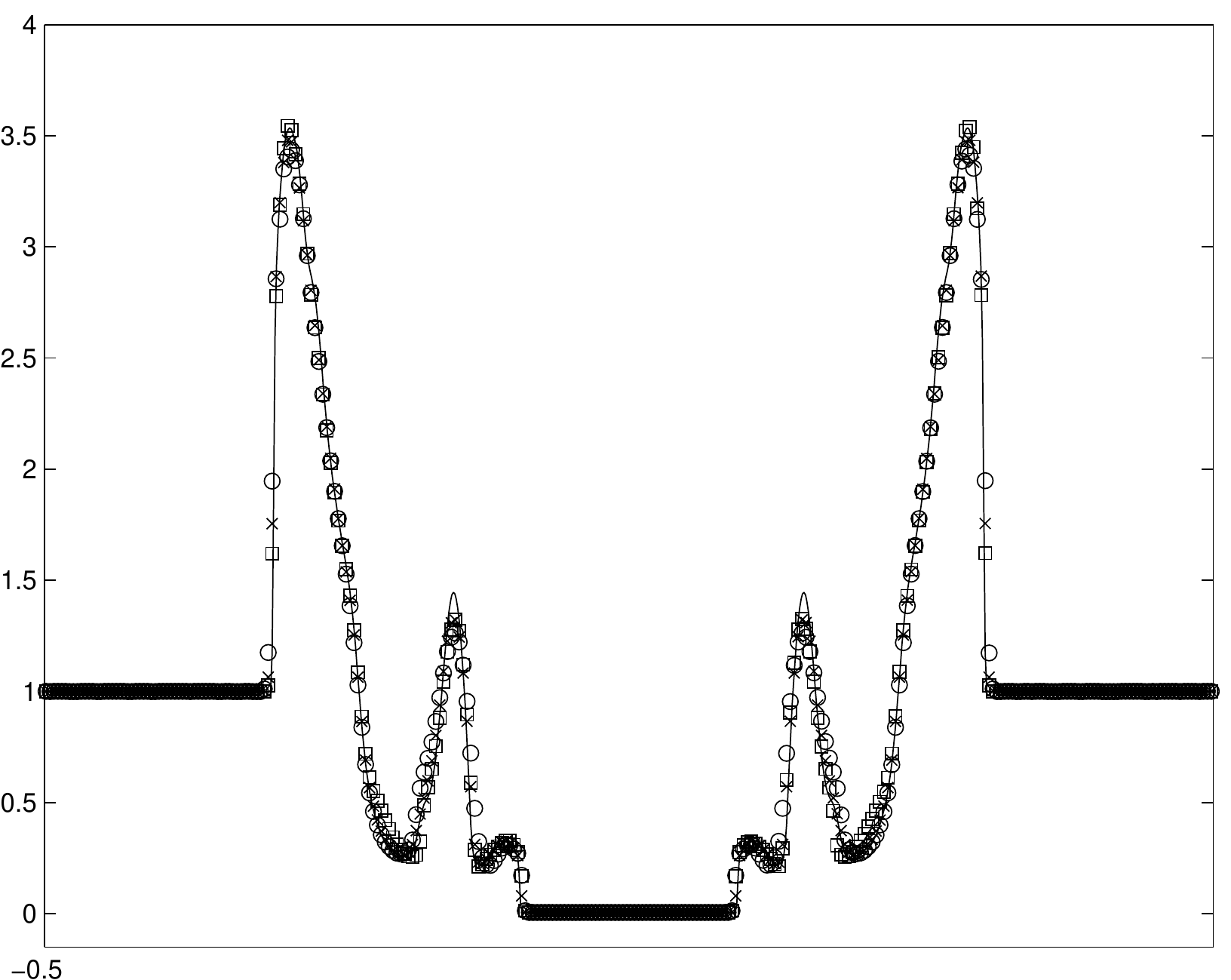}
\end{tabular}
    \caption{Example \ref{exRMHDRotor}:
The densities (top) and gas pressure $p$ (bottom)
at $t=0.4$ along the line $y=x$.
The solid line denotes the reference solution obtained by using
the MUSCL scheme with $800\times 800$ uniform cells, while the symbols
``$\circ$'', ``$\times$'', and ``$\square$'' denotes
      the numerical solutions obtained by the $P^1$-, $P^2$-, and $P^3$-based methods with
      $300\times 300$ cells, respectively. }
    \label{fig:Rocmprhopg}
  \end{figure}

  \begin{Example}[Rotor problem]\label{exRMHDRotor}\rm
  This problem  is a relativistic extension of the classical MHD rotor
test problem \cite{Zanna:2003} and  has been considered in the literature,
  The computational domain is
    $[-0.5,0.5]\times [-0.5,0.5]$ with four outflow conditions.
    Initially, the gas pressure and the
magnetic field are uniform, and there is a disk of fluid centered at
    $(0,0)$ and with high density rotating in a anti-clockwise direction at a
high relativistic speed. The disk radius  is 0.1.
The ambient fluid is homogeneous for $r > 0.115$ and
changing linearly for $0.1\le r\le 0.115$, where $r=\sqrt{x^2+y^2}$.
Specifically, the initial data
are taken as
    $$(\rho,v_x,v_y,v_z,B_x,B_y,B_z,p)=\begin{cases} (10,-\alpha y, \alpha x,
      0,1,0,0,1), & r<0.1,\\
      (1+9\delta, -\alpha y \delta/r, \alpha x \delta/r, 0,1,0,0,1), &
      0.1\le r \le 0.115,\\
      (1,0,0,0,1,0,0,1), & r>0.115,\end{cases}$$
    where  $\alpha=9.95$, $\delta(r):=(0.115-r)/0.015$ is a ``taper''
    function helping to reduce initial transients.
 The adiabatic index is taken as $\Gamma=5/3$.

Figs.~\ref{fig:RMHDRotorrho}, \ref{fig:RMHDRotorpre}, \ref{fig:RMHDRotorpm},
and \ref{fig:RMHDRotorgam} show respectively
the contour plots of density $\rho$, gas pressure $p$, magnetic pressure $p_m$,
and Lorentz factor $\gamma$ at $t=0.4$ obtained by using DG methods with $300\times 300$ cells.
 The CFL numbers of  $P^1$-, $P^2$-, $P^3$-based non-central DG
methods are $0.2,~0.15,~0.1$, respectively, while
those of corresponding central DG are $0.3/\theta,~0.25/\theta,~0.2/\theta$ with
$\theta=0.3$, that is,
    $\Delta t_n=0.3\tau_n$.
    The parameter $M$ in the modified TVB minmod function is taken as $M=500$.
    From those results, one can see that as
the time increases, the winding magnetic field lines are formed and
decelerated the disk speed, and the initial high density at the center
is swept away completely and a oblong-shaped shell is formed.
At $t = 0.4$, the central magnetic field
lines are rotated almost 90$^\circ$.
Those computed results agree quite well with other published
solutions.
Fig.~\ref{fig:Rocmprhopg}
gives a comparison of the density and
gas pressure at $t=0.4$ along the line $y=x$ to
the reference solutions, which are obtained by using the MUSCL
scheme with $800\times 800$ uniform cells.
It is seen that the resolution of central DG methods
is slightly better than corresponding non-central DG methods,
and the numerical solution of the high order methods
is in good agreement with the reference solution.
Moreover, high order methods identify  more ``troubled'' cells, see
Table \ref{tab:cellperRMHDSc} for the percentage of ¡°troubled¡± cells.
Corresponding CPU times given in Table \ref{tab:RMHDSccmpRC}
show that the CPU time used in central DG methods
is significantly more than that of the same order non-central DG methods.
    \end{Example}

  {}

\begin{figure}[!htbp]
    \centering{}
  \begin{tabular}{cc}
    \includegraphics[width=0.35\textwidth]{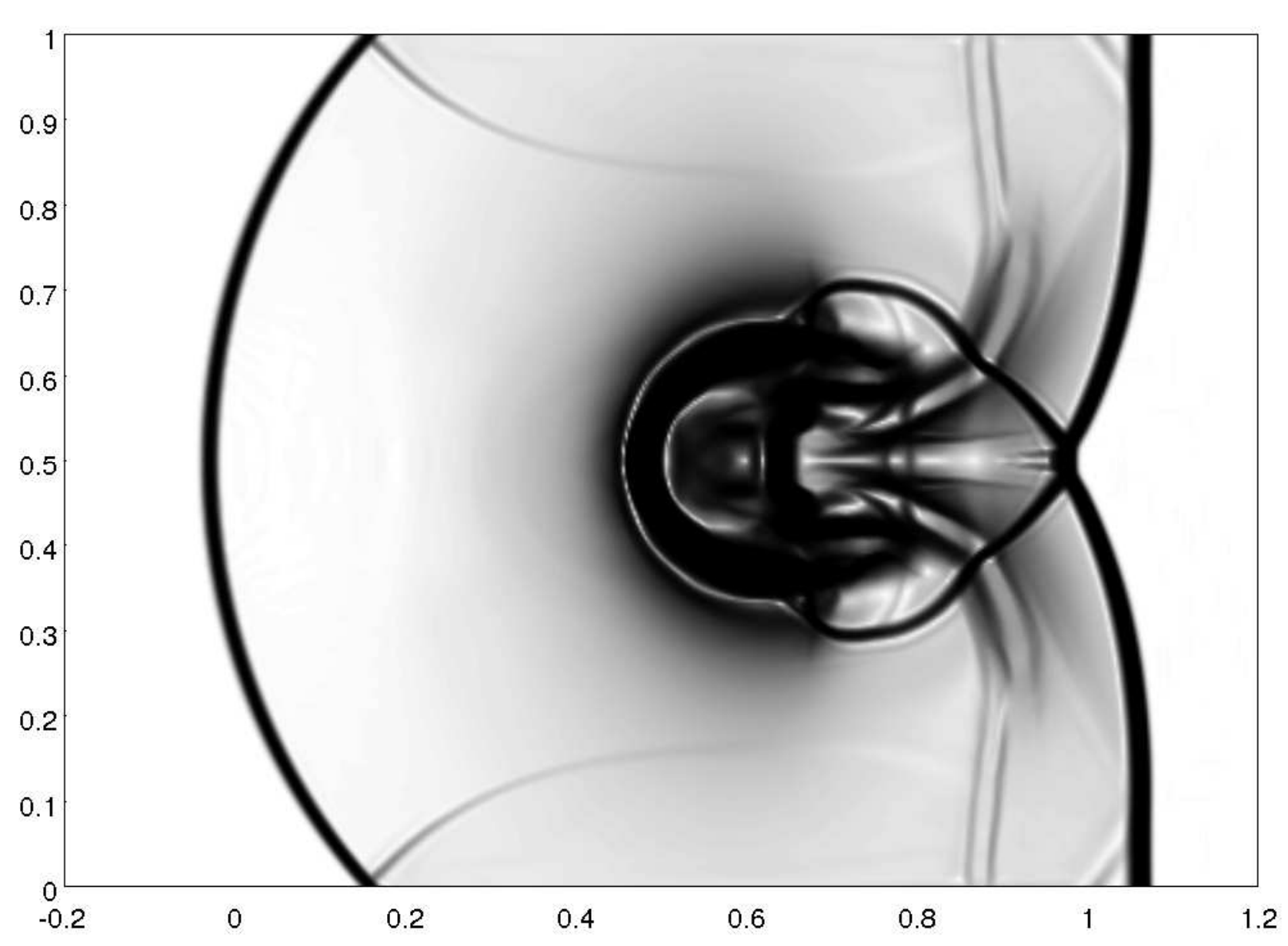}&
 \includegraphics[width=0.35\textwidth]{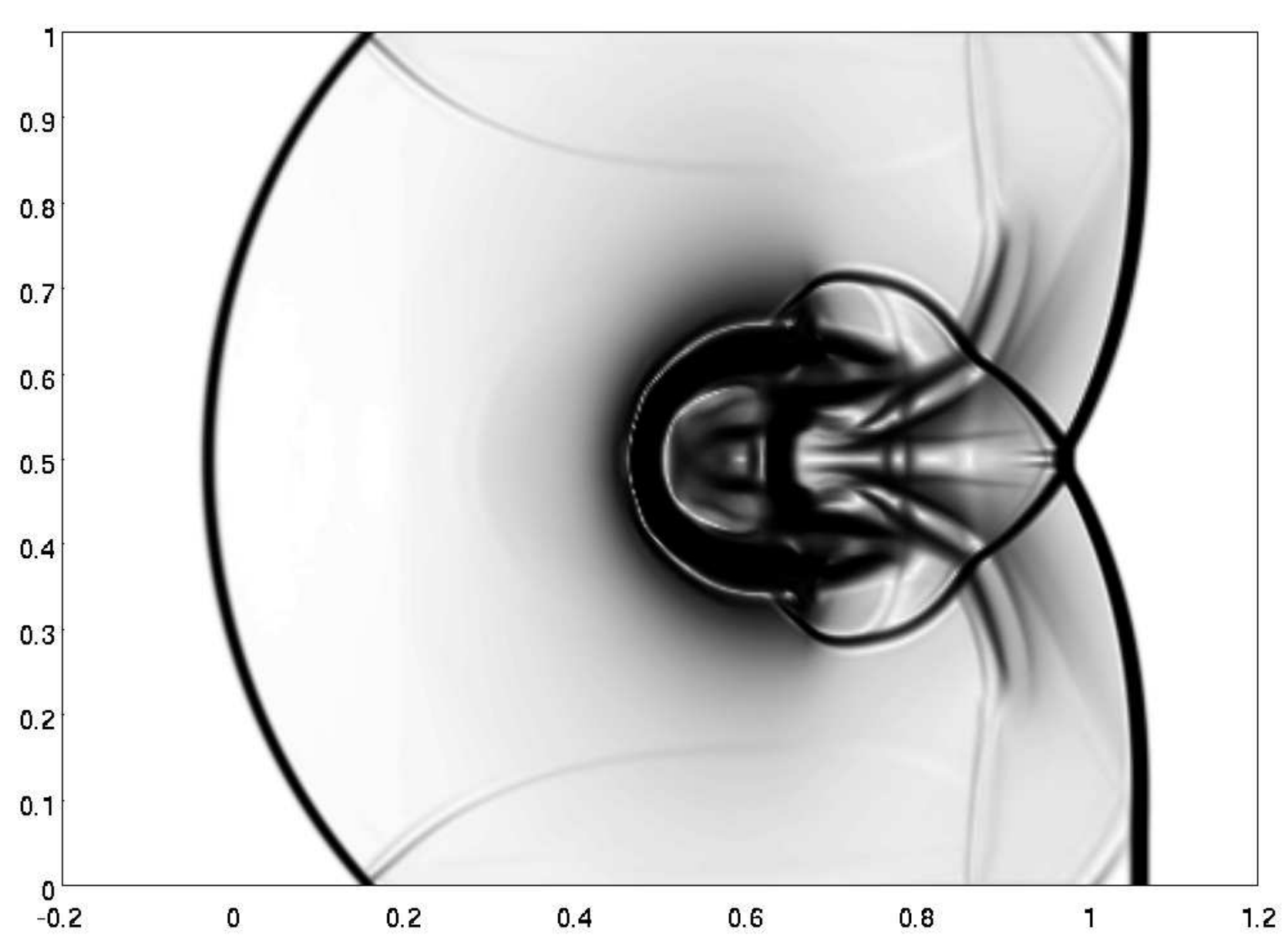}\\
    \includegraphics[width=0.35\textwidth]{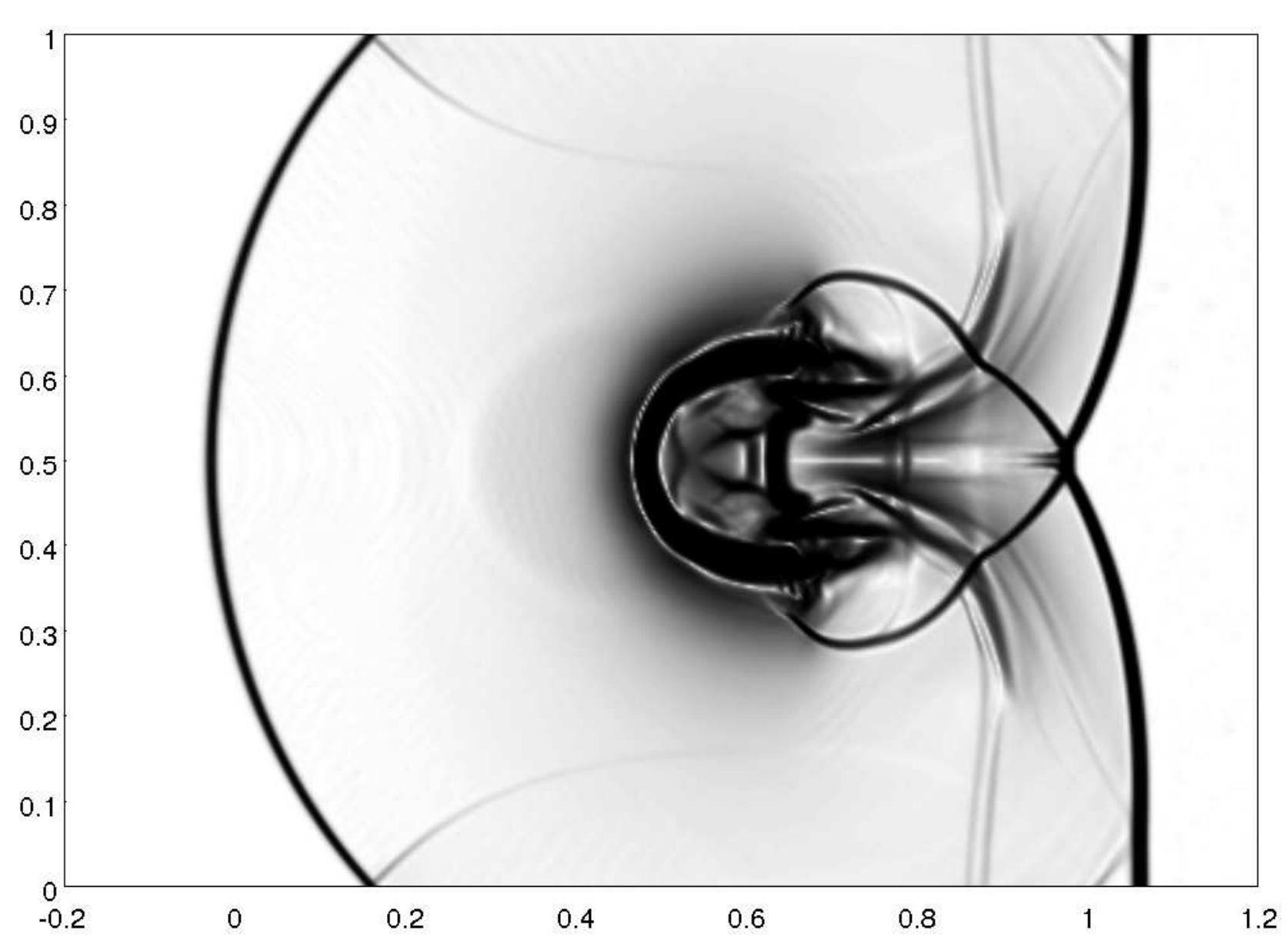}&
 \includegraphics[width=0.35\textwidth]{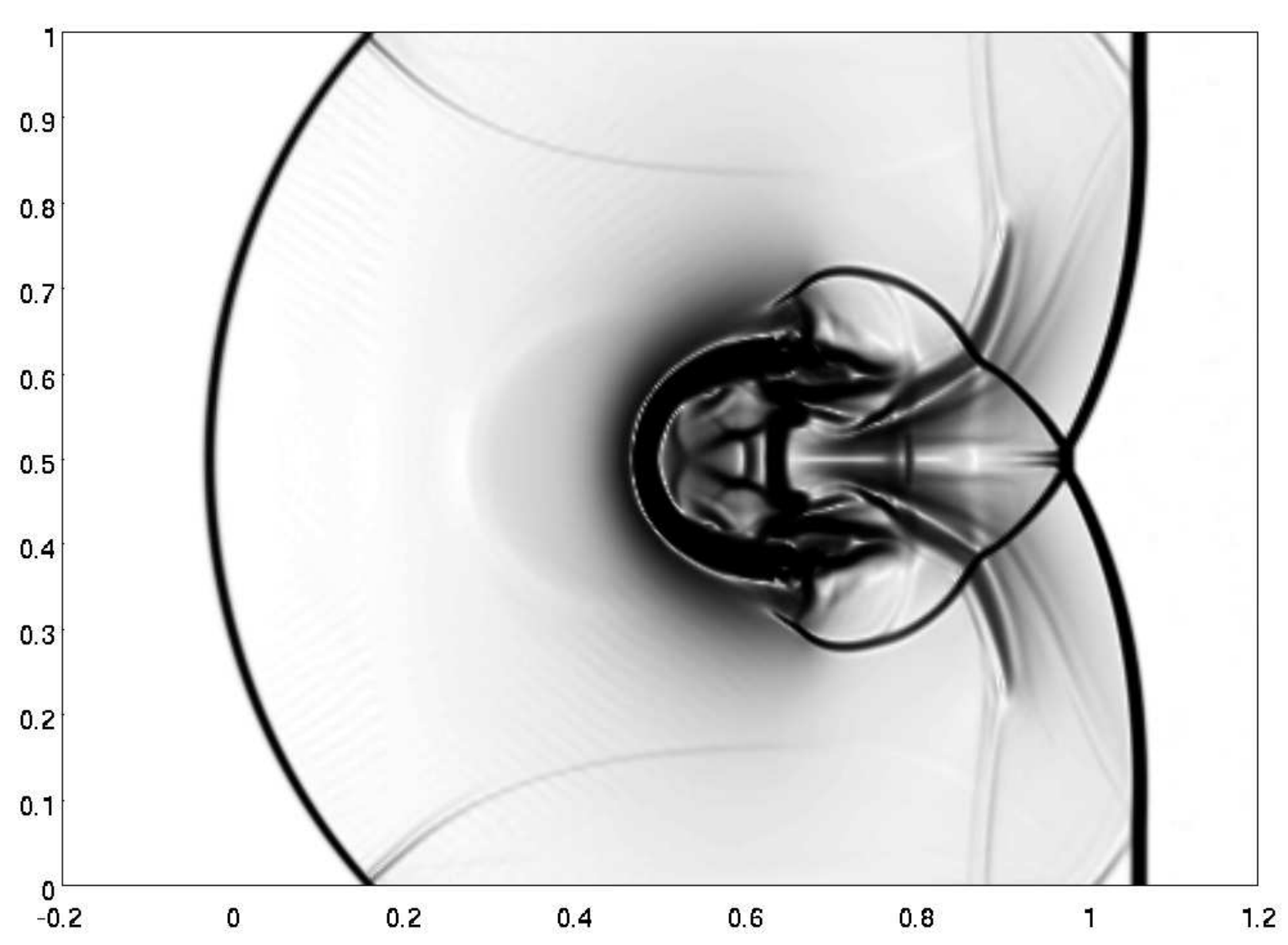}\\
 \includegraphics[width=0.35\textwidth]{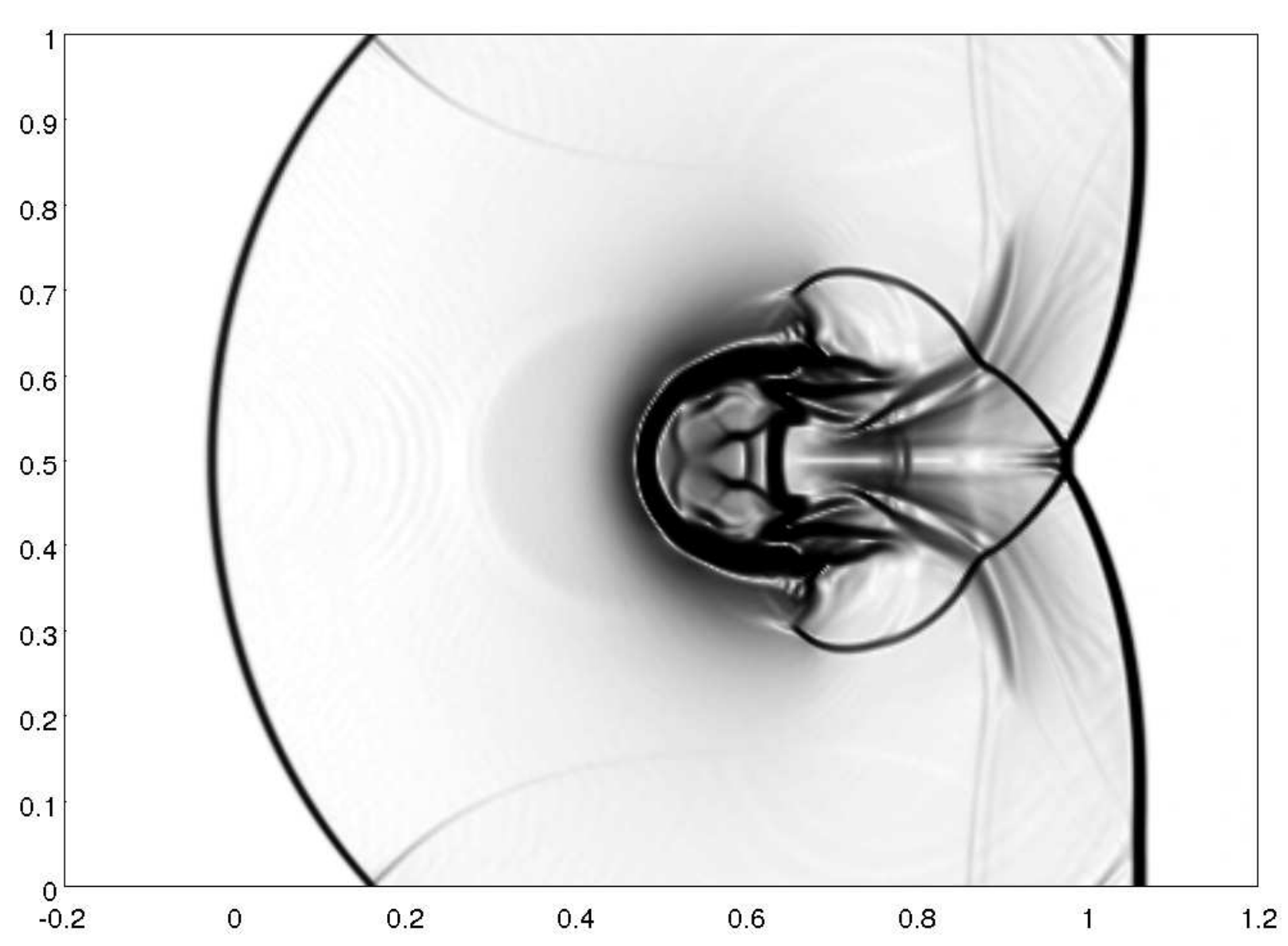}&
  \includegraphics[width=0.35\textwidth]{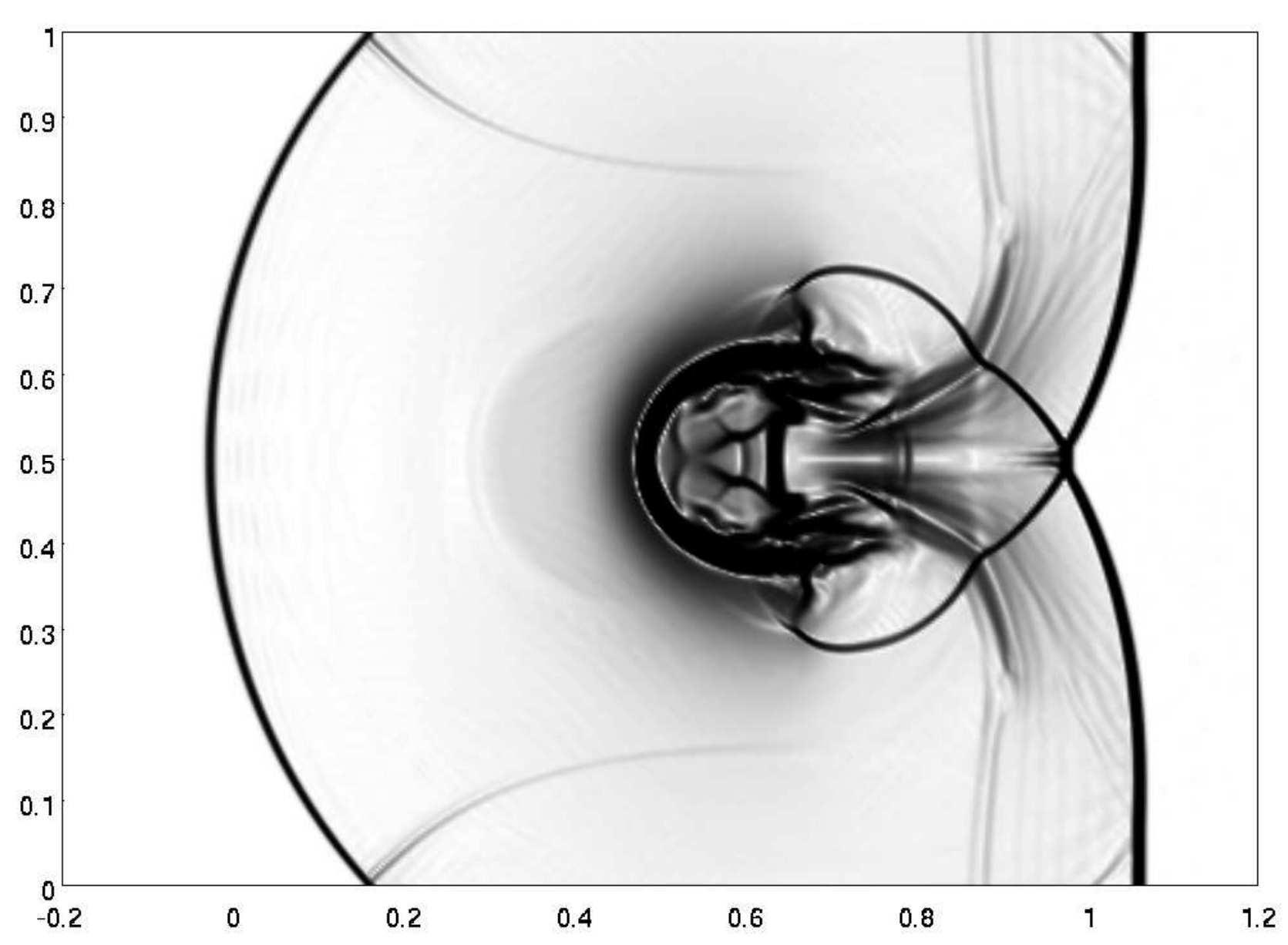}\\
    \end{tabular}
    \caption{Example~\ref{exRMHDSC}:   Schlieren images of density $\rho$ at $t=1.2$ obtained with $420\times 300$ cells.
Left: $P^K$-based \DG{}; right: $P^K$-based \CDG{}.
From top to bottom: $K=1,~2,~3$.  }
    \label{fig:RMHDSCrho}
  \end{figure}

\begin{figure}[!htbp]
    \centering{}
  \begin{tabular}{cc}
    \includegraphics[width=0.35\textwidth]{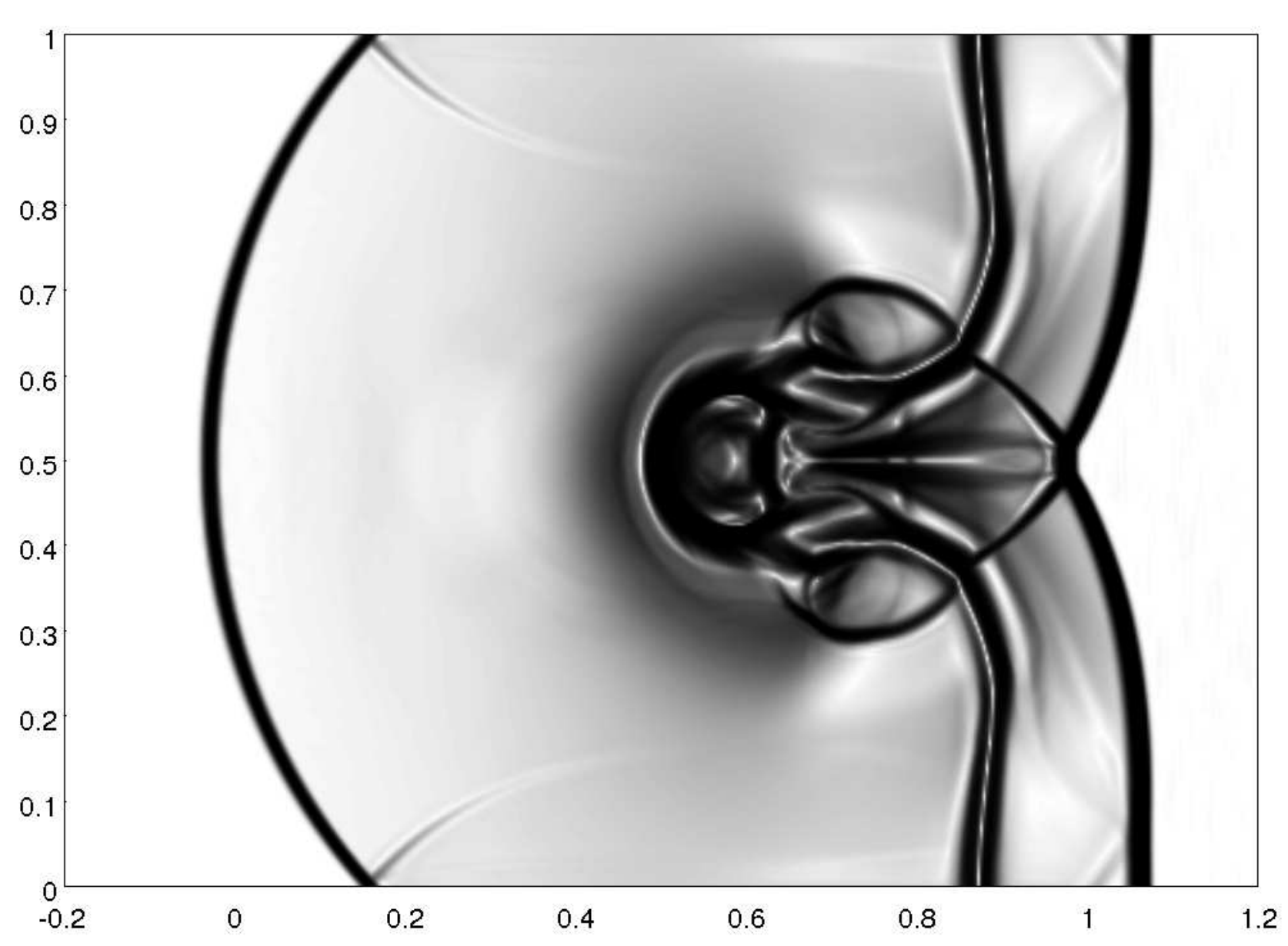}&
 \includegraphics[width=0.35\textwidth]{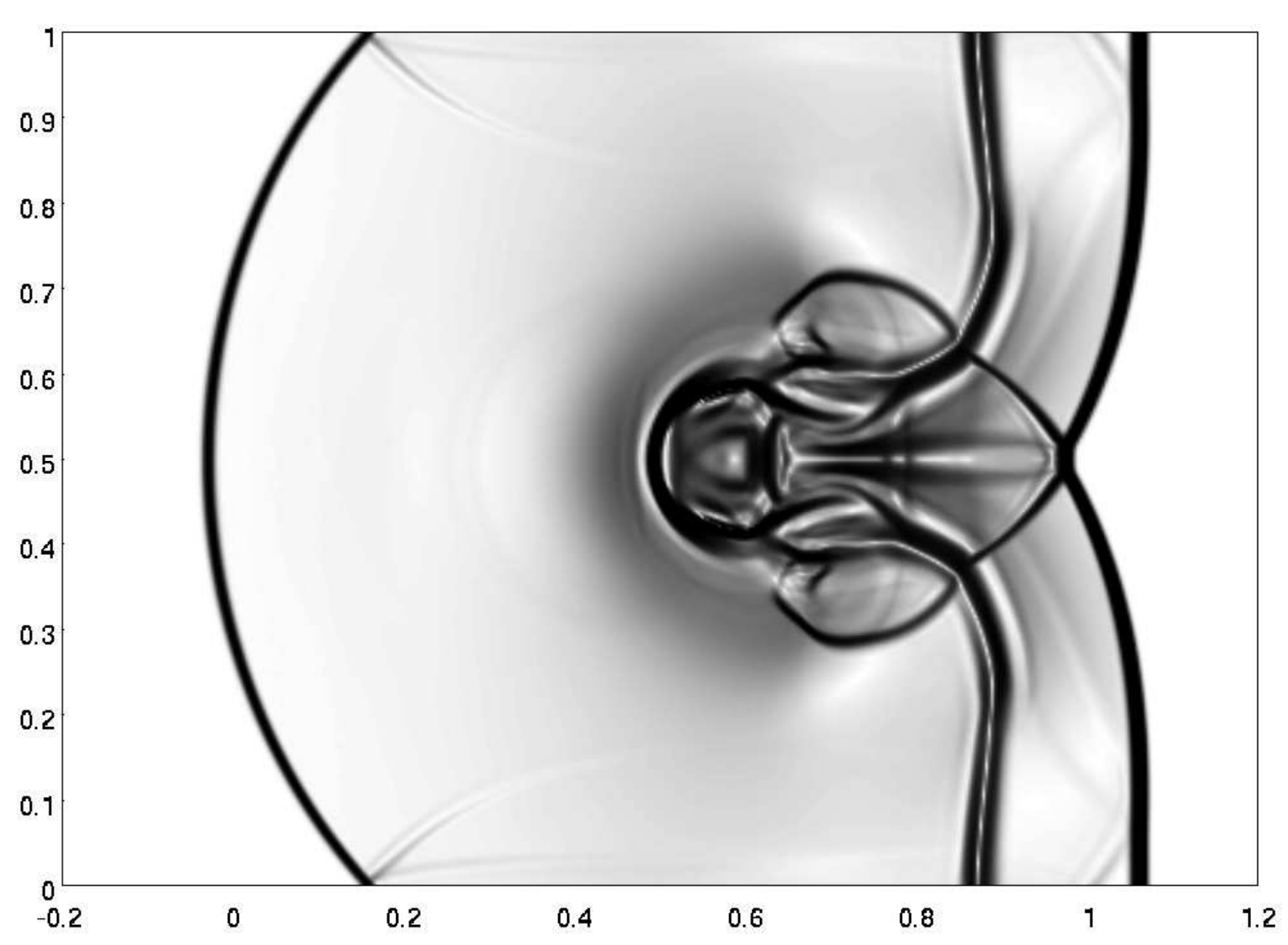}\\
     \includegraphics[width=0.35\textwidth]{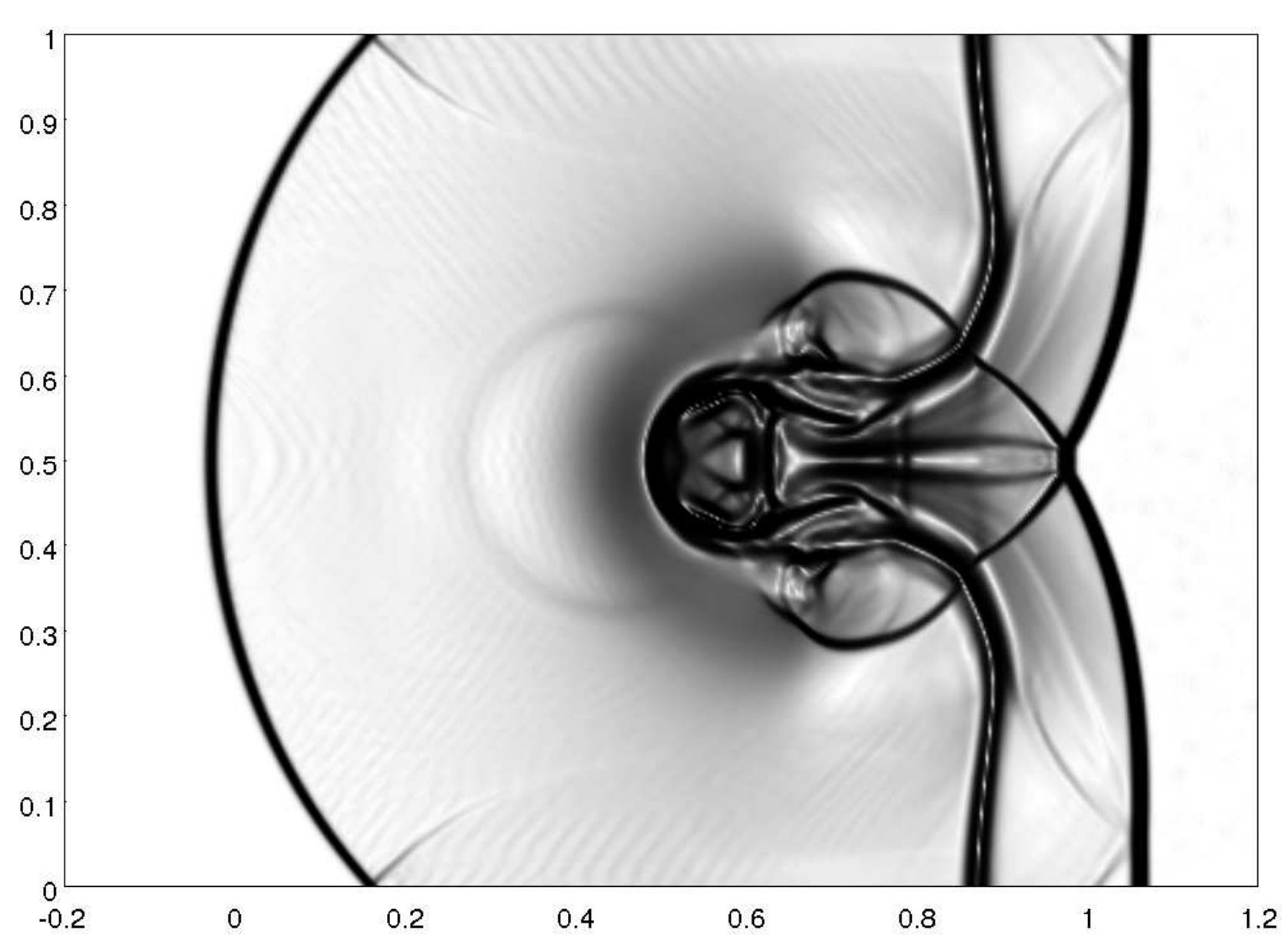}&
 \includegraphics[width=0.35\textwidth]{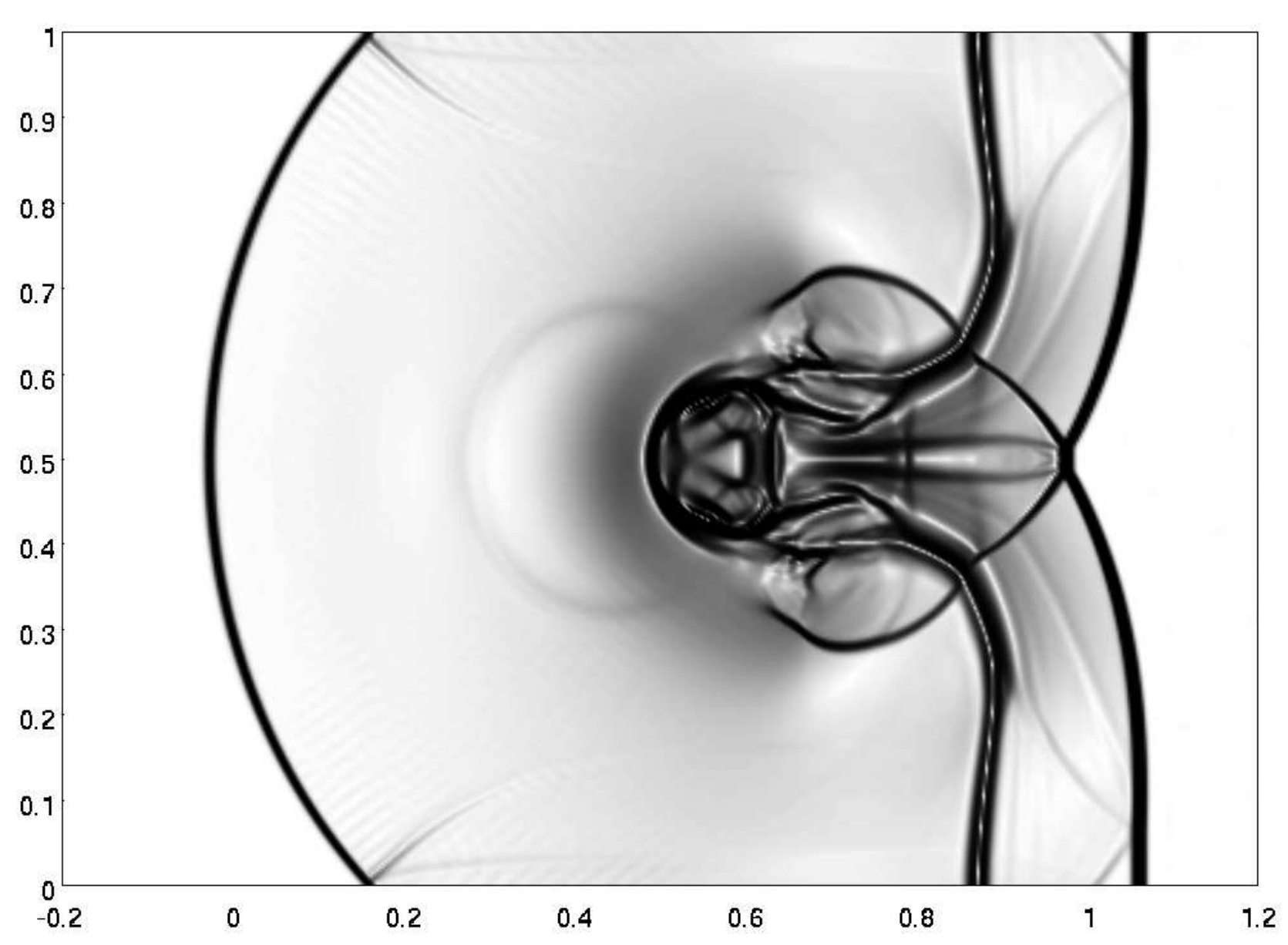}\\
 \includegraphics[width=0.35\textwidth]{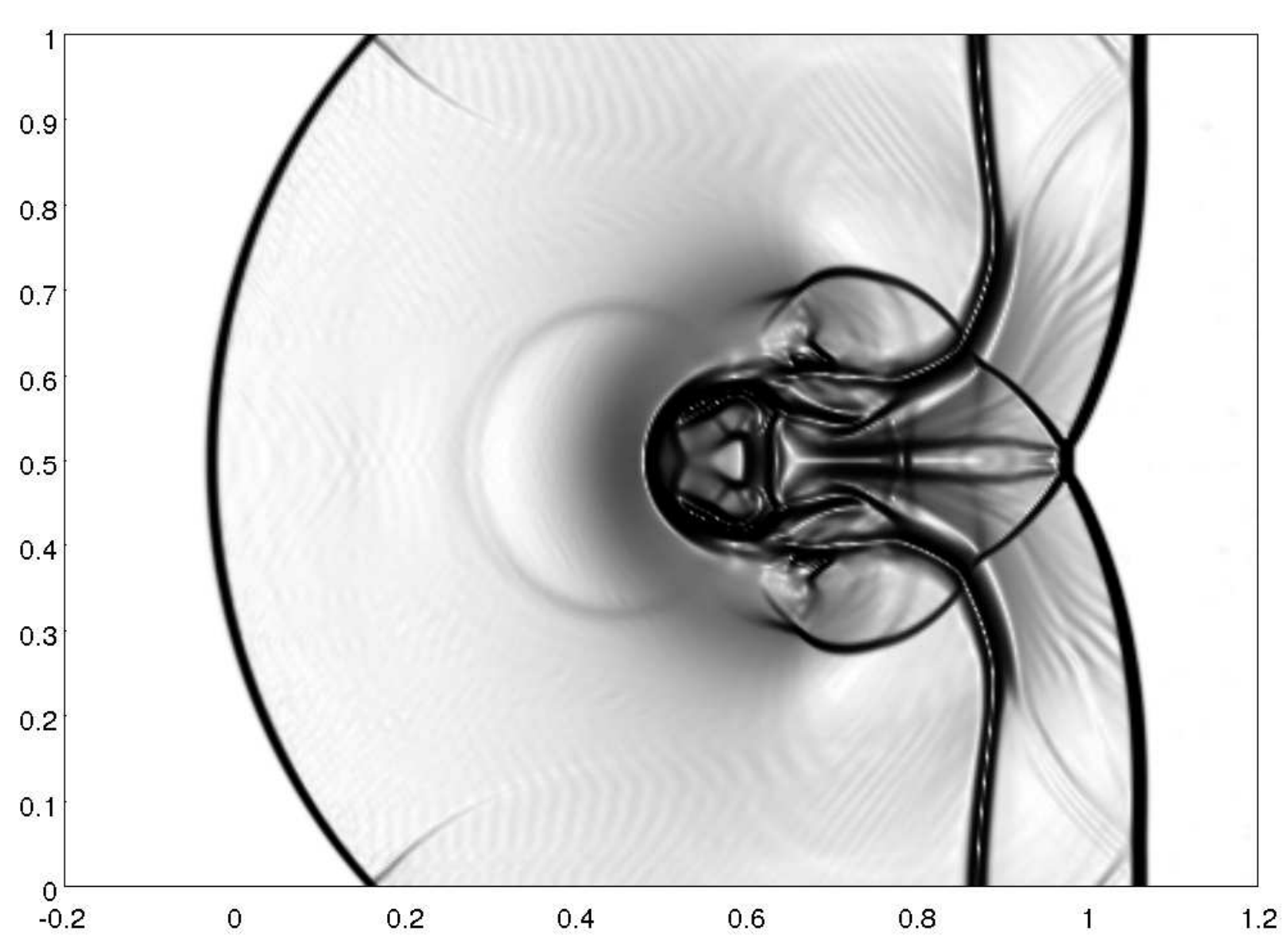}&
 \includegraphics[width=0.35\textwidth]{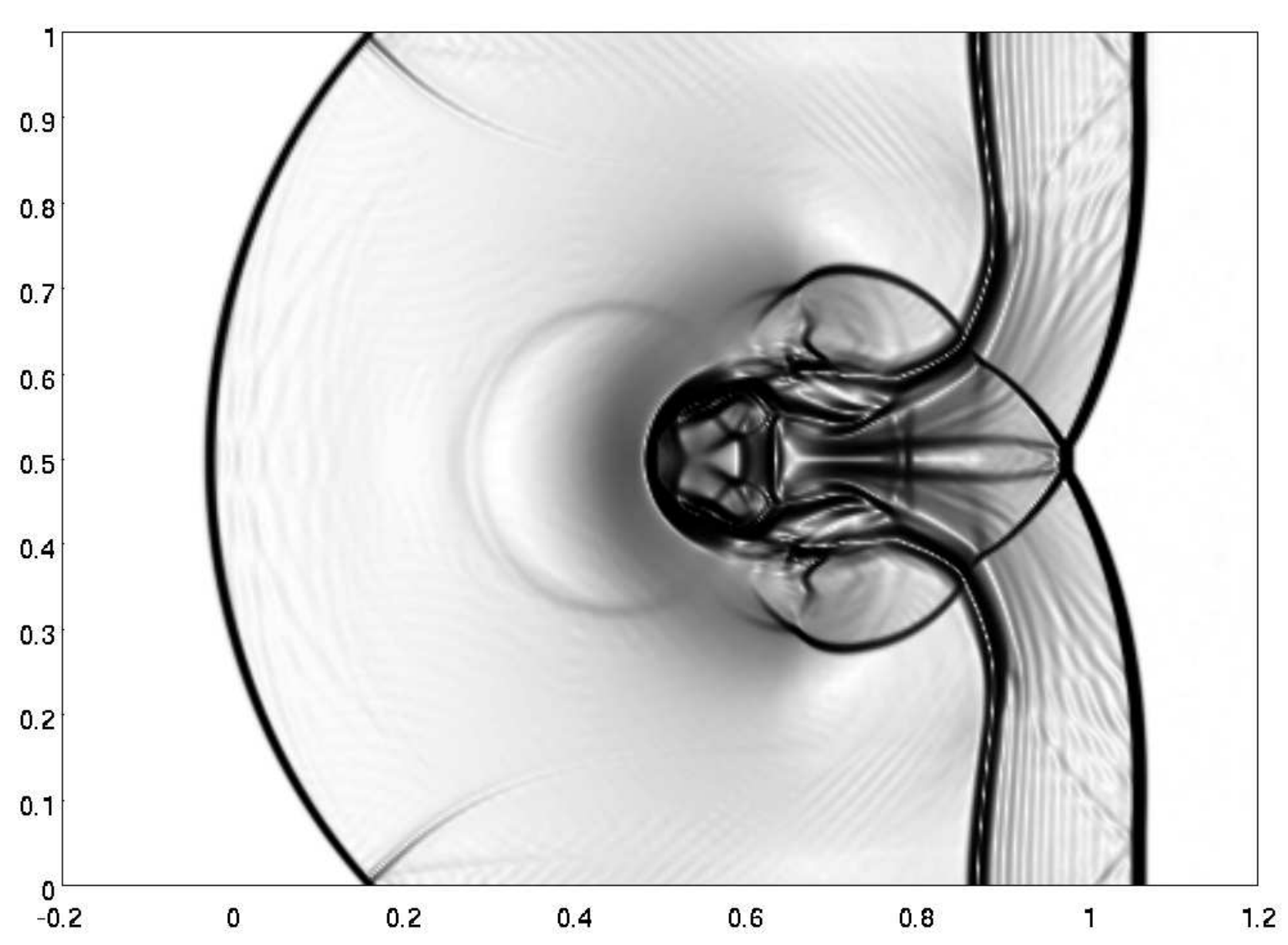}\\
    \end{tabular}
    \caption{Same as Fig.~\ref{fig:RMHDSCrho} except for the magnetic pressure $p_m$.}
    \label{fig:RMHDSCpm}
  \end{figure}

\begin{figure}[!htbp]
    \centering{}
    \begin{tabular}{cc}
    \includegraphics[width=0.35\textwidth]{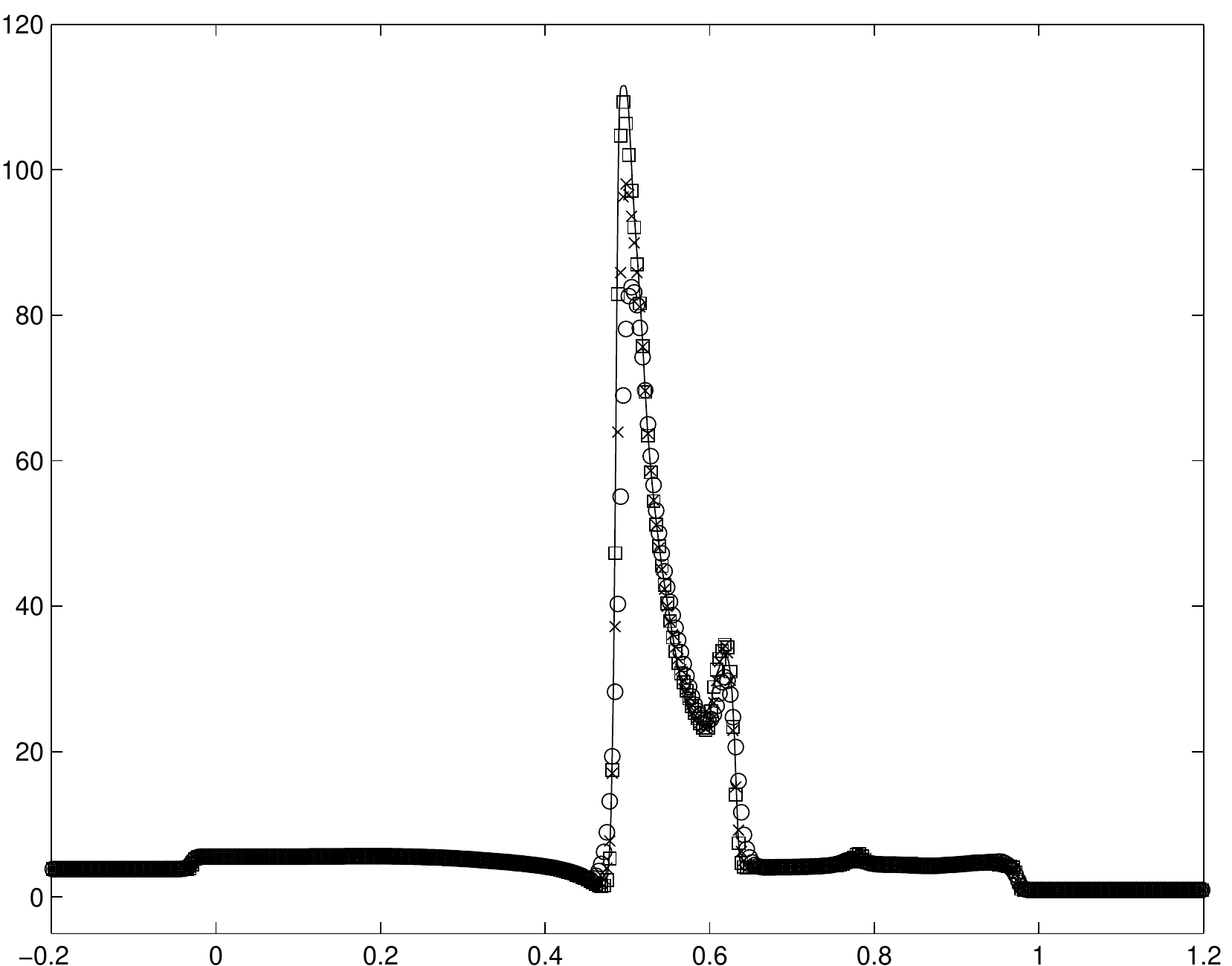}&
        \includegraphics[width=0.35\textwidth]{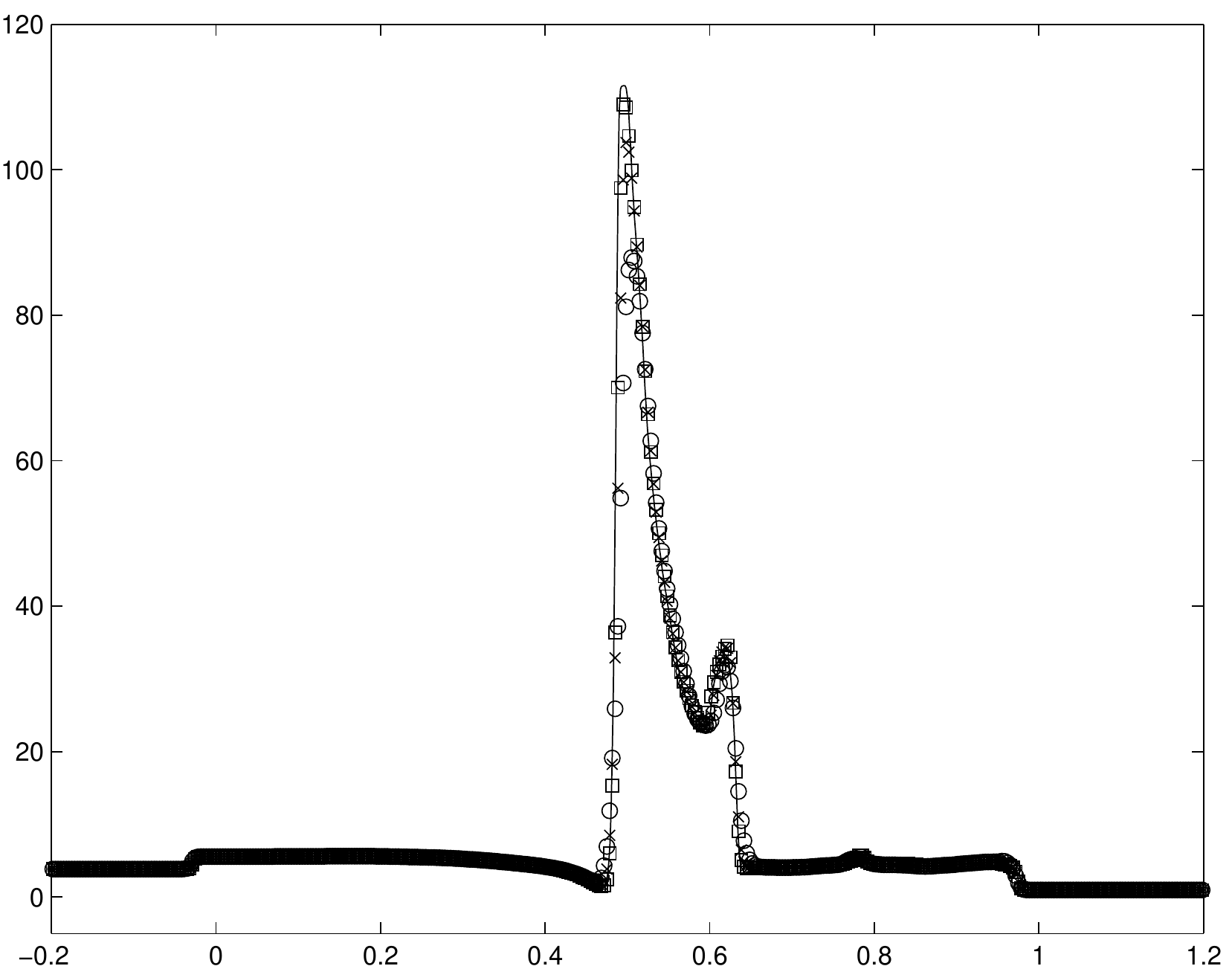}\\
         \includegraphics[width=0.35\textwidth]{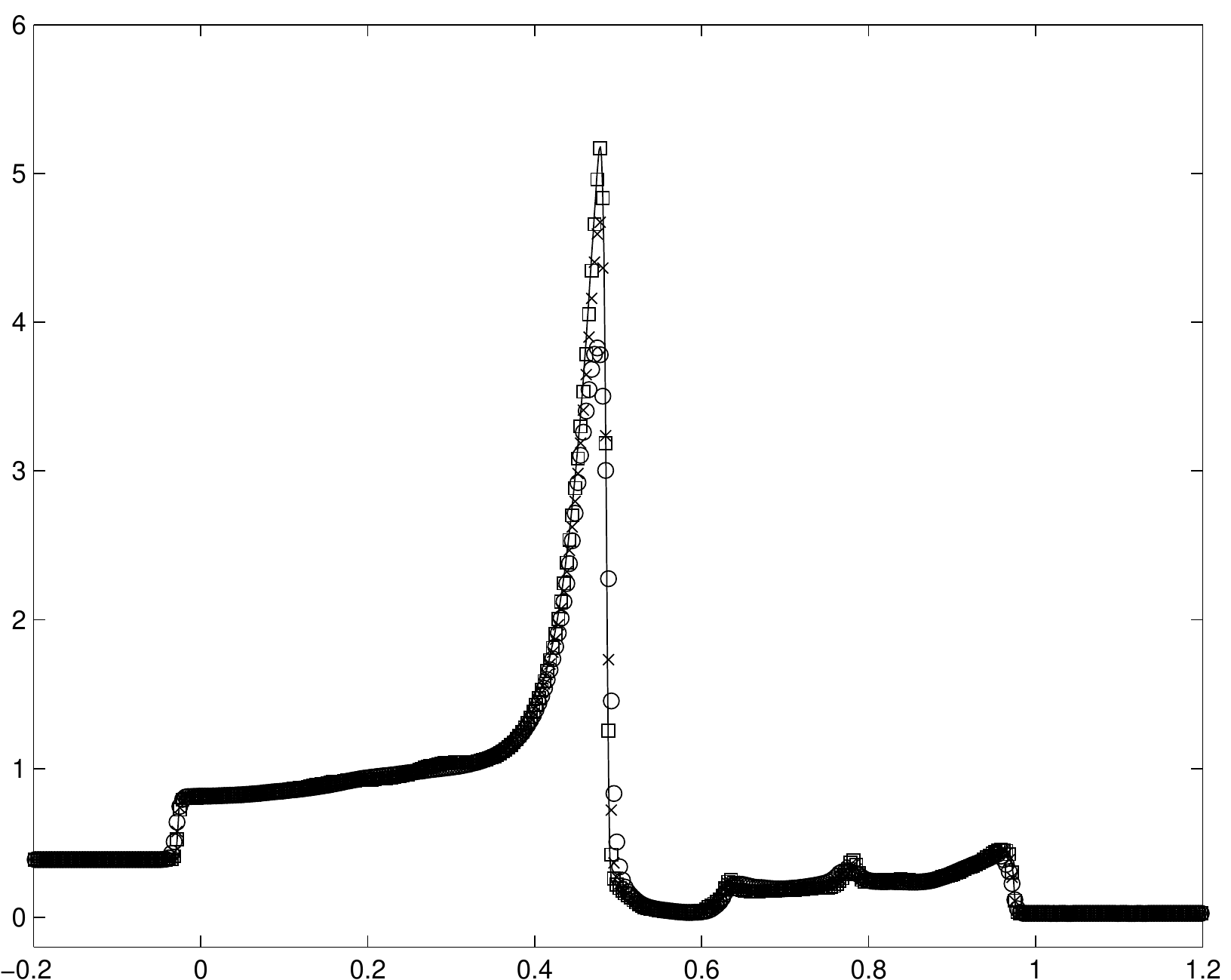}&
        \includegraphics[width=0.35\textwidth]{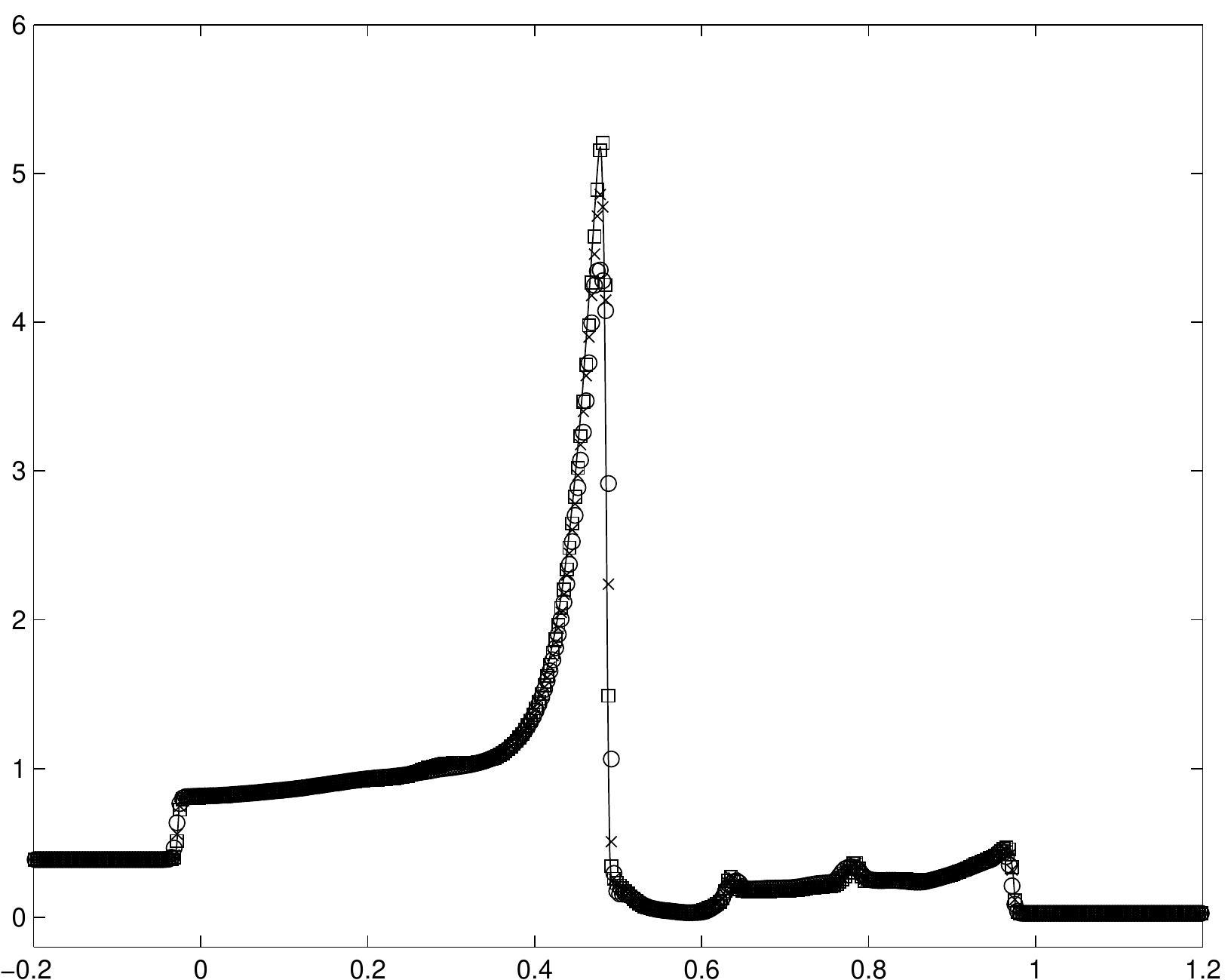}\\
        \end{tabular}
    \caption{Example~\ref{exRMHDSC}:
    The density $\rho$ (top) and the magnetic pressure $p_m$ (bottom) at
    $t=1.2$ along the line $y=0.5$. The solid line denotes the reference solution obtained by using the MUSCL scheme with $980\times 700$ cells, while the symbols ``$\circ$'', ``$\times$'', and ``$\square$''
    denote the numerical solutions obtained by the $P^1$-, $P^2$-, and $P^3$-based methods with $420\times 300$ cells, respectively. Left: \DG{}; right: \CDG{}. }
    \label{fig:Sccmprhopm}
  \end{figure}

     \begin{Example}[Shock and cloud interaction]\label{exRMHDSC}\rm
     The shock-cloud interaction
problem modeling the disruption of a high density cloud
by a strong shock wave has been widely used to test the classical
MHD codes, see e.g. \cite{Dai:1998,Ross:2004,Toth:2000}. It is extended to
the  relativistic case with the magnetic field orthogonal
to the slab plane in \cite{MignoneHLLCRMHD}
so that any magnetic divergence-free treatment
is not needed.
 Here we consider a different extension  of this
problem \cite{HeAdaptiveRMHD}, in which the magnetic field is not orthogonal to the
slab plane so that the magnetic divergence-free treatment has
to be imposed.
 The adiabatic index  $\Gamma$ and  computational domain are $5/3$ and $[-0.2,1.2]\times [0,1]$,
 respectively.
 The inflow condition is specified on the left boundary, while the outflow conditions
are  on other boundaries.
Initially, there is a right-moving shock wave located at $x=0.05$
    and with the left and right states
    \begin{align*}
    &(\rho,v_x,v_y,v_z,B_x,B_y,B_z,p)(x,y,0)
    \\
   & =\begin{cases} (3.86859,0.68, 0,
      0,0,0.84981,-0.84981,1.25115), & x<0.05,\\
      (1,0,0,0,0,0.16106,0.16106,0.05), & x>0.05.\end{cases}
    \end{align*}
    There exists a rest circular cloud which is in magneto-hydrostatic balance
with the surrounding fluid and centered at the point (0.25, 0.5) with a
high density $\rho = 30$ and radius 0.15.

        Figs.~\ref{fig:RMHDSCrho} and \ref{fig:RMHDSCpm}
        show the the densities $\rho$ and magnetic pressures $p_m$ at $t=1.2$ obtained
        by the proposed DG methods with $420\times 300$ cells.
    The CFL numbers of $P^1$-, $P^2$-, $P^3$-based non-central DG methods
    are $0.2,~0.15,~0.1$, respectively, while those of corresponding central
    DG methods are $0.3/\theta,~0.25/\theta,~0.2/\theta$, with $\theta=0.3$, that is,
    $\Delta t_n=0.3\tau_n$.
One can see from those plots that compared to the low order methods,
high order methods  better capture the reflected wave structure and  resolved the complex wave structure generated due to  the interaction between the shock wave and cloud,
while the number of ``troubled'' cells identified by the higher order methods
is larger, see
Table \ref{tab:cellperRMHDSc}. Fig.~\ref{fig:Sccmprhopm} displays
the densities $\rho$ and magnetic pressures $p_m$
at $t=1.2$ along $y=0.5$. It finds
the good agreement between the numerical solutions  obtained by the higher order methods
and the reference solutions, which are obtained by using the MUSCL scheme with $980\times 700$ uniform cells.
The solutions of   higher order methods  are obviously better than the low order methods.

        \end{Example}

  \begin{figure}[!htbp]
    \centering{}
  \begin{tabular}{cc}
    \includegraphics[width=0.35\textwidth]{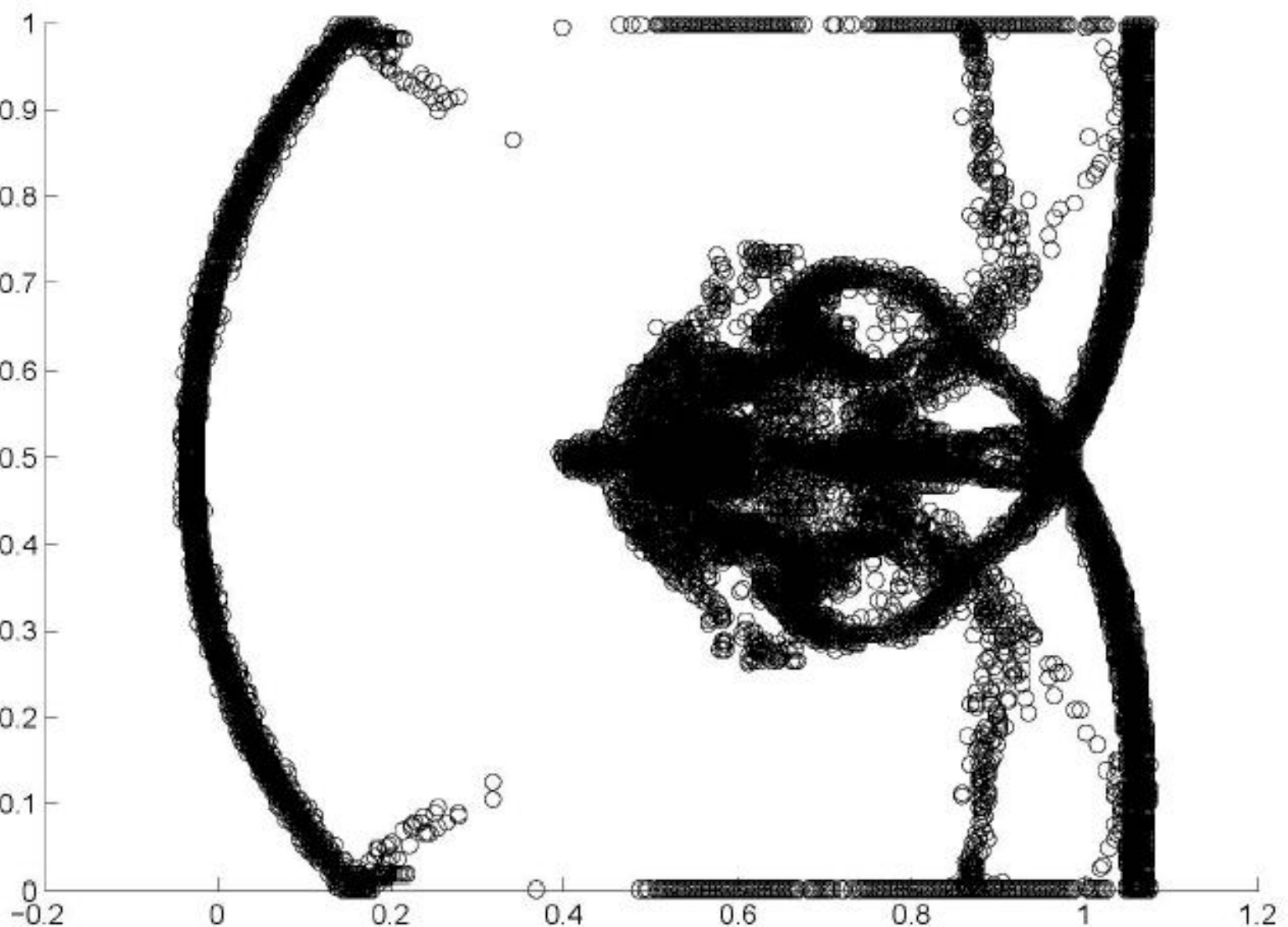}&
 \includegraphics[width=0.35\textwidth]{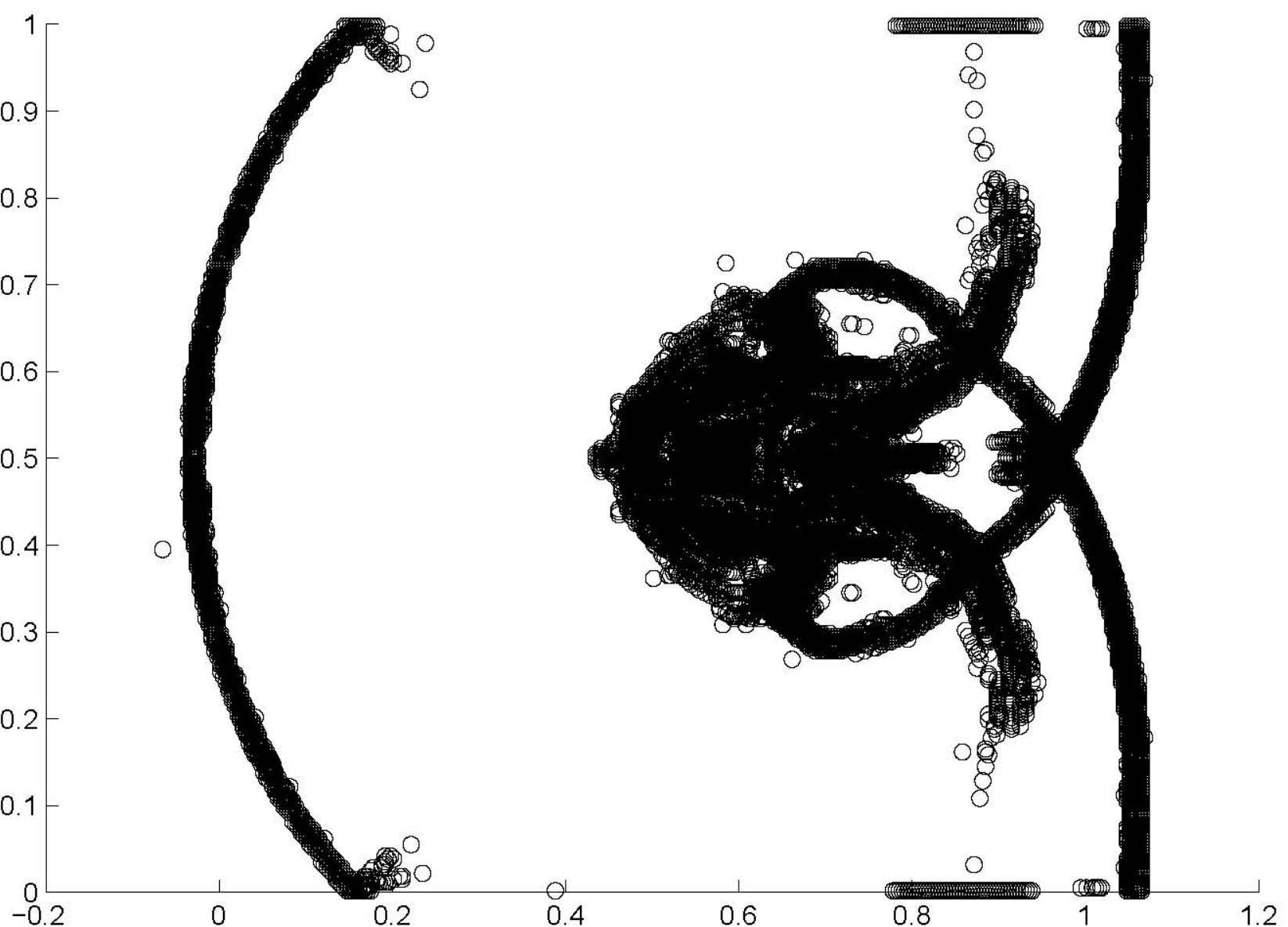}\\
     \includegraphics[width=0.35\textwidth]{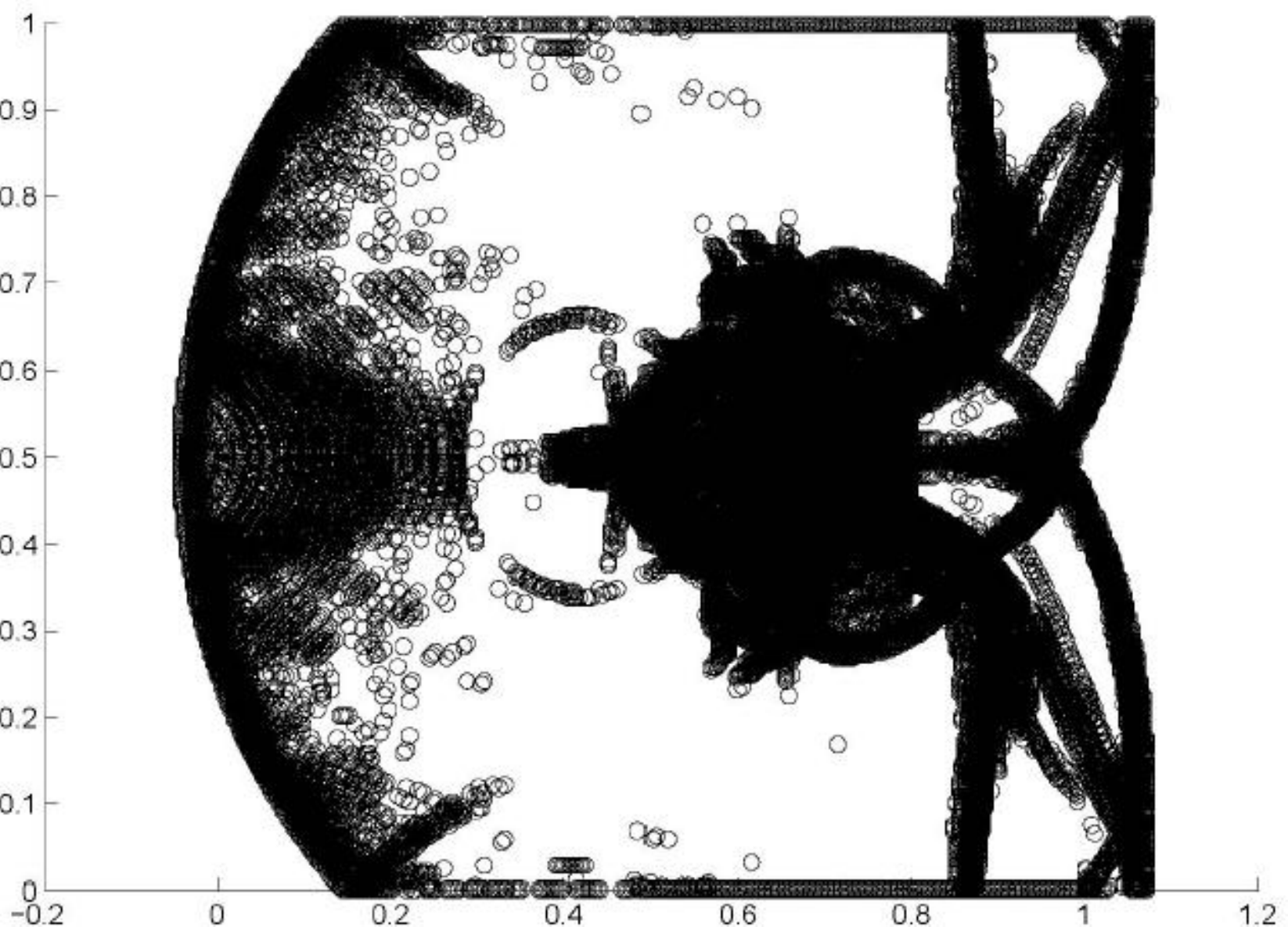}&
 \includegraphics[width=0.35\textwidth]{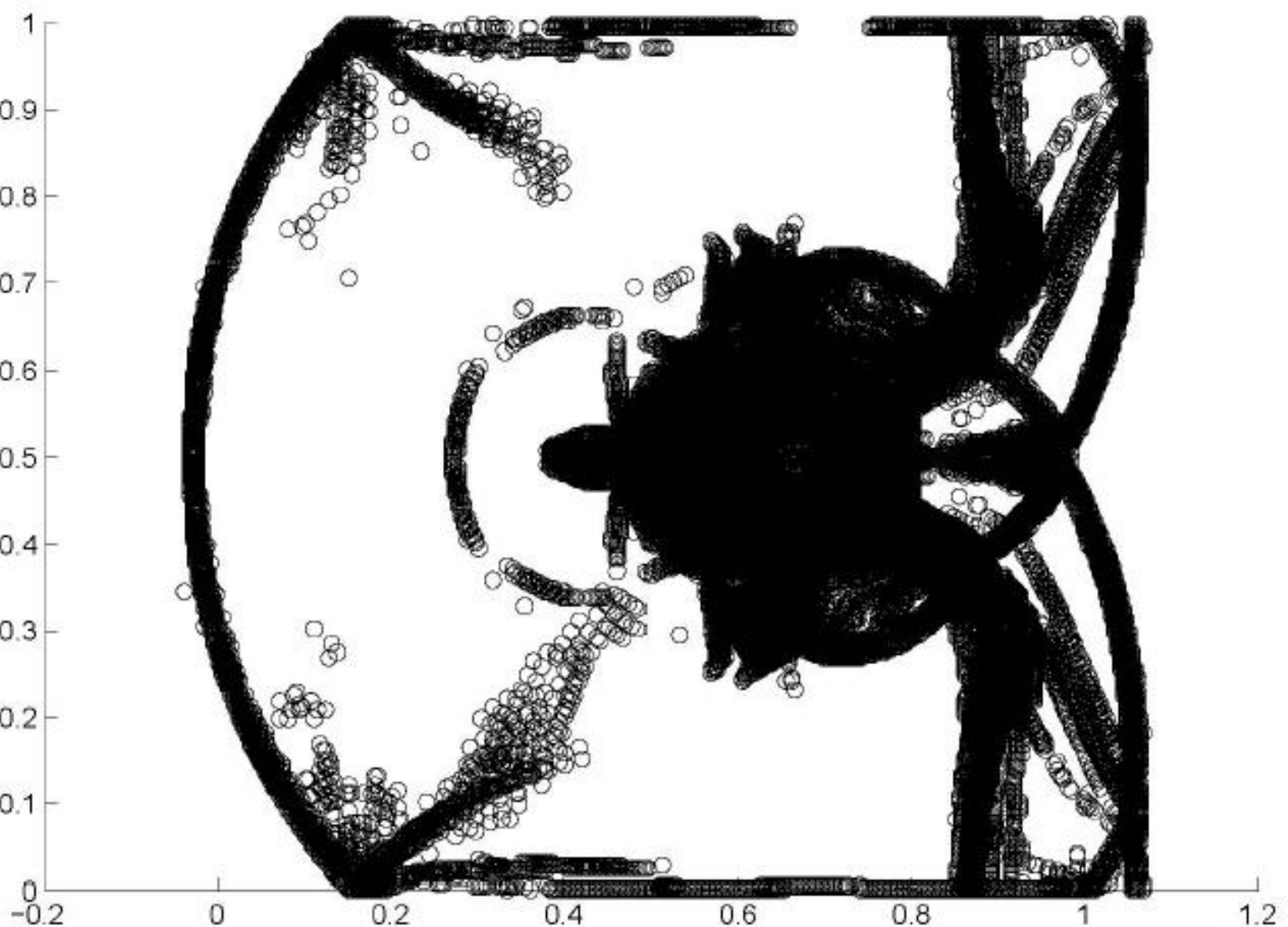}\\
  \includegraphics[width=0.35\textwidth]{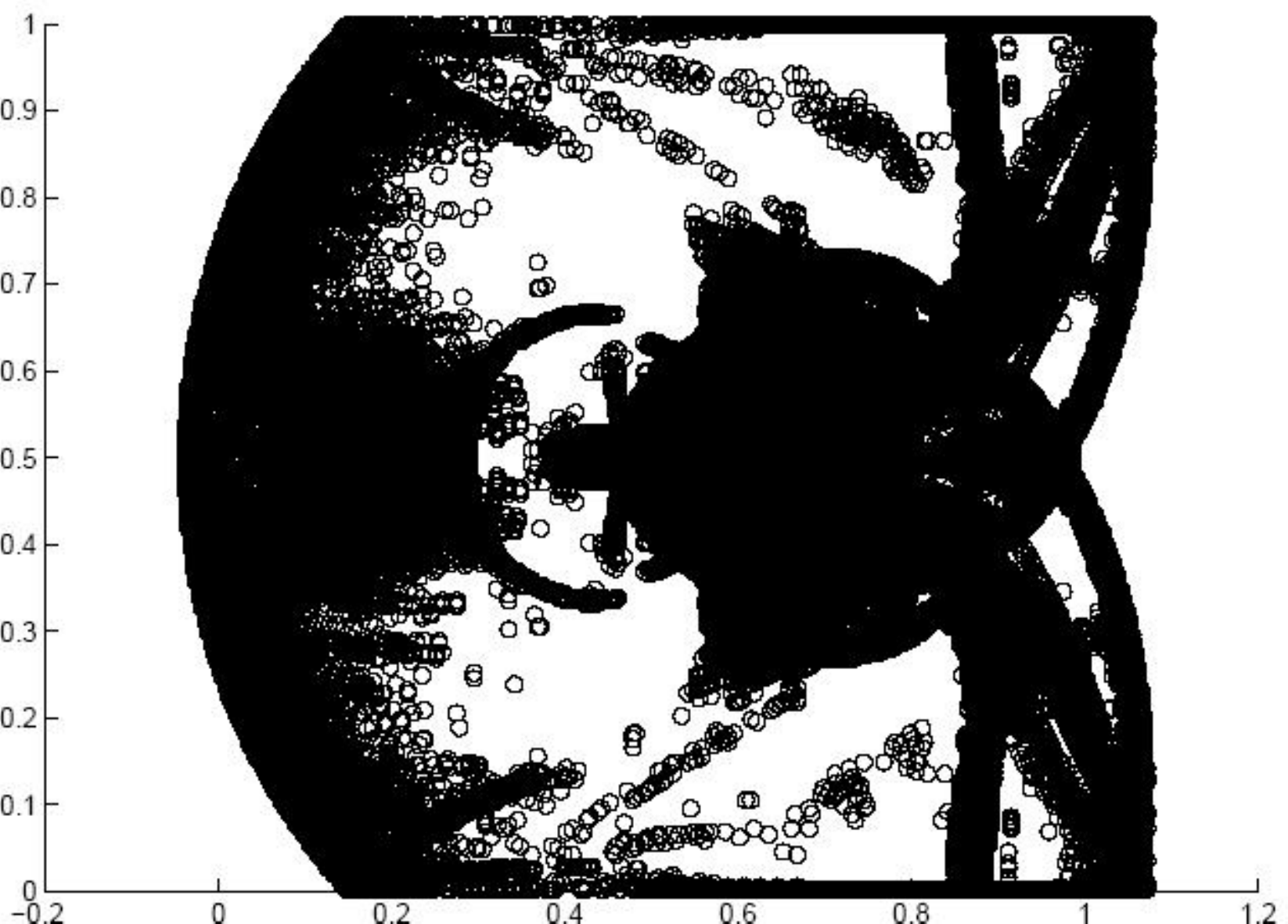}&
  \includegraphics[width=0.35\textwidth]{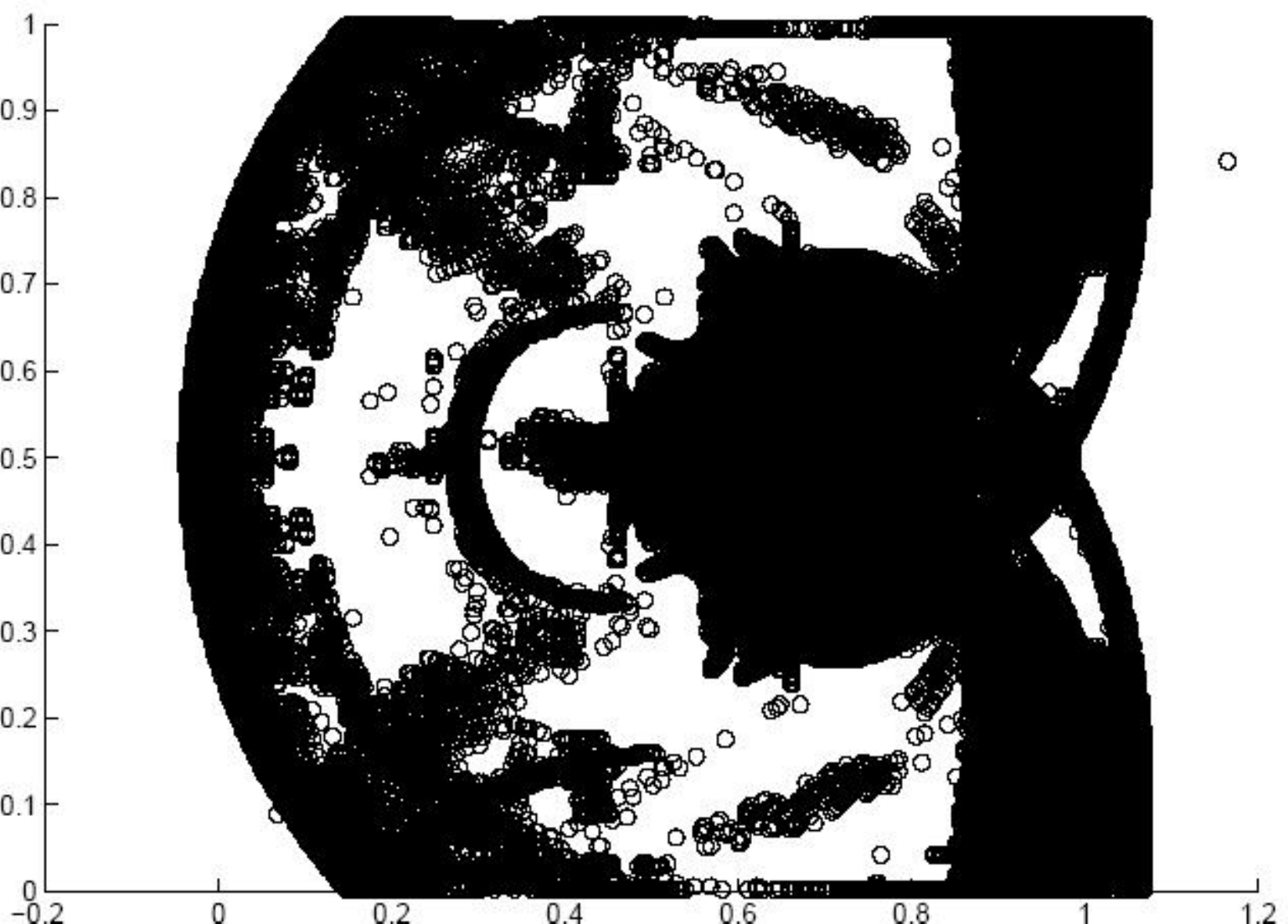}\\
    \end{tabular}
    \caption{Same as Fig.~\ref{fig:RMHDSCrho} except for the ``troubled'' cells.}
    \label{fig:RMHDSCcell}
  \end{figure}

\section{Conclusions}
\label{Section-conclusion}

 The paper developed  Runge-Kutta $P^K$-based non-central and central  DG methods with WENO limiter to the one- and two-dimensional special relativistic magnetohydrodynamical (RMHD) equations, $K=1,2,3$. The non-central DG methods were locally divergence-free, while the central DG methods were ``exactly'' divergence-free  but had to find
two approximate solutions defined on mutually dual meshes.
For each mesh, the central DG approximate solutions on its dual mesh were used to calculate the flux values in the cell and on the cell boundary so that the approximate solutions on mutually dual meshes were coupled with
each other, and the use of numerical flux might be avoided. In addition the central DG methods allowed
the use of a larger  CFL number.
The adaptive WENO limiter  was directly implemented for the physical variables
$\vec R=(D,m_x,m_y,m_z,B_z,E)^T$
by two steps: the ``troubled'' cells were first identified by using a modified TVB minmod function,
and then  new  polynomials of degree $(2K+1)$ inside the ``troubled'' cells were
locally reconstructed to replace the non-central or central DG solutions
by using the  WENO technique based on the cell average values of the DG solutions
in the neighboring cells as well as the original cell averages of the ``troubled'' cells.

However, in order that the WENO limiting procedure did not destroy  the locally or ``exactly'' divergence-free property of magnetic field, it should be
specially implemented for the magnetic field $\vec Q=(B_x,B_y)^T$.
In view of the fact that the non-central DG methods  used
a piecewise polynomial approximation of the magnetic field which satisfied the divergence free property locally, our WENO limiting procedure
also used the polynomial
satisfying the divergence free property to give a new approximation
of the magnetic field so that  the approximate magnetic field of non-central DG methods with
 WENO limiter is still locally divergence free.
In the ``exactly'' divergence-free  central DG methods, the approximated
normal magnetic field was first obtained by solving the governing equation
of magnetic field on the cell boundary, and then
used to reconstruct the new  magnetic field within the cell,
which is locally divergence-free in the cell and whose normal component
was continuous across the cell boundary.
Thus the WENO limiting procedure
was first applied to the approximated normal magnetic field on the cell boundary
and then  the ``exactly'' divergence-free WENO magnetic field is reconstructed.
Because the WENO limiter was only employed for finite ``troubled'' cells,
the computational cost can be as little as possible.

Several test problems in one and two dimensions were solved by using
our locally and ``exactly'' divergence-free DG methods with WENO limiter.
The numerical results demonstrated
that our methods were stable, accurate, and robust in resolving complex wave structures.
From the results of two-dimensioanl examples, it was seen that
the resolution of   ``exactly'' divergence-free central DG methods was
   slightly higher
    than the local divergence free non-central DG methods, but the CPU times of central DG methods were longer.
In solving RMHD problems with large Lorentz factor, or strong discontinuities,
or low rest-mass density or pressure etc., it is still possible for the $P^K$-based
non-central and central  DG methods
to give nonphysical solutions. To cure such difficulty, the $P^0$-based methods
may be locally used to replace  the $P^K$-based.
The genuinely effective way is to employ the physical-constraints preserving methods, see e.g.
\cite{Wu-Tang-RMHD2016}.

\section*{Acknowledgements}
This work was partially supported by the National Natural Science
Foundation of China (Nos. 91330205 \& 11421101).

\bibliography{ref/journalname,ref/pkuth}
           \bibliographystyle{plain}

\end{document}